%% file: thesis.tex
\documentclass[10pt]{report}

\usepackage{amsmath,amssymb,latexsym,amsthm,epsfig,stmaryrd}

\newtheorem{thm}{Theorem}[chapter]
\newtheorem{definition}{Definition}[chapter]
\newtheorem{prp}[thm]{Proposition}
\newtheorem{lem}[thm]{Lemma}
\newtheorem{cor}[thm]{Corollary}

\newenvironment{rmk}{\begin{trivlist}\item[]{\bf Remark:}}
{\end{trivlist}}
\newenvironment{ex}{\begin{trivlist}\item[]{\bf Example:}}
{\end{trivlist}}
\newenvironment{prf}{\begin{trivlist}\item[]{\bf Proof:}}
{\hfill $\blacksquare$ \end{trivlist}}
\newenvironment{lemprf}{\begin{trivlist}\item[]{\bf Proof:}}
{\hfill$\square$\end{trivlist}}

\def\part#1{\frac{\partial\phantom{q}}{\partial#1}}
\def\End{\mathop{\rm End}\nolimits}

\def\ker{\mathop{\rm ker}\nolimits}

\def\im{\mathop{\rm Im}\nolimits}
\def\deg{\mathop{\rm deg}\nolimits}

\def\cliff{\mbox{\sl Cliff}}
\def\ric{\mathop{\rm Ric}\nolimits}

\def\mod{\mathop{\rm mod}\nolimits}

\def\im{\mathop{\rm Im}\nolimits}

\newcommand{\norm}[1]{\parallel\!\!{#1}\!\!\parallel}
\newcommand{\pin}{\hat{\circ}}
\newcommand{\htimes}{\hat{\otimes}}
\newcommand{\R}{\mathbb{R}}
\newcommand{\C}{\mathbb{C}}
\newcommand{\K}{\mathbb{H}}
\newcommand{\Z}{\mathbb{Z}}
\newcommand{\T}{\mathbb{T}}
\newcommand{\N}{\mathbb{N}}
\newcommand{\Oc}{\mathbb{O}}

\newcommand{\mf}{\mathfrak}
\newcommand{\mc}{\mathcal}
\newcommand{\dstar}{d\!\star\!}

\begin{document}

\title{\huge Special metric structures and closed forms}
\author{Frederik Witt\\
New College\\
[5pt] University of Oxford\\
[8.5cm] A thesis submitted for the degree of\\
[5pt] {\em Doctor of Philosophy}\\
[1cm] Trinity Term 2004}
\date{}

\maketitle

\newpage

\thispagestyle{empty}

\newpage

\thispagestyle{empty}

\begin{abstract}
In recent work, N.~Hitchin described special geometries in terms of
a variational problem for closed generic $p$-forms. In particular,
he introduced on 8-manifolds the notion of an integrable
$PSU(3)$-structure which is defined by a closed and co-closed
3-form.

In this thesis, we first investigate this $PSU(3)$-geometry further.
We give necessary conditions for the existence of a topological
$PSU(3)$-structure (that is, a reduction of the structure group to
$PSU(3)$ acting through its adjoint representation). We derive
various obstructions for the existence of a topological reduction to
$PSU(3)$. For compact manifolds, we also find sufficient conditions
if the $PSU(3)$-structure lifts to an $SU(3)$-structure. We find
non-trivial, (compact) examples of integrable $PSU(3)$-structures.
Moreover, we give a Riemannian characterisation of topological
$PSU(3)$-structures through an invariant spinor valued 1-form and
show that the $PSU(3)$-structure is integrable if and only if the
spinor valued 1-form is harmonic with respect to the twisted Dirac
operator.

Secondly, we define new generalisations of integrable $G_2$- and
$Spin(7)$-manifolds which can be transformed by the action of both
diffeomorphisms and 2-forms. These are defined by special closed
even or odd forms. Contraction on the vector bundle $T\oplus T^*$
defines an inner product of signature $(n,n)$, and even or odd forms
can then be naturally interpreted as spinors for a spin structure on
$T\oplus T^*$. As such, the special forms we consider induce
reductions from $Spin(7,7)$ or $Spin(8,8)$ to a stabiliser subgroup
conjugate to $G_2\times G_2$ or $Spin(7)\times Spin(7)$. They also
induce a natural Riemannian metric for which we can choose a spin
structure. Again we state necessary and sufficient conditions for
the existence of such a reduction by means of spinors for a spin
structure on $T$. We classify topological $G_2\times G_2$-structures
up to vertical homotopy. Forms stabilised by $G_2\times G_2$ are
generic and an integrable structure arises as the critical point of
a generalised variational principle. We prove that the integrability
conditions on forms imply the existence of two linear metric
connections whose torsion is skew, closed and adds to 0. In
particular we show these integrability conditions to be equivalent
to the supersymmetry equations on spinors in supergravity theory of
type~ IIA/B with NS-NS background fields. We explicitly determine
the Ricci-tensor and show that over compact manifolds, only trivial
solutions exist. Using the variational approach we derive weaker
integrability conditions analogous to weak holonomy $G_2$. Examples
of generalised $G_2$- and $Spin(7)$-structures are constructed by
the device of T-duality.
\end{abstract}

\newpage

\tableofcontents

\newpage

\begin{center}
\textbf{\large Acknowledgements}
\end{center}

I want to acknowledge the Engineering and Physical Sciences Research
Council (EPSRC), the Deutscher Akademischer Austausch Dienst (DAAD),
the University of Oxford and the Studienstiftung des deutschen
Volkes for various financial support over the past three years. I
also want to thank the SFB 288 ``Differentialgeometrie und
Quantenphysik", the international Graduiertenkolleg ``Arithmetic and
Geometry" and the Graduiertenkolleg ``Strukturuntersuchungen,
Pr\"azisionstests und Erweiterungen des Standardmodells der
Elementarteilchenphysik" for several kind invitations.

I want to thank my fellow DPhil students Gil Cavalcanti and Marco
Gualtieri for many helpful discussions. Special thanks go to Wilson
Sutherland and Michael Crabb (University of Aberdeen) for helping me
out with my numerous questions in topology, and to Claus Jeschek
(MPI M\"unchen) for explaining to me a great deal of String- and
M-theory, and for so many enlightening discussions we had together.

I want to express my sincere thankfulness to my former supervisor
Helga Baum (Humboldt Universit\"at zu Berlin) for introducing me to
differential geometry through her brilliant lectures and for her
continuous interest in and support of my work.

Finally, I wish to express my deep gratitude to my supervisor Nigel
Hitchin who made these three years at Oxford so enjoyable and who
proved to be a never failing source of motivation and inspiration.

\setcounter{chapter}{-1}


\chapter{Introduction}


The present thesis is part of a programme initiated by N.~Hitchin in
a series of papers~\cite{hi00},~\cite{hi01} and~\cite{hi02} which
aims to characterise geometries through closed differential forms
enjoying special algebraic properties. The algebraic property which
plays the prime r\^ole in this thesis is {\em stability}. A form
$\rho\in\Omega^p(M)$ is said to be {\em stable} if at any point, its
orbit under the natural action of $GL(T_xM)$ is open. Such a stable
form induces a reduction of the structure group of $M$ to the
stabiliser subgroup of $\rho$. Analysing this reduction process
closely yields in all interesting cases a natural volume form which
over a compact manifold can be integrated to give a
diffeomorphism-invariant function on the space of stable forms.
Since nearby forms are also stable, we can look for critical points
for the variation of this function restricted to closed forms in a
given de Rham cohomology class, thereby performing a non-linear
version of Hodge theory, and adopt the resulting equations as
integrability conditions. The appeal of this method is that it
unravels in a straightforward way the special geometrical structure
of the local moduli space for the resulting geometry.

In this thesis, we will deal with those stable forms whose
stabiliser induces a natural Riemannian metric. Now stability is
clearly a low-dimensional phenomenon since for large $n$ and $p>2$,
$$
\dim GL(n)=n^2\ll\dim\Lambda^pT^*=n!/p!(n-p)!
$$
(for $p=2$ and $n$ even, any symplectic form is stable). The
complete list of complex representations of reductive Lie groups
admitting an open orbit was given by Sato and Kimura~\cite{saki77}.
It follows that except for the well-known $G_2$-geometry, there is
only one more case of interest in dimension 8. It is associated with
the structure group $PSU(3)$ induced by a stable 3-form, where
$PSU(3)$ acts through its adjoint representation ${\rm
Ad}:PSU(3)\hookrightarrow SO(8)$. The condition to define a critical
point for the variational problem is for both $G_2$- and
$PSU(3)$-structures that the defining 3-form is closed and
co-closed.

The first task of this research project was to take further the
investigation of $PSU(3)$-structures initiated in~\cite{hi01}. First
we derive necessary topological conditions for their existence
(Section \ref{obstrucpsu3}). The embedding of $PSU(3)$ into $SO(8)$
naturally lifts to $Spin(8)$ and restricted to this lift, the vector
representation $\Lambda^1$ and the two irreducible spin
representations $\Delta_+$ and $\Delta_-$ of $Spin(8)$ coincide
(\cite{hi01} and Proposition~\ref{equivalence}). As a result, the
tangent space is isomorphic to the spinor bundles which imposes
quite severe constraints on the topology of the underlying manifold.
In particular, we establish the existence of four linearly
independent vector fields (Proposition \ref{obpsu3}). With any
principal $PSU(3)$-fibre bundle we also associate a characteristic
class in $H^2(M,\Z_3)$ which we call the {\em triality class}
(Section~\ref{trialob}). It is the obstruction to lift a
$PSU(3)$-structure to an $SU(3)$-structure and
Theorem~\ref{psu3_structure} completely settles the existence
problem for those $PSU(3)$-structures whose triality class vanishes.

Unlike $G_2$ the group $PSU(3)$ does not act transitively on a
sphere. In particular, it follows from the Berger-Simons theorem
that $PSU(3)$ does not appear as a possible holonomy group for the
Levi-Civita connection except for the symmetric space $SU(3)$ and
its non-compact dual. As a consequence, the 3-form of an integrable
$PSU(3)$-structure is not parallel with respect to the Levi-Civita
connection. Explicit examples of such integrable $PSU(3)$-structures
are described in Section~\ref{geomexclass}. Continuing with the
integrable case, we characterise reductions to $PSU(3)$-structures
in Riemannian terms through a spinor-valued 1-form which we show to
be harmonic with respect to the twisted Dirac operator on
$T^*M\otimes(\Delta_+\oplus\Delta_-)$ if and only if the structure
is integrable (Theorem~\ref{raritaschwinger}), thereby completing
and correcting an argument given in~\cite{hi01}. In particular, this
implies the existence of a Rarita-Schwinger field which is present
in various supergravity and string theories. The projection of this
spinor-valued 1-form onto the modules $T^*M\otimes\Delta_+$ and
$T^*M\otimes\Delta_-$ induces two orientation-{\em preserving}
isometries $TM\to\Delta_{\pm}$. More generally, we show that the
existence of a {\em supersymmetric map}
(Definition~\ref{susymapdef}), that is, an isometry
$TM\to\Delta_{\pm}$, induces a reductive Lie algebra structure on
$\Delta_{\mp}$ whose adjoint group preserves the metric structure
(Theorem \ref{orbit}). As a result, the structure group reduces to
the intersection of $SO(\Delta_{\pm})$ with the automorphism group
of the Lie algebra. By Cartan's classification we obtain, other than
the group $PSU(3)$ associated with the Lie algebra $\mf{su}(3)$, the
two groups $SO(3)\times SO(3)\times SO(2)$ and $SO(3)\times SO(5)$
which arise as automorphism groups inside $SO(\Delta_{\pm})$ of the
remaining Lie algebra structures $\mf{su}(2)\oplus\R^5$ and
$\mf{su}(2)\oplus\mf{su}(2)\oplus\R^2$. The corresponding
supersymmetric maps are orientation-{\em reversing}. These groups
also arise as stabilisers of certain 3-forms (albeit non-stable, cf.
Theorem \ref{orbitdecomposition}). We will briefly analyse the
resulting geometries where we impose the same integrability
conditions as in the $PSU(3)$-case (see Sections \ref{relgeo}
and~\ref{class_sol}). In particular, we shall see that an
orientation-reversing isometry $\gamma_{\pm}:TM\to\Delta_{\pm}$
associated with $SO(3)\times SO(5)$ induces an almost quaternionic
structure on the tangent space. If the $Sp(1)\cdot Sp(2)$-invariant
4-form is closed, then $\gamma_{\pm}$ is harmonic with respect to
the twisted Dirac-operator. Moreover, we will construct an example
showing that integrability does not imply the defining 3-form to be
covariant constant with respect to the induced Levi-Civita
connection.

The aforementioned classification established by Sato and Kimura
also indicated the existence of stable spinors in either of the spin
representations $\Lambda^{ev,od}$ of $Spin(n,n)$ for $n=6$ and $7$,
where the stabiliser subgroup is conjugate to a real form of
$SL(6,\C)$ or $G_2(\C)\times G_2(\C)$ respectively. If we consider
the spin structure associated with the orthogonal vector bundle
$T\oplus T^*$ together with contraction as inner product of
signature $(n,n)$, then spinors can be naturally identified with
differential forms. In~\cite{hi02}, Hitchin dealt with the
6-dimensional case and introduced the notion of a {\em generalised
Calabi-Yau manifold} which is defined by a stable form with
stabiliser $SU(3,3)$. From the algebra of the $Spin(12,\C)$-action
he constructed a diffeomorphism-invariant function defined on the
open orbit and a similar variational principle as for forms of pure
degree applies. Moreover, this function is also invariant under the
natural action of closed 2-forms. At any point of the manifold, a
2-form (or a {\em $B$-field} in accordance with physicists' jargon)
can be regarded as an element of $\mf{so}(6,6)=\Lambda^2(T_x\oplus
T^*_x)$ and exponentiation to $Spin(6,6)$ yields an action on all
geometric objects living in a representation space of this group. In
particular, this action preserves the orbit structure of
$\Lambda^{ev,od}$ and thus the induced geometry. This observation
subsequently led to the broader concept of a {\em geometry with
$B$-fields} or {\em generalised geometry} associated with the
$SO(n,n)$-structure on $T\oplus T^*$. A ``special" generalised
geometry designates a $G$-structure inside this $SO(n,n)$-bundle, a
concept which appears all over this thesis.

Proposition~\ref{speccase} suggested the existence of stable spinors
whose stabiliser subgroup is conjugate to $G_2\times G_2$ and
$G_2(\C)$. The first one induces a natural metric giving rise to
what we call a {\em generalised $G_2$-structure} and the
investigation of those constitutes the second part of the thesis.
There is also an invariant function which was first explicitly
described by Gyoja in purely algebraic terms~\cite{gy90}. However,
with a view towards setting up the variational principle for these
structures, this formulation proved to be rather cumbersome for our
purposes. We therefore approached these structures from a
group-theoretic point of view and discussed the geometry in terms of
$G$-structures on the manifold. The inclusion $G_2\times
G_2\hookrightarrow Spin(7)\times Spin(7)$ prescribes a $G_2\times
G_2$-structure of the vector bundle $T\oplus T^*$. First the
inclusion $G_2\times G_2\hookrightarrow SO(7)\times SO(7)$ induces
what we call a generalised metric structure (Section
\ref{gen_met_struc}), i.e. an orthogonal splitting
$$
T\oplus T^*=V_+\oplus V_-
$$
into a positive and negative definite subbundle $V_+$ and $V_-$ with
respect to the inner product on $T\oplus T^*$. Equivalently, this
splitting can be encoded in a linear map $P:T\to T^*$, and taking
the symmetric and skew-symmetric part of its dual $P\in T^*\otimes
T^*$ induces an honest Riemannian metric $g$ and a 2-form $b$ on the
manifold. The $G_2\times G_2$-structure then defines a
$G_2$-structure on $V_+$ and $V_-$ which can be pulled back to yield
reductions to two $G_2$-subbundles inside the $SO(7)$-bundle
associated with $g$. Taking the induced spin structure, this can be
rephrased by the existence of two unit spinors $\Psi_+$ and
$\Psi_-$. Interpreting the square $\Delta\otimes\Delta$ of the
irreducible spin representation space $\Delta$ of $Spin(7)$ as the
module of either even or odd forms (and thus as spinors for
$Spin(7,7)$), we show that any stable form whose stabiliser is
conjugate to $G_2\times G_2$ can be written as
\begin{equation}\label{decomposablerep}
\rho=e^{-F}\cdot e^{b/2}(\Psi_+\otimes\Psi_-)
\end{equation}
(Corollary~\ref{lmapcor}). The scalar $F$ is the so--called {\em
dilaton}. Here it appears as a scaling factor but as becomes
apparent in Section~\ref{gengeom}, it is assuredly the same object
as the well-known scalar field which occurs in various string
theories. The identification of the tensor product of spinors with
forms is classical, and a special case of the
formula~(\ref{decomposablerep}) was already considered
in~\cite{gmpw02}. What is new in our approach is to interpret the
tensor product $\Delta\otimes\Delta$ not as a $Spin(7)$-module but
as a module for $Spin(7)\times Spin(7)$ which also acts on even or
odd forms through the natural map $Spin(7)\times Spin(7)\to
Spin(7,7)$.

The same approach also makes sense for $Spin(7)$-structures over
8-manifolds and an analogous formula to~(\ref{decomposablerep})
yields an even or odd form whose stabiliser is conjugate to
$Spin(7)\times Spin(7)$ (Section~\ref{ges}). However, as in the
classical case, such a form is not generic. The parity of $\rho$
depends on the chirality of the spinors $\Psi_+$ and $\Psi_-$ and in
this way, we obtain the natural notion of a {\em generalised
$Spin(7)$-structure of even} or {\em odd type}. More generally it
makes sense to consider {\em generalised $G$-structures} associated
with $G\times G\hookrightarrow SO(n,n)$ for any subgroup $G\leqslant
SO(n)$. For instance, Gualtieri introduced generalised K\"ahler
structures associated with $U(n)\times U(n)$~\cite{gu03}.

Any classical $G_2$- and $Spin(7)$-structure induces a canonical
generalised $G_2$- and $Spin(7)$-structure where the two $G_2$- or
$Spin(7)$-bundles inside the principal $SO(7)$- or $SO(8)$-bundle
coincide. From this, necessary and sufficient conditions for the
existence of reductions to $G_2\times G_2$ and $Spin(7)\times
Spin(7)$ readily follow (Proposition~\ref{gentopob}). We will also
prove the existence of {\em exotic} structures, that is the two
reductions to a $G_2$- and $Spin(7)$-subbundle inside the principal
$Spin(7)$- or $Spin(8)$-bundle cannot be transformed into each other
by a vertical, that is, fibrewise homotopy. In the generalised
$G_2$-case we will actually classify these structures up to vertical
homotopy by an integer invariant which essentially counts (with an
appropriate sign convention) the number of points where the two
$G_2$-structures coincide (Theorem~\ref{intersection}).

We move on to set up the variational formalism to derive various
notions of integrability (Sections~\ref{uvpgen},~\ref{cvpgen} and
\ref{tvpgen}). The first one is that of {\em strong} integrability
which means that $\rho$ is closed and co-closed in an appropriate
sense. Trivial examples are $G_2$- or $Spin(7)$-structures whose
defining form is closed and co-closed. We reformulate this condition
in terms of the right-hand side data of~(\ref{decomposablerep})
which leads to the main result of this thesis.
Theorem~\ref{integrability} characterises integrable generalised
$G_2$- and $Spin(7)$-structures in terms of two linear metric
connections whose torsion is skew, closed and adds to 0, and a
further equation which links the differential of the dilaton with
the torsion 3-form $T$,
$$
(dF\pm\frac{1}{2}T)\cdot\Psi_{\pm}=0.
$$
For the first variational principle we set up in
Section~\ref{uvpgen}, the resulting torsion is just $db$ which we
can think of as ``internal" or ``intrinsic" torsion coming from the
omnipresent twist by the B-field $b$. In presence of an additional
closed 3-form $H$ (a bosonic background field in physicists'
terminology), we can also consider a twisted version of this
principle taking place over a $d_H$-cohomology class, where $d_H$ is
the twisted differential operator on forms defined by
$d_H\alpha=d\alpha+H\wedge\alpha$. This 3-form then appears as
``external torsion", that is $T=db+H$. Connections with
skew-symmetric torsion gained a lot of attention in the recent
mathematical literature (see, for
instance,~\cite{agfr04},~\cite{chsw04},~\cite{friv02},~\cite{friv03}
and~\cite{iv01}) due to their importance in string and M-theory, and
eventually, our reformulation yields the supersymmetry equations in
supergravity of type~ IIA/B with bosonic background
fields~\cite{gmpw02}. Therefore we drop any distinction between
``internal" and ``external" torsion and work with $T$ rather than
with $db$ or $H$. A further benefit of this approach is an explicit
formula for the Ricci tensor of an integrable generalised $G_2$- and
$Spin(7)$-structure (Theorem~\ref{riccitensor} and
\ref{riccitensorspin7}). As a striking consequence of the closedness
of the torsion we obtain the following version of a no-go theorem,
by using arguments going back to~\cite{ivpa01} (the result is wrong
for non-closed $T$). Over a compact manifold, any integrable
generalised $G_2$- or $Spin(7)$-structure is induced by two
classical $G_2$- or $Spin(7)$-structures whose defining spinor is
parallel with respect to the Levi-Civita connection. In particular,
no exotic structure can be strongly integrable
(Corollary~\ref{vanishingg2} and Corollary~\ref{vanishing}). This
provides a full solution to the variational problem which in this
sense does not give rise to any truly new geometry, contrasting
sharply with the $U(n)\times U(n)$-geometry mentioned above, for
which compact examples that are not defined by two integrable
K\"ahler structures do exist~\cite{gu03}. This motivated us to
consider a constrained variational problem (as introduced
in~\cite{hi01}) which subsequently led to the notion of a {\em
weakly} integrable generalised $G_2$-structure of {\em even} or {\em
odd type} (there is no meaningful $Spin(7)$-equivalent). This notion
corresponds to weak holonomy $G_2$, but as the discussion in
Section~\ref{exgeneral} reveals, it defines a new type of geometry
in the sense that a weak holonomy $G_2$-manifold does not induce a
weakly integrable generalised $G_2$-structure.

Setting up all these notions would, of course, prove fruitless if
they were not supplied by non-trivial examples. In the context of
generalised geometry, we dispose of the powerful device of T-duality
to construct solutions. This duality was known for some time to
relate the string theories of type~ IIA and IIB, but the first
rigorous mathematical formulation was given only recently in a paper
by Bouwknegt et al.~\cite{bem03} (see also~\cite{ca05}). Starting
with a principal $S^1$-fibre bundle which carries an integrable
generalised structure, we define a second topological generalised
structure by exchanging the $S^1$-fibre with another one, but
without destroying integrability. With this tool at hand, we can
easily construct local examples of strongly integrable structures
with closed torsion out of classical $G_2$- or $Spin(7)$-structures
(Section~\ref{exgeneral}). T-duality also applies to weakly
integrable structures. However, the lack of a classical counterpart
makes a straightforward application difficult. If examples of weakly
integrable structures exist remains a question to be settled.

\bigskip

As we have seen, all the geometries we introduce give rise to a
metric and can be defined by closed forms. Consequently,
$G$-structures play a predominant r\^ole, and this requires a
thorough understanding of the representation-theoretic aspects of
the groups involved. The recurrent theme of this thesis is the
equivalent description of these $G$-structures by means of spinors.
We therefore chose to deal simultaneously with both {\em classical}
structures (which are defined by a form of {\em pure} degree) and
{\em generalised} structures (which are defined by an even or odd
form which we regard as a spinor for $T\oplus T^*$). As a result, we
organised the material as follows.

In the first chapter, we discuss the linear algebra of
supersymmetric maps which identify vectors isometrically with
spinors. Particular emphasis will be given to the {\em triality
principle} which states that in dimension 8, spinors and vectors
have the same internal structure. Since the decomposition into
irreducible $Spin(n)$-modules of the space of linear maps from the
tangent space to the spinor space yields two pieces of ``spin 1/2"
and ``spin 3/2", we can distinguish supersymmetric maps accordingly,
each type giving rise to different geometries. From this point of
view, we can characterise classical and generalised structures by
the existence of one or two supersymmetric maps of according spin.

We are then in a position to deal with global existence issues of
topological $G$-structures. In the second chapter, we first analyse
topological reductions to $PSU(3)$ before moving on to generalised
$G_2$- and $Spin(7)$-structures.

The third and final chapter describes several notions of
integrability. If the defining form is stable, we can set up various
variational principles which provide us with a natural set of
integrability conditions. We then adopt these for the non-stable
cases, too. The discussion of the representation-theoretic aspects
laid down in the first chapter is taken further for both classical
and generalised structures in order to derive geometrical properties
of those. Finally, we construct concrete examples.


\chapter{Linear algebra}\label{linearalgebra}


In this chapter we deal with the representation-theoretic aspects of
some groups which give rise to a natural metric $g$ on a vector
space $T$ and also stabilise a form of either pure degree (inducing
what we call a {\em classical} structure for reasons which become
apparent in Section~\ref{ges}) or an even or odd form (subsequently
referred to as a {\em generalised} structure). The primary aim is
first to gain a better grasp of the linear algebra of generic or
{\em stable} forms (Section~\ref{stableforms}) associated with the
groups $PSU(3)$ and $G_2\times G_2$. The latter case provides an
example of a generalised structure, where even or odd forms can be
naturally identified with spinors for the orthogonal vector space
$T\oplus T^*$ with contraction as inner product of signature $(n,n)$
(Section~\ref{so(n,n)}).

In the context of these ``metric" structures, spinors associated
with the metric $g$ on $T$ naturally appear and motivate an approach
in Riemannian terms. Moreover, stability is essentially a
low-dimensional phenomenon and as such, special isomorphisms between
vectors and spinors exist. This idea will be formalised by the
concept of a {\em supersymmetric map} (Definition~\ref{susymapdef}).
Such a map identifies the tangent space isometrically with the
spinor space and completely characterises the geometry. Since the
decomposition into irreducible $Spin(n)$-modules of the space of
linear maps from the tangent space to the spinor space yields two
pieces of ``spin 1/2" and ``spin 3/2" we can distinguish
supersymmetric maps accordingly, each type giving rise to different
geometries (Sections~\ref{class_susy} and~\ref{gravitino}).
Particular emphasis will be given to the {\em triality principle}
which states that in dimension 8, spinors and vectors have the same
internal structure (Section~\ref{triality}). Here, a supersymmetric
map of spin 1/2 leads to the well-known $Spin(7)$-structures, but
for 3/2, we get something new. Theorem~\ref{orbit} states that a
supersymmetric map $T\to\Delta_{\pm}$ induces a reductive Lie
algebra structure on $\Delta_{\mp}$ whose adjoint group preserves
the metric structure. Other than the group $PSU(3)$ associated with
the Lie algebra $\mf{su}(3)$, the groups $SO(3)\times SO(3)\times
SO(2)$ and $SO(3)\times SO(5)$ associated with
$\mf{su}(2)\oplus\mf{su}(2)\oplus\R^2$ and $\mf{su}(2)\oplus\R^5$
appear. These groups can be equivalently defined by a form which,
however, is not stable. We will investigate these structures in
Section~\ref{relgeo}.

We then move on to generalised structures. These share the special
feature of being invariant not only under linear isomorphisms, but
also under a natural action of a 2-form or {\em $B$-field} in
physicists' language. We spend some time on describing this setup
(Sections~\ref{so(n,n)} and ~\ref{gen_met_struc}) before introducing
what we call {\em generalised exceptional structures}
(Section~\ref{ges}) which are associated with the groups $G_2\times
G_2$ and $Spin(7)\times Spin(7)$. These stabilise a spinor for
$Spin(7,7)$ and $Spin(8,8)$ which, as we have just explained above,
can be naturally identified with a form. Propositions~\ref{linalg}
and~\ref{lmap} show how these structures fit into the spinorial
picture which will be a key theme of this thesis.

For a more detailed discussion and proofs on the general facts we
use we suggest Harvey's book~\cite{ha91}.

\section{Clifford algebras and their representations}\label{clifford}

In this section we recall some fundamental properties of Clifford
algebras to fix our notation.

\subsection{Clifford algebras}

For any inner product space $(T,g)$ of dimension $n$, we can
consider the {\em Clifford algebra} $\cliff(T,g)$ which is generated
by $T$ and a unit $\mathbf{1}$ subject to the relations
$$
v\cdot v=-g(v,v)\mathbf{1}
$$
for any $v\in T$. To ease notation, we shall drop the inner product
$g$ and write simply $\cliff(T)$ whenever the underlying metric is
clear from the context. Note that as a vector space, the map
\begin{equation}\label{Jmap}
J:v_{i_1}\cdot\ldots\cdot v_{i_k}\mapsto v_{i_1}\wedge\ldots\wedge
v_{i_k}
\end{equation}
identifies $\cliff(T)$ and $\Lambda^*T^*$. We will usually abuse
notation and write informally $v_{i_1\ldots i_k}$ for either
$v_{i_1}\cdot\ldots\cdot v_{i_k}$ or $v_{i_1}\wedge\ldots\wedge
v_{i_k}$. In the same vein, we speak of the $\Z$- and $\Z_2$-grading
on Clifford algebras or exterior forms. For instance, if
$J(x)\in\Lambda^pT$ then we shall refer to $p$ as the {\em degree}
of $x$ and denote it by $\deg(x)$. If $x\in\cliff(T)$ and
$v\in T$ then
$$
J(v\cdot x)=v\wedge x-v\llcorner x,\quad J(x\cdot
v)=(-1)^{\deg(x)}(v\wedge x+v\llcorner x).
$$
A Clifford algebra comes along with several canonical involutions.
We denote by $\sigma$ the most useful for our purposes. It is an
anti-automorphism of $\cliff(T)$ defined on elements $x$ of degree
$p$ by
\begin{equation}\label{sigma}
\sigma(x)=\epsilon(p)x
\end{equation}
where $\epsilon(p)=1$ for $p\equiv 0,3\: {\rm mod}\:4$ and
$\epsilon(p)=-1$ for $p\equiv 1,2\: {\rm mod}\: 4$. The Hodge star
operator $\star$ on $\Lambda^*$ then translates into
\begin{equation}\label{cliffhodge}
\star J(x)= J(\sigma(x)\cdot vol),
\end{equation}
where $vol$ now denotes the Riemannian volume element seen as an
element in $\cliff(T)$.

Finally, we briefly recall some aspects of the representation theory
of Clifford algebras. These can be identified with (sometimes a
direct sum of) endomorphism algebras over the so-called pinor space
$P$ which is unique up to isomorphism. We shall usually denote by
$\cdot$ the action of the Clifford algebra on $P$ which is also
called {\em Clifford multiplication}. The pinor spaces come equipped
with a bilinear form $q$ (this letter being henceforth reserved for
it) satisfying
\begin{equation}\label{q_prod}
q(x\cdot\phi,\psi)=q(\phi,\sigma(x)\cdot\psi).
\end{equation}
This property renders $q$ unique up to a scalar. If the signature of
$g$ is $0\mod 8$, then $q$ is symmetric if and only if $p=[n/2]$
equals $0$ or $3\mod 4$ (which, for us, will always be fulfilled),
and skew-symmetric otherwise. If we restrict Clifford multiplication
to even elements, then the pinor space $P$ can be decomposed into
irreducible components referred to as {\em spin representations}. In
particular, we obtain representations for the so--called {\em spin}
group
$$
Spin(T,g)=\{v_1\cdot\ldots\cdot v_{2k}\:|\:g(v_i,v_i)=\pm 1\}.
$$
For $a\in Spin(T,g)$ we define $\pi_0(a)\in SO(T,g)$ by
$\pi_0(a)v=a\cdot v\cdot a^{-1}$. The map
\begin{equation}\label{cover}
\pi_0: Spin(T,g)\to SO(T,g)
\end{equation}
is in fact a double covering. The induced Lie algebra isomorphism is
given on basis elements $e_i\cdot e_j$ by
$$
\pi_{0*}:e_i\cdot e_j\in\mf{so}(T,g)\leqslant\cliff(T,g)\to
2e_i\wedge e_j\in\mf{so}(T,g)=\Lambda^2T^*.
$$
The construction of explicit representations in the cases we will
subsequently consider will occupy us next.

\subsection{The triality principle}\label{triality}

Roughly speaking, the triality principle asserts that in dimension 8, spinors and vectors have the same
internal structure. It therefore plays a
central r\^ole if we wish to approach special geometries in
dimensions 7 and 8 from a ``supersymmetric'' point of view.

Recall that the vector representation $\Lambda^1$ and the two
irreducible spin representations of $Spin(8)$ are all 8-dimensional
and real. A convenient way to think of these spaces is to adopt the
octonions $\Oc$ as the underlying vector space. More concretely, let
us fix an orthonormal basis $e_1,\ldots,e_8$ in $\Lambda^1$ and
identify these vectors with the standard basis $1,i,\ldots,e\cdot k$
of $(\Oc,\norm{\cdot})$ according to the prescription
\begin{equation}\label{octonions}
e_0\equiv 1,\:e_1\equiv i,\:e_2\equiv j,\:e_3\equiv k,\:e_4\equiv
e,\:e_5\equiv e\cdot i,\:e_6\equiv e\cdot j,\:e_7\equiv e\cdot k.
\end{equation}
Here, $\norm{\cdot}$ denotes the norm induced by the standard inner
product $(x,y)={\rm Re}(x\cdot\bar{y})/2$. If $R_u$ denotes right
multiplication by $u\in\Oc$, then the map
\begin{equation}\label{cliff(8)}
u\in\Oc\mapsto \left(\begin{array}{cc} 0 & R_u\\ -R_{\bar{u}} & 0
\end{array}\right)\in\End(\Oc\oplus\Oc)
\end{equation}
extends to an isomorphism $\cliff(\Oc)\cong\End(\Oc\oplus\Oc)$ where
$\Delta=\Oc\oplus\Oc$ is the (reducible) space of spinors for
$Spin(8)$. These two summands can be distinguished after fixing an
orientation, since the Clifford volume element acts on those by
$\pm{\rm Id}$. Consequently, we will denote these two spin
representation spaces by $\Delta_+$ and $\Delta_-$. Then, if
$E_{ij}$ is the basis of skew-symmetric matrices $\Lambda^2$ given
by
\begin{eqnarray*}
E_{ij} & = & \left(\begin{array}{ccccc} 0 & \ldots & \ldots & \ldots & 0\\
\ldots & \ldots & \ldots & -1 & \ldots\\ \ldots & \ldots & \ldots & \ldots & \ldots\\
\ldots & 1 & \ldots & \ldots & \ldots \\ 0 & \ldots & \ldots &
\ldots & 0\end{array}\right)
\begin{array}{c} \\ \ldots i \\ \\ \ldots j\\ \\\end{array}\\
& &\qquad\qquad\!\vdots \qquad \quad \; \; \;\vdots\\
& &\qquad\qquad\! i \qquad \quad \; \; \; j
\end{eqnarray*}
it is straightforward to check that the resulting matrices of the
inclusion
$\Lambda^1\cong\Oc\hookrightarrow\End(\Delta_+\oplus\Delta_-)$ are
\begin{equation}\label{spin8_rep}
\begin{array}{lcr}
e_0 & \equiv & -E_{1,9}-E_{2,10}-E_{3,11}-E_{4,12}-E_{5,13}-E_{6,14}-E_{7,15}-E_{8,16},\\
e_1 & \equiv & \phantom{-}E_{1,10}-E_{2,9}-E_{3,12}+E_{4,11}-E_{5,14}+E_{6,13}+E_{7,16}-E_{8,15},\\
e_2 & \equiv & \phantom{-}E_{1,11}+E_{2,12}-E_{3,9}-E_{4,10}-E_{5,15}-E_{6,16}+E_{7,13}+E_{8,14},\\
e_3 & \equiv & \phantom{-}E_{1,12}-E_{2,11}+E_{3,10}-E_{4,9}-E_{5,16}+E_{6,15}-E_{7,14}+E_{8,13},\\
e_4 & \equiv & \phantom{-}E_{1,13}+E_{2,14}+E_{3,15}+E_{4,16}-E_{5,9}-E_{6,10}-E_{7,11}-E_{8,12},\\
e_5 & \equiv & \phantom{-}E_{1,14}-E_{2,13}+E_{3,16}-E_{4,15}+E_{5,10}-E_{6,9}+E_{7,12}-E_{8,11},\\
e_6 & \equiv & \phantom{-}E_{1,15}-E_{2,16}-E_{3,13}+E_{4,14}+E_{5,11}-E_{6,12}-E_{7,9}+E_{8,10},\\
e_7 & \equiv &
\phantom{-}E_{1,16}+E_{2,15}-E_{3,14}-E_{4,13}+E_{5,12}+E_{6,11}-E_{7,10}-E_{8,9}.
\end{array}
\end{equation}
Moreover, the inner product on $\Oc$ can be adopted as the
$Spin(8)$-invariant inner product on $\Delta_+$ and $\Delta_-$.
Consequently the three irreducible representations $\Lambda^1$,
$\Delta_+$ and $\Delta_-$, albeit non-equivalent as representation
spaces of $Spin(8)$, coincide as Euclidean vector spaces. The representations are distinguished by the action of the volume element. With the three representations $\pi_0$, $\pi_+$ and $\pi_-$ at hand,
the triality principle states that
\begin{equation}\label{trialprin}
\pi_0(g)(x\cdot y)=\pi_-(g)(x)\cdot\pi_+(g)(y)\quad \mbox{for all
}x,y\in\Oc.
\end{equation}
Another way of formulating triality is this. The group of symmetries
of the Dynkin diagram of $Spin(8)$ coincides with the outer
isomorphism group. The symmetries are given by a reflection $\kappa$
which maps the vertex representing the fundamental representation
$\Delta_+$ to the vertex of $\Delta_-$ and a rotation $\lambda$ of
order 3. Hence $Out(Spin(8))=\mf{S}_3=$
the permutation group of three elements, and if we regard $\lambda$
and $\kappa$ acting as outer automorphisms of $Spin(8)$, we have
$$
\pi_0=\pi_+\circ\kappa\circ\lambda\mbox{ and
}\pi_-=\pi_+\circ\lambda^2.
$$
Morally this means that we can exchange by an outer automorphism any
two of the representations $\Lambda^1$, $\Delta_+$ and $\Delta_-$
while the remaining third one is fixed. We will use these outer
automorphisms later on to derive necessary topological conditions
for the existence of $PSU(3)$-structures (cf.
Section~\ref{obstrucpsu3})\label{3al1}.

\subsection{The Lie algebra $\mathfrak{so}(n,n)$ and its spin
representations}\label{so(n,n)}

Let $T^n$ be an $n$-dimensional real vector space and consider the
direct sum with its dual $T\oplus T^*$. This direct sum carries a
natural orientation by requiring an oriented basis to be of the form
$(v_1,\ldots,v_n,\xi_1,\ldots,\xi_n)$ where $\{\xi_j\}$ is the dual
basis to $\{v_j\}$. The choice of such a basis identifies the group
$SL(T\oplus T^*)$ - the orientation preserving endomorphisms of
$T\oplus T^*$ - with $SL(2n,\R)$. Furthermore, $T\oplus T^*$ carries
the natural inner product
$$
(v+\xi,v+\xi)=-\frac{1}{2}\xi(v)
$$
of signature $(n,n)$, where $v\in T$ and $\xi\in T^*$. The factor is
chosen for convenience and has no geometrical significance. To pick
this inner product means singling out a group conjugate to $O(n,n)$
inside $GL(T\oplus T^*)=GL(2n,\R)$, and requiring compatibility with the
natural orientation brings us down to the group $SO(T\oplus
T^*)=SO(n,n)$. We can then consider the Lie algebra $\mf{so}(n,n)$
inside $\mf{gl}(2n)$, the set of endomorphisms of $T\oplus T^*$.
Note that $GL(T)\leqslant SO(T\oplus T^*)$ and as a $GL(T)$-space,
we have
$$
\mf{so}(n,n)=\Lambda^2(T\oplus
T^*)=End(T)\oplus\Lambda^2T^*\oplus\Lambda^2T.
$$
Hence any $X\in\mf{so}(T\oplus T^*)$ can be written in the form
\begin{equation}\label{elem_so}
X=A\oplus b\oplus \beta=(\sum A^i_jv_i\otimes\xi_j,\sum
b_{ij}\xi_i\wedge\xi_j,\sum\beta^{ij}v_i\wedge v_j).
\end{equation}
In particular, any 2-form $b$ defines an element in the Lie
algebra $\mf{so}(n,n)$. We will often refer to such a 2-form as a
{\em B-field}. In order to understand its action on $T\oplus T^*$
through exponentiation we represent $\mf{so}(T\oplus T^*)$ as
matrices with respect to the splitting $T\oplus T^*$, i.e.
$$
\left(\begin{array}{cc} A & \beta \\ B & -A^{tr}\end{array}\right),
$$
where $A\in\mf{gl}(T)$, $b\in\Lambda^2T^*$ and $\beta\in\Lambda^2T$.
Now the action of $e^b=\exp_{SO(n,n)}(b)$ on $v\oplus\xi\in T\oplus
T^*$ is given by
\begin{equation}\label{ebaction}
e^b(v\oplus\xi)=v\oplus (v\llcorner b +\xi),
\end{equation}
which written in matrix form is
\begin{equation}\label{ebactionmatrix}
\left(\begin{array}{cc} 1 & 0\\ b &
1\end{array}\right)\cdot\left(\begin{array}{c}v\\
\xi\end{array}\right)=\left(\begin{array}{c}v\\b(v)+\xi\end{array}\right).
\end{equation}
Next we consider the Clifford algebra $\cliff(T\oplus T^*)$
associated with the inner product $(\cdot\,,\cdot)$ and construct a
pinor representation space. We define an action of $T\oplus T^*$ on
forms by
$$
(v+\xi)\bullet\tau=v\llcorner\tau+\xi\wedge\tau.
$$
As this squares to minus the identity it gives rise to an
isomorphism
$$
\cliff(T\oplus T^*)\cong End(\Lambda^*T^*).
$$
We shall reserve the symbol $\bullet$ for this particular Clifford
multiplication. The exterior algebra $\Lambda^*T^*$ becomes thus the
pinor representation space of $\cliff(T\oplus T^*)$. It splits into
the irreducible representation spaces
$$
S^{\pm}=\Lambda^{ev,od}T^*
$$
of $Spin(T\oplus T^*)$.

\begin{rmk}\label{glvs.u}
There is an embedding
$$
GL_+(T)\hookrightarrow Spin(T\oplus T^*)
$$
of the identity component of $GL(T)$ into the spin group of $T\oplus
T^*$. From the $GL_+(T)$-representation theoretic point of view we
have to twist the spin representations by $(\Lambda^nT)^{1/2}$ so
that
$$
S^{\pm}=\Lambda^{ev,od}\otimes(\Lambda^nT)^{1/2}.
$$
This is analogous to the complex case
$$
U(n)\hookrightarrow Spin^{\C}(2n)=Spin(2n)\times_{\Z_2}S^1,
$$
where the even and odd forms get twisted with the square root of the
canonical bundle $\kappa$,
$$
S^{\pm}=\Lambda_{\C}^{ev,od}T^*\otimes\kappa^{1/2}.
$$
As long as we are doing linear algebra this is a mere notational
issue but in the global situation we cannot trivialise $\Lambda^nT$
if the manifold is non-orientable. However, this will always be the
case for the class of manifolds we consider. Consequently, we omit
the twist to ease notation, relying on the context to make it clear
whether we consider the standard $GL(T)$-action on exterior forms or
the induced action on spinors.
\end{rmk}

We can define the non-degenerate bilinear form $q$ in (\ref{q_prod}) by taking the top degree of the wedge product
$\alpha\wedge\sigma(\beta)$, i.e.
$$
q(\alpha,\beta)=(\alpha\wedge\sigma(\beta))_n.
$$
It is seen to be invariant under the action of the spin group and
from the previous remark we gather that from the $GL(T)$ point of
view, $q$ takes values in $\R$ rather than in the volume forms. This
form is non-degenerate and symmetric if $n\equiv 0,3\: (4)$ and skew
if $n\equiv 1,2\: (4)$. Moreover, $S^+$ and $S^-$ are nondegenerate
and orthogonal if $n$ is even and totally isotropic if $n$ is odd.

Finally, we note that the action of a B-field $b\in\mf{so}(n,n)$ on
the form $\tau\in S^{\pm}$ exponentiated to $Spin(T\oplus T^*)$ is
given by
$$
\exp_{Spin(T\oplus
T^*)}(b)\bullet\tau=e^{b}\bullet\tau=(1+b+\frac{1}{2}b\wedge
b+\ldots)\wedge\tau.
$$
We shall therefore also write
$\exp(b)\bullet\tau=\exp(b)\wedge\tau$.

\section{Supersymmetric maps of spin 1/2 and related geometries}\label{class_susy}

As pointed out in the introduction, the purpose of this work is to
investigate the Riemannian geometries which arise in Hitchin's
variational approach. These geometries fit into a larger picture
based on the notion of what we shall call a {\em supersymmetric
map}. In physicists' language a supersymmetry is supposed to
exchange bosons (particles that transmit forces, mathematically
described as elements lying in a vector representation of the
structure group) with fermions (particles that make up matter
represented by elements in a spin representation). To formalise this
concept, we make the

\begin{definition}\label{susymapdef}
Let $T=\Lambda^1$ be a Euclidean vector space of dimension 7 or 8
and let $\Delta$ denote an irreducible spin representation of the
associated spin group $Spin(T,g)$. Then a {\em supersymmetric map}
is an isometry
$$
\gamma:\Lambda^1\to\Delta
$$
onto its image.
\end{definition}

A basic instance of this is a unit spinor $\psi$ in the irreducible
spin representation $\Delta$ of $Spin(7)$ which by Clifford
multiplication induces an isometry
$$
x\in\Lambda^1=\R^7\to x\cdot\psi\in\Delta.
$$
Picking a unit spinor is equivalent to the choice of a copy of $G_2$
inside $Spin(7)$ and it is a general fact that supersymmetric maps
are tied to special geometries. Since it is induced by a spinor we say that this
isometry is a {\em supersymmetric map of spin 1/2}. The
$1/2$ refers to the highest coefficient in the dominant weight of
the spin representation. If $x_1,\ldots,x_m$ are the coefficients of
a maximal torus of $\mf{so}(2m)$ or $\mf{so}(2m+1)$, then the spin
representations $\Delta_{\pm}$ or $\Delta$ are of highest weight
$$
\frac{1}{2}x_1+\frac{1}{2}x_2+\ldots\pm\frac{1}{2}x_m
$$
(where the $-$ occurs for $\Delta_-$). There are also
representations of higher spin and in the Section~\ref{gravitino},
we shall also consider {\em gravitinos}, that is particles of spin
$3/2$. In this section, however, we will consider geometries which
can be characterised by a spinor and give rise to a spin $1/2$
supersymmetric map in dimensions 7 or 8. Our treatment is far from
being exhaustive. For details and notations, we
suggest~\cite{chsa01} for $SU(3)$-structures and~\cite{br87} for
$G_2$- and $Spin(7)$-structures. An exposition which comes
particulary close to the spirit of this work and covers several
other interesting geometries is~\cite{gmw03}. In this article, the
authors also deal with structures which are described by two
supersymmetric maps. Such geometries will be investigated in
Sections~\ref{gen_met_struc} and~\ref{ges} where we introduce the
``generalised" setup.

\subsection{The group $SU(3)$}\label{su3}

In this subsection, we consider $T^6=\R^6$ together with an
$SU(3)$-structure. To see what this means in terms of invariant
tensors, we first endow $T$ with an orthogonal complex structure
$I$. Recall that if a real vector space $T^{2n}$ admits an almost
complex structure with complex coordinates $z_1,\ldots,z_n$, the
complexification of $\Lambda^nT^*\otimes\C$ can be decomposed into
the direct sum of complex subspaces $\bigoplus_{p+q=n}\Lambda^{p,q}$
where $\Lambda^{p,q}$ is the space spanned by $dz_{i_1}\wedge
\ldots\wedge dz_{i_p}\wedge d\bar{z}_{i_1}\wedge\ldots\wedge
d\bar{z}_{i_q}$. We now introduce an important piece of notation
(following~\cite{sa89})\label{salnot} which we shall use throughout
this text. Whenever we deal with a complex representation space $V$,
we can also consider its conjugate $\overline{V}$ where a complex
scalar $z$ acts by $\bar{z}$. For example, we have
$\overline{\wedge^{p,q}}=\wedge^{q,p}$. If $V$ arises as the
complexification of a real vector space, we denote that space by
$[V]$ which means that
$$
[V]\otimes\C=V.
$$
In particular, $V\cong \overline{V}$ as {\em complex}
representations and ${\rm dim}_{\R}[V]={\rm dim}_{\C}V$. An instance
of this are the spaces $\wedge^{p,p}$. If $V$ is not equivalent to
$\overline{V}$ then the only way to produce a real vector space is
to forget the almost complex structure and to obtain the underlying
real vector space $\llbracket V\rrbracket=\llbracket
\overline{V}\rrbracket$. Its complexification is
$$
\llbracket V\rrbracket\otimes\C=V\oplus\overline{V}
$$
and its real dimension is twice the complex dimension of $V$. For
example, we have $T\otimes\C=\Lambda^{1,0}\oplus\Lambda^{0,1}$ so
that $T=\llbracket\Lambda^{1,0}\rrbracket$.

To return to the mainstream of our development we note that the
choice of an almost complex structure $I$ brings us down to the
structure group $GL(T,I)\cong GL(3,\C)$. On the other hand, a
non-degenerate 2-form $\omega$ yields the structure group
$Sp(T,\omega)\cong Sp(6,\R)$. If this form lies in
$[\Lambda^{1,1}]$, then the structure group reduces to the
intersection of $GL(3,\C)$ and $Sp(6,\R)$ which can be $U(2,1)$ or
$U(3)$. The latter occurs if the metric
$$
g(x,y)=\omega(Ix,y)
$$
is positive definite. This metric also distinguishes a unit circle
inside $\llbracket\Lambda^{3,0}\rrbracket$. A reduction to $SU(3)$
is then achieved by the choice of a 3-form $\psi_+$ lying in that
circle. This form determines a holomorphic volume form
$$
\psi=\psi_++i\psi_-\in\Lambda^{3,0},
$$
where $\psi_-=I\psi_+$. Choosing a suitable orthonormal basis
$e_2,\ldots,e_7$ we can write
\begin{eqnarray*}
\omega & = & e_{32}+e_{54}+e_{67}\\
\Psi & = & (e_3+ie_2)\wedge(e_5+ie_4)\wedge(e_6+ie_7)
\end{eqnarray*}
and decomposing $\Psi$ into the real and imaginary part yields
\begin{eqnarray*}
\psi_+ & = & e_{356}-e_{347}-e_{257}-e_{246}\\
\psi_- & = & e_{357}+e_{346}+e_{256}-e_{247}
\end{eqnarray*}
where as usual we identified $e_i$ with its metric dual $e^i$ and
$e_{ijk}$ is shorthand for $e_i\wedge e_j\wedge e_k$. The rather
unorthodox choice of the normal form for $\omega$ and $\psi_{\pm}$
in comparison to the standard notation~\cite{chsa01} will be
justified in Proposition~\ref{normalform}. Note also the identities
\begin{eqnarray*}
\omega\wedge\psi_{\pm} & = & 0\\
\psi_+\wedge\psi_- & = & \frac{2}{3}\omega^3.
\end{eqnarray*}

Under the action of $SU(3)$, the exterior powers $\Lambda^p$
decompose as follows.

\begin{prp}$\!\!\!${\rm~\cite{chsa01}}\hfill
$$
\begin{array}{rcl}
T^* & = & \llbracket\Lambda^1\rrbracket\\
\Lambda^2T^* & = &
\llbracket\Lambda^{1,0}\rrbracket\oplus[\Lambda^{1,1}_0]\oplus\R\omega\\
\Lambda^3T^* & = &
\R\psi_+\oplus\R\psi_-\oplus\llbracket\odot^{2,0}\rrbracket\oplus\llbracket\Lambda^{2,1}_0\rrbracket\oplus\llbracket\Lambda^1\rrbracket\\
\end{array}
$$
\end{prp}
Note that the induced Hodge star operator determines an isomorphism
$\Lambda^p\cong\Lambda^{6-p}$.

Next we want to give a characterisation of $SU(3)$-structures in
Riemannian terms, that is we start with a metric on $T$ and ask for
the additional invariants to achieve a reduction from $SO(6)$ to
$SU(3)$. To this effect, we will discuss $SU(3)$-structures from a
spinorial point of view. The pinor space of $\cliff(T,g)$ decomposes
into two irreducible spin representation spaces $S$ and
$\overline{S}$ of $Spin(6)$ which are complex conjugate to each other. In fact
there is a special isomorphism $SU(4)=Spin(6)$ and under this
identification, the spin representation spaces become the vector
representations $S=\C^4$ and $\overline{S}=\overline{\C}^4$. Since
$$
SU(4)/SU(3)\cong S^7,
$$
we see that $SU(3)$ is the stabiliser of a spinor in $S$ and
$\overline{S}$. If we denote these spinors by $\psi$ and
$\overline{\psi}$ and choose a representation of the Clifford
algebra $\cliff(\R^6,g)$, we can reconstruct the forms $\omega$ and
$\Psi$ by
\begin{eqnarray*}
\omega(e_i,e_j) & = & -iq(e_{ij}\cdot\psi,\psi)\\
\Psi(e_i,e_j,e_k) & = & q(e_{ijk}\cdot\psi,\overline{\psi}).
\end{eqnarray*}
There is clearly some choice involved in this process -- we followed
the convention in~\cite{gmw03}.

\subsection{The group $G_2$}\label{g2}

Next we move up one dimension and consider $T=\R^7$ together with
the structure group $G_2$. It is one of the two so-called
exceptional holonomy groups among Berger's list of special holonomy
groups. There are several closely interrelated ways of approaching
$G_2$. We can see $G_2$ inside $GL(8)$ as the automorphism group of
the product structure of the octonions. This product structure
induces a vector cross product on the imaginary
octonions~\cite{brgr67} which is also invariant under the action of
$G_2$. From this orthogonal product we can build a 3-form which is
generic and stabilised by $G_2$ and this is the definition of $G_2$
we shall adopt. More formally, consider the natural action of
$GL(7)$ on $\Lambda^3$. It has two open orbits one of which contains
the aforementioned 3-form and is thus isomorphic to $GL(7)/G_2$.
This orbit is commonly denoted by $\Lambda^3_+$ and its elements are
called {\em stable}, following the language of~\cite{hi01}. Stable
forms will be treated in more detail in Section~\ref{stableforms}.
We will usually denote a given $G_2$-invariant 3-form by $\varphi$.
Since $G_2$ sits inside $SO(7)$, the choice of a $G_2$ structure
also induces an orientation and a metric $g$. The exterior algebra
decomposes as follows.

\begin{prp}$\!\!\!${\rm~\cite{br87}}\label{g2_form}\hfill\newline
{\rm (i)} $\Lambda^1=\Lambda^1_7=\R^7$

{\rm (ii)} $\Lambda^2=\Lambda^2_7\oplus\Lambda^2_{14}$ where
$$
\begin{array}{lcl}
\Lambda^2_7 & = &
\{\alpha\in\Lambda^2\:|\:\star(\varphi\wedge\alpha)=2\alpha\}=\{X\llcorner\varphi\:|\:X\in
T\}\\
\Lambda^2_{14} & = &
\{\alpha\in\Lambda^2\:|\:\star(\varphi\wedge\alpha)=-\alpha\}=\mf{g}_2
\end{array}
$$

{\rm (iii)}
$\Lambda^3=\Lambda^3_1\oplus\Lambda^3_7\oplus\Lambda^3_{27}$,
where
$$
\begin{array}{lcl}
\Lambda^3_1 & = & \{t\varphi\:|\:t\in\R\}\\
\Lambda^3_7 & = & \{\star(\varphi\wedge\alpha)\:|\:\alpha\in\Lambda^1_7\}\\
\Lambda^3_{27} & = & \{\alpha\in\Lambda^3\:|\:\alpha\wedge\varphi=0,\:\alpha\wedge\star\varphi=0\}.\\
\end{array}
$$
\end{prp}
The notation is standard and the subscript indicates the dimension
of the module. The decomposition of the higher order powers follows
again from the isomorphism $\Lambda^{7-p}\cong\Lambda^p$ induced by
the Hodge $\star$-operator, and modules of the same dimension are
equivalent. We will use the notation of Proposition~\ref{g2_form}
throughout this work.

As in the case of $SU(3)$, a spinorial formulation of
$G_2$-structures is available. The group $G_2$ is simply-connected and the
inclusion $G_2\hookrightarrow SO(7)$ lifts to $Spin(7)$ whose
irreducible spin representation $\Delta$ is real and 8-dimensional.
We have
$$
Spin(7)/G_2\cong S^7,
$$
that is, $G_2$ is the stabiliser of a unit spinor $\psi_0$ which
induces an isometry
$$
X\in\Lambda^1\mapsto X\cdot\psi_0\in\psi_0^{\perp}\leqslant\Delta.
$$
By fixing a representation of the Clifford algebra $\cliff(\R^7,g)$
we obtain an invariant 3-form by the formula
\begin{equation}\label{fierz}
\varphi(e_i,e_j,e_k)=q(e_{ijk}\cdot\psi,\psi).
\end{equation}
The representation we shall use throughout this work is constructed
as follows. Identify $T^7=\R^7$ with the imaginary octonions by
$e_1=i,\ldots,e_7=e\cdot k$ as in~(\ref{octonions}) and embed $T$
into
$$
\End_{\R}(\Oc)\oplus\End_{\R}(\Oc)\cong\cliff(T,g)
$$
by
$$
u\in T\cong{\rm Im}\,\Oc\mapsto\left(\begin{array}{cc} R_u & 0\\ 0 &
-R_u\end{array}\right),
$$
where $R_u$ denotes right multiplication as in~(\ref{cliff(8)}).
Projection on the first summand $\End_{\R}(\Oc)$ yields the explicit
matrix representation given by
$$
\begin{array}{lcr}
e_1 & \equiv & E_{1,2}-E_{3,4}-E_{5,6}+E_{7,8},\\
e_2 & \equiv & E_{1,3}+E_{2,4}-E_{5,7}-E_{6,8},\\
e_3 & \equiv & E_{1,4}-E_{2,3}-E_{5,8}+E_{6,7},\\
e_4 & \equiv & E_{1,5}+E_{2,6}+E_{3,7}+E_{4,8},\\
e_5 & \equiv & E_{1,6}-E_{2,5}+E_{3,8}-E_{4,7},\\
e_6 & \equiv & E_{1,7}-E_{2,8}-E_{3,5}+E_{4,6},\\
e_7 & \equiv & E_{1,8}+E_{2,7}-E_{3,6}-E_{4,5}.
\end{array}
$$
Using this, expression~(\ref{fierz}) becomes
\begin{equation}\label{normalg2}
\varphi=e_{123}+e_{145}-e_{167}+e_{246}+e_{257}+e_{347}-e_{356}
\end{equation}
which we shall adopt as a normal form of $\varphi$. Consequently, we
obtain
\begin{equation}\label{normalg2*}
\star\varphi=-e_{1247}+e_{1256}+e_{1346}+e_{1357}-e_{2345}+e_{2367}+e_{4567}.
\end{equation}

Finally, we want to investigate the relationship between $SU(3)$-
and $G_2$-structures. We start with the form point of view. Choose a
unit vector $\alpha\in T$ and decompose $T=\T\oplus\R\alpha$ where
$\T=\R^6$. Then~(\ref{normalg2}) and~(\ref{normalg2*}) can be
written as
\begin{equation}\label{g2su3}
\varphi=\psi_++\omega\wedge\alpha,\quad
\star\varphi=\psi_-\wedge\alpha+\frac{1}{2}\omega^2
\end{equation}
where the forms $\omega$ and $\psi_{\pm}$ are pulled back from $\T$.
As usual we identify vectors with 1-forms in the presence of a
metric. The data $g_{|\T}$, $\omega$ and $\psi_+$ defines an
$SU(3)$-structure on $\T$. Conversely, take the invariant forms
$\omega$ and $\psi_+$ and the metric $g_0$ of an $SU(3)$-structure
and introduce a 1-form $\alpha$. Define $T=\T\oplus\R\alpha$ with
the metric $g=g_0\oplus\alpha\otimes\alpha$. Then the form $\varphi$
given as in~(\ref{g2su3}) defines a reduction to $G_2$.

In the spinorial picture, let $\psi$ be the invariant spinor under
$G_2$. The choice of a unit vector $\alpha$ induces an almost
complex structure on $\Delta$, since
$\alpha\cdot\alpha\cdot\psi=-\psi$. Hence the complexification of
$\Delta$ becomes
$$
\Delta\otimes\C=\Delta^{1,0}\oplus\Delta^{0,1}.
$$
In terms of the almost complex structure $I(\psi)=\alpha\cdot\psi$
induced by $\alpha$, these spaces are
$$
\Delta^{1,0}=\{\phi-i\alpha\cdot\phi\:|\:\phi\in\Delta\}\mbox{ and
} \Delta^{0,1}=\{\phi+i\alpha\cdot\phi\:|\:\phi\in\Delta\}.
$$
Since $\Delta^{1,0}$ and $\Delta^{0,1}$ are 4-dimensional complex
vector spaces we have a natural $SU(4)=Spin(6)$ action induced by
the cover of $Spin(6)\hookrightarrow Spin(7)$ of the reduction from
$SO(7)$ to $SO(6)$ associated with the vector $\alpha$. Now
$$
\Delta=\llbracket S\rrbracket=\llbracket \overline{S}\rrbracket
$$
so that the choice of an inclusion $SU(3)\hookrightarrow SU(4)$
yields a 2-dimensional invariant real subspace spanned by $\psi_+$
and $\psi_-=\alpha\cdot\psi_+$. The group $SU(3)$ thus sits inside
the two copies of $G_{2\pm}$ stabilising the spinors $\psi_+$ and
$\psi_-$ in $\Delta$. The corresponding 3-forms $\varphi_{\pm}$ are
\begin{equation}\label{varphi+-}
\varphi_{\pm}=\omega\wedge\alpha\pm\psi_+.
\end{equation}

\subsection{The group $Spin(7)$}

Finally we discuss the geometry associated with $Spin(7)$ over
$T=\R^8$. This group is obtained as the reduction from $GL(8)$ to
the stabiliser of a certain 4-form which we denote by $\Omega$.
Again it can be constructed out of the algebraic data we have on the
octonions, but unlike the $G_2$-invariant 3-form $\varphi$ it is not
generic. However, we still have an inclusion $Spin(7)\hookrightarrow
SO(8)$ which induces a metric and an orientation, and $\Omega$ is
self-dual with respect to the corresponding Hodge $\star$-operator.
Since $Spin(7)$ is simply-connected it also lifts to $Spin(8)$.

\begin{prp}$\!\!\!${\rm~\cite{br87}}\label{spin7_form}\hfill\newline
{\rm (i)} $\Lambda^1=\Lambda^1_8=\R^8$

{\rm (ii)} $\Lambda^2=\Lambda^2_7\oplus\Lambda^2_{21}$ where
$$
\begin{array}{lcl}
\Lambda^2_7 & = &
\{\alpha\in\Lambda^2\:|\:\star(\Omega\wedge\alpha)=3\alpha\}\\
\Lambda^2_{21} & = &
\{\alpha\in\Lambda^2\:|\:\star(\Omega\wedge\alpha)=-\alpha\}=\mf{spin}(7)
\end{array}
$$

{\rm (iii)} $\Lambda^3=\Lambda^3_8\oplus\Lambda^3_{48}$, where
$$
\begin{array}{lcl}
\Lambda^3_8 & = & \{\star(\Omega\wedge\alpha)\:|\:\alpha\in T\}\\
\Lambda^3_{48} & = &
\{\alpha\in\Lambda^3\:|\:\alpha\wedge\Omega=0\}.
\end{array}
$$

{\rm (iv)}
$\Lambda^4=\Lambda^4_1\oplus\Lambda^4_7\oplus\Lambda^4_{27}\oplus\Lambda^4_{35}$
where
$$
\begin{array}{lcl}
\Lambda^4_1 & = & \{t\Omega\:|\:t\in\R\}\\
\Lambda^4_7 & = &
\{\alpha^*\Omega\:|\:\alpha\in\Lambda^2_7\}\mbox{ $($where $\alpha^*$ denotes the action of }\Lambda^2=\mf{so}(8))\\
\Lambda^4_{27} & = &
\{\alpha\in\Lambda^4\:|\:\Omega\wedge\alpha=0,\:\star\alpha=\alpha,
\:\beta\wedge\alpha=0\;{\rm for\;all }\;\beta\in\Lambda^4_7\},\\
\Lambda^4_{35} & = &
\{\alpha\in\Lambda^4\:|\:\star\alpha=-\alpha\}.
\end{array}
$$
\end{prp}

We saw earlier in our discussion of triality that vectors and
spinors are equivalent under an outer automorphism. An invariant
vector in $\Lambda^1$ reduces the structure group from $SO(8)$ to
$SO(7)$ and this is covered by an inclusion of $Spin(7)$ to
$Spin(8)$. We denote the image by $Spin(7)_0$ to indicate that it
fixes an element under the vector representation of $Spin(8)$
induced by $\pi_0$ (cf.~(\ref{cover})). The triality automorphisms
$\lambda$ and $\kappa$ permute the representation spaces and
consequently, the stabiliser of a unit spinor in $\Delta_{\pm}$ are
two copies of $Spin(7)$ which we denote by $Spin(7)_{\pm}$. Hence we
have
$$
Spin(8)/Spin(7)_{\pm}\cong S^7\subset\Delta_{\pm}
$$
and $Spin(7)_{\pm}$ acts irreducibly on $\Delta_{\mp}$.

As in the previous cases we can compute a normal form of $\Omega$ by
the formula
$$
\Omega(e_i,e_j,e_k,e_l)=q(e_{ijkl}\cdot\psi,\psi).
$$
If we fix the representation given by~(\ref{spin8_rep}) we find
\begin{eqnarray*}
\Omega & = &
-e_{0123}-e_{0145}+e_{0167}-e_{0246}-e_{0257}-e_{0347}+e_{0356}+e_{1247}-e_{1256}-e_{1346}\\
& & -e_{1357}+e_{2345}-e_{2367}-e_{4567}.
\end{eqnarray*}
To link $G_2$- with $Spin(7)$-structures we choose a unit vector
$\gamma\in T$ and decompose $T=\R\gamma\oplus\T$ where we endow
$\T=\R^7$ with a $G_2$-structure $\varphi\in\Lambda^3\T^*$. Then the
4-form
$$
\Omega=\gamma\wedge\varphi+\star\varphi
$$
is stabilised by $Spin(7)$. The spinor description follows from the
decomposition
$$
\Delta_+\cong\Delta_-\cong\Lambda^1\oplus\R\psi_{\pm}
$$
into $G_2$-modules. In fact, $Spin(8)$ acts transitively on
$S^7\times S^7$ inside $\Delta_+\times\Delta_-$ with stabiliser
$G_2$,
$$
Spin(8)/G_2=S^7\times S^7\subset\Delta_+\times\Delta_-
$$
and consequently, $G_2$ sits in and actually equals the intersection
of $Spin(7)_+$ with $Spin(7)_-$. Moreover, we can show the

\begin{prp}\label{g2ds}
Restricted to $G_2\subset Spin(8)$, the representations $\Lambda^1$, $\Delta_+$ and $\Delta_-$ coincide.
\end{prp}

\begin{prf}
We denote the four fundamental weights of $Spin(8)$ by $\omega_1$, $\omega_2$,
$\omega_3$ and $\omega_4$. The weights of the representation
$\Lambda^1=[1,0,0,0]$ are
\begin{equation}\label{weights0}
\pm\omega_1,\:\pm(\omega_1-\omega_2),\:\pm(\omega_2-\omega_3-\omega_4),\:\pm(\omega_3-\omega_4).
\end{equation}
For the representations $\Delta_+=[0,0,1,0]$ and
$\Delta_-=[0,0,0,1]$ we find the weights
\begin{equation}\label{weights+}
\pm(\omega_1-\omega_4),\:\pm(\omega_1-\omega_2+\omega_4),\:\pm(\omega_2-\omega_3),\:\pm\omega_3
\end{equation}
and
\begin{equation}\label{weights-}
\pm(\omega_1-\omega_3),\:\pm(\omega_1-\omega_2+\omega_3),\:\pm(\omega_2-\omega_4),\:\pm\omega_4.
\end{equation}
Now $\Delta_+=\Lambda^1_7\oplus\R\psi_+$ as a $G_2$-space where $G_2$ acts with weights
$$
0,\,\pm(\sigma-\tau),\,\pm\sigma,\,\pm(2\sigma-\tau),
$$
where $\sigma$ and $\tau$ are the fundamental weights of $G_2$. Substituting these into~\ref{weights0}, \ref{weights+} and \ref{weights-} gives in all three cases the same $G_2$-weights, i.e. they are equivalent.
\end{prf}

We will also consider the case where we have to orthogonal unit
spinors $\psi_{\pm}\in\Delta_+$. The spinor $\psi_+$ gives rise to
the stabiliser $Spin(7)_+$, so that $\Delta_+=\R\psi_+\oplus \R^7$.
The group $Spin(7)_+$ acts on $\R^7$ in its vector representation
and thus transitively on the sphere $S^6$ with stabiliser $Spin(6)$.
As a $Spin(6)$-space, $\Lambda^1=\R^8=\llbracket\C^4\rrbracket$,
where $\C^4$ denotes the spin representation under the
identification $Spin(6)=SU(4)$ as above. Note that $Spin(6)$
stabilises a 2-form in $\Lambda^2$, given by
$$
\varpi=e_{01}-e_{23}-e_{45}+e_{67}
$$
and the three self-dual 4-forms
\begin{eqnarray*}
\Omega_1 & = &
e_{0246}+e_{0257}+e_{0347}-e_{0356}-e_{1247}+e_{1256}+e_{1346}+e_{1357}\\
\Omega_2 & = &
e_{0247}-e_{0256}-e_{0346}-e_{0357}+e_{1246}+e_{1257}+e_{1347}-e_{1356}\\
\Omega_3 & = &
-e_{0123}-e_{0145}+e_{0167}+e_{2345}-e_{2367}-e_{4567},
\end{eqnarray*}
where $\Omega_3=\varpi\wedge\varpi/2$.

We summarised the relationship between the groups $SU(3)$, $G_2$ and
$Spin(7)$ in the next table. Note that the decomposition of the
vector and spin spaces in the first and second row are taken with
respect to the group $G_2$ and $SU(3)$.
\begin{table}[hbt]
\begin{center}
\begin{tabular}{|c|c|c|}
\hline Group & vector space & spinor space\\
\hline $Spin(8)$ & $T=\R\gamma\oplus\R^7$ &
$\Delta_{\pm}=\R\psi_{\pm}\oplus\R^7$\\
$\phantom{G_2}\cup$ & & \\
$G_2=Spin(7)_+\cap Spin(7)_-$ & $T=\R\alpha\oplus\R^6$ &
$\Delta=\R\psi_+\oplus\R\psi_-\oplus\R^6$\\
$\phantom{G_2}\cup$ & & \\
$SU(3)=G_{2+}\cap G_{2-}$ & $T=\R^6$ &
$S=\C^3\oplus\C\psi,\,\bar{S}=\overline{\C^3}\oplus\overline{\psi}$\\
\hline
\end{tabular}
\end{center}
\caption{The groups $SU(3)$, $G_2$ and $Spin(7)$}\label{overview}
\end{table}

\section{Supersymmetric maps of spin $3/2$ and related
geometries}\label{gravitino}

In dimension 8, a supersymmetric map sits inside the module
$$
\Lambda^1\otimes\Delta_{\pm}=\Delta_{\mp}\oplus\ker\mu_{\mp}
$$
where $\ker\mu_{(\pm)}$
denotes the kernel of the Clifford multiplication
$$
\mu_{\pm}:\Lambda^1\otimes\Delta_{\pm}\to\Delta_{\mp}.
$$
As a spin module, $\ker\mu$ is an irreducible space of highest
weight
$$
\frac{3}{2}x_1+\frac{1}{2}x_2+\ldots\pm\frac{1}{2}x_n
$$
and consequently, an element of this space will be referred to as a
spin $3/2$ particle. In the same vein as in the previous section, we
try to analyse the geometries associated with a supersymmetric map,
this time of spin $3/2$. As an example of such a geometry, we first
consider the group $PSU(3)$.

\subsection{The group $PSU(3)$}\label{psu3}

We shall first adopt a definition of the group $PSU(3)$ which
follows closely the approach of~\cite{br87} for $G_2$, emphasising
the form point of view. Let $T$ be a real vector space of dimension
8. We fix a basis $e_1,\ldots,e_8$ whose dual basis we denote by
$e^1,\ldots,e^8$ and choose $e^{12345678}$ as orientation.

\begin{definition}
Let
\begin{equation}\label{rho}
\rho=e^{123}+\frac{1}{2}e^1(e^{47}-e^{56})+\frac{1}{2}e^2(e^{46}+e^{57})+\frac{1}{2}e^3(e^{45}-e^{67})+\frac{\sqrt{3}}{2}e^8(e^{45}+e^{67}).
\end{equation}
Then we define
$$
PSU(3)=\{a\in GL_+(T)\:|\: a^*\rho=\rho\}.
$$
\end{definition}

The coefficients $c_{ijk}$ of~(\ref{rho}) are taken
from~\cite{gmn64}, p.50. We shall give a more detailed description
of this 3-form in a moment. First we note that as in the case of
$G_2$, we may regard $PSU(3)$ as (the connected component of) an
automorphism group.

\begin{prp}
The subgroup $PSU(3)\leqslant GL(T)=GL(8)$ is the identity component
of the automorphism group of the Lie algebra $\mf{su}(3)$, i.e.
$$
PSU(3)\cong Ad(SU(3))=SU(3)/Z(SU(3)),
$$
where $Z(SU(3))$ denotes the center of $SU(3)$.
\end{prp}

\begin{prf}\label{psu3_stab}
We define a bracket $[\cdot\,,\cdot]$ on $T\times T$ by linearly
extending the relations
$$
[e_i,e_j]=\sum_k\rho(e_i,e_j,e_k)e_k,
$$
that is, the structure constants $c^k_{ij}$ are just the
coefficients of $\rho$ and hence are totally skew. We will show that
this bracket endows $T$ with a Lie algebra structure isomorphic to
$\mf{su}(3)$. It is immediate to verify that we actually have
defined a Lie bracket, that is the Jacobi identity holds. Next,
consider the Killing form $B(x,y)={\rm Tr}(ad(x)\circ ad(y))$. It is
straightforward to check that the basis $e_1,\ldots,e_8$ verifies
$$
B(e_i,e_j)=-3\delta_{ij}
$$
so that $B$ is negative definite. Hence, $(T,[\cdot\,,\cdot])$ is a
compact semi-simple Lie algebra. By Cartan's classification theorem
we conclude it must be simple on dimensional grounds and thus
isomorphic to $\mf{su}(3)$. Hence $PSU(3)$ is contained in
$Aut(\mf{su}(3))=Ad(SU(3))\times \Z_2$. Since the orbit of $\rho$ is
of dimension less than or equal to 56, the dimension of $\Lambda^3$,
it follows that $PSU(3)$ is of dimension greater than or equal to 8
and hence equals 8 since this is the dimension of $Ad(SU(3))$.
Consequently, the identity component of the automorphism group is
contained in $PSU(3)$. Since the other component of
$Aut(\mf{su}(3))$ reverses the orientation the assertion follows.
\end{prf}

Invariantly formulated, we have
$$
\rho(x,y,z)=-\frac{1}{3}B([x,y],z)
$$
where $B$ and $[\cdot\,,\cdot]$ denotes the Killing form and the
bracket $[\cdot\,,\cdot]$ of $\mf{su}(3)$. The vectors
$e_1,\,e_2,\,e_3,\,e_8$ span a $\mf{u}(2)=\mf{su}(2)\oplus\R e_8$
subalgebra. As an $\mf{su}(2)$-space, we have
$\mf{su}(3)=\mf{su}(2)\oplus\R e_8\oplus\llbracket\C^2\rrbracket$,
where $\C^2$ denotes the vector representation of $\mf{su}(2)$ which
we view as a real space $\R^4$ spanned by $e_4,\ldots,e_7$. The
space of 2-forms $\Lambda^2\R^4$ decomposes into the direct sum of
self-dual and anti-self-dual forms (the $\pm1$-eigenspaces of the
Hodge $\star$-operator). By choosing a suitable orientation, the
anti-self-dual forms $\omega_{1-}=e^{47}-e^{56}$,
$\omega_{2-}=e^{46}+e^{57}$ and $\omega_{3-}=e^{45}-e^{67}$ are
acted on trivially by $\mf{su}(2)$, while
$\omega_{3+}=e^{45}+e^{67}$ belongs to the space of self-dual forms
which is just $\mf{su}(2)$. We can then write
\begin{equation}\label{rho1}
\rho=
e_{123}+\frac{1}{2}e_1\wedge\omega_{1-}+\frac{1}{2}e_2\wedge\omega_{2-}+\frac{1}{2}e_3\wedge\omega_{3-}
+\frac{\sqrt{3}}{2}e_8\wedge\omega_{3+}.
\end{equation}
and
\begin{equation}\label{starrho1}
\star\rho=e_{45678}-\frac{1}{2} e_{238}\wedge\omega_{-1}+\frac{1}{2}
e_{138}\wedge\omega_{-2}-\frac{1}{2}e_{128}\wedge\omega_{-3}
+\frac{\sqrt{3}}{2}e_{123}\wedge\omega_{3+}.
\end{equation}
The next corollary displays yet another feature that $PSU(3)$ has in
common with $G_2$ and which will be crucial in the sequel.

\begin{cor}
Under the natural action of $GL_+(8)$, the orbit of the 3-form
$\rho$ is diffeomorphic to $GL_+(8)/PSU(3)$ and therefore open.
\end{cor}

Again, we shall refer to a form which lies in this open orbit as
{\em stable}. Expressed in a suitable frame, a stable element has
the canonical form~(\ref{rho}). We refer to such a frame as a {\em
$PSU(3)$-frame}.

As a further corollary of Proposition~\ref{psu3_stab} we note the

\begin{cor}\label{psu3lift}
The Lie group $PSU(3)$ is a compact, connected Lie group of dimension 8
with fundamental group $\pi_1(PSU(3))=\Z_3$. It acts irreducibly on
$T$. The adjoint representation $Ad:SU(3)\to SO(8)$ descends to an
embedding $PSU(3)\hookrightarrow SO(8)$ that lifts to an embedding
$PSU(3)\hookrightarrow Spin(8)$.
\end{cor}

Hence the choice of an orientation and a stable 3-form induces a
metric $g$, a Hodge $\star$-operator $\star_g$, a Lie bracket
$[\cdot,\cdot]$ and a spin structure. In particular, we will drop
any distinction between $\{e_i\}$ and its dual basis $\{e^i\}$.

Next, we will discuss some elements of the representation theory for
$PSU(3)$. We will always suppose that we work with a $PSU(3)$-frame
$e_1,\ldots,e_8$. We will label irreducible representations by their
highest weight expressed in the basis of fundamental weights which
in the case of $\mf{su}(3)$ is provided by two elements, $\sigma_1$
and $\sigma_2$. For example, $T\otimes\C$, the complexification of
the adjoint representation of $SU(3)$, is
$T\otimes\C=V(\sigma_1+\sigma_2)$ which we more succinctly write as
$(1,1)$, the i-th component referring to $\sigma_i$. This can be
seen by fixing the Cartan subalgebra spanned by $e_3$ and $e_8$. We
obtain the following decomposition of the complexified Lie algebra
$\mf{su}(3)\otimes\C$ into root spaces:
$$
\mf{su}(3)\otimes\C=\mf{t}\oplus V(\alpha) \oplus
V(\bar{\alpha})\oplus V(\beta)\oplus V(\bar{\beta})\oplus
V(\alpha+\beta)\oplus V(\bar{\alpha}+\bar{\beta}).
$$
We use the notation $T=[1,1]$ for the {\em real} representation in
accordance with the convention introduced on page~\pageref{salnot}.
For subsequent computations we will adopt the following convention.
The fundamental roots $\alpha$ and $\beta$ will be labeled in such a
way that the root spaces $V(\alpha)$ and $V(\beta)$ are spanned by
$$
x_{\alpha}=e_4+ie_5\mbox{ and }x_{\beta}=e_6-ie_7.
$$
This corresponds to the expressions
\begin{equation}\label{cartanmatrix}
\alpha=-\frac{i}{2}(e^3+\sqrt{3}e^8)\mbox{ and
}\beta=\frac{i}{2}(-e^3+\sqrt{3}e^8)
\end{equation}
for the fundamental roots $\alpha$ and $\beta$. It follows that
$V(\alpha+\beta)$ is generated by $x_{\alpha+\beta}=e_1+ie_2$.
Moreover, the fundamental weights are
\begin{equation}\label{cmat2}
\begin{array}{lcl}
\sigma_1 & = & \frac{2}{3}\alpha+\frac{1}{3}\beta\\
\sigma_2 & = & \frac{1}{3}\alpha+\frac{2}{3}\beta.
\end{array}
\end{equation}
Note that the Lie structure on $T=\Lambda^1$ induces a
$PSU(3)$-invariant operator $d:\Lambda^k\to\Lambda^{k+1}$ which is
just the exterior differential restricted to left-invariant forms on
the group $SU(3)$,
$$
de_i=\sum_{j<k}c_{ijk}e_j\wedge e_k.
$$
The following proposition gives the decomposition of the exterior
algebra into irreducible $PSU(3)$-modules.

\begin{prp}\label{psu3decomp}\hfill
\begin{enumerate}
\item $\Lambda^1=T$ is irreducible and of real type.
\item
$\Lambda^2=\Lambda^2_8\oplus\Lambda^2_{20}$, where
$$
\begin{array}{l}
\Lambda^2_{8\phantom{0}}=\{X\llcorner\rho\:|\:X\in \Lambda^1\}=\mf{su}(3)\\
\Lambda^2_{20}=\{\tau\in\Lambda^2\:|\:\tau\wedge\star\rho=0\}
\end{array}
$$
$\Lambda^2_{20}$ is a representation of complex type and its
complexification decomposes as
$$
\Lambda^2_{20}\otimes\C=\Lambda^{2\C}_{10+}\oplus\Lambda^{2\C}_{10-}
$$
where
$$
\Lambda^{2\C}_{10\pm}=\{\alpha\in\Lambda^2\otimes\C\:|\:\star
(\rho\wedge \alpha)=\pm\sqrt{3}i\cdot\alpha^*\rho\}.
$$
Expressed in terms of dominant weights with respect to the basis provided by the fundamental weights $\sigma_1$ and $\sigma_2$, these spaces are
$$
\begin{array}{l}
\Lambda^2_{8\phantom{0+}}=[1,1]\\
\Lambda^{2\C}_{10+}=(0,3)\\
\Lambda^{2\C}_{10-}=(3,0).
\end{array}
$$
\item
$\Lambda^3=\Lambda^3_1\oplus\Lambda^3_8\oplus\Lambda^3_{20}\oplus\Lambda^3_{27}$,
where
$$
\begin{array}{l}
\Lambda^3_1=\{t\rho\:|\:t\in\R\}\\
\Lambda^3_8=\{\star((X\llcorner\rho)\wedge\rho)\:|\:X\in \Lambda^1\}\cong\Lambda^2_8\\
\Lambda^3_{20}=d\Lambda^2_{20}\\
\Lambda^3_{27}=\{\alpha\in\Lambda^3\:|\:\alpha\wedge\rho=0,\:\alpha\wedge\star\rho=0\}
\end{array}
$$
Moreover,
$$
\Lambda^3_{27}=[2,2].
$$
\item
$\Lambda^4=\{\alpha\in\Lambda^4\:|\:\star\alpha=\alpha\}\oplus
\{\alpha\in\Lambda^4\:|\:\star\alpha=-\alpha\}=\Lambda^{4+}_8\oplus\Lambda^{4+}_{27}\oplus\star\Lambda^{4-}_8\oplus\star\Lambda^{4-}_{27}$,
where $\Lambda^{4\pm}_8$ and $\Lambda^{4\pm}_{27}$ are obtained as
projections of the spaces
$$
\begin{array}{l}
\{X\llcorner\star\rho\:|\:X\in \Lambda^1\}\cong\Lambda^1\\
d\Lambda^3_{27}\cong\Lambda^3_{27}
\end{array}
$$
on the $\pm1$-eigenspaces of $\star$.
\end{enumerate}
\end{prp}

The proof of the proposition is a routine application of Schur's
lemma. For instance, $\tau=x_{\beta}\wedge x_{\alpha+\beta}$ is
contained in $\Lambda^{2\C}_{10+}=(0,3)$. It is straightforward to
check that $\star(\rho\wedge\tau)=\sqrt{3}i\cdot\tau^*\rho$ and thus holds for any form in $\Lambda^{2\C}_{10+}$. Note that the spaces
$\Lambda^{2\C}_{10\pm}\cong\Lambda^{3\C}_{10\pm}$ are interchanged
if we swap orientations. The forms $\sqrt{3}e_{123}+e_{458}+e_{678}$
and $-\sqrt{3}(e_{145}+e_{167})/2+e_{238}+e_{478}/2-e_{568}/2$ lie
in $\Lambda^3_8$ and $\Lambda^3_{20}$ and wedged with $\rho$ these
get mapped to non-zero elements in $\Lambda^6$ etc..

We can approach this decomposition also from a cohomological point of view well suited for our later purposes. The Lie algebra structure on $\Lambda^1=\mf{su}(3)$ induces a $PSU(3)$-invariant operator $b_k:\Lambda^k\to\Lambda^{k+1}$ by extension of
$$
be_i=\sum_{j<k}c_{ijk}e_j\wedge e_k.
$$
Since $b$ is built out of the structure constants, it is just the exterior
differential operator restricted to the left-invariant differential forms of $SU(3)$ with adjoint $b^*=d^*=-\star\,d\,\star$. The resulting elliptic complex is isomorphic to the de Rham cohomology $H^*(SU(3),\R)$ which is trivial except for the Betti numbers $b_0=b_3=1=b_5=b_8$. Hence, $\im b_k=\ker b_{k+1}$ for $k=0,\,1,\,3,\,5,\,6$ and $\im b_k=\ker b_{k+1}\oplus \R$ for $k=-1,\,2,\,4,\,7$. Schematically, we have
\begin{equation}\label{b_decomp}
\begin{array}{lclclclclclclclcl}
\Lambda^0_1 & \phantom{\stackrel{b}{\longrightarrow}} & & & & & \Lambda^3_1 & & & & \Lambda^5_1 & & & & & \phantom{\stackrel{b}{\longrightarrow}} & \Lambda^8_1\\[2pt]
& & \Lambda^1_8 & \stackrel{b}{\longrightarrow} & \Lambda^2_8 & & \Lambda^3_8 & \stackrel{b}{\longrightarrow} & \Lambda^4_8 & & & & \Lambda^6_8 & \stackrel{b}{\longrightarrow} & \Lambda^7_8 & & \\[2pt]
& & & & & &  & & \Lambda^4_8 &  \stackrel{b}{\longrightarrow} & \Lambda^5_8 & & & & & & \\ 
& & & & \Lambda^2_{20} & \stackrel{b}{\longrightarrow} & \Lambda^3_{20} & & & & \Lambda^5_{20} & \stackrel{b}{\longrightarrow} & \Lambda^6_{20} & & & & \\[2pt]
& & & & & & \Lambda^3_{27} & \stackrel{b}{\longrightarrow} & \Lambda^4_{27} & & & & & & & & \\[2pt]
& & & & & & & & \Lambda^4_{27} & \stackrel{b}{\longrightarrow} & \Lambda^5_{27} & & & & &
\end{array}
\end{equation}
In particular, we will use the more natural splitting of $\Lambda^4$ into $\Lambda^4_o=\ker b_3$ and $\Lambda^4_i=\im b^*$ instead of the $SO(8)$-equivariant splitting into self- and anti-self-dual forms. 

Using the $b$-operator and its co-differential, we can easily construct the projection operators $\pi^p_q:\Lambda^p\to\Lambda^p_q$. In particular, we find for $p=2$:

\begin{prp}\label{liealgsu3} For any $\alpha\in\Lambda^2$ we have $b(\alpha)=-\alpha^*\rho$. Moreover, $\Lambda^2_8=\ker b_2$ and the projection operator on the complement is $\pi^2_{20}(\alpha)=\frac{4}{3}b^*_3b_2(\alpha)$. For the complexified modules $\Lambda^2_{10\pm}$, the projection operators are $\pi^2_{10\pm}(\alpha)= \frac{2}{3}b^*_3b_2(\alpha)\mp\frac{8\sqrt{3}}{9}i\star(b_2(\alpha)\wedge\rho)$.
\end{prp}

As in the case of $G_2$-structures, there is also a half-spin
$PSU(3)$-invariant.

First, we can decompose
$$
\Delta_+\otimes\Delta_-=\Lambda^1\oplus\Lambda^3
$$
into $Spin(8)$-spaces and observe that the 3-form stabilised by
$PSU(3)$ and divided by its norm induces an isometry
$\Delta_+\to\Delta_-$. The triality principle permutes any two of
the representations $\Lambda^1$, $\Delta_+$ and $\Delta_-$ and since
these spaces are isomorphic as Euclidean vector spaces, we conclude
that we get two corresponding isometries in
$\Lambda^3\Delta_{\pm}=\ker\mu_{\pm}$ which consequently define
supersymmetric maps of spin $3/2$.

Alternatively, we can establish the existence of such an isomorphism
by a weight argument as in Proposition~\ref{g2ds} or~\cite{hi01}.
We restrict the $Spin(8)$-representations $\Lambda^1$, $\Delta_+$ and $\Delta_-$ to the embedding of $PSU(3)$ which acts on $\Lambda^1$ through the adjoint
representation of $SU(3)$ with weights
$$
0,\:\pm(2\sigma_1-\sigma_2),\:\pm(-\sigma_1+2\sigma_2),\:\pm(\sigma_1+\sigma_2).
$$
Hence, restricting the weights~(\ref{weights0}) to $PSU(3)$ we yield
$$
\omega_1=0,\:\omega_2=2\sigma_1-\sigma_2,\:\omega_3=\sigma_1+\sigma_2,\:\omega_4=0.
$$
Substituting this into~(\ref{weights+}) and~(\ref{weights-}) we see
that $PSU(3)$ acts with equal weights on $\Lambda^1$, $\Delta_+$,
$\Delta_-$. We showed the

\begin{prp}\label{equivalence}$\!\!\!${\rm~\cite{hi01}}\hspace{2pt} Restricted to $PSU(3)$, these three
representation spaces are equivalent
$$
\Lambda^1=\Delta_+=\Delta_-=\mf{su}(3).
$$
\end{prp}

\begin{rmk}
In particular Clifford multiplication
$\mu:\Lambda^1\otimes\Delta_{\pm}\to\Delta_{\mp}$ induces an
orthogonal product
$$
\times:\Lambda^1\otimes \Delta_+\cong \Lambda^1\otimes\Lambda^1\to
\Delta_-\cong\Lambda^1
$$
in $\mf{su}(3)$. If we think of $\mf{su}(3)$ as the algebra of skew
hermitian $3\times 3$ matrices, this product is explicitly given by
\begin{equation}\label{orthprod}
A\times B=\omega AB-\overline{\omega}BA-\frac{i}{\sqrt{3}}tr(AB)I,
\end{equation}
where $\omega=(1+i\sqrt{3})/2$~\cite{hi01}.
\end{rmk}

We conclude that vectors and spinors are essentially the same objects when regarded as being acted on by $PSU(3)$. If we denote the isomorphisms by $\gamma_{\pm}:\Lambda^1\to\Delta_{\pm}$, then these are characterised (up to a scalar) by the equations
$$
x\llcorner\rho(\gamma_{\pm})\,=\,\frac{1}{2}\kappa(x\llcorner\rho)\cdot\gamma_{\pm}-\gamma_{\pm}\circ x\llcorner\rho\,=\,0
$$
valid for all $x\in\Lambda^1$, i.e. $\mf{su}(3)$ acts trivially on $\gamma_{\pm}$. Their matrices with respect to a $PSU(3)$-frame and a fixed orthonormal basis of $\Delta_{\pm}$ are given by
$$
\gamma_+=\left(\begin{array}{rrrrrrrr}
\scriptstyle-\frac{1}{2} &\scriptstyle 0 &\scriptstyle \frac{1}{2} &\scriptstyle \frac{\sqrt{3}}{4} &\scriptstyle -\frac{1}{4} &\scriptstyle -\frac{\sqrt{3}}{4} &\scriptstyle \frac{1}{4} &\scriptstyle 0\\[3pt]

\scriptstyle0 &\scriptstyle -\frac{1}{2} &\scriptstyle 0 &\scriptstyle \frac{1}{4} &\scriptstyle \frac{\sqrt{3}}{4} &\scriptstyle \frac{1}{4} &\scriptstyle \frac{\sqrt{3}}{4} &\scriptstyle \frac{1}{2}\\[3pt]

\scriptstyle-\frac{1}{2} &\scriptstyle 0 &\scriptstyle -\frac{1}{2} &\scriptstyle -\frac{\sqrt{3}}{4} &\scriptstyle \frac{1}{4} &\scriptstyle -\frac{\sqrt{3}}{4} &\scriptstyle \frac{1}{4} &\scriptstyle 0\\[3pt]

\scriptstyle0 &\scriptstyle \frac{1}{2} &\scriptstyle 0 &\scriptstyle \frac{1}{4} &\scriptstyle \frac{\sqrt{3}}{4} &\scriptstyle -\frac{1}{4} &\scriptstyle -\frac{\sqrt{3}}{4} &\scriptstyle \frac{1}{2}\\[3pt]

\scriptstyle-\frac{1}{2} &\scriptstyle 0 &\scriptstyle \frac{1}{2} &\scriptstyle -\frac{\sqrt{3}}{4} &\scriptstyle \frac{1}{4} &\scriptstyle \frac{\sqrt{3}}{4} &\scriptstyle -\frac{1}{4} &\scriptstyle 0\\[3pt]

\scriptstyle0 &\scriptstyle \frac{1}{2} &\scriptstyle 0 &\scriptstyle \frac{1}{4} &\scriptstyle \frac{\sqrt{3}}{4} &\scriptstyle \frac{1}{4} &\scriptstyle \frac{\sqrt{3}}{4} &\scriptstyle -\frac{1}{2}\\[3pt]

\scriptstyle\frac{1}{2} &\scriptstyle 0 &\scriptstyle \frac{1}{2} &\scriptstyle -\frac{\sqrt{3}}{4} &\scriptstyle \frac{1}{4} &\scriptstyle -\frac{\sqrt{3}}{4} &\scriptstyle \frac{1}{4} &\scriptstyle 0\\[3pt]

\scriptstyle0 &\scriptstyle \frac{1}{2} &\scriptstyle 0 &\scriptstyle-\frac{1}{4} &\scriptstyle -\frac{\sqrt{3}}{4} &\scriptstyle \frac{1}{4} &\scriptstyle \frac{\sqrt{3}}{4} &\scriptstyle \frac{1}{2} 
\end{array}\right)
$$
and
$$
\gamma_-=\left(\begin{array}{rrrrrrrr}
\scriptstyle0 &\scriptstyle -\frac{1}{2} &\scriptstyle 0 &\scriptstyle \frac{1}{4} &\scriptstyle -\frac{\sqrt{3}}{4} &\scriptstyle -\frac{1}{4} &\scriptstyle \frac{\sqrt{3}}{4} &\scriptstyle \frac{1}{2}\\[3pt]

\scriptstyle-\frac{1}{2} &\scriptstyle 0 &\scriptstyle -\frac{1}{2} &\scriptstyle -\frac{\sqrt{3}}{4} &\scriptstyle -\frac{1}{4} &\scriptstyle -\frac{\sqrt{3}}{4} &\scriptstyle -\frac{1}{4} &\scriptstyle 0\\[3pt]

\scriptstyle0 &\scriptstyle \frac{1}{2} &\scriptstyle 0 &\scriptstyle \frac{1}{4} &\scriptstyle -\frac{\sqrt{3}}{4} &\scriptstyle \frac{1}{4} &\scriptstyle -\frac{\sqrt{3}}{4} &\scriptstyle \frac{1}{2}\\[3pt]

\scriptstyle-\frac{1}{2} &\scriptstyle 0 &\scriptstyle \frac{1}{2} &\scriptstyle \frac{\sqrt{3}}{4} &\scriptstyle \frac{1}{4} &\scriptstyle -\frac{\sqrt{3}}{4} &\scriptstyle -\frac{1}{4} &\scriptstyle 0\\[3pt]

\scriptstyle0 &\scriptstyle -\frac{1}{2} &\scriptstyle 0 &\scriptstyle -\frac{1}{4} &\scriptstyle \frac{\sqrt{3}}{4} &\scriptstyle \frac{1}{4} &\scriptstyle -\frac{\sqrt{3}}{4} &\scriptstyle \frac{1}{2}\\[3pt]

\scriptstyle\frac{1}{2} &\scriptstyle 0 &\scriptstyle \frac{1}{2} &\scriptstyle -\frac{\sqrt{3}}{4} &\scriptstyle -\frac{1}{4} &\scriptstyle -\frac{\sqrt{3}}{4} &\scriptstyle -\frac{1}{4} &\scriptstyle 0\\[3pt]

\scriptstyle0 &\scriptstyle -\frac{1}{2} &\scriptstyle 0 &\scriptstyle \frac{1}{4} &\scriptstyle -\frac{\sqrt{3}}{4} &\scriptstyle \frac{1}{4} &\scriptstyle -\frac{\sqrt{3}}{4} &\scriptstyle -\frac{1}{2}\\[3pt] 

\scriptstyle-\frac{1}{2} &\scriptstyle 0 &\scriptstyle \frac{1}{2} &\scriptstyle -\frac{\sqrt{3}}{4} &\scriptstyle -\frac{1}{4} &\scriptstyle \frac{\sqrt{3}}{4} &\scriptstyle \frac{1}{4} &\scriptstyle 0
\end{array}\right).
$$
To what extent we can characterise $PSU(3)$-structures through the existence of a
supersymmetric map of spin $3/2$ will be the topic of the following section.

We close our discussion of the group $PSU(3)$ with a remark about
special $PSU(3)$-orbits in the Grassmannians
$\tilde{G}_3(\Lambda^1)$ and $\tilde{G}_5(\Lambda^1)$, the spaces of
oriented 3- and 5-dimensional subspaces of $\Lambda^1$. This links
into another interesting feature of exceptional geometry, namely the
existence of special submanifolds.

In general, let $(V,g)$ be an oriented vector space with a metric
$g$ and $\tau\in\Lambda^p$ a $p$-form. We say that $\tau$ is a {\em
calibration}~\cite{hala82} if for every oriented $k$-plane
$\xi=e_1\wedge\ldots\wedge e_k$ in $\Lambda^1$ the inequality
$$
\tau(e_1,\ldots,e_k)\le 1
$$
holds. In other words any oriented $k$-plane $i:E\hookrightarrow
\Lambda^1$ satisfies $i^*_E\tau\le vol_E$, where $vol_E$ is the
volume element in $\Lambda^pE^*$ induced by the metric $g$
restricted to $E$. $E$ is {\em calibrated} by $\tau$ if equality
actually holds. The {\em contact set} of a calibration is the set of
calibrated planes.

A classical example are the so-called associative and co-associative
planes in $\R^7$ which are calibrated by the $G_2$-invariant form
$\varphi$ and $\star\varphi$. The contact set is parametrised by
$G_2/SO(4)$.

Let $\rho$ be a stable form associated with $PSU(3)$.

\begin{prp}
We have $\rho(\xi)\le 1$ with equality if and only if
$\xi=Ad(g)\mf{h}$ for $g\in SU(3)$ where $\mf{h}$ is an
$\mf{su}(2)$-subalgebra coming from a highest root. Furthermore,
$\star\rho(\xi)\le 1$ with equality if and only if $\xi$ is
perpendicular to a 3-plane calibrated by $\rho$. In particular,
$PSU(3)$ acts transitively on either contact set.
\end{prp}

\begin{prf}
We adapt the proof from~\cite{ta85}. Let $e_1,\ldots,e_8$ be a
$PSU(3)$-frame inducing the Euclidean norm $\norm{\cdot}$ and fix
the Cartan subalgebra $\mf{t}$ spanned by $e_3$ and $e_8$. Define
$$
\begin{array}{lll}
\alpha_1=\alpha=\frac{1}{2}e^3+\frac{\sqrt{3}}{2}e^8, & E_{\alpha_1}=e_5, & F_{\alpha_1}=-e_4\\
\alpha_2=\beta=\frac{1}{2}e^3-\frac{\sqrt{3}}{2}e^8, & E_{\alpha_2}=-e_6, & F_{\alpha_2}=e_7\\
\alpha_3=\alpha+\beta=e^3, & E_{\alpha_3}=e_1, & F_{\alpha_3}=e_2.
\end{array}
$$
Then $\norm{\alpha_i}\le 1$ and we immediately verify the relations
\begin{equation}\label{su2rel}
[H,E_{\alpha_i}]=\alpha_i(H)F_{\alpha_i},\quad
[H,F_{\alpha_i}]=-\alpha_i(H)E_{\alpha_i}\quad\mbox{and}\quad
[E_{\alpha_i},F_{\alpha_i}]=\alpha_i,
\end{equation}
for $H\in\mf{t}$ and $i=1,2,3$. Let $\xi\in\tilde{G}(3,\mf{su}(3))$.
Since $\mf{t}$ is a Cartan subalgebra,
$Ad(SU(3))A\cap\mf{t}\not=\emptyset$ for any $0\not=A\in\mf{su}(3)$.
Moreover, $\rho$ is $Ad$-invariant, so we may assume that
$\xi\cap\mf{t}\not=\emptyset$ up to the action of a $g\in SU(3)$.
Pick $T\in\xi\cap\mf{t}$ and extend it to a positively oriented
basis $\{T,X,Y\}$ of $\xi$. Then
$$
X=T_0+\sum\limits_{i=1}^3s_iE_{\alpha_i}+\sum\limits_{i=3}^3t_iF_{\alpha_i},
$$
where $T_0\in\mf{t}$. Hence,
$$
\norm{X}^2=\norm{T_0}^2+\sum\limits_{i=1}^3s^2_i+t^2_i
$$
which implies
\begin{eqnarray*}
\norm{[T,X]}^2 & = & \norm{\sum\limits_{i=1}^3s_i\alpha_i(T)F_{\alpha_i}-t_i\alpha_i(T)E_{\alpha_i}}^2\\
& \le & \sum\limits_{i=1}^3\norm{\alpha_i}^2\norm{T}^2(s_i^2+t^2_i)\\
& \le & \norm{T}^2\norm{X}^2,
\end{eqnarray*}
and finally
$$
\norm{[T,X]}\le\norm{T}\norm{X}.
$$
The result follows now by applying the Cauchy-Schwarz inequality
$$
|\rho(T,X,Y)|=\norm{g([T,X],Y)}\le\norm{[T,X]}\norm{Y}\le\norm{T}\norm{X}\norm{Y}=1.
$$
Furthermore, equality holds if and only if (i) $T\in\R\alpha_i$ and
$X\in\R E_{\alpha_i}\oplus\R F_{\alpha_i}$ for an $i\in\{1,2,3\}$,
(ii) $Y$ is a multiple of $[T,X]$ and (iii) $\rho(T,X,Y)>0$. Then
$Y\in\R E_{\alpha_i}\oplus\R F_{\alpha_i}$ and because of
(\ref{su2rel}), $\xi$ is an $\mf{su}(2)$-algebra.

Since $(\star\!\,\rho)_{|\xi^{\perp}}=\star\!\,(\rho_{|\xi})$ any
calibrated 5-plane is the orthogonal complement of an $\mf{su}(2)$
algebra. Moreover, any two subalgebras of highest root are conjugate.
\end{prf}

\begin{rmk}\hfill\newline
By a similar argument, Tasaki showed in~\cite{ta85} that for any
simple Lie group $G$ with Killing form $B$, the 3-form
$$
\rho(X,Y,Z)=\frac{1}{\norm{\delta}}B(X,[Y,Z])
$$
(where $\delta$ is the highest root of the Lie algebra) defines a
calibration on $G$ and that any calibrated submanifold is a
translate of $SU(2)$.
\end{rmk}

\subsection{Supersymmetric maps in dimension 8}\label{susy}

In dimension 8 triality implies
$\ker\mu_{\pm}=\Lambda^3\Delta_{\pm}$ so that
\begin{eqnarray*}
\Lambda^1\otimes\Delta_+ & = & \Delta_-\oplus \Lambda^3\Delta_-\\
\Lambda^1\otimes\Delta_- & = & \Delta_+\oplus \Lambda^3\Delta_+.
\end{eqnarray*}
We continue in the vein of Section~\ref{class_susy} and ask what
kind of geometry is induced by a supersymmetric map
$\gamma_{\pm}:\Lambda^1\to\Delta_{\pm}$ of spin $3/2$, that is
$\gamma_{\pm}\in\Lambda^3\Delta_{\mp}$. Since this takes only the
metric structure into account, by permuting the representations
$\Lambda^1$, $\Delta_+$ and $\Delta_-$ with the triality
automorphisms $\lambda$ and $\kappa$ (see Page~\pageref{3al1}), we
can rephrase our question as follows. What are the stabilisers of a
3-form that induces an isometry $\Delta_+\to\Delta_-$?

\begin{rmk}
Note that there is clearly a choice in this identification of
$\Delta_+\otimes\Delta_-$ with $\Lambda^1\oplus\Lambda^3$ since the
Hodge $\star$-operator induces an isomorphism between $\Lambda^p$
and $\Lambda^{8-p}$.
\end{rmk}

In what follows, it will be convenient to consider the full module
of spinors $\Delta=\Delta_+\oplus\Delta_-$ endowed with its spin
invariant metric $q=q_+\oplus q_-$, and to think of an odd form as a
linear map $\Delta\to\Delta$ whose matrix with respect to this
splitting is given by
\begin{equation}\label{rho_matrix}
\left(\begin{array}{cc} 0 & A_{\rho}\\ B_{\rho} &
0\end{array}\right).
\end{equation}
In the case of $\rho$ being a 3-form, the matrix has to be symmetric
as Clifford multiplication is skew for~$q$, hence
$B_{\rho}=A_{\rho}^{tr}$. Let $\mf{I}_g$ denote the set of 3-forms
in $\Delta\otimes\Delta$ that induce an isometry. It is
characterised by the next theorem.

\begin{thm}\label{orbit}
If a 3-form $\rho$ lies in $\mf{I}_g$, then $\rho$ is of unit length
and there exists a Lie bracket $[\cdot\,,\cdot]$ on $\Lambda^1$ such
that
\begin{equation}\label{bracket}
\rho(x,y,z)=g([x,y],z).
\end{equation}
In particular, the adjoint group of the resulting Lie algebra acts
as a group of isometries on $\Lambda^1$.

Conversely, if there exists a Lie algebra structure on $\Lambda^1$
whose adjoint group leaves $g$ invariant, then the 3-form defined by
(\ref{bracket}) and divided by its norm lies in $\mf{I}_g$.
\end{thm}

\begin{prf}
Inducing an isometry and defining a Lie bracket through
(\ref{bracket}) are both quadratic conditions on the coefficients of
$\rho$ which we show to coincide.

To this effect we define the linear map
$$
\mathop{{\rm Jac}}:\Lambda^3\otimes\Lambda^3\to\Lambda^4
$$
by skew-symmetrising the contraction to $\Lambda^2\otimes\Lambda^2$.
This is clearly $SO(8)$-equivariant and written in indices with
respect to some orthonormal basis $\{e_i\}$, we have
\begin{eqnarray*}
\mathop{\rm Jac}(c_{ijk}d_{lmk}) & = &
c_{[ij}^{\phantom{[ij}k}d_{lm]k}\\
& = &
c_{ij}^{\phantom{ij}k}d_{lmk}+c_{il}^{\phantom{il}k}d_{mjk}+c_{im}^{\phantom{im}k}d_{jlk}+c_{jl}^{\phantom{jl}k}d_{imk}+c_{jm}^{\phantom{jm}k}d_{lik}+c_{lm}^{\phantom{lm}k}d_{ijk}
\end{eqnarray*}
and in particular
\begin{equation}\label{jacid}
\mathop{\rm
Jac}(c_{ijk}c_{lmk})=2(c_{ij}^{\phantom{ij}k}c_{klm}+c_{li}^{\phantom{li}^k}c_{kjm}+c_{jl}^{\phantom{jl}k}c_{kim}).
\end{equation}
If we consider the skew-symmetric map
$[\,\cdot\,,\cdot\,]:\Lambda^2\to\Lambda^1$ defined through
$\rho=\sum_{i<j<k}c_{ijk}e_{ijk}$ by (\ref{bracket}), then
$$
[e_i,e_j]=c_{ij}^{\phantom{ij}k}e_k
$$
and~(\ref{jacid}) entails that this defines a Lie bracket if and
only if
$$
\mathop{\rm Jac}(\rho\otimes\rho)=0.
$$

Next we analyse the conditions for $\rho$ to induce an isometry.
Represent the map $\rho:\Delta\to\Delta$ by a mtarix $A_{\rho}$ as in~(\ref{rho_matrix}). Then $\rho$ defines an isometry if and only if
\begin{eqnarray}
q(\psi,\phi) & = & q(\rho\cdot \psi,\rho\cdot
\phi)\nonumber\\
& = & q(\rho\cdot\rho\cdot
\psi,\phi)\nonumber\\
& = &
q_+(A_{\rho}A_{\rho}^{tr}\psi_+,\phi_+)+q_-(A_{\rho}^{tr}A_{\rho}\psi_-,\phi_-)\nonumber\\
& = & q_+(\psi_+,\phi_+)+q_-(\psi_-,\phi_-).\label{isomeq}
\end{eqnarray}
This motivates considering the $Spin(8)$-equivariant linear maps
$$
\Gamma_{\pm}:\rho\otimes\tau\in\Lambda^3\otimes\Lambda^3\mapsto
pr_{\Delta_{\pm}}(\rho\cdot\tau)\in\Delta_{\pm}\otimes\Delta_{\pm}.
$$
The condition~(\ref{isomeq}) then reads
$$
\rho\in\mf{I}_g\:\mbox{ if and only if
}\:\Gamma_{\pm}(\rho\otimes\rho)={\rm Id}_{\Delta_{\pm}}.
$$
If we decompose both the domain and the target space into
irreducible components using the algorithm in~\cite{sa89}, we find
\begin{eqnarray*}
\Lambda^3\otimes\Lambda^3 & \cong &
\phantom{\oplus}\mathbf{1}\oplus2\Lambda^2\oplus\Lambda^4_+\oplus\Lambda^4_-\oplus
[0,0,2,2]\oplus[0,1,0,2]\oplus[0,1,2,0]\oplus\\
& & \oplus[0,2,0,0]\oplus[1,0,1,1]\oplus[2,0,0,0]
\end{eqnarray*}
and
$$
\Delta_+\otimes\Delta_+=\Lambda^2\Delta_+\oplus\odot^2\Delta_+\cong\Lambda^2\oplus\mathbf{1}\oplus\Lambda^4_+.
$$
Note that we labeled some irreducible representations by their
highest weight expressed in the basis of fundamental weights of
$\mf{so}(8)$ as we did earlier in this section. The modules
$\Lambda^4_+=[0,0,2,0]$ and $\Lambda^4_-=[0,0,0,2]$ are the spaces
of self-dual and anti-self-dual 4-forms respectively. Note that
$$
\Gamma_+(\rho\otimes\tau)^{tr}=\Gamma_-(\tau\otimes\rho)
$$
and consequently, it will suffice to consider the map $\Gamma_+$
only. Since the map induced by $\rho$ is symmetric,
$\Gamma_+(\rho\otimes\rho)\in\odot^2\Delta_+=\mathbf{1}\oplus\wedge^4_+$
follows. Moreover $\Gamma_+$ maps $\Lambda^4_+$ non-trivially into
$\Delta_+\otimes\Delta_+$. To see this, complexify the
representations and consider the root vectors $x_{\alpha}$,
$x_{\beta}$, $x_{\gamma}$ and $x_{\delta}$ of $\mf{so}(8)\otimes\C$.
Then
$$
x_{\alpha}\wedge x_{\beta}\wedge x_{\gamma}\otimes x_{\gamma}\wedge
x_{\beta}\wedge x_{\delta}\in \Lambda^4_+=
V(\alpha+2\beta+2\gamma+\delta)\leqslant\Lambda^3\otimes\Lambda^3,
$$
where we labeled $\Lambda^4$ by its highest weight expressed this
time in terms of the fundamental roots $\alpha$, $\beta$, $\gamma$,
and $\delta$. This vector gets mapped to
$$
\Gamma_+(x_{\alpha}\wedge x_{\beta}\wedge x_{\gamma}\otimes
x_{\gamma}\wedge x_{\beta}\wedge
x_{\delta})=\norm{x_{\beta}}^2\norm{x_{\gamma}}^2 x_{\alpha}\cdot
x_{\delta}.
$$
Hence a necessary condition for $\rho$ to lie in $\mf{I}_g$ is
$$
\Gamma_+(\rho\otimes\rho)_{\bigodot^2_0\Delta_+}=0.
$$
This obstruction in $\Lambda^4_+$ can be identified with the
obstruction ${\rm Jac}(\rho\otimes\rho)$ for $\rho$ to induce an
adapted Lie algebra as we shall show next. In fact, we claim that
\begin{equation}\label{isom++}
\Gamma_+(\rho\otimes\rho)=-\frac{1}{2}\mathop{\rm
Jac}(\rho\otimes\rho)_{\Lambda^4_+}+\norm{\rho}^2\mathop{\rm Id}.
\end{equation}
This yields directly the assertion of the theorem.

In order to derive~(\ref{isom++}), we first note that Clifford
multiplication induces a map
$$
\rho\otimes\tau\in\Lambda^3\otimes\Lambda^3\mapsto\rho\cdot\tau\in\Lambda^0\oplus\Lambda^2\oplus\Lambda^4\oplus\Lambda^6
$$
if we regard the product $\rho\cdot\tau$ as an element of
$\cliff(\Lambda^1,g)=\Lambda^*$ rather than as a map
$\Delta_+\otimes\Delta_+\oplus\Delta_-\otimes\Delta_-$. The various
components of $\rho\cdot\tau$ under this identification come from
the number of ``coinciding pairs" in the expression $e_{ijklmn}$.
This means that if we have, say, three coinciding pairs, then $i=l$, $j=m$
and $k=n$, hence $e_{ijklmn}=1$. The 2-form component is obtained
from two such pairs etc. (e.g. $e_{ijklmn}=e_{im}$ if, say, $j=l$
and $k=n$). Since the map induced by $\rho$ is symmetric,
$\Gamma_+(\rho\otimes\rho)$ takes values in the components
$\bigodot^2\Delta_+=\Lambda^0\oplus\Lambda^4_+$ only. Hence
$$
\rho=\sum\limits_{i<j<k}c_{ijk}e_{ijk}
$$
gets mapped to
\begin{eqnarray*}
\rho\cdot\rho & = &\sum\limits_{\begin{array}{c}\scriptstyle
i<j<k,\, l<m<n\end{array}}c_{ijk}c_{lmn}e_{lmnijk}\\
& = & \sum\limits_{\begin{array}{c} \scriptstyle i<j<k,\,
l<m<n\\\mbox{\scriptsize 3 pairs
coincide}\end{array}}c_{ijk}c_{lmn}e_{lmnijk}
+\sum\limits_{\begin{array}{c} \scriptstyle
i<j<k,\,l<m<n\\\mbox{\scriptsize 1 pair
coincides}\end{array}}c_{ijk}c_{lmn}e_{lmnijk}.
\end{eqnarray*}
Now the first sum is just
$$
\sum\limits_{i<j<k}c^2_{ijk}\mathbf{1}=\norm{\rho}^2\mathbf{1}
$$
which leaves us with the contribution of the sum with one pair of
equal indices. No matter which indices of the two triples
$\{i<j<k\}$ and $\{l<m<n\}$ coincide, the skew-symmetry of the
$c_{ijk}$ and $e_{ijk}$ allows us to rearrange and rename the
indices in such a way that
\begin{eqnarray*}
\sum\limits_{\begin{array}{c} \scriptstyle i<j<k,\,
l<m<n\\\mbox{\scriptsize 1 pair
coincides}\end{array}}c_{ijk}c_{lmn}e_{lmnijk} & = &\phantom{-}
\sum\limits_a\mathop{\sum}\limits_{\begin{array}{c} \scriptstyle
j<k,\,m<n\\\mbox{\scriptsize{\it j,k,m,n} distinct}\end{array}}c_{ajk}c_{amn}e_{amnajk}\\
& = & -\sum\limits_a\sum\limits_{\begin{array}{c} \scriptstyle
j<k,\, m<n\\\mbox{\scriptsize{\it j,k,m,n}
distinct}\end{array}}c_{ajk}c_{amn}e_{mnjk}\\
& = & -\frac{1}{2}\mathop{\rm Jac}(\rho\otimes\rho).
\end{eqnarray*}
whence~(\ref{isom++}).

Finally note that the metric $g$ is necessarily $ad$-invariant
because of the skew-symmetry of $\rho$.
\end{prf}

The 3-forms in $\mf{I}_g$ thus encapsulate the data of a Lie algebra
structure whose adjoint action preserves the metric on $\Lambda^1$.
We also say that the Lie structure is {\em adapted} to the metric
$g$. We shall write $\mf{k}$ if we think of $\Lambda^1$ as a Lie
algebra. The classification of the resulting structures shall occupy
us next.

According to the theorem of Levi-Malcev, we can decompose
$\mf{k}=\mf{m}\oplus\mf{r}$ into a semi-simple sub-algebra $\mf{m}$
and the so-called {\em radical} $\mf{r}$, i.e. the maximal solvable
ideal of $\mf{k}$. Since a semi-simple Lie algebra is the direct sum
of simple ones (which are completely classified), we are left with
the solvable part in order to determine the Lie algebra structure of
$\mf{k}$ completely. By restricting the metric to $\mf{r}$, we see
that $\mf{r}$ is adapted to $g$.

\begin{prp}
Let $\mf{s}$ be a solvable Lie algebra which is adapted to some
metric $g$. Then $\mf{s}$ is abelian.
\end{prp}

\begin{prf}
We shall proceed by induction over $n$, the dimension of $\mf{s}$.
If $n=1$, then $\mf{s}$ is abelian and the assertion is trivial. Now
assume that the assumption holds for all $1\leq m<n$. Let $\mf{a}$
be an {\em abelian ideal} of $\mf{s}$. Then, for all $A\in\mf{a}$
and $X,\,Y\in\mf{s}$, $ad$-invariance of $g$ implies
$$
g(X,[Y,A])=g(Y,[A,X])=0,
$$
for if $Y\in\mf{a}$, then $[Y,A]=0$ and if $Y\in\mf{a}^{\perp}$,
then $[A,X]\in\mf{a}\perp Y$. Hence $\mf{a}\subset\mf{z}$, the
center of $\mf{s}$. If $\mf{z}$ were to be trivial, then any abelian
ideal of $\mf{s}$ would also be trivial which is equivalent to
saying that $\mf{s}$ is semi-simple, contradicting our assumption.
Hence we can write $\mf{s}=\mf{z}\oplus\mf{h}$, where $\mf{h}$ is an
orthogonal complement to $\mf{z}$. Now for all $X\in\mf{s}$,
$Z\in\mf{z}$ and $H\in\mf{h}$ we have
$$
g(Z,[X,H])=g(H,[Z,X])=0
$$
so that $[X,H]\in\mf{z}^{\perp}=\mf{h}$. Equivalently, $\mf{h}$ is
an ideal of $\mf{s}$. As such, it is adapted and solvable since
$\mf{s}$ is adapted and solvable. Since the dimension of $\mf{h}$ is
strictly less than $n$, our induction hypothesis applies and we
deduce that $\mf{s}$ is abelian.
\end{prf}

\begin{cor}
An adapted Lie algebra $\mf{k}$ is reductive.
\end{cor}

\begin{prf}
We need to show that the center $\mf{z}\subset\mf{r}$ actually
equals the radical of $\mf{k}$. But for any $X\in\mf{k}$ and
$R_1,\,R_2\in\mf{r}$, the previous proposition implies that
$$
g(R_1,[X,R_2])=g(X,[R_2,R_1])=0.
$$
Hence $[X,R_2]\in\mf{r}\cap\mf{r}^{\perp}=\{0\}$ for any
$X\in\mf{k}$ and consequently, $\mf{r}=\mf{z}$.
\end{prf}

As a result, any Lie algebra which is adapted to $g$ is of the form
$\mf{m}\oplus\mf{z}$ with $\mf{m}$ semi-simple. Cartan's
classification of simple Lie algebras then asserts that $\Lambda^1$
is isomorphic to one of the following Lie algebras.
\begin{enumerate}\label{lie_types}
\item $\mf{k}_1=\mf{su}(3)$
\item $\mf{k}_2=\mf{su}(2)\oplus\mf{su}(2)\oplus E_2$
\item $\mf{k}_3=\mf{su}(2)\oplus E_5$.
\end{enumerate}
Here, $E_5$ denotes the abelian Lie algebra $\R^5$ which we think of
as a Euclidean vector space. We refer to a 3-form in $\mf{I}_g$ to
be of {\em type}~1, 2 or 3 according to the Lie algebra it induces.
Hence there is a disjoint decomposition of $\mf{I}_g$ into the sets
$\mf{I}_{g1}$, $\mf{I}_{g2}$ and $\mf{I}_{g3}$ acted on by $SO(8)$
and pooling the forms of type~1, 2 or 3 respectively. Since the
automorphism groups $Aut(\mf{k}_i)$ of $\mf{k}_1$, $\mf{k}_2$ and
$\mf{k}_3$ are $PSU(3)\times\Z_2$, $SO(3)\times GL(5)$ and
$SO(3)\times SO(3)\times GL(2)$ (recall that
$Aut(\mf{su}(2))=SU(2)/Z_2=SO(3)$), the stabiliser of $\rho$
in~(\ref{bracket}) inside $SO(8)$ is $PSU(3)$, $SO(3)\times SO(5)$
and $SO(3)\times SO(3)\times SO(2)$ respectively. On the other hand,
there is a natural decomposition of $\mf{I}_g$ into the subsets
$\mf{I}_{g\pm}$ of forms whose induced map
$A_{\rho}:\Delta_-\to\Delta_+$ has positive or negative determinant.
The orbit structure of the $SO(8)$-action on $\mf{I}_g$ and
$\mf{I}_{g\pm}$ is given by the next proposition.

\begin{prp}\label{orbitdecomposition}
The sets $\mf{I}_{g1}$, $\mf{I}_{g2}$ and $\mf{I}_{g3}$ can be
described as follows.
\begin{enumerate}
\item $\mf{I}_{g1}=SO(8)/PSU(3)$
\item $\mf{I}_{g2}=(S^1-\{pt\})\times SO(8)/(SO(3)\times SO(3)\times SO(2))$
\item $\mf{I}_{g3}=SO(8)/(SO(3)\times SO(5))$
\end{enumerate}
Furthermore,
$$
\mf{I}_{g-}=SO(8)/PSU(3)
$$
and $\mf{I}_{g+}$ foliates over the circle $S^1$ with principal
orbits $SO(8)/(SO(3)\times SO(3)\times SO(2))$ over $S^1-\{pt\}$ and
a degenerate orbit $SO(8)/(SO(3)\times SO(5))$ at $\{pt\}$.
\end{prp}

\begin{prf}
First we note that since the 3-forms under consideration are of the
form $g(\cdot,[\cdot\,,\cdot])$, these forms are built out of the
structure constants of the corresponding Lie algebra. To begin with,
let us show that $\mf{I}_{g1}=SO(8)/PSU(3)$ and that this is
contained in $\mf{I}_{g-}$. We will choose a certain
$\rho_1\in\mf{I}_{g1}$ whose orbit under $SO(8)$ is the whole set
$\mf{I}_{g1}$. Fix an orthonormal basis $\{e_i\}$ and define
\begin{equation}\label{nf0}
\rho_1=\frac{1}{2}e_{123}+\frac{1}{4}e_{147}-\frac{1}{4}e_{156}+\frac{1}{4}e_{246}+\frac{1}{4}e_{257}+\frac{1}{4}e_{345}-\frac{1}{4}e_{367}+\frac{\sqrt{3}}{4}e_{458}+\frac{\sqrt{3}}{4}e_{678}.
\end{equation}
The form $\rho_1$ is of unit length and its coefficients are $1/2$
of the structure constants $c_{ijk}$ of $\mf{su}(3)$ as given in
Section~\ref{psu3}. Hence $\rho_1\in\mf{I}_{g1}$ and its stabiliser
is $PSU(3)$. To establish transitivity, take any
$\rho\in\mf{I}_{g1}$ and endow $\Lambda^1$ with the induced
$\mf{su}(3)$ structure. We have to find a transformation in $SO(8)$
which puts $\rho$ into the normal form of~(\ref{nf0}). Since $g$ is
$ad$-invariant and positive definite, there exist a strictly
negative real number $\lambda$ such that
$$
\lambda g=B_{\mf{su}(3)},
$$
where $B_{\mf{k}}(x,y)={\rm Tr}(ad_x\circ ad_y)$ designates the
Killing form of the Lie algebra $\mf{k}$. Let
$\tilde{f}_1,\ldots,\tilde{f}_8$ be a basis for $\mf{su}(3)$ such
that $[\tilde{f}_i,\tilde{f}_j]=c_{ijk}\tilde{f}_k$. Then
$$
\lambda
g(\tilde{f}_i,\tilde{f}_j)=B_{\mf{su}(3)}(\tilde{f}_i,\tilde{f}_j)=-3\delta_{ij}.
$$
Consequently, the basis $f_i=\sqrt{\frac{|\lambda |}{3}}\tilde{f}_i$
is orthonormal and has structure constants
$c_{ijk}/\sqrt{3|\lambda|}$. Hence we have shown the existence of an
orthonormal basis $f_1,\ldots,f_8$ of $g$ such that
$$
\rho(f_i,f_j,f_k)=g([f_i,f_j],f_k)=-\frac{1}{3}\sqrt{\frac{|\lambda
|}{3}}B_{\mf{su}(3)}([\tilde{f}_i,\tilde{f}_j],\tilde{f}_k)=\sqrt{\frac{|\lambda
|}{3}}c_{ijk}
$$
whence
$$
\rho=\sqrt{\frac{|\lambda |}{3}}\sum\limits_{i<j<k}c_{ijk}f_{ijk}.
$$
Since $\rho$ is of unit norm we have $\lambda =-3/4$. Mapping the
basis $e_i$ into $f_i$ yields the desired transformation in $SO(8)$.
Next we compute the determinant of $A_{\rho_1}$ by identifying
${e_i}$ with the octonions as in~(\ref{octonions}) and fixing the
representation~(\ref{spin8_rep}). This gives
\begin{equation}\label{arho1}
A_{\rho_1}=\frac{1}{4}\left (\begin {array}{cccccccc}
\sqrt{3}&0&0&3&-\sqrt {3}&0&0&1\\
2 &-\sqrt{3}&-1&0&2&-\sqrt{3}&-1&0\\
0&3&-\sqrt{3}&0&0&-1&-\sqrt {3}&0\\
-1&0&2&\sqrt{3}&1&0&-2 &-\sqrt {3}\\
-\sqrt{3}&0&0&1&\sqrt {3}&0&0&3\\
-2&-\sqrt{3}&-1&0&-2& -\sqrt {3}&-1&0\\
0&-1&-\sqrt{3}&0&0&3&-\sqrt {3}&0\\
1&0&2&-\sqrt{3}&-1 &0&-2&\sqrt {3}
\end {array}\right ),
\end{equation}
hence $\det(A_{\rho_0})=-1$. Moreover, we have $\det\pi_{\pm}(a)=1$
for any $a\in Spin(8)$. This follows from the fact that its
generators $e_i\cdot e_j$ square to $-Id$, therefore they are of
determinant 1. The $Spin(8)$-equivariance of the embedding
$$
\rho\in\Lambda^3\hookrightarrow\left(\begin{array}{cc} 0 &
A_{\rho}\\ A^{\rm tr}_{\rho} &
0\end{array}\right)\in\Delta\otimes\Delta
$$
entails
$$
A_{\pi_0(a)^*\rho_1}=\pi_+(a)\circ A_{\rho_1}\circ\pi_-(a)^{-1},
$$
whence $\mf{I}_1\subset\mf{I}_-$.

Next we turn to structures of type~2 and 3 which are characterised
by a non-trivial center. We shall see in a moment that these two
cases are intimately related as type~3 can be conceived as a limit
case of type~2. First we assume that $\Lambda^1$ carries a Lie
algebra structure isomorphic to
$\mf{su}(2)_1\oplus\mf{su}(2)_2\oplus E_2$ where for the sake of
clarity we labeled the two subalgebras isomorphic to $\mf{su}(2)$ by
the subscripts 1 and 2. Hence there exist strictly-negative
constants $\lambda_1$ and $\lambda_2$ such that
$$
B_{\mf{k}|\mf{su}(2)_1}=\lambda_1g,\quad
B_{\mf{k}|\mf{su}(2)_2}=\lambda_2g.
$$
Since $\mf{su}(2)_1$ and $\mf{su}(2)_2$ are ideals of $\mf{k}$, the
restriction of $B_{\mf{k}}$ to these is just the Killing form of
$B_{\mf{su}(2)}$. We claim that the set $\{\lambda_1,\lambda_2\}$ is
invariant under the action of $SO(8)$. More concretely, if we
transform $\rho$ by $A\in SO(8)$, then $A^*\rho$ gives rise to a new
Lie bracket $[\cdot\,,\cdot]_A$ which is just
$$
[v,w]_A=A[A^{-1}v,A^{-1}w].
$$
We write $A^*\mf{k}$ for the resulting Lie algebra which has the
same type as $\mf{k}$. We obtain again a set
$\{\lambda_1^A,\lambda_2^A\}$ of negative numbers which we claim to
be $\{\lambda_1,\lambda_2\}$. To show this we pick two non-zero
vectors $v_{1}$ and $v_{2}$ in $\mf{su}(2)_1$ and $\mf{su}(2)_2$
respectively so that
$$
\lambda_1=\frac{B_{\mf{su}(2)_1}(v_1,v_1)}{g(v_1,v_1)},\quad
\lambda_2=\frac{B_{\mf{su}(2)_2}(v_2,v_2)}{g(v_2,v_2)}.
$$
The Killing form of $A^*\mf{k}$ can be evaluated by the formula
$$
B_{A^*\mf{k}}(Av,Av)=B_{\mf{k}}(v,v),
$$
since for an orthonormal basis $\{e_i\}$ of $\Lambda^1$ we have
\begin{eqnarray*}
B_{A^*\mf{k}}(Av,Av) & = & \mathop{\rm Tr}({\rm
ad}_{Av}^{A^*\mf{k}}\circ{\rm ad}^{A^*\mf{k}}_{Av})\\
& = & \sum\limits_ig([Av,[Av,e_i]_A]_A,e_i)\\
& = & \sum\limits_ig(A[v,[v,A^{-1}e_i]],e_i)\\
& = & \mathop{\rm Tr}({\rm ad}^{\mf{k}}_v\circ{\rm ad}^{\mf{k}}_v)\\
& = & B_{\mf{k}}(v,v).
\end{eqnarray*}
Hence the set $\{\lambda_1^A,\lambda_2^A\}$ of the transformed Lie
algebra $A^*\mf{k}$ is just
$$
\lambda_1^A=\frac{B_{A^*\mf{su}(2)_1}(Av,Av)}{g(Av,Av)}=\frac{B_{\mf{su}(2)}(v,v)}{g(v,v)}=\lambda_1
$$
and
$$
\lambda_2^A=\frac{B_{A^*\mf{su}(2)_2}(Av,Av)}{g(Av,Av)}=\frac{B_{\mf{su}(2)^{\perp}}(v,v)}{g(v,v)}=\lambda_2
$$
which proves our claim.

For $\lambda_1+\lambda_2=-2$, an example of a form in $\mf{I}_{g2}$
is given by
$$
\rho_2=\sqrt{\frac{|\lambda_1|}{2}}e_{123}+\sqrt{\frac{|\lambda_2|}{2}}e_{456}.
$$
Now pick a $\rho\in\mf{I}_{g2}$. In each ideal $\mf{su}(2)_i$ we can
find a basis $\tilde{f}_{i_1},\tilde{f}_{i_2},\tilde{f}_{i_3}$ such
that
$$
[\tilde{f}_{i_l},\tilde{f}_{i_m}]=\epsilon_{lmn}\tilde{f}_{i_n}
$$
(where $\epsilon_{ijk}$ is totally skew in its indices with
$\epsilon_{123}=1$) and
$$
B_{\mf{su}(2)}(\tilde{f}_{i_l},\tilde{f}_{i_m})=-2\delta_{lm}.
$$
Hence for $i=1,\:2$ the basis
$f_{i_l}:=\sqrt{\frac{|\lambda_i|}{2}}\tilde{e}_{i_l}$ is
orthonormal for $g$ which we extend to a basis of $\Lambda^1$. The
only non-trivial coefficients of $\rho$ are
$$
\rho(f_{i_1},f_{i_2},f_{i_3})=g([f_{i_1},f_{i_2}],f_{i_3})=-\frac{1}{2}\sqrt{\frac{|\lambda_i|}{2}}B([\tilde{f}_{i_1},\tilde{f}_{i_2}],\tilde{f}_{i_3})=\sqrt{\frac{|\lambda_i|}{2}}.
$$
Any mixed term like $\rho(f_{1_1},f_{1_2},f_{2_1})$ etc. has to
vanish since the two sub-algebras $\mf{su}(2)$ are ideals in
$\mf{k}_1$. Moreover, $\rho$ is of norm 1 which implies
$\lambda_1+\lambda_2=-2$ so that $\lambda_2=-2-\lambda_1<0$ for $-2<
\lambda_1<0$. If $\rho$ induces the same pair
$\{\lambda_1,\lambda_2\}$, then
$$
\rho=\sqrt{\frac{|\lambda_1|}{2}}f_{1_11_21_3}+\sqrt{\frac{|\lambda_2|}{2}}f_{2_12_22_3}=A^*\rho_2
$$
for the transformation $A\in SO(8)$ which maps $e_i$ to $f_i$ (after
choosing a suitable labeling for the basis $\{f_i\}$).

The matrix of $\rho_2$ is given by
\begin{equation}\label{matrix2}
\frac{1}{\sqrt{2}}\left (\begin {array}{cccccccc}
0&0&\sqrt{2+\lambda_1}&\sqrt{-\lambda_1}&0&0&0&0\\
0&0&-\sqrt{-\lambda_1}&\sqrt{2+\lambda_1}&0&0&0&0\\
\sqrt{2+\lambda_1}&\sqrt{-\lambda_1}&0&0&0&0&0&0\\
-\sqrt{-\lambda_1}&\sqrt{2+\lambda_1}&0&0&0&0&0&0\\
0&0&0&0&0&0&\sqrt{2+\lambda_1}&\sqrt{-\lambda_1}\\
0&0&0&0&0&0&-\sqrt{-\lambda_1}&\sqrt{2+\lambda_1}\\
0&0&0&0&\sqrt{2+\lambda_1}&\sqrt{-\lambda_1}&0&0\\
0&0&0&0&-\sqrt{-\lambda_1}&\sqrt{2+\lambda_1}&0&0\\
\end{array}\right )
\end{equation}
and therefore $\det(A_{\rho_2})=1$, that is
$\mf{I}_{g2}\subset\mf{I}_{g+}$.

The limit case $\{-2,0\}$ corresponds precisely to structures of
type~3. Any form in $\mf{I}_{g3}$ can be written as $\rho_3=e_{123}$
for a suitably chosen orthonormal basis. Its associated matrix is
$$
A_{\rho_3}=\left(\begin {array}{cccccccc} 0&0&0&1&0&0&0&0\\0&0
&-1&0&0&0&0&0\\0&1&0&0&0&0&0&0\\-1
&0&0&0&0&0&0&0\\0&0&0&0&0&0&0&1\\0
&0&0&0&0&0&-1&0\\0&0&0&0&0&1&0&0\\0 &0&0&0&-1&0&0&0\end
{array}\right)
$$
whose determinant equals 1, hence $\mf{I}_{g3}\subset\mf{I}_{g+}$.

It follows that after rescaling, $\mf{I}_{g2}=(S^1-\{pt\})\times
SO(8)/(SO(3)\times SO(3)\times SO(2)),$ and
$\mf{I}_{g+}=\mf{I}_{g2}\cup\mf{I}_{g3}$ foliates over the circle
$S^1$ with principal orbits $SO(8)/(SO(3)\times SO(3)\times SO(2))$
over $S^1-\{pt\}$ and a degenerate orbit $SO(8)/(SO(3)\times SO(5))$
at $pt$.
\end{prf}

Since Theorem~\ref{orbit} only uses the Euclidean structure of
$\Lambda^1$ which by triality coincides with that of $\Delta_+$ and
$\Delta_-$, we derive the same conclusion for isometries in
$\Lambda^1\otimes\Delta_+$ and $\Lambda^1\otimes\Delta_-$. Note however that the  characterisation of $\mf{I}_{g\pm}$ does depend on the module under consideration as the outer triality morphisms reverse the orientation. By abuse of language, any element sitting inside $\Lambda^3$,
$\Lambda^3\Delta_+$ or $\Lambda^3\Delta_-$ will be referred to as a
map of spin $3/2$. If such an element defines an isometry, it
corresponds to the reduction of the structure group $SO(\Lambda^1)$,
$SO(\Delta_+)$ or $SO(\Delta_-)$ to the group $K_1=PSU(3)$,
$K_2=SO(3)\times SO(3)\times SO(5)$ and $K_3=SO(3)\times SO(5)$. We
obtain then a subgroup $\widetilde{K}_i$ of $Spin(8)$ through the
lift by $\pi_0$, $\pi_+$ and $\pi_-$. We will also refer to the
sub-types $i_0$, $i_+$ and $i_-$ to render the provenance of
$\widetilde{K}_i$ explicit. Unlike for 3-forms of type~2 and 3,
$\widetilde{K}_1$ is not connected. Table~\ref{coverK} displays all
three possibilities where we have to divide by $\Z_2$ in case of
$i=2$ or $3$. Note that $Sp(1)\times_{\Z_2} Sp(2)$ is usually
denoted by $Sp(1)Sp(2)$ or $Sp(1)\cdot Sp(2)$, giving rise to what
is called an {\em almost quaternionic} structure.

\begin{table}[hbt]
\begin{center}
\begin{tabular}{|c|c|c|c|}
\hline type $i$ & 1 & 2 & 3\\ \hline $K_i$ & $PSU(3)$ & $SO(3)\times
SO(3)\times SO(2)$ &
$SO(3)\times SO(5)$\\
[4pt] $\widetilde{K}_i$ & $PSU(3)\times\Z_2$ & $SU(2)\cdot
SU(2)\times U(1)$ & $Sp(1)\cdot Sp(2)$\\
\hline
\end{tabular}
\end{center}
\caption{The stabilisers and their covers in $Spin(8)$}
\label{coverK}
\end{table}

\subsection{The groups $SO(3)\times SO(3)\times SO(2)$ and $SO(3)\times SO(5)$}\label{relgeo}

A group $\widetilde{K}_i$ coming, say, from an $i_+$-structure might
intersect the kernel of $\pi_0$ and $\pi_-$ trivially and therefore
induce a geometry on $\Lambda^1$ and $\Delta_-$ (the two being
identified through the isometry $\Lambda^1\to\Delta_-$) which is
different from that of $\Delta_+$. In
the case of $PSU(3)$, we have already seen that this is not the case
as we get isometries which identify $\Lambda^1$, $\Delta_+$ and
$\Delta_-$ as $PSU(3)$-representations spaces.

Next, we want to carry out a similar analysis for structures of
type~2 and 3. Again, we may restrict ourselves to structures of
type~$2_0$ and $3_0$ -- the remainder is a matter of applying
triality.

To this effect, we recall that $\mf{su}(2)$ has one fundamental
weight $\sigma$ and that the irreducible representations are
$$
[n]=\bigodot^n\C^2
$$
of dimension $n+1$. In particular, we get the standard
representation $\C^2$ for $n=1$ and the adjoint representation
$\mf{su}(2)$ for $n=2$. They are real for $n$ even and quaternionic
for $n$ odd. Consequently, the irreducible representations of
$\widetilde{K}_2$ can be labeled by
$[n_1,n_2,m]=[n_1]\otimes[n_2]\otimes[m]$ where the first two
factors correspond to the $SU(2)$-representations $n_1\sigma_1$ and
$n_2\sigma_2$. The third factor is the irreducible
$S^1$-representation $S_m:\theta(z)\mapsto e^{im\theta}\cdot z$
which is one dimensional and complex. According to the discussion in
the previous section we obtain
$$
\Lambda^1=\mf{su}(2)\oplus\mf{su}(2)\oplus
E_2=[2,0,0]\oplus[0,2,0]\oplus\llbracket0,0,2\rrbracket.
$$
Hence, $\widetilde{K}_2$ acts with weights (cf.~(\ref{weights0}))
$$
0,\,2\sigma_l,\,2\sigma_r,\,2m.
$$
As an $SO(8)$-space, $\Lambda^1=[1,0,0,0]$ and the following
substitution of weights
$$
\omega_1=0,\;-\omega_1+\omega_2=2m,\;-\omega_2+\omega_3+\omega_4=2\sigma_l,\;\omega_3-\omega_4=2\sigma_r
$$
yields
\begin{equation}\label{g1weights}
\omega_1=0,\;\omega_2=2m,\;\omega_3=\sigma_l+\sigma_r+m,\;\omega_4=\sigma_l-\sigma_r+m
\end{equation}
for the restrictions of $\omega_i$ to $\tilde{G}_2$. Substituting
(\ref{g1weights}) into~(\ref{weights+}) and~(\ref{weights-}) gives
$$
\pm(-\sigma_l+\sigma_r-m),\;\pm(\sigma_l-\sigma_r-m),\;\pm(-\sigma_l-\sigma_r+m),\;\pm(\sigma_l+\sigma_r+m).
$$
We conclude
$$
\Delta_+=\Delta_-=\llbracket1,1,1\rrbracket.
$$
In particular, the action of $K_2$ on $\Delta_{\pm}$ preserves an
almost complex structure. Note however that this structure does not
reduce to $SU(4)$ as the torus component acts non-trivially on
$\lambda^{4,0}\Delta_{\pm}$.

Permuting with the triality automorphisms yields analogous results
for structures of type~$2_+$ and $2_-$. We summarise our results in
the Table~\ref{g2action}.

\begin{table}[hbt]
\begin{center}
\begin{tabular}{|c|c|c|c|}
\hline $\rho\in$ & $\Lambda^3$ & $\Lambda^3\Delta_+$ &
$\Lambda^3\Delta_-$\\
\hline type & 0 & -- & +\\
\hline $\Lambda^1$ &
$\mf{su}(2)\oplus\mf{su}(2)\oplus E_2$ &
$\llbracket\C^2\otimes\C^2\otimes S_1\rrbracket$ &
$\llbracket\C^2\otimes\C^2\otimes S_1\rrbracket$\\
[4pt] $\Delta_+$ & $\llbracket\C^2\otimes\C^2\otimes S_1\rrbracket$
& $\mf{su}(2)\oplus\mf{su}(2)\oplus E_2$ &
$\llbracket\C^2\otimes\C^2\otimes S_1\rrbracket$\\
[4pt] $\Delta_-$ & $\llbracket\C^2\otimes\C^2\otimes S_1\rrbracket$
& $\llbracket\C^2\otimes\C^2\otimes S_1\rrbracket$ &
$\mf{su}(2)\oplus\mf{su}(2)\oplus E_2$\\
\hline
\end{tabular}
\end{center}
\caption{The action of $\widetilde{K}_2=SU(2)\cdot SU(2)\times
U(1)$ on $\Lambda^1$, $\Delta_+$ and $\Delta_-$}
\label{g2action}
\end{table}

Next we turn to structures of type~3. Again we consider
$\rho\in\Lambda^3$ and assume that as a $\widetilde{K}_3=Sp(1)\cdot
Sp(2)$ representation space, we have
$$
\Lambda^1=\mf{su}(2)\oplus E_5=[2,0,0]\oplus[0,1,0].
$$
The first index $n_1$ refers to the $SU(2)=Sp(1)$-representation labeled by
$n\sigma$, while the last two indices $(m_1,m_2)$ designate the
irreducible $Sp(2)=Spin(5)$-representation with respect to the basis of
fundamental weights $\tau_1$ and $\tau_2$. The weights of the action
on $\Lambda^1$ are
$$
0,\:2\sigma,\:\tau_1,\:\tau_1-2\tau_2.
$$
Substituting
$$
\omega_1=0,\:\omega_1-\omega_2=2\sigma,\:\omega_2-\omega_3-\omega_4=\tau_1,\:\omega_3-\omega_4=\tau_1-2\tau_2
$$
we obtain as weights on $\Delta_+$ and $\Delta_-$
$$
\pm(-\sigma-\tau_1+\tau_2),\:\pm(\sigma-\tau_1+\tau_2),\:\pm(\sigma-\tau_2),\:\pm(\sigma+\tau_2)
$$
and thus
$$
\Delta_+=\Delta_-=[1,0,1]=[\C^2\otimes\K^2].
$$

Our conclusions are displayed in the Table~\ref{g3action}. Choosing the orientation suitably, we see in particular that an almost quaternionic structure on $\Lambda^1$ is equivalent to a Lie algebra structure of type 3 on $\Delta_+$.

\begin{table}[hbt]
\begin{center}
\begin{tabular}{|c|c|c|c|}
\hline $\rho\in$ & $\Lambda^3$ & $\Lambda^3\Delta_+$ &
$\Lambda^3\Delta_-$\\
\hline type & 0 & -- & +\\
\hline $\Lambda^1$ & $\mf{su}(2)\oplus E_5$ & $[\C^2\otimes\K^2]$ &
$[\C^2\otimes\K^2]$\\
[4pt] $\Delta_+$ & $[\C^2\otimes\K^2]$ & $\mf{su}(2)\oplus E_5$ &
$[\C^2\otimes\K^2]$\\
[4pt] $\Delta_-$ & $[\C^2\otimes\K^2]$ & $[\C^2\otimes\K^2]$ &
$\mf{su}(2)\oplus E_5$\\
\hline
\end{tabular}
\end{center}
\caption{The action of $\widetilde{K}_3=Sp(1)\cdot Sp(2)$ on
$\Lambda^1$, $\Delta_+$ and $\Delta_-$} \label{g3action}
\end{table}

\section{Generalised metric structures}\label{gen_met_struc}

So far we have discussed ``classical" supersymmetric structures.
This means that we are given one supersymmetric map which forces the
structure group to reduce. Equivalently, we could consider a {\em
homogeneous} form, i.e. a form of pure degree. Next we want to
describe geometries which are induced by the existence of two
supersymmetric maps. For this we use the ``generalised" setup as
introduced in~\cite{hi02} which describes geometrical structures by
even or odd forms. The idea is to work with the so-called {\em
Narain} or {\em generalised $T$-duality group} $O(n,n)$
(cf.~\cite{ha00} and references therein) which is naturally
associated with the vector space $T\oplus T^*$ and contains
subgroups of the form $G\times G$. For example, $G_2\times G_2$ and
$Spin(7)\times Spin(7)$ sit naturally inside in $O(7,7)$ and
$O(8,8)$. In this section, we investigate the geometry associated
with $O(n)\times O(n)$ giving rise to what we call a {\em
generalised metric structure}.

Recall that the Euclidean structures on a real vector space $T=\R^n$
are parametrised by the homogeneous space $GL(n)/O(n)$ where the
subgroup $O(n)$ determines a certain Euclidean inner product $g$ up
to a scale. In the generalised case we are then looking for a tensor
which yields a reduction from $O(n,n)$ to $O(n)\times O(n)$, so we
are interested in the space
$$
O(n,n)/O(n)\times O(n).
$$
To this effect let us analyse the invariants attached to $O(n)\times
O(n)$. First, picking a particular copy of $O(n)\times O(n)$ inside
$O(n,n)$ defines an orthogonal splitting
$$
(T\oplus T^*,(\cdot,\cdot))=(V_+\oplus V_-,g_+\oplus g_-)
$$
into a positive and a negative definite subspace $V_+$ and $V_-$. We
refer to such a splitting also as a {\em metric splitting}.
Conversely, any choice of such a metric splitting defines the
subgroup
$$
O(V_+,g_+)\times O(V_-,g_-)\cong O(n)\times O(n)
$$
inside $O(n,n)$.

\begin{definition}
A {\em generalised metric structure} is a reduction from $O(n,n)$ to
$O(n)\times O(n)$.
\end{definition}

\begin{figure}[ht]
\begin{center}
\input{fig1.pstex_t}
\end{center}
\caption{Metric splitting of $T\oplus T^*$}\label{lightcone}
\end{figure}
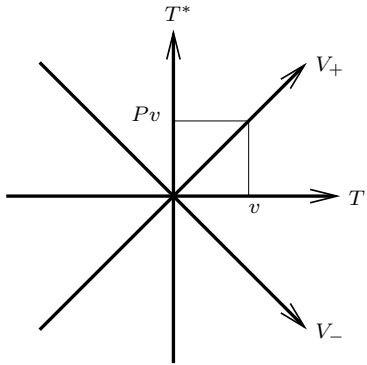

Figure~\ref{lightcone} suggests how to characterise a metric
splitting algebraically.
If we think of the coordinate axes $T$ and $T^*$ as a lightcone,
choosing a subgroup conjugate to $O(n)\times O(n)$ inside $O(n,n)$
boils down to the choice of a spacelike $V_+$ with a timelike $V_-$
as orthogonal complement. Interpreting $V_+$ as the graph of a map
$P:T\to T^*$ yields a metric and a 2-form as the symmetric and the
skew part of the dual $P\in T^*\otimes T^*$. More formally, we prove
the

\begin{prp}\label{genmetric}
The choice of an equivalence class in the space $O(n,n)/O(n)\times
O(n)$ is equivalent to either set of the following data:

{\rm (i)} A metric splitting
$$
T\oplus T^*=V_+\oplus V_-
$$
into subbundles $(V_+,g_+)$ and $(V_-,g_-)$ with positive and
negative definite metrics $g_{\pm}=(\cdot\,,\cdot)_{|V_{\pm}}$
respectively.

{\rm (ii)} A Riemannian metric $g$ and a 2-form $b$ on $T$.
\end{prp}

\begin{prf}
The equivalence between a choice of a group isomorphic to
$O(n)\times O(n)$ and (i) is true almost by definition. Consequently
we are left to show that a metric splitting gives rise to the second
set of data on $T$ and that conversely, this data allows us to
reconstruct the original metric splitting.

Suppose then we are given a metric splitting $T\oplus
T^*=V_+\oplus V_-$. Since the spaces $V_{\pm}$ are definite and
$T$ and $T^*$ are isotropic, $V_{\pm}$ can be graphed over $T$,
that is there are linear maps
$$
P_{\pm}:T\to T^*
$$
such that
$$
V_{\pm}=\{t\oplus P_{\pm}t\:|\:t\in T\}.
$$
Define the map $\widetilde{P}_+:T\times T\to \R$ by
$$
\widetilde{P}_+(s,t)=(s,P_+(t)).
$$
Then the antisymmetric part
$$
b(s,t)=\widetilde{P}_+(s,t)/2-\widetilde{P}_+(t,s)/2
$$
gives rise to a 2-form and the symmetric part defines a metric
$$
g(s,t)=\widetilde{P}_+(s,t)/2+\widetilde{P}_+(t,s)/2.
$$
Indeed, $g$ is positive definite. If $t\in T$ is non-zero, then
$$
g(t,t)=(t,P_+t)=(t\oplus P_+t,t\oplus P_+t)/2=g_+(t\oplus
P_+t,t\oplus P_+t)/2>0.
$$
Since $V_+$ and $V_-$ are orthogonal to each other an analogous
decomposition of $P_-$ yields a metric $g_-$ and a 2-form $b_-$
that satisfy
$$
g_-=-g\quad\mbox{and}\quad b_-=b.
$$

Conversely, assume we are given a metric $g$ and a 2-form $b$ on
$T$. Define the ``metric diagonal" by $D_{\pm}=\{t\oplus \mp
t\llcorner g\:|\: t\in T \}$ (the unnatural sign in the definition
of $D_{\pm}$ stems from our definition of the inner product on
$T\oplus T^*$). Then $D_+\oplus D_-$ is a metric splitting of
$T\oplus T^*$ which induces $g$ but yields a vanishing 2-form.
Twisting by $e^b$ (cf.~(\ref{ebaction}) and~(\ref{ebactionmatrix}))
gives
$$
V_{\pm}=\{t\oplus t\llcorner (\mp g+b)\:|\:t\in T\}=e^b(D_{\pm}),
$$
which induces $b$ as well as the metric $g$.
\end{prf}

\begin{rmk}
Since we are working with a fixed inner product $(\cdot\,,\cdot)$ on
$T\oplus T^*$, the scale of the induced tensor $\mc{G}$ is also
fixed, that is, $g$ and $b$ are uniquely determined. Note also that
the left and the right hand side $O(n)$ in $O(n)\times O(n)$ are
essentially the same in the sense that they induce (up to a sign)
the same metric on $T$ and thus the same $O(T)$. The picture to bear
in mind is then the following. Endowing $T$ with a metric $g$
reduces the structure group of $T$ to $O(T)$. This allows us to
define the ``trivial" generalised metric defined by the metric
diagonal $D_+\oplus D_-$. To achieve full generality we have to
twist this splitting with a $B$-field.
\end{rmk}

As a corollary, we obtain

\begin{cor}
$$
O(n,n)/O(n)\times O(n)=\{P:T\to T^*\:|\:\langle Pt,t\rangle>0\mbox{
for all }t\not= 0\}.
$$
\end{cor}

In the same vein we can also consider oriented structures. The space
$T\oplus T^*$ carries a natural orientation , and an orientation of
$V_+$ fixes thus the orientation on $V_-$ by requiring the
orientation of $V_+\oplus V_-$ to coincide with the natural one on
$T\oplus T^*$ and vice versa. This determines then an orientation on
$T$.

\begin{prp}
The choice of an equivalence class in the space $SO(n,n)/SO(n)\times
SO(n)$ is equivalent to the choice of
\begin{itemize}
\item an orientation of $T$
\item a metric $g$ on $T$
\item a $B$-field $b\in\Lambda^2T^*$.
\end{itemize}
\end{prp}

\section{Generalised exceptional structures}\label{ges}

Having introduced the notion of a generalised metric structure we
are now in a position to discuss generalised exceptional structures,
that is, structures associated with the groups $G_2\times G_2$ and
$Spin(7)\times Spin(7)$. As we saw earlier, the ``classical" $G_2$-
and $Spin(7)$-structures are associated with special orbits in an
irreducible spin representation space of $Spin(7)$ and $Spin(8)$. If
we want to continue along these lines to obtain a meaningful notion
of generalised exceptional structures we have to exhibit a special
orbit of the $Spin(7,7)$- and $Spin(8,8)$-action on even and odd
forms. The next result shows that at least locally, the groups
$G_2\times G_2$ and $Spin(7)\times Spin(7)$ arise as the stabiliser
of a spinor.

\begin{prp}\label{speccase}\hfill\newline
{\rm (i)} Assume we are given a $G_2$-structure on $T=\R^7$ with
stable 3-form $\varphi$. Then the stabiliser of the even form
$$
\rho_{\pm}=1\pm \star\varphi\in S^+=\Lambda^{ev}T^*
$$
under the spin action of $\mf{so}(7,7)$ is isomorphic to
$\mf{g}_2\otimes\C$ and $\mf{g}_2\oplus\mf{g}_2$ respectively. The
same is true for the odd form
$$
\hat{\rho}_{\pm}=\pm\varphi+vol_{\varphi}\in S^-=\Lambda^{od}T^*.
$$

{\rm (ii)} Assume we are given a $Spin(7)$-structure on $T=\R^8$
with invariant 4-form $\Omega$. Then the stabiliser of the even form
$$
\rho_{\pm}=1\pm \Omega+vol_{\Omega}\in S^+=\Lambda^{ev}T^*
$$
under the spin action of $\mf{so}(8,8)$ is isomorphic to
$\mf{so}(7)\otimes\C$ and $\mf{so}(7)\oplus \mf{so}(7)$
respectively.
\end{prp}

\begin{prf}
We consider the action of the Lie algebra $\mf{so}(7,7)$ on
$S^{\pm}$ and show that the stabiliser in either case is an
irreducible symmetric Lie algebra which are dual to each other. We
then deduce the result by appealing to Cartan's complete
classification of those.

Let $T$ be endowed with a $G_2$-structure. Then it comes equipped
with a Riemannian metric $g$ and we subsequently identify vectors
with their duals. We denote the stabiliser of $\rho_{\pm}$ by
$\mf{g}_{\pm}$.

First we show that $\mf{g}_{\pm}=\mf{g}_2\oplus \mf{g}_2$ as a
$G_2$-module. To this effect, pick a generic element in $A\oplus
b\oplus\beta\in\mf{so}(T\oplus T^*)$ (cf.~(\ref{elem_so})) where we
now see $\beta$ as 2-form. By decomposing the equation
\begin{equation}\label{eq1}
(A\oplus b\oplus\beta)\bullet\rho_{\pm}=0
\end{equation}
into homogeneous components,~(\ref{eq1}) is equivalent to the
conditions
\begin{equation}\label{set}
\begin{array}{lcl}
A^*\star\varphi=0 &\mbox{ and } & {\rm Tr}(A)=0\\
b\wedge \star\varphi=0 &\mbox{ and }&\beta\llcorner\star\varphi\pm
b=0
\end{array}
\end{equation}
The symbol $A^*$ denotes the usual action of $\mf{gl}$ on forms so
that the set of equations in the first row simply states that
$A\in\mf{g}_2(T)$ -- the Lie algebra of the $G_2$-structure on $T$.
Moreover, $b\wedge\star\varphi=0$ implies that $b\in\Lambda^2_{14}$
(cf. Proposition~\ref{spin7_form}). Here and in the sequel, we will
make intensive use of the general formulae
\begin{eqnarray*}
v\llcorner\star\xi & = & (-1)^{\deg(\xi)}\star(v\wedge\xi)\\
v\wedge\star\xi & = & (-1)^{\deg(\xi)}\star(v\llcorner\xi)
\end{eqnarray*}
which hold for an arbitrary form $\xi$ of pure degree. If we
decompose $\beta$ in its components $\beta_7\oplus\beta_{14}$ with
respect to the $G_2$-decomposition, Proposition~\ref{g2_form} gives
$$
-2\beta^{\flat}_{7}+\beta^{\flat}_{14}=\mp b_{14}
$$
and thus $\beta_{14}=\mp b$ and $\beta_{7}=0$. We define
$$
\mf{m}_{\pm}=\{b\oplus\mp b\in\Lambda^2\oplus\Lambda^2\:|\:
b\in\Lambda^2_{14}\}
$$
which is isomorphic to $\Lambda^2_{14}=\mf{g}_2$. Hence, as a
$G_2$-module, we get
$$
\mf{g}_{\pm}=\mf{g}_2\oplus\mf{g}_2.
$$
It remains to understand the Lie algebra structure of
$\mf{g}_{\pm}$. We claim that $\mf{g}_{\pm}$ is an orthogonal
symmetric irreducible Lie algebra. If we use the matrix picture of
$\mf{so}(T\oplus T^*)$ then the stabiliser $\mf{g}_{\pm}$ consists
of matrices of the form
$$
\left(\begin{array}{cc} A & \mp b\\b & A\end{array}\right),
$$
where $A\in\mf{g}_2(T)$ and $b\in\mf{m}_{\pm}$ are skew-symmetric
matrices. Put $\mf{h}=\mf{g}_2(T)$. The commutator $[C,D]=C\circ
D-D\circ C$ is the usual one for matrix algebras and we immediately
deduce the relations
$$
[\mf{h},\mf{h}]\subset\mf{h},\;[\mf{h},\mf{m}_{\pm}]\subset\mf{m}_{\pm},\;[\mf{m}_{\pm},\mf{m}_{\pm}]\subset\mf{h}.
$$
This amounts to saying that the Lie algebra $\mf{g}_{\pm}$ is
symmetric. In particular, $\mf{m}_{\pm}$ is an $\mf{h}$-irreducible
representation space. We conclude from Cartan's classification
theorem that we are left with only two possibilities for the Lie
algebra structure of $\mf{g}_{\pm}$, either the direct sum of Lie
algebras $\mf{g}_2\oplus\mf{g}_2$ or the complex Lie algebra
$\mf{g}_2\otimes\C$. We show that both cases occur and
$\mf{g}_+=\mf{g}_-^*=\mbox{the dual of }\mf{g}_-$. Recall that the
dual of $\mf{g}_-$ is defined to be the real subalgebra
$\mf{h}\oplus i\mf{m}_-$ inside the complexification
$\mf{g}_-\otimes\C$. Hence the dual is represented by matrices of
the form
$$
\left(\begin{array}{cc} A& ib\\ ib & A\end{array}\right).
$$
We easily check that
$$
\left(\begin{array}{cc} A& ib\\ ib &
A\end{array}\right)\in\mf{g}_-^*\mapsto\left(\begin{array}{cc} A&
-b\\ b & A\end{array}\right)\in\mf{g}_+
$$
sets up a (real) Lie algebra isomorphism.

Finally we have to determine which of the two Lie algebras
$\mf{g}_{\pm}$ is compact. To that end we will show that the Killing
form $B_{\mf{g}_-}$ restricted to $\mf{m}_-$ is negative definite,
hence $\mf{g}_-=\mf{g}_2\oplus\mf{g}_2$. First, we introduce a
suitable inner product on $\mf{g}_-$ to compute the Killing form.
For this we remark that the matrices in $\mf{g}_-$ are
skew-symmetric, since the matrices $A\in\mf{h}$ and $b\in\mf{m}_-$
are skew. Therefore, for any non-zero matrix $M\in\mf{g}_-$, we have
$$
-{\rm Tr}(M\circ M)=-\sum_{i,j} M_{ij}M_{ji}=\sum_{i,j} M_{ij}^2>0.
$$
It defines therefore a positive definite, $ad$-invariant inner
product on $\mf{g}_-$ which we denote by
$\langle\cdot,\cdot\rangle$. Fix an orthonormal basis $\{X_i\}$ for
this inner product and pick a nonzero $X\in\mf{m}$. Then
\begin{eqnarray*}
B_{\mf{g_-}}(X,X)&=&{\rm Tr}(ad_X\circ ad_X)=\sum_i\langle[X,[X,X_i]],X_i\rangle\\
&=&-\sum_i\langle[X,X_i],[X,X_i]\rangle<0
\end{eqnarray*}
The last equality follows from $ad$-invariance. This completes the
proof for $\rho_{\pm}$ in the first case. For $\hat{\rho}$ the
equations~(\ref{set}) become
$$
\begin{array}{lcl}
A^*\varphi=0 &\mbox{ and } & {\rm Tr}(A)=0\\
\beta\wedge\varphi=0 & \mbox{ and } & b\llcorner\varphi\pm \beta=0.
\end{array}
$$
Thus we deduce again $\beta_{14}=\mp b_{14}$ and hence the same
result holds for $\hat{\rho}_{\pm}$. Finally in (ii), we obtain the
equations
$$
\begin{array}{lcl}
A^*\Omega=0 &\mbox{ and } & {\rm Tr}(A)=0\\
b\pm\beta\llcorner\Omega=0 &\mbox{ and }&\pm
b\wedge\Omega+\beta\llcorner vol_{\Omega}=0.
\end{array}
$$
The first row states that $A\in\mf{spin}(7)$. As
$\star\Omega=\Omega$, decomposing $b$ and $\beta$ into irreducible
$Spin(7)$-components (cf. Proposition~\ref{spin7_form}) yields from
the first equation in the second row
$$
b_7=\mp3\beta_7\mbox{ and }b_{21}=\mp\beta_{21}.
$$
The second equation, however, implies
$$
b_7=\mp\frac{1}{3}\beta_7\mbox{ and }b_{21}=\mp\beta_{21}.
$$
Hence $b_7=\beta_7=0$ and $b_{21}=\mp\beta_{21}$ and the proof of
(i) carries over word for word.
\end{prf}

\begin{rmk}
Note that $\hat{\rho}_{\pm}=\star\sigma(\rho_{\pm})$ and
$\rho_{\pm}=\star\sigma(\rho_{\pm})$, where $\sigma$ is the
involution defined in~(\ref{sigma}).
\end{rmk}

Motivated by the previous proposition we define

\begin{definition}\hfill\newline
{\rm (i)} We call a reduction from the structure group $SO(7,7)$ of
$T\oplus T^*$ to $G_2\times G_2$ a {\em generalised $G_2$-
structure}.

{\rm (ii)} We call a reduction from the structure group $SO(8,8)$ of
$T\oplus T^*$ to $Spin(7)\times Spin(7)$ a {\em generalised
$Spin(7)$-structure}.

We refer to both of these structures also as {\em generalised
exceptional structures}.
\end{definition}

The previous proposition suggests that these structures can be, like
in the classical case, defined by a form which, however, is of mixed
degree. To get things rolling, we first look at the tensorial
invariants on $T$ which are induced by such a reduction along the
lines of Proposition~\ref{genmetric}. Since $G_2\times G_2$
determines some group $SO(V_+)\times SO(V_-)$ conjugate to
$SO(7)\times SO(7)$, such a structure induces a generalised oriented
metric structure $(g,b)$. We then consider the irreducible spin
representation $\Delta$ of $Spin(7)=Spin(T,g)$. The group $G_2\times
G_2$ acts on $\Delta\otimes\Delta$ and stabilises a decomposable
line element $\Psi_+\otimes\Psi_-$ which we choose to be of unit
norm. This gives rise to two groups $G_{2\pm}$ conjugated to $G_2$
inside $Spin(7)$ associated with the spinors $\Psi_+$ and $\Psi_-$.
However, these are only determined up to a scalar since
$e^{-F}\Psi_+\otimes e^F\Psi_-=\Psi_+\otimes\Psi_-$. To remove this
ambiguity up to a simultaneous sign flip in $\Psi_+$ and $\Psi_-$ we
let $\Psi_+$ and $\Psi_-$ be of unit norm and introduce the
additional scalar $F$. We refer to it as the {\em dilaton}.

We can carry out the same analysis for the group $Spin(7)\times
Spin(7)$ which sits inside some $SO(8)\times SO(8)$. Note that the
stabilised spinors can lie in either of the
$Spin(T,g)$-representation spaces $\Delta_+$ or $\Delta_-$ according
to whether or not the specific copy of $Spin(7)$ inside $Spin(T,g)$
contains plus or minus the volume element of $T=\R^8$. By choosing a
suitable global orientation on $T$, we can always assume that the
spinor $\psi_+$ stabilised by the left-hand side factor of
$Spin(7)\times Spin(7)$ lives in $\Delta_+$.

We summarise this discussion in the next proposition.

\begin{prp}\label{linalg}\hfill\newline
{\rm (i)} A generalised $G_2$-structure induces the following data
on $T=\R^7$:
\begin{itemize}
\item an orientation
\item a metric $g$
\item a $2$-form $b$
\item two unit spinors $\psi_{\pm}\in \Delta$.
\item a scalar $f$.
\end{itemize}
Here, $\Delta$ denotes the irreducible spin representation space of
$Spin(T,g)=Spin(7)$.

{\rm (ii)} A generalised $Spin(7)$-structure induces the following
data on $T=\R^8$:
\begin{itemize}
\item an orientation
\item a metric $g$
\item a $2$-form $b$
\item a unit spinor $\psi_+\in\Delta_+$ and a unit spinor $\psi_-$ in either $\Delta_+$ or $\Delta_-$.
\item a scalar $f$
\end{itemize}
Here, $\Delta_+$ and $\Delta_-$ denote the irreducible spin
representation spaces of $Spin(T,g)=Spin(8)$.
\end{prp}

Since for the generalised $Spin(7)$-case we have to take into
account the chirality of the spinors, we make the

\begin{definition}
A generalised $Spin(7)$-structure which induces two spinors of the
same or opposite chirality is said to be of {\em even} or {\em odd}
type.
\end{definition}

\begin{rmk}
The reason for this nomenclature will become apparent in a moment.
Note that the subscript $\pm$ of the spinors does not indicate the
chirality of the spinor. Its meaning will become clear in
Theorem~\ref{integrability}.
\end{rmk}

The key idea for establishing a converse of Proposition~\ref{linalg}
is to interpret the exterior algebra as the spin representation
space for $Spin(7)\times Spin(7)$ and $Spin(8)\times Spin(8)$ and to
relate it to the tensor product $\Delta\otimes\Delta$ of the
``classical" spin representations of $Spin(7)$ and $Spin(8)$. To
discuss both cases in parallel we let $(T,g)$ denote an inner
product space of dimension 7 or 8 and consider the groups
$Spin(7)=Spin(T^7,g)$ or $Spin(8)=Spin(T^8,g)$. Again, let $\Delta$
be the spin representation space of $Spin(7)$ or $Spin(8)$ which in
the latter case can be decomposed into the irreducible
$Spin(8)$-modules $\Delta_+$ and $\Delta_-$. We regard the tensor
product $\Delta\otimes\Delta$ as a $Spin(7)\times Spin(7)$- or
$Spin(8)\times Spin(8)$-representation space. We can construct a
further representation space of these groups out of the following
embedding into $\cliff(T\oplus T^*,(\cdot,\cdot))$ defined through a
metric splitting of $T\oplus T^*=V_+\oplus V_-$. Write the
associated Clifford algebra as a $\Z_2$-graded tensor product
$$
\cliff(V_+,g_+)\htimes\cliff(V_-,g_-)\cong\cliff(T\oplus T^*),
$$
where the isomorphism is given by extension of
$$
v_+\htimes v_-\mapsto v_+\cdot v_-.
$$
In particular, this maps $Spin(V_+,g_+)\times Spin(V_-,g_-)$ into
$Spin(T\oplus T^*)$. On the other hand, the maps
$$
\pi_{\pm}:x\in (T,\pm g)\mapsto -x\oplus x\llcorner(\pm g - b)
$$
induce isomorphisms
$$
\cliff(T,\pm g)\cong\cliff(V_{\pm},g_{\pm})
$$
which map $Spin(T,g)=Spin(T,-g)$ isomorphically onto
$Spin(V_{\pm},g_{\pm})$. The compounded algebra isomorphism
$$
\begin{array}{ccl}
\iota_b:\cliff(T,g)\htimes\cliff(T,-g) & \to &
\cliff(T\oplus T^*)\\
x\htimes y & \mapsto & (-x\oplus x\llcorner(g-b))\bullet(-y\oplus
y\llcorner(-g-b))
\end{array}
$$
maps
$$
Spin(T)\times Spin(T)\to Spin(V_+)\times Spin(V_-)\leqslant
Spin(T\oplus T^*).
$$
We will suppress the subscript $b$ and simply write $\iota$ if
$b\equiv 0$. Recall also our convention to denote by $\bullet$ the
product and the Clifford action of $\cliff(T\oplus T^*)$. Now we are
in a position to compare the actions on $\Delta\otimes\Delta$ and
$\Lambda^{ev,od}$. The so-called {\em pinor product}
$\psi_0\pin\psi_1$ for two spinors $\psi_0$ and $\psi_1$ is the
endomorphism
$$
\psi_0\pin\psi_1(\phi)=q_{\Delta}(\psi_1,\phi)\psi_0
$$
(where $q_{\Delta}$ denotes the spin-invariant inner product on
$\Delta$). Clifford multiplication identifies any element in the
Clifford algebra with an endomorphism of $\Delta$. It relates to the
pinor product by~\cite{ha91}
\begin{equation}\label{cliffordid}
(x\cdot\psi_0)\pin\psi_1=x\circ(\psi_0\pin\psi_1)\mbox{ and
}\psi_0\pin(x\cdot\psi_1)=(\psi_0\pin\psi_1)\circ\sigma(x),
\end{equation}
where $\sigma$ is the Clifford involution of~(\ref{sigma}). If we
regard $\psi_0\pin\psi_1$ as an element in $\cliff(T\oplus T^*)$,
then
$$
\sigma(\psi_0\pin\psi_1)=\psi_1\pin\psi_0.
$$
Let $R$ denote the representations of $\cliff(T^8)$ or $\cliff(T^7)$
as defined in Sections~\ref{triality} and~\ref{g2}. Recall that $J$
was the vector space isomorphism between $\cliff(T)$ and
$\Lambda^*T^*$ defined in~(\ref{Jmap}). Then we define the following
maps
$$
L^8_b:\Delta\otimes\Delta\stackrel{\pin}{\longrightarrow}
\End(\Delta_+\oplus\Delta_-)\stackrel{R^{-1}}{\longrightarrow}\cliff(T^8,g)\stackrel{J}{\longrightarrow}
\Lambda^*T^*\stackrel{e^{b/2}\wedge}{\longrightarrow}\Lambda^*T^*
$$
and
$$
L^7_b:\Delta\otimes\Delta\stackrel{\pin\oplus0}{\longrightarrow}
\End(\Delta)\oplus\End(\Delta)\stackrel{R^{-1}}{\longrightarrow}\cliff(T^7,g)\stackrel{J}{\longrightarrow}
\Lambda^*T^*\stackrel{e^{b/2}\wedge}{\longrightarrow}\Lambda^*T^*.
$$
As the 2-form $b$ is to play a crucial part in the latter
development, it is worthwhile emphasising it in the notation. If
$b=0$ we shall simply write $L_0=L$. Here and from now on we will
drop the superscripts 7 and 8, relying on the context to determine
which particular map is under consideration. We denote by
$L^{ev,od}_b$ the map $L_b$ followed by the projections on the
even or odd part of the exterior algebra, i.e.
$$
L_b=L^{ev}_b\oplus L^{od}_b.
$$
Note that
$$
\begin{array}{lcl}
L_b(\Delta_+\otimes\Delta_+\oplus\Delta_-\otimes\Delta_-)
& = & \Lambda^{ev}\\
L_b(\Delta_+\otimes\Delta_-\oplus\Delta_-\otimes\Delta_+)
& = & \Lambda^{od},\\
\end{array}
$$
which justifies our notion of a generalised ``even" or ``odd"
$Spin(7)$-structure. Whenever it makes sense, we regard the maps
$L^{ev,od}_b$ as the restriction of $L_b$ to the spinors of equal or
opposite chirality. Moreover, we denote by $L^{ev}_{\pm}$ and
$L^{od}_{b,\pm}$ the restrictions to $\Delta_+\otimes\Delta_+$,
$\Delta_-\otimes\Delta_-$, $\Delta_-\otimes\Delta_+$ and
$\Delta_+\otimes\Delta_-$ respectively.\label{Lnot} The next
proposition says that up to a sign twist, the action of $T$ on
$\Delta\otimes\Delta$ and $\Lambda^{ev,od}$ commute.

\begin{prp}\label{lmap}
For any $x\in T$ we have
\begin{eqnarray*}
L_b^{ev,od}(x\cdot\varphi\otimes\psi) & = &
\phantom{\pm}\iota_b(x\htimes1)\bullet L_b^{od,ev}(\varphi\otimes\psi)\\
L_b^{ev,od}(\varphi\otimes x\cdot\psi) & = & \pm\iota_b(1\htimes
x)\bullet L_b^{od,ev}(\varphi\otimes\psi).
\end{eqnarray*}
\end{prp}

We postpone the proof of Proposition~\ref{lmap} to draw several corollaries first.

\begin{cor}\label{lmapcor}
For unit vectors $x_1\ldots,x_k\in (T,g)$ and $y_1,\ldots,y_l\in
(T,-g)$ with $k+l$ even we have
$$
L^{ev,od}_b(x_1\cdot\ldots \cdot x_k\cdot\varphi\otimes
y_1\cdot\ldots \cdot y_l\cdot\psi)=(-1)^ls(l)\iota_b(x_1\cdot\ldots\cdot
x_k\htimes y_1\cdot\ldots\cdot y_l)\bullet
L^{ev,od}_b(\varphi\otimes\psi),
$$
where $s(l)=1$ for $l\equiv 0,1 \mod 4$ and $s(l)=-1$ for $l\equiv
2,3 \mod 4$. For $k+l$ odd we obtain
$$
L^{ev,od}_b(x_1\cdot\ldots\cdot x_k\cdot\varphi\otimes
y_1\cdot\ldots \cdot y_l\cdot\psi)=(-1)^ls(l)\iota_b(x_1\cdot\ldots\cdot
x_k\htimes y_1\cdot\ldots\cdot y_l)\bullet
L^{od,ev}_b(\varphi\otimes\psi).
$$
Moreover, $L^{ev,od}_b$ preserves any invariant space of a subgroup
of $Spin(7)\times Spin(7)$ or $Spin(8)\times Spin(8)$. In
particular, we get the following decomposition into irreducible
k-dimensional representations $V_k$ of $G_2\times G_2$ and
$Spin(7)\times Spin(7)$.

{\rm (i)} $\Lambda^{ev,od}T^{7*}=V_1\oplus V_{l,7}\oplus
V_{r,7}\oplus V_{49}$

{\rm (ii)} $\Lambda^{ev,od}T^{8*}=V_1\oplus V_{l,7}\oplus
V_{r,7}\oplus V_{49}\oplus V_{64}.$

$V_1$ corresponds to the trivial representation spanned by
$\psi_+\otimes\psi_-$, $V_{l,7}$ to $\psi_+^{\perp}\otimes\psi_-$,
$V_{r,7}$ to $\psi_-\otimes\psi_+^{\perp}$ and $V_{49}$ to
$\psi_+^{\perp}\otimes\psi_-^{\perp}$.
\end{cor}

\begin{prf}
If $k+l$ is even, even- or oddness of $k$ implies that of $l$. We
abbreviate $x_1\cdot\ldots\cdot x_l$ and $y_1\cdot\ldots\cdot y_k$
by $\underline{x}$ and $\underline{y}$. Now suppose that $k$ is
even, then
$$
L^{ev}_b(\underline{x}\cdot\varphi\otimes\underline{y}\cdot\psi)=\iota_b(\underline{x}\htimes1)\bullet
L^{ev}_b(\varphi\otimes\underline{y}\cdot\psi).
$$
Moreover, $l$ is even and because of the sign flip for $L^{od}_b$,
we get an extra minus sign whenever $l$ is of the form $4m+2$ while
this sign cancels for $l=4m$. Hence the sign is given by $s(l)$. The
other cases are discussed in the same way where the parity flips if
$k+l$ is odd. The remaining claim is straightforward.
\end{prf}

\begin{cor}
Restricted to the action of $a\in Spin(T)\mapsto(a,a)\in
Spin(T)\times Spin(T)$, we have
$$
L^{ev,od}(a\cdot\varphi\otimes
a\cdot\psi)=\pi_0(a)^*L^{ev,od}(\varphi\otimes\psi)
$$
for any $a\in Spin(T,g)$, where $\pi_0(a)^*$ denotes the action on
forms induced by the vector representation $\pi_0:Spin(T,g)\to
SO(T,g)$. Hence, $L^{ev,od}$ defines a $Spin(T)$-equivalence between
$\Delta\otimes\Delta$ and $\Lambda^{ev,od}$ or $\Lambda^*$ (which is
unique up to a scalar).
\end{cor}

\begin{prf}
To derive the assertion we look at the action of the Lie algebra,
that is we want to show
$$
L^{ev,od}(A\cdot\varphi\otimes\psi+\varphi\otimes
A\cdot\psi)=\pi_{0*}(A)^*L^{ev,od}(\varphi\otimes\psi)
$$
for any $A\in\mf{spin}(T,g)$. We fix an orthonormal basis $\{e_i\}$.
Since $e_i\wedge(e_j\llcorner\tau)=-e_j\llcorner(e_i\wedge\tau)$ for
any form $\tau$ we get for $e_{ij}=e_i\cdot e_j$ with $i\not =j$
that
\begin{eqnarray*}
L^{ev,od}(e_{ij}\cdot\varphi\otimes\psi+\varphi\otimes
e_{ij}\cdot\psi) & = &
\phantom{-}\iota(e_{ij}\oplus0)\bullet L^{ev,od}(\varphi\otimes\psi)-\\
& &-\iota(0\oplus e_{ij})\bullet L^{ev,od}(\varphi\otimes\psi)\\
& = & \phantom{-}2(e_j\llcorner(e_i\wedge L^{ev,od}(\varphi\otimes\psi))+\\
& & +e_j\wedge(e_i\llcorner L^{ev,od}(\varphi\otimes\psi))\\
& = & \phantom{-}\pi_{0*}(e_i\cdot
e_j)^*L^{ev,od}(\varphi\otimes\psi),
\end{eqnarray*}
hence the result.
\end{prf}

\begin{cor}\hfill\newline
{\rm (i)} Generalised $G_2$-structures are in 1-1 correspondence
with lines of forms $\rho$ in $\Lambda^{ev}$ or $\Lambda^{od}$ whose
stabiliser under the action of $Spin(7,7)$ is isomorphic to
$G_2\times G_2$. We refer to $\rho$ as the {\em structure form} of
the generalised $G_2$-structure. This form can be uniquely written
(modulo a simultaneous sign change for $\Psi_+$ and $\Psi_-$) as
$$
\rho=e^{-f}L^{ev}(\psi_+\otimes\psi_-)\mbox{ or
}\rho=e^{-f}L^{od}(\psi_+\otimes\psi_-).
$$

{\rm (ii)} Generalised $Spin(7)$-structures are in 1-1
correspondence with lines of forms $\rho$ in $\Lambda^{ev}$ or
$\Lambda^{od}$ whose stabiliser under the action of $Spin(8,8)$ is
isomorphic to $Spin(7)\times Spin(7)$. We refer to $\rho$ as the
{\em structure form} of the generalised $Spin(7)$-structure. We
speak of an {\em even} or {\em odd} structure according to whether
$\rho$ is even or odd. This form can be uniquely written (modulo a
simultaneous sign change for $\Psi_+$ and $\Psi_-$) as
$$
\rho=e^{-f}L(\psi_+\otimes\psi_-).
$$
\end{cor}

Now we turn to the proof for Proposition~\ref{lmap}.

\begin{prf} (of proposition~\ref{lmap}) First we assume that $b\equiv 0$.
Starting with the 8-dimensional case, we see that
\begin{eqnarray*}
L(x\cdot\varphi\otimes\psi) & = &
\phantom{-}J(R^{-1}(R(x)(\varphi)\pin\psi))\\
& = & \phantom{-}J(R^{-1}(R(x)\circ(\varphi\pin\psi)))\\
& = & \phantom{-}J(x\cdot R^{-1}(\varphi\pin\psi))\\
& = & -x\llcorner L(\varphi\otimes\psi)+x\wedge
L(\varphi\otimes\psi)\\
& = & \phantom{-}\iota(x\htimes 1)\bullet L(\varphi\otimes\psi),
\end{eqnarray*}
where we have used~(\ref{cliffordid}). Hence we get
$$
L^{ev}(x\cdot\varphi\otimes\psi)=\iota(x\htimes 1)\bullet
L^{od}(\varphi\otimes\psi),\quad
L^{od}(x\cdot\varphi\otimes\psi)=\iota(x\htimes 1)\bullet
L^{ev}(\varphi\otimes\psi).
$$
Similarly, we find
\begin{eqnarray*}
L(\varphi\otimes x\cdot\psi) & = &
\phantom{-}J(R^{-1}(\varphi\pin R(x)(\psi)))\\
& = & \phantom{-}J(R^{-1}(\varphi\pin\psi\circ R(\hat{x})))\\
& = & -J( R^{-1}(\varphi\pin\psi)\cdot x)\\
& = & \phantom{+}x\llcorner L^{ev}(\varphi\otimes\psi)-x\wedge
L^{ev}(\varphi\otimes\psi)+\\
&  & +x\llcorner L^{od}(\varphi\otimes\psi)+x\wedge
L^{od}(\varphi\otimes\psi)\\
& = & -\iota(1\htimes x)\bullet
L^{ev}(\varphi\otimes\psi)+\iota(1\htimes x)\bullet
L^{od}(\varphi\otimes\psi),
\end{eqnarray*}
so that
$$
L^{ev}(\varphi\otimes x\cdot\psi)=\iota(1\htimes x)\bullet
L^{od}(\varphi\otimes\psi),\quad L^{od}(\varphi\otimes
x\cdot\psi)=-\iota(1\htimes x)\bullet L^{ev}(\varphi\otimes\psi).
$$
Next we turn to the case $n=7$. By convention, we let $T^7$ act
through the inclusion
$$
T^7\hookrightarrow\cliff(T^7,g)\stackrel{
R}{\cong}\End(\Delta)\oplus\End(\Delta)
$$
followed by projection on the first summand. Thus
\begin{eqnarray*}
L(x\cdot\varphi\otimes\psi) & = &
J( R^{-1}((x\cdot\varphi)\pin\psi\oplus0))\\
& = & J(x\cdot R^{-1}(\varphi\pin\psi\oplus0))
\end{eqnarray*}
and we can argue as above. The same applies to $L(\varphi\otimes
x\cdot\psi)$ so that the case $b=0$ is shown.

Now let $b$ be an arbitrary B-field. For the sake of clarity we will
temporarily denote by $\widetilde{\exp}$ the exponential map from
$\mf{so}(7)$ to $Spin(7)$, while the untilded exponential takes
values in $SO(7)$. The adjoint representation $Ad$ of the group of
units inside a general Clifford algebra $\cliff(V)$ restricts to the
double cover $Spin(V)\to SO(V)$ still denoted by $Ad$. As a
transformation in $SO(V)$ we then have
$Ad\circ\widetilde{\exp}=\exp\circ ad$. Since $ad(v\wedge
w)=[v,w]/4$ we obtain
$$
e^b=Ad(\widetilde{e}^{\sum b_{ij}e_i\cdot e_j/2})
$$
for the $B$-field $b=\sum_{i<j}b_{ij}e_i\wedge e_j$. Hence
\begin{eqnarray*}
L^{ev,od}_b(x\cdot\varphi\otimes\psi) & = &
\widetilde{e}^{b/2}\bullet
L^{ev,od}(x\cdot\varphi\otimes\psi)\\
& = &
(\widetilde{e}^{b/2}\bullet\iota(x\htimes 1))\bullet L^{od,ev}(\varphi\otimes\psi)\\
& = & (\widetilde{e}^{b/2}\bullet\iota(x\htimes
1)\bullet \widetilde{e}^{-b/2})\bullet L^{od,ev}_b(\varphi\otimes\psi)\\
& = & Ad(\widetilde{e}^{b/2})(\iota(x\htimes 1))\bullet
L^{od,ev}_b(\varphi\otimes\psi)\\
& = & e^b(\iota(x\htimes 1))\bullet
L^{od,ev}_b(\varphi\otimes\psi)\\
& = & \iota_b(x\htimes 1)\bullet L^{od,ev}_b(\varphi\otimes\psi).
\end{eqnarray*}
Similarly, the claim is checked for $L^{ev,od}_b(\varphi\otimes
x\cdot\psi)$ which completes the proof.
\end{prf}

The next corollary works out how a 3-form acts by wedging and
contraction. This will become important in Section
\ref{spinsolgencase} where 3-forms naturally appear as the torsion
form of a linear metric connection.

\begin{cor}\label{baction}
Let $\tau\in\Lambda^3$. Then the identities
\begin{eqnarray*}
\tau\llcorner L^{ev,od}(\varphi\otimes\psi) & =
&
\frac{1}{8}L^{od,ev}(-\tau\cdot\varphi\otimes\psi\pm\varphi\otimes\tau\cdot\psi\mp\\
& & \phantom{\frac{1}{8}L^{od,ev}(}\mp
\sum\limits_i(e_i\llcorner\tau)\cdot\varphi\otimes
e_i\cdot\psi+\sum\limits_ie_i\cdot\varphi\otimes(e_i\llcorner\tau)\cdot\psi)
\end{eqnarray*}
and
\begin{eqnarray*}
\tau\wedge L^{ev,od}(\varphi\otimes\psi) & = &
\frac{1}{8}L^{od,ev}(\phantom{\pm}\tau\cdot\varphi\otimes\psi\pm\varphi\otimes\tau\cdot\psi\mp\\
& & \phantom{\frac{1}{8}L^{od,ev}(}\mp
\sum\limits_i(e_i\llcorner\tau)\cdot\varphi\otimes
e_i\cdot\psi-\sum\limits_ie_i\cdot\varphi\otimes(e_i\llcorner\tau)\cdot\psi)
\end{eqnarray*}
hold.
\end{cor}

\begin{prf}
From Proposition~\ref{lmap} we derive
\begin{eqnarray*}
L^{ev,od}(x\cdot\varphi\otimes\psi) & = & \phantom{\pm}x\wedge
L^{od,ev}(\varphi\otimes\psi)-x\llcorner
L^{od,ev}(\varphi\otimes\psi)\\
L^{ev,od}(\varphi\otimes x\cdot\psi) & = & \pm x\wedge
L^{od,ev}(\varphi\otimes\psi)\pm x\llcorner
L^{od,ev}(\varphi\otimes\psi)
\end{eqnarray*}
which implies
\begin{eqnarray*}
x\wedge L^{ev}(\varphi\otimes\psi) & = &
\frac{1}{2}L^{od}(x\cdot\varphi\otimes\psi-\varphi\otimes x\cdot\psi)\\
x\wedge L^{od}(\varphi\otimes\psi) & = & \frac{1}{2}L^{ev}(x\cdot\varphi\otimes\psi+\varphi\otimes x\cdot\psi)\\
\end{eqnarray*}
and
\begin{eqnarray*}
x\llcorner L^{ev}(\varphi\otimes\psi) & = &
\frac{1}{2}L^{od}(-x\cdot\varphi\otimes\psi-\varphi\otimes
x\cdot\psi)\\
x\llcorner L^{od}(\varphi\otimes\psi) & = &
\frac{1}{2}L^{ev}(-x\cdot\varphi\otimes\psi+\varphi\otimes
x\cdot\psi).
\end{eqnarray*}
Fix an orthonormal basis $e_1,\ldots,e_7$. By using the identities
above we gather that
\begin{eqnarray*}
e_{ijk}\llcorner L^{ev,od}(\varphi\otimes\psi) & =
&\phantom{\pm}\frac{1}{8}L^{od,ev}(-e_{ijk}\cdot\varphi\otimes\psi\pm\varphi\otimes e_{ijk}\cdot\psi\pm\\
& &
\pm e_{ik}\cdot\varphi \otimes e_j\cdot\psi - e_i\cdot\varphi\otimes e_{kj}\cdot\psi\pm\\
& & \pm e_{ji}\cdot\varphi\otimes e_k\cdot\psi-e_j\cdot\varphi\otimes e_{ik}\psi\pm\\
& & \pm e_{kj}\cdot\varphi\otimes
e_i\cdot\psi-e_k\cdot\varphi\otimes e_{ji}\cdot\psi)
\end{eqnarray*}
and
\begin{eqnarray*}
e_{ijk}\wedge L^{ev,od}(\varphi\otimes\psi) & =
&\phantom{\pm}\frac{1}{8}L^{od,ev}(e_{ijk}\cdot\varphi\otimes\psi\pm\varphi\otimes e_{ijk}\cdot\psi\pm\\
& &
\pm e_{ik}\cdot\varphi \otimes e_j\cdot\psi+e_i\cdot\varphi\otimes e_{kj}\cdot\psi\pm \\
& & \pm e_{ji}\cdot\varphi\otimes e_k\cdot\psi+e_j\cdot\varphi\otimes e_{ik}\psi\pm\\
& & \pm e_{kj}\cdot\varphi\otimes e_i\cdot\psi+ e_k\cdot\varphi
\otimes e_{ji}\cdot\psi).
\end{eqnarray*}
For instance, we have
\begin{eqnarray*}
e_{ijk}\wedge L^{ev}(\varphi\otimes\psi) & = &
\phantom{-}\frac{1}{2}e_{ij}\wedge
L^{od}(e_k\cdot\varphi\otimes\psi-\varphi\otimes e_k\cdot\psi)\\
& = & \phantom{-}\frac{1}{4}e_i\wedge
L^{ev}(e_{jk}\cdot\varphi\otimes\psi+e_k\cdot\varphi\otimes
e_j\cdot\psi-\\
& & -e_j\cdot\varphi\otimes e_k\cdot\psi-\varphi\otimes e_{jk}\cdot\psi)\\
& = &
\phantom{-}\frac{1}{8}L^{od}(e_{ijk}\cdot\varphi\otimes\psi-e_{jk}\cdot\varphi\otimes
e_i\cdot\psi+e_{ik}\cdot\varphi\otimes
e_j\cdot\psi-\\
& & {}-e_k\cdot\varphi\otimes
e_{ij}\cdot\psi-e_{ij}\cdot\varphi\otimes
e_k\cdot\psi+e_j\cdot\varphi\otimes
e_{ik}\cdot\psi-\\
& & {}-e_i\cdot\varphi\otimes e_{jk}\cdot\psi+\varphi\otimes
e_{ijk}\cdot\psi).
\end{eqnarray*}
The remaining assertions are established in the same way. If
$\tau=\sum\limits_{i<j<k}\tau_{ijk}e_{ijk}$ is a general 3-form,
then
\begin{eqnarray*}
\tau\llcorner L^{ev,od}(\varphi\otimes\psi) & = &
\phantom{\pm}\frac{1}{6}\sum\limits_{i,j,k}\tau_{ijk}e_{ijk}\llcorner
L^{ev,od}(\varphi\otimes\psi)\\
& = &
\phantom{\pm}\frac{1}{48}L^{od,ev}(-6\tau\cdot\varphi\otimes\psi\pm\varphi\otimes6\tau\cdot\psi\mp\\
& & \mp 3\sum\limits_{ijk}\tau_{ijk}e_{jk}\cdot\varphi\otimes e_i\cdot\psi+3\sum\limits_{ijk}\tau_{ijk}e_i\cdot\varphi\otimes e_{jk}\cdot\tau)\\
& = &
\phantom{\pm}\frac{1}{8}L^{od,ev}(-\tau\cdot\varphi\otimes\psi\pm\varphi\otimes\tau\cdot\psi\mp\\
& & \mp \sum\limits_i(e_i\llcorner\tau)\cdot\varphi\otimes
e_i\cdot\psi+\sum\limits_ie_i\cdot\varphi\otimes(e_i\llcorner\tau)\cdot\psi),
\end{eqnarray*}
since $e_i\llcorner\tau=\sum_{j,k}\tau_{ijk}e_{jk}/2$. The second
identity is derived in a similar fashion.
\end{prf}

Choosing a generalised metric structure and an orientation on $T$
also sets up the isomorphism $L_b$. Since the spin representation
$\Delta$ of $Spin(7)$ is of real type so is the tensor product. We
remark that up to isomorphism, there is only one 49-dimensional
irreducible representation space of $Spin(7)\times Spin(7)$, hence
$\Lambda^{ev}$ and $\Lambda^{od}$ are isomorphic as $Spin(7)\times
Spin(7)$-spaces. Consequently, there is (up to a scalar) a unique
$Spin(7)\times Spin(7)$-invariant in
$\Lambda^{ev,od}\otimes\Lambda^{od,ev}$, that is we obtain an
equivariant map $\Lambda^{ev,od}\to\Lambda^{od,ev}$. The description
of this invariant will occupy us next. Morally it is the Hodge
$\star$-operator twisted with a $B$-field and the anti-automorphism
$\sigma$ defined in (\ref{sigma}).

\begin{definition}
The {\em box operator} or {\em generalised Hodge $\star$-operator}
$\Box_{g,b}:\Lambda^{ev,od}T^*\to\Lambda^{od,ev}T^*$ ($n$ odd) or
$\Box_{g,b}:\Lambda^{ev,od}T^*\to\Lambda^{ev,od}T^*$ ($n$ even)
associated with $g$ and $b$ is defined by
$$
\Box_{g,b}\rho=e^{b/2}\wedge\star_g\sigma(e^{-b/2}\wedge\rho).
$$
\end{definition}

We will need the following technical lemma which is an immediate
consequence from the definition of $\sigma$.

\begin{lem}\label{hat}
Let $\rho^{ev,od}\in\Lambda^{ev,od}T^*$ be an even or an odd form
and let $x\in T^*$. Then
\begin{enumerate}
\item $\begin{array}{lcll}\star\sigma(\rho^{ev,od}) & = & \phantom{\pm}\sigma(\star\rho^{ev,od}), & n=7\\
\star\sigma(\rho^{ev,od}) & = & \pm\sigma(\star\rho^{ev,od}), &
n=8\end{array}$
\item $\sigma(e^b\wedge\rho^{ev,od})=e^{-b}\wedge\sigma(\rho^{ev,od})$
\end{enumerate}
\end{lem}

Now we can prove the following

\begin{prp}\label{selfdual}\hfill\newline
{\rm (i)} If $n=7$, then
$$
\Box_{g,b}L_b=L_b
$$
or equivalently,
$$
\Box_{g,b}L_b^{ev,od}=L_b^{od,ev}.
$$

{\rm (ii)} If $n=8$, then
$$
\Box_{g,b}L_{b,\pm}^{ev,od}=\pm L^{ev,od}_{b,\pm}.
$$

In particular, $\Box_{g,b}$ is $Spin(7)\times Spin(7)$- and
$Spin(8)\times Spin(8)$-equivariant.
\end{prp}

\begin{prf}
For any $\varphi,\psi\in\Delta\otimes\Delta$ ($n=7$ or $8$) we have
according to~(\ref{cliffhodge}) and Lemma~\ref{hat}
\begin{eqnarray*}
\Box_{g,b}L_b(\varphi\otimes\psi) & = &
e^{b/2}\wedge\star\sigma(e^{-b/2}\wedge L_b(\varphi\otimes\psi))\\
& = & e^{b/2}\wedge\star\sigma(L(\varphi\otimes\psi))\\
& = & e^{b/2}\wedge J( R^{-1}(\varphi\pin\psi)\cdot\lambda).
\end{eqnarray*}
Now for $n=7$, $ R^{-1}((\varphi\pin\psi\oplus0)\cdot vol)$ is just
$ R^{-1}(\varphi\pin\psi\oplus0)$ while for $n=8$ the sign changes
according to the chirality of the spinor, since the volume element
acts on $\Delta_{\pm}$ by $\pm {\rm id}$.
\end{prf}

If $g$ and $b$ are induced by $\rho$ we will also use the sloppier
notation $\Box_{\rho}$ or drop the subscript altogether.
Nevertheless it is important to bear in mind that the $\Box$
operator is really induced by the choice of a generalised oriented
metric structure. Consequently, if this generalised metric gets
conjugated by an element in $A\in O(7,7)$ or $O(8,8)$, then the
$\Box$-operator will transform naturally for the lift
$\widetilde{A}\in Pin(7,7)$ or $Pin(8,8)$ which means that
\begin{equation}\label{boxtrans}
\Box_{\widetilde{A}\bullet\rho}\widetilde{A}\bullet\rho=\widetilde{A}\bullet\Box_{\rho}\rho.
\end{equation}

For concrete computations it is useful to have an explicit
description of the structure forms. Taking our standard
representation $R$ of $\cliff(\R^7)$ or $\cliff(\R^8)$, we can
readily compute the coefficients from the description as a tensor
product by using the following formula.

\begin{lem}
Extend $g$ to an inner product on $\Lambda^*$ (which we still denote
by $g$) and let $I$ be a multi-index. Then
$$
g(L(\varphi\otimes\psi),e_I)=\frac{1}{16}q(
R(e_I)\cdot\psi,\varphi).
$$
\end{lem}

This lemma is a direct consequence of Theorem 13.73 in~\cite{ha91}.
The process of passing from the tensor product to the exterior
algebra is also referred to as {\em fierzing} in the physics
literature -- we saw an instance of this procedure in
Section~\ref{class_susy} where we related the square of a unit
spinor with a form. To ease the notation we assume our map $L$ to be
rescaled by $1/16$ so that we can henceforth discard this factor.
Now suppose we are given two unit spinors $\psi_+$ and $\psi_-$, say
in the 7-dimensional case. There is a pair of orthonormal spinors
$\widetilde{\psi}_0,\widetilde{\psi}_1$ with
$\psi_+=\widetilde{\psi}_0$ and
$\psi_-=c\widetilde{\psi}_0+s\widetilde{\psi}_1$. As the notation
suggests, this is merely the projection onto $\widetilde{\psi}_0$
and $\widetilde{\psi}_1$ whose coefficients are given by the sine
and cosine of the angle $\sphericalangle(\psi_+,\psi_-)$. By Theorem
14.69 in~\cite{ha91}, the action of $Spin(7)$ on the Stiefel variety
$V_2(\Delta)$, the set of pairs of orthonormal spinors, is
transitive. Hence, there is an $w\in Spin(T,g)$ such that
$\widetilde{\psi}_0=w\cdot\psi_0$ and
$\widetilde{\psi}_1=w\cdot\psi_1$, where,
$\psi_0=1,\psi_1=i,\ldots,\psi_7=e\cdot k$ is the standard basis of
$\Oc$ which we use in our explicit matrix representation. Since
$L^{ev,od}(w\cdot\psi_0\otimes
w\cdot\psi_1)=\pi_0(w)^*L^{ev,od}(\psi_0\otimes \psi_1)$, we may
assume that $\psi_+=\psi_0$ and $\psi_-=c\psi_0+s\psi_1$. In the
case $n=8$, we have to distinguish between the odd and even case.
$Spin(8)$ acts transitively on $V_2(\Delta_+)$ and hence for
structures of even type we can apply the same argument as before. In
the odd case we appeal again to Theorem 14.69 in~\cite{ha91} to see
that the action of $Spin(8)$ on $\Delta_+\times\Delta_-$ is
transitive on the product of spheres $S^7\times S^7$. Hence we may
assume that $\psi_+$ and $\psi_-$ are both unit spinors which are
orthogonal to each other. The computation of the normal form is now
straightforward.

\begin{prp}\label{normalform}\hfill\newline
{\rm (i)} Let $\rho^{ev,od}\in\Lambda^{ev,od}T^{7*}$ be a form
stabilised by $G_2\times G_2$. Then there exists a unique
$b\in\Lambda^2$ and $f\in\R$ with
$$
\rho=e^{-f}e^{b/2}\wedge \rho_0^{ev,od}
$$
and an orthonormal basis $e_1,\ldots,e_7$ such that
\begin{eqnarray*}
\rho^{ev}_0 & = &
\phantom{+}\cos(a)+\sin(a)(-e_{23}-e_{45}+e_{67})+\\
& &
+\cos(a)(-e_{1247}+e_{1256}+e_{1346}+e_{1357}-e_{2345}+e_{2367}+e_{4567})+\\
& & +\sin(a)(e_{1246}+e_{1257}+e_{1347}-e_{1356})-\sin(a)e_{234567}
\end{eqnarray*}
and
\begin{eqnarray*}
\rho^{od}_0 & = &
\phantom{+}\sin(a)e_1+\sin(a)(e_{247}-e_{256}-e_{346}-e_{357})+\\
& &
+\cos(a)(e_{123}+e_{145}-e_{167}+e_{246}+e_{257}+e_{347}-e_{356})+\\
& &
+\sin(a)(-e_{12345}+e_{12367}+e_{14567})+\cos(\alpha)e_{1234567},
\end{eqnarray*}
where $a=\sphericalangle(\psi_+,\psi_-)$.

{\rm (ii)} Let $\rho\in\Lambda^{ev}T^{8*}$ be an even form
stabilised by $Spin(7)\times Spin(7)$. Then there exists a unique
$b\in\Lambda^2$ and $f\in\R$ with
$$
\rho=e^{-f}e^{b/2}\wedge \rho^{ev}_0
$$
and an orthonormal basis $e_0,\ldots,e_7$ such that
\begin{eqnarray*}
\rho_0^{ev} & = &
\phantom{+}\cos(a)+\sin(a)(e_{01}-e_{23}-e_{45}+e_{67})+\\
& &
+\cos(a)(e_{0123}+e_{0145}-e_{0167}+e_{0246}+e_{0257}+e_{0347}-e_{0356}-\\
& &
-e_{1247}+e_{1256}+e_{1346}+e_{1357}-e_{2345}+e_{2367}+e_{4567})+\\
& &
+\sin(a)(e_{0247}-e_{0256}-e_{0346}-e_{0357}+e_{1246}+e_{1257}+e_{1347}-e_{1356})+\\
& &
+\sin(a)(-e_{012345}+e_{012367}+e_{014567}-e_{234567})+\cos(a)e_{01234567},
\end{eqnarray*}
where $a=\sphericalangle(\psi_+,\psi_-)$.

{\rm (iii)} Let $\rho\in\Lambda^{od}T^{8*}$ be an odd form
stabilised by $Spin(7)\times Spin(7)$. Then there exists a unique
$b\in\Lambda^2$ and $f\in\R$ with
$$
\rho=e^{-f}e^{b/2}\wedge\rho^{od}_0.
$$
There exists an orthonormal basis $e_0,\ldots,e_7$ such that
\begin{eqnarray*}
\rho_0^{od} & = & -e_0+
e_{123}+e_{145}-e_{167}+e_{246}+e_{257}+e_{347}-e_{356}+\\
& &
{}+e_{01247}-e_{01256}-e_{01346}-e_{01357}+e_{02345}-e_{02367}-e_{04567}+e_{1234567}.
\end{eqnarray*}
\end{prp}

If the spinors are not parallel, we can express the homogeneous
components of $\rho$ in terms of the invariant forms of the
intersection $G_{2+}\cap G_{2-}=SU(3)$, $Spin(7,\psi_+)_+\cap
Spin(7,\psi_-)_+=Spin(6)$ or $Spin(7)_+\cap Spin(7)_-=G_2$ (cf.
Section~\ref{class_susy}).

\begin{cor}\label{normalformcor1}
In terms of the underlying $SU(3)$- and $G_2$-invariants (cf.
Section~\ref{class_susy}, in particular Table~\ref{overview}), we
find for these 3 cases (letting $s$ and $c$ be shorthand for the
(co-)sine of $a$)

{\rm (i)} generalised $G_2$
$$
\rho_0^{ev}=c+s\omega-c(\psi_-\wedge\alpha+\frac{1}{2}\omega^2)+s\psi_+\wedge\alpha-\frac{1}{6}s\omega^3
$$
and
$$
\rho_0^{od}=s\alpha-c(\psi_++\omega\wedge\alpha)-s\psi_--\frac{1}{2}s\omega^2\wedge\alpha+c
vol_g.
$$

{\rm (ii)} generalised $Spin(7)$, even type
\begin{eqnarray*}
\rho_0^{ev} & = &
\phantom{-}c+s\varpi+c(\Omega_1-\frac{1}{2}\varpi^2)+s\Omega_2-\frac{s}{3}\varpi^3+cvol_g\\
& = & \phantom{-}\gamma\wedge\rho_0^{od,G_2}+\rho_0^{ev,G_2},
\end{eqnarray*}
where $\rho_0^{ev,G_2}$ and $\rho_0^{od,G_2}$ are the forms defined
in (i).

{\rm (iii)} generalised $Spin(7)$, odd type
\begin{eqnarray*}
\rho_0^{od} & = &
-\gamma-\varphi+\gamma\wedge\star_7\varphi+\frac{1}{7}\varphi\wedge\star_7\varphi\\
& = & \phantom{-}\gamma\wedge(-1+\star_7\varphi)-\varphi+vol_7
\end{eqnarray*}
where $\star$ and $vol_7$ are taken with respect to the metric $g$
restricted to $\gamma^{\perp}$.
\end{cor}

\begin{rmk}
The 1-form $s\alpha$ of the generalised $G_2$-structure form is the
dual of the vector field induced by $\psi_+$ and the projection of
$\psi_-$ onto the line $\R\widetilde{\psi}_1$. Similarly the 1-form
$-\gamma$ in (iii) is the dual of the vector field $x$ such that
$x\cdot\psi_+=\psi_-$.
\end{rmk}

Since $B$-fields act as isometries on the spin spaces
$\Lambda^{ev,od}$ we obtain the following relation between the
dilaton $f$, the induced Riemannian volume form $vol_g$ and the
``norm" $q(\rho,\Box_{\rho}\rho)$ of the spinor $\rho$.

\begin{cor}\label{normalformcor}\hfill\newline
{\rm (i)} If $\rho=e^{-f}L^{ev}_b(\psi_+\otimes\psi_-)$ defines a
generalised $G_2$-structure, then
$$
q(\rho,\Box_{\rho}\rho)=8e^{-2f}vol_g.
$$
{\rm (ii)} If $\rho=e^{-f}L_b(\psi_+\otimes\psi_-)$ defines a
generalised $Spin(7)$-structure of even or odd type, then
$$
q(\rho,\Box_{\rho}\rho)=16e^{-2f}vol_g.
$$
\end{cor}

\section{Stable forms}\label{stableforms}

Let $T$ be an $n$-dimensional real vector space and consider a form
$\rho$ in $\Lambda^pT^*$ or $\Lambda^{ev,od}T^*$. We regard these
spaces as $GL(n)$- or $\R^{\times}\times Spin(n,n)$-modules under
their natural action respectively and refer to the $GL(n)$-case as
{\em classical} and to the $\R_{>0}\times Spin(n,n)$-case as {\em
generalised}.

\begin{definition}
$\rho$ in $\Lambda^pT^*$ or $\Lambda^{ev,od}T^*$ is said to be {\em
stable} if it lies in an open orbit under the action of $GL(n)$ or
$\R_{>0}\times Spin(n,n)$.
\end{definition}

We refer to a geometrical structure defined by a stable form also as
a {\em variational structure} for reasons that become apparent in
Chapter~\ref{geometry}. As pointed out in~\cite{hi01}, save for the
obvious cases $(n\in\N,p=1)$ and $(n\in2\N,p=2)$ (the orbit of a
{\em symplectic} form), stability is an exceptional phenomenon
restricted to low dimensions since $\dim GL(T)=n^2$ or $\dim
\R_{>0}\times Spin(n,n)=n(n-1)/2+1$ is usually much smaller than
$\dim \Lambda^pT^*=n!/p!(n-p)!$ or $\dim\Lambda^{ev,od}T^*=2^{n-1}$.
A complete classification of representation spaces that admit an
open orbit under the action of a connected reductive complex Lie
group was given in~\cite{saki77}. Note also that if there is an open
orbit under the $GL(T)$ action in $\Lambda^pT^*$ then there is also
one in the dual space
$\Lambda^pT\cong\Lambda^{n-p}T^*\otimes\Lambda^nT$ and consequently
in $\Lambda^{n-p}T^*$~\cite{hi01}.

Let us consider some examples that are relevant to us.

\begin{ex}

(i) Fix a $G_2$-structure over $T=\R^7$ defined by the 3-form
$\varphi$. The dimension of its orbit is
$$
\dim GL(7)-\dim G_2=35=\dim\Lambda^3T^*,
$$
and hence the form is stable in accordance to our earlier use of
this notion.

(ii) Take the case of a spin 3/2 supersymmetric map induced by a
$\mf{su}(3)$-structure on $T=\R^8$. The stabiliser of the 3-form
$\rho$ inside $GL_+(8)$ is conjugate to $PSU(3)$ and so the orbit is
open:
$$
\dim GL(8)-\dim PSU(3)=56=\dim\Lambda^3T^*.
$$
Note, however, that 3-forms inducing an orientation-preserving
isometry are not stable.

(iii) Finally, consider a generalised $G_2$-form
$\rho\in\Lambda^{ev,od}T^*$ over $T=\R^7$. Its stabiliser is
conjugate to $G_2\times G_2$ inside $Spin(7,7)$ and thus
$$
\dim \R_{>0}\times Spin(7,7)-2\dim G_2=64=\dim\Lambda^{ev}T^*.
$$
Hence, if we let the real numbers act by rescaling then we get an
open orbit under the action of $\R_{>0}\times Spin(7,7)$.
\end{ex}

\begin{rmk}
Introducing the additional action of $\R_{>0}$ to obtain an open
orbit seems {\em ad hoc} at first glance. It appears, however,
already in the classical cases, albeit in disguise. In the classical
cases the forms are homogeneous and the scalar matrices in $GL(T)$
act by rescaling.
\end{rmk}

The important feature all of the previous examples share is that the
respective structure groups give rise to a natural metric. There are
other examples -- symplectic forms, a special 3-form on $\R^6$ with
stabiliser $SL(3,\C)$ and a generalised version with stabiliser
$SU(3,3)\leqslant Spin(6,6)$ which lead to non-metrical
structures~\cite{hi00},~\cite{hi02}.

Moreover all these structures we have mentioned (metrical or not)
give rise to a canonically defined volume form and thus fit into a
general formalism set up in~\cite{hi01}. We shall explain the
classical case first. Let $U$ denote the open orbit under the
$GL(T)$ action. Then we can construct a smooth, $GL(T)$-equivariant
map
$$
\phi:U\subset\Lambda^pT^*\to\Lambda^nT^*.
$$
Taking a scalar matrix in $GL(T)$, the equivariance of $\phi$
implies that
$$
\phi(\lambda^p\rho)=\lambda^n\phi(\rho).
$$
Therefore $\phi$ is necessarily homogeneous of degree $n/p$.

Since $\phi$ is smooth we can consider its derivative $D\phi_{\rho}$
at $\rho\in U$ which is an invariant element in
$(\Lambda^pT^*)^*\otimes\Lambda^nT^*$. Now
$(\Lambda^pT^*)^*\cong\Lambda^{n-p}T^*\otimes\Lambda^nT$ so that the
derivative lives in $\Lambda^{n-p}T^*$. It follows that there exists
a unique element $\hat{\rho}\in\Lambda^{n-p}T^*$ for which
$$
D\phi_{\rho}(\dot{\rho})=\hat{\rho}\wedge\dot{\rho}.
$$
We call the form $\hat{\rho}$ the {\em companion} of $\rho$. To
determine what $\hat{\rho}$ is we can take $\dot{\rho}=\rho$ and
apply Euler's formula to get
$$
\hat{\rho}\wedge\rho=\frac{n}{p}\phi(\rho).
$$
In particular, we see that $\hat{\rho}$ is also invariant under
the stabiliser of $\rho$.

\begin{ex}

(i) In the case of a symplectic form $\omega$ over $\R^{2m}$ one can
directly compute the differential $D\phi$ to find
$\hat{\omega}=\omega^{m-1}/(m-1)!$~\cite{hi01}.

(ii) In the $G_2$- and $PSU(3)$-case we have to look at the elements
in $\Lambda^4$ and $\Lambda^5$ which are stabilised by $G_2$ and
$PSU(3)$ -- up to a factor this is just $\star\rho$. As we can
always rescale the volume map $\phi$, we may assume that in these
cases $\hat{\rho}=\star\rho$~\cite{hi01}.
\end{ex}

Next we consider the action of $\R_{>0}\times Spin(n,n)$ on the spin
modules $\Lambda^{ev,od}$. Here, we are left with the cases $n=6$
and $n=7$. The case $n=6$ with stabiliser isomorphic to $SU(3,3)$
leads to a so-called {\em generalised Calabi-Yau} structure and was
dealt with in~\cite{hi02}. We will focus on generalised
$G_2$-structures. Let $U$ be the open orbit $U$ in $\Lambda^{ev}$ or
$\Lambda^{od}$. We define
$$
\phi:\rho\in U\mapsto q(\Box_{\rho}\rho,\rho)\in\Lambda^7T
$$
and since $\Box$ transforms naturally under the action of
$Spin(7,7)$ (cf.~(\ref{boxtrans})), we immediately conclude

\begin{prp}\label{volformg2g2}
$\phi$ is homogeneous of degree 2 and $Spin(7,7)$-invariant.
\end{prp}

\begin{rmk}
The existence of such invariants for complex Lie groups acting with
an open orbit holds in general~\cite{saki77} and follows from purely
algebraic considerations. In this context, an explicit formula of
the invariant for $Spin(14,\C)$ was given in~\cite{gy90}. However,
in view of the variational formalism (see Section
\ref{unc_var_prob}), this description proves to be rather cumbersome
for our purposes which motivated our approach in terms of
$G$-structures.
\end{rmk}

The differential of this map is obtained in the same way as for
classical $G_2$-structures. Since the form $q$ is non-degenerate, we
can write
$$
D\phi_{\rho}(\dot{\rho})=q(\hat{\rho},\dot{\rho}).
$$
for a unique $\hat{\rho}\in\Lambda^{od}$. The $Spin(7,7)$-invariance
of $\phi$ then implies that $\hat{\rho}$ is a $G_2\times G_2$
invariant and by rescaling $\phi$ appropriately we conclude that
$\hat{\rho}=\Box_{\rho}\rho$.


\chapter{Topology}


This chapter deals with global issues of classical and generalised
supersymmetric structures.

First we consider the case of $PSU(3)$-structures and exhibit
obstructions to their existence. $PSU(3)$-manifolds are spinnable
and the key feature we use is the natural identification between the
tangent and the spinor bundle which imposes quite severe constraints
on the topology of the underlying manifold. As a consequence, the
topological constraints we derive hold in the more general situation
where the tangent bundle coincides with the spinor bundles
$\Delta_+$ and $\Delta_-$. We call such a structure {\em doubly
supersymmetric}. The $PSU(3)$-case is special in the sense that it
admits an orthogonal product and in particular, we can establish the
existence of four linearly independent vector fields
(Proposition~\ref{obpsu3}). Basic examples of $PSU(3)$-manifolds are
obtained from topological $SU(3)$-structures ($SU(3)$ acting in its
adjoint representation) through the projection $SU(3)\to PSU(3)$.
Conversely, a $PSU(3)$-structure induces such an $SU(3)$-structure
if a certain cohomology class -- the so-called {\em triality class}
-- vanishes (Section \ref{trialob}). We will give necessary and
sufficient conditions for a $PSU(3)$-structure with vanishing
triality class to exist over a compact manifold
(Theorem~\ref{psu3_structure}).

In the second part of the chapter we investigate the topology of
generalised structures. The two unit spinors of these structures
induce reductions to their respective stabilisers, resulting in two
subbundles associated with $G_2$ or $Spin(7)$ inside the orthonormal
frame bundle. Necessary and sufficient conditions for their
existence easily follow from the classical $G_2$- or $Spin(7)$-case
(Proposition~\ref{gentopob}). Generalised $G_2$-structures are
classified by the top cohomology module $H^7(M,\Z)$
(Theorem~\ref{intersection}). If $M$ is compact, this yields an
integer invariant which essentially counts the number of points
where the two $G_2$-subbundles inside the orthonormal frame bundle
coincide. We also discuss generalised $Spin(7)$-structures, but
these are harder to analyse on dimensional grounds. We conclude with
some examples.

\section{The topology of $PSU(3)$-structures}

\subsection{Obstructions to the existence of $PSU(3)$-structures}\label{obstrucpsu3}

\begin{definition}
A Riemannian 8-dimensional spin manifold is said to be {\em doubly
supersymmetric} if and only if the tangent bundle $T=\Lambda^1$ and
the spinor bundles $\Delta_+$ and $\Delta_-$ are associated with a
principal $G$-fibre bundle such that there exist $G$-invariant
isomorphisms between any two of the three bundles $\Lambda^1$,
$\Delta_+$ or $\Delta_-$, i.e. $T=\Delta_+=\Delta_-$.
\end{definition}

As an example, consider an orientable 8-manifold $M$ which comes
equipped with a stable 3-form $\rho$. This means that at any point
$x\in M$ the form $\rho_x\in\Lambda^3T^*_x$ lies in the open orbit
isomorphic to $GL_+(8)/PSU(3)$. The existence of such a stable form
requires the structure group of $M$ to reduce to $PSU(3)$. As we
have seen in the first chapter, this implies the existence of a
(then globally) defined $PSU(3)$-invariant isometry
$\gamma_{\pm}:T\to\Delta_{\pm}$ which identifies $T$, $\Delta_+$ and
$\Delta_-$. Another example which will be important in the sequel
are generalised $Spin(7)$-structures of {\em odd} type to be defined
in Section~\ref{topges}. As pointed out in Section~\ref{g2ds}, the
existence of unit spinors in $\Delta_+$ and $\Delta_-$ implies a
reduction to $G_2$ for which $T$, $\Delta_+$ and $\Delta_-$ are
isomorphic as $G_2$-representation spaces. In this section we are
primarily interested in the $PSU(3)$-case, but the results hold in
general for any doubly supersymmetric structure.

\begin{definition}
A {\em topological $PSU(3)$-structure} on $M$ is a reduction from
the frame bundle to a principal $PSU(3)$-fibre bundle via the
inclusion $Ad:PSU(3)\hookrightarrow GL(8)$.
\end{definition}

The Lie algebra $\mf{su}(3)$, the symmetric space $SU(3)=SU(3)\times
SU(3)/SU(3)$ and its non-compact symmetric dual $SL(3,\C)/SU(3)$
provide trivial examples of topological $PSU(3)$-structures defined
by the bi-invariant 3-form $ \rho(X,Y,Z)=B(X,[Y,Z])$, where $B$
denotes again the Killing form.

\begin{rmk}
Since $PSU(3)$ is a subgroup of $SO(8)$ and $Spin(8)$ (cf. Corollary
\ref{psu3lift}), every topological $PSU(3)$-structure induces
\begin{itemize}
\item an orientation
\item a metric $g$
\item an associated Hodge $\star$-operator
\item a spin structure $P_{Spin(8)}=P_{PSU(3)}\times Spin(8)$
\item an orthogonal multiplication $\times:T\otimes T\to T$ (cf.~(\ref{orthprod})).
\end{itemize}
\end{rmk}

In this section, we will reserve the letter $P$ for a $G$-principal
fibre bundle which induces a doubly supersymmetric structure. The
existence of a reduction to $G$-principal fibre bundle is
topologically obstructed. Generally speaking, if we have a
$G$-structure over $M$ with classifying map $f:M\to BG$, then there
exists a reduction to a subgroup $i:H\to G$ if and only if there
exists a map $g:M\to BH$ such that the diagram in Figure~\ref{fig4}
commutes.

\begin{figure} [ht]
\begin{center}

\begin{minipage}{11.5cm}\unitlength0.4cm
\begin{picture}(2,5)
\put(11.15,0){$BG$}

\put(18.5,4.4){$EG\times_HG$}

\put(18.5,0){$EG$}

\put(11.3,4.4){$BH$}

\put(4.5,0){$M$}

\put(7.9,3.3){$g$}

\put(7.9,0.6){$f$}

\put(12.3,2.2){$i$}

\put(5.7,0.7){\vector(4,3){5.0}}

\put(5.7,0.25){\vector(1,0){5.1}}

\put(18.2,4.60){\vector(-1,0){5.1}}

\put(18.2,0.25){\vector(-1,0){5.1}}

\put(11.9,4.1){\vector(0,-1){3.1}}

\put(19.0,4.1){\vector(0,-1){3.1}}

\end{picture}
\end{minipage}

\end{center}
\caption{Reductions from $G$ to $H$} \label{fig4}
\end{figure}

Now assume that $C$ is a characteristic class in $H^*(BG,R)$ such
that $i^*C=0$ in $H^*(BH,R)$. Then $f^*C=g^*i^*C$ has to vanish in
$H^*(M,R)$. A classical example is the reduction of $SO(n)$ to
$SO(n-1)$. In this case, the Euler class of $M$ has to be zero.

It is known~\cite{qu71},~\cite{cava94} that the cohomology rings for
$BSpin(8)$ are
$$
H^*(BSpin(8),\Z_2)\cong\Z_2[w_4(E),w_6(E),w_7(E),w_8(E),\epsilon(E)]
$$
and
$$
H^*(BSpin(8),\Z)\cong\Z[q_1(E),q_2(E),e(E),\delta w_6(E)]/\langle
2\delta w_6(E)\rangle
$$
where $\epsilon$ and the so-called {\em spin characteristic classes}
$q_1(E)$ and $q_2(E)$ are defined through the relations
$$
p_1(E)=2q_1(E),\; p_2(E)=q_1^2(E)+2e(E)+4q_2(E)\mbox{ and
}\rho_2(q_2(E))=\epsilon(E).
$$
Here, $\rho_2:H^*(BSpin(8),\Z)\to H^*(BSpin(8),\Z_2)$ is the
morphism induced from the reduction $\mod2$.

In our situation, we also have the outer automorphisms $\kappa$ and
$\lambda$ to play with. Let $E=ESpin(8)$
denote the universal $Spin(8)$-bundle. The automorphisms $\kappa$
and $\lambda$ define further $Spin(8)$-bundles over $BSpin(8)$ by
$E_{\kappa}=E\times_{\kappa}Spin(8)$ and
$E_{\lambda}=E\times_{\lambda}Spin(8)$. They induce two classifying
maps $BSpin(8)\to BSpin(8)$ which we still denote by $\kappa$ and
$\lambda$, that is, we have $E_{\kappa}\cong \kappa^*E$ and
$E_{\lambda}\cong \lambda^*E$. Let $p_i(E)$, $w_i(E)$ and $e(E)$
denote the i-th Pontrjagin, the i-th Stiefel-Whitney and the Euler
class of the universal spin bundle $E$. The corresponding
characteristic classes for $M$ will be denoted by $p_i$, $w_i$ and
$e$ which are the pull-back of the universal classes under the
classifying map $f:M\to BSpin(8)$. Now
$\pi_+=\pi_0\circ\lambda^2\circ\kappa$ and
$\pi_-=\pi_0\circ\lambda^2\circ\kappa\circ\lambda^2$
by~(\ref{trialprin}), so we have
$$
C(\Delta_+)=f^*\circ(\lambda^2\circ\kappa)^*C(E),\quad
C(\Delta_-)=f^*\circ(\lambda^2\circ\kappa\circ\lambda^2)^*C(E),
$$
where $C$ denotes any characteristic class.

\begin{lem}$\!\!\!${\rm~\cite{cava97}}\hspace{2pt}
For $\kappa,\,\lambda:BSpin(8)\to BSpin(8)$, we have
$$
\begin{array}{ll}
\kappa^*(q_2(E))=\phantom{-}q_2(E)+e(E) & \lambda^*(q_2(E))=-e(E)-q_2(E)\\
\kappa^*(e(E))=-e(E) & \lambda^*(e(E))=\phantom{-}q_2(E).
\end{array}
$$
\end{lem}

Since the bundles $T$, $\Delta_+$ and $\Delta_-$ coincide, we derive
$$
e=f^*e(E)=f^*\circ(\lambda^2\circ\kappa)^*e(E)=-f^*q_2(E)=-q_2
$$
and
$$
e=f^*e(E)=f^*\circ(\lambda^2\circ\kappa\circ\lambda^2)^*e(E)=-f^*e(E)=-e.
$$
Since $H^8(M,\Z)$ has no torsion, $e$ and thus $q_2$ must be zero.
Furthermore~\cite{cava94}, $q_2=q_2(\Delta_+)$ is uniquely
determined by the relation
$$
16q_2(\Delta_+)=4p_2(\Delta_+)-p_1^2(\Delta_+)-8e(\Delta_+)
$$
and therefore we obtain as a

\begin{cor}\label{primob}
If $M^8$ is doubly supersymmetric, then
$$
\begin{array}{c}
w_1=w_2=0\\
e=0\\
4p_2=p_1^2.
\end{array}
$$
\end{cor}

Since the degree of the first non-vanishing Stiefel-Whitney class
must always be a power of 2 (cf. 8-B in~\cite{mist74}), we also get
that $w_3=0$.

\begin{rmk}
Corollary~\ref{primob} fails to hold if we merely assume a ``simply"
supersymmetric structure. For example, a $Spin(7)$-structure which
stabilises a unit spinor $\Psi_+\in\Delta_+$ yields a
$Spin(7)$-equivariant isometry $X\in T\mapsto
X\cdot\Psi_+\in\Delta_-$, but this does not imply the vanishing of
the Euler class (see also~\cite{lami89} Subsection IV.10 or
Proposition~\ref{gentopob}).
\end{rmk}

As an application of Corollary~\ref{primob}, we prove the following

\begin{prp}
Let $G$ be a simple compact Lie group and $(G/H,g)$ a Riemannian
homogeneous space which is doubly supersymmetric. If $M=G/H$ has
vanishing Euler class, then $G/H$ is diffeomorphic to $SU(3)$.
\end{prp}

\begin{prf}
Since $G$ sits inside the isometry group of $(M,g)$, its dimension
is less than or equal to $9\cdot 8/2=36$. If we had equality, then
$M$ would be diffeomorphic to a torus or -- up to a finite covering
-- to the 8-sphere (\cite{wo67}, Corollary 11.6.3). While the first
case is excluded for $G$ has to be simple, the second case is
impossible since $e(S^8)\not= 0$. Hence $G$ must be -- up to a
covering -- a group of type $A_1,\ldots,A_5$, $B_2,B_3$, $C_3$,
$D_4$ or $G_2$. As a closed subgroup of $G$, $H$ is compact and
hence reductive. Therefore $H$ is covered by a direct product of
simple Lie groups and a torus, that is the Lie algebra of $H$ is
isomorphic to
$$ \mf{h}=\mf{g}_1\oplus\ldots\oplus\mf{g}_k\oplus\mf{t}^l.$$
If we denote by $rk$ the rank of a Lie group, we get the following
necessary conditions.
$$
\begin{array}{l}
k\le rk(G) \\
l+\sum rk(\mf{g}_i)\le rk(G)\\
l+\sum{dim(\mf{g}_i)}=dim(G)-8,
\end{array}
$$
which yields the possibilities displayed in Table~\ref{tab1}.
\begin{table}[hbt]
\begin{center}
\begin{tabular}{c|c|c|c}
$G$ & $H$ up to a covering & $dim(H)$ & $rk(H)$\\
\hline
$A_2$ & $\{1\}$ & $0$ & $0$\\
$A_3$ & $A_1\times A_1\times S^1$ & $7$ & $3$\\
$A_4$ & $A_3\times S^1,\; G_2\times S^1\times S^1,\; A_2\times A_2$ & 16 & 4\\
$A_5$ & $A_1\times A_4,\; A_1\times A_1\times B_3,\; A_1\times
A_1\times C_3$ & $27$ & $5$\\
$B_2$ & $S^1\times S^1$ & $2$ & $2$\\
$B_3$ & $A_1\times A_2$ & $13$ & $3$\\
$C_3$ & $A_1\times A_2$ & $13$ & $3$\\
$D_4$ & $A_1\times A_1\times G_2$ & $20$ & $4$\\
$G_2$ & $A_1\times A_1$ & $6$ & $2$\\
\end{tabular}
\end{center}
\caption{} \label{tab1}
\end{table}
It follows that $H$ is of maximal rank, that is $rk(H)=rk(G)$,
unless $G=SU(3)$ and $H=\{1\}$. But the first case would imply that
$e(G/H)\not= 0$~\cite{sa58} which is impossible by Corollary
\ref{primob}.
\end{prf}

Since the Euler class is zero, there exists a nowhere vanishing
vector field. In fact we are going to see that if $M$ is closed and
carries a $PSU(3)$-structure, then there exist four linearly
independent vector fields. Let $[M]$ denote the fundamental class of
$M$. We start with the

\begin{lem}
If $M^8$ is a compact spin manifold such that $p_1^2=4p_2$, then
$$
sgn(M)=16\hat{A}[M].
$$
In particular, $sgn(M)\equiv 0\:\mod 16$.
\end{lem}

\begin{prf}
We have $sgn(M)=(7p_2-p_1^2)[M]/45=p_2[M]/15$. On the other hand,
$$
\hat{A}[M]=\frac{1}{2^7\cdot 45}(7p_1^2-4p_2)[M]=\frac{1}{2^4\cdot
15}p_2[M]=\frac{1}{16}sgn(M).
$$
Since $M^8$ is spin, $\hat{A}[M]$ is the index of the Dirac operator
and is therefore an integer. This implies the second assertion.
\end{prf}

If $(b_4^+,b_4^-)$ denotes the signature of the Poincar\'e pairing
on $H^4$, then $sgn(M)=b_4^+-b_4^-$. As a result of this and the
previous lemma we obtain

\begin{cor}
Let $M^8$ be a compact simply-connected doubly supersymmetric
manifold. If $\hat{A}[M]=0$ (e.g. if there exists a metric with
strictly positive scalar curvature), then
$$
1+b_2+b_4^+=b_3.
$$
\end{cor}

For example, $H^*(SU(3),\R)$ is isomorphic to the space of
invariants in $\Lambda^*\mf{su}(3)$ since $SU(3)$ is symmetric.
Therefore, $b_2=0$, $b_3=1$ and $b_4^+=0$ in accordance with the
corollary.

As $e(M)=0$ and $sgn(M)\equiv 0\mod4$, we can assert the existence
of two linearly independent vector fields $X$ and $Y$ (\cite{th69}.
If the doubly-supersymmetric structure is induced by a
$PSU(3)$-structure, then the orthogonal product $\times$
(\ref{orthprod}) produces a third non-vanishing vector field. In
particular, this causes the sixth Stiefel-Whitney class of $M^8$ to
be zero. We are then in a position to apply the following
Proposition

\begin{prp}$\!\!\!${\rm~\cite{cava98a}}\hspace{2pt}
Let $M$ be a closed connected smooth spin manifold of dimension 8.
If $w_6(M)=0$, $e(M)=0$ and $(4p_2(M)-p_1^2(M))[M]\equiv 0 \mod 128$
and there is a $k\in\Z$ such that $4p_2(M)=(2k-1)^2p_1^2(M)$ and
$\:k(k+2)p_2(M)[M]\equiv 0 \mod 3$, then $M$ has 4 linearly
independent vector fields.
\end{prp}

Taking $k=0$ establishes the existence of four linearly independent
vector fields. As a consequence, the fifth Stiefel-Whitney class has
to vanish.

\begin{prp}
We have $w_4^2=0$. In particular, all the Stiefel-Whitney numbers
have to vanish.
\end{prp}

\begin{prf}
Note that by Wu's formula, we get
$$
{\rm
Sq}^k(w_m)=w_kw_m+\binom{k-m}{1}w_{k-1}w_{m+1}+\ldots+\binom{k-m}{k}w_0w_{m+k},
$$
\label{wu}
where
$$
\binom{x}{i}=\frac{x(x-1)\cdot\ldots\cdot(x-i+1)}{i!}.
$$
A further theorem of Wu asserts that
$$
w_k=\sum\limits_{i+j=k}{\rm Sq}^i(v_j),
$$
with elements $v_k\in H^k(M,\Z_2)$ defined through the identity
$v_k\cup x[M]={\rm Sq}^k(x)[M]$ which has to hold for any $x\in
H^{n-k}(M,\Z_2)$. In particular, we have $v_i=0$ for $i>4$. It
follows that $v_1=v_2=v_3=0$, $w_4=v_4$ and $w_8={\rm
Sq}^4w_4=w^2_4=0.$
\end{prf}

We summarise our results in the following theorem.

\begin{thm}\label{obpsu3}
For a topological $PSU(3)$-structure over a closed and oriented
8-manifold $M$ to exist, the following conditions have to hold.
$$
\begin{array}{c}
w_i=0\mbox{ except for }i=4\mbox{ and }
w_4^2=0 \\
e=0\\
p_1^2=4p_2.
\end{array}
$$
In particular, there exist four linearly independent vector fields
on $M$ and all Stiefel-Whitney numbers vanish.
\end{thm}

We close this section by a brief remark about supersymmetric
structures of type 3. As follows from the discussion
in~\ref{relgeo}, their existence is tied to almost quaternionic
structures (i.e. a reduction from $SO(8)$ to $Sp(1)\cdot Sp(2)$).
This problem was settled in~\cite{cava98d}. In particular, the
authors proved the following result.

\begin{prp}$\!\!\!${\rm~\cite{cava98d}}\hspace{2pt}
Let $\xi$ be an oriented 8-dimensional vector bundle over a closed connected smooth spin manifold $M$. If there is $R\in H^4(M,\Z)$ such that the conditions
\begin{itemize}
\item $Sq^2\rho_2R=0$
\item $(Rp_1-2R^2)[M]\equiv 0\mod 16$
\item $w_2(\xi)=0$
\item $w_6(\xi)=0$
\item $4p_2(\xi)-p_1^2(\xi)-8e(\xi)=0$
\item $(p_1^2(\xi)-p_1p_1(\xi)-8e(\xi)+8R^2+4Rp_1(\xi)+4Rp_1)[M]\equiv 0\mod 32$
\end{itemize}
are satisfied, then the structure group of $\xi$ can be reduced to $Sp(1)\cdot Sp(2)$. If $H^2(M,\Z_2)=0$, then all the previous conditions are also necessary.
\end{prp}

If $\xi=T$, the result can be refined to

\begin{prp}$\!\!\!${\rm~\cite{cava98d}}\hspace{2pt} Let $M$ be an oriented closed connected smooth manifold of dimension 8. If
\begin{itemize}
\item $w_2=0$
\item $w_6=0$
\item $4p_2-p_1^2-8e=0$
\end{itemize}
and there exists an $R\in H^4(M,\Z)$ such that
\begin{itemize}
\item $Sq^2\rho_2R=0$
\item $(Rp_1-2R^2)[M]\equiv0\mod 16$
\item $(R^2+Rp_1-e)[M]\equiv 0\mod 4$,
\end{itemize}
then $M$ carries an almost quaternionic K\"ahler structure. The first and third condition are always necessary while the remaining ones are necessary if $H^2(M,\Z_2)=0$.
\end{prp}

These propositions give similar necessary conditions as Theorem~\ref{obpsu3}. Sufficient conditions for $PSU(3)$-structures shall occupy us next.

\subsection{$PSU(3)$-structures and the triality
class}\label{trialob}

Let $M^8$ be again a connected, closed and spinnable 8-manifold.
Consider a complex rank 3 vector bundle $E$ with vanishing first
Chern class, i.e. $E$ is associated with a principal $SU(3)$-fibre
bundle. If the adjoint bundle $\mf{su}(E)=P_{SU(3)}\times\mf{su}(3)$
is isomorphic to the tangent bundle, then $T$ is associated with a
$PSU(3)$-structure coming from the projection $p:SU(3)\to PSU(3)$.
However, not every $PSU(3)$-structure arises in this way. A basic
\v{C}ech cohomology argument implies that principal $G$-fibre
bundles over $M$ are classified by $H^1(M,G)$ (see, for
example,~\cite{lami89} Appendix A). The exact sequence (where $\Z_3$
is central)
$$
1\to\Z_3\to SU(3)\stackrel{p}{\to} PSU(3)\to 1
$$
gives rise to an exact sequence
$$
\ldots\to H^1(M,\Z_3)\to {\rm Prin}_{SU(3)}(M) \stackrel{p_*}{\to}
{\rm Prin}_{PSU(3)}(M)\stackrel{t}{\to} H^2(M,\Z_3).
$$
(where ${\rm Prin}_G(M)$ denotes the set of $G$-principal bundles
over $M$). Hence, a principal $PSU(3)$-bundle $P$ is induced by an
$SU(3)$-bundle if and only if the obstruction class $t(P)\in
H^2(M,\Z_3)$ vanishes. Following~\cite{avis79} where the authors
consider $PSU(3)$-structures over 4-manifolds, we call this class
the {\em triality class}. By the universal coefficients theorem this
obstruction vanishes trivially if $H^2(M,\Z)=0$ and $H^3(M,\Z)$ has
no torsion elements of order divisible by three.

If $f:M\to BPSU(3)$ is a classifying map for $P$, then $t(P)=f^*t$
for the {\em universal triality class} $t\in H^2(BPSU(3),\Z_3)$. It
is induced by $c_1(E_{U(3)})$, the first Chern class of the
universal $U(3)$-bundle $E_{U(3)}$~\cite{wo82a}. Let
$\overline{p}:U(3)\to PU(3)$ be the natural projection. Note that
the inclusion $SU(3)\subset U(3)$ induces an isomorphism between
$PSU(3)$ and $PU(3)$ and therefore identifies $BPSU(3)$ with
$BPU(3)$. Since $BPU(3)$ is simply connected and
$\pi_2(BU(3))=\Z\to\pi_2(BPU(3))=\Z_3$ is the reduction $\mod3$ map
$\rho_3:Z\to\Z_3$, the Hurewicz isomorphism theorem and the
universal coefficients theorem imply that $H^2(BPU(3),\Z_3)=\Z_3$
and that
$$
B\overline{p}^*:H^2(BPU(3),\Z_3)\to H^2(BU(3),\Z_3)
$$
is an isomorphism. Then
$$
t=(B\overline{p}^*)^{-1}\rho_{3*}c_1(E_{U(3)}).
$$

If the triality class vanishes the problem of finding sufficient
conditions for the existence of a $PSU(3)$-structure reduces to the
problem for the existence of a complex rank 3 vector bundle $E$ with
$\mf{su}(E)=T$. This question is settled in the next theorem.

\begin{thm}\label{psu3_structure}
Suppose that $M$ is a connected, closed and spinnable 8-manifold.
Then the frame bundle reduces to a principal $PSU(3)$-fibre bundle
$P$ with $t(P)=0$ if and only if $e=0$, $4p_2=p_1^2$, $w_6=0$, $p_1$
is divisible by 6 and $p_1^2[M]\in 216\Z$.
\end{thm}

\begin{prf}
Let us first assume that there exists a $PSU(3)$-fibre bundle $P$
coming from a reduction $Ad:SU(3)\to SO(8)$ with principal bundle
$\widetilde{P}$. In virtue of Theorem~\ref{obpsu3} we only need to
show the last two conditions. We define the complex vector bundle
$E=\widetilde{P}\times_{SU(3)}\C^3$ so that $\mf{su}(E)=T$, and
compute the Pontrjagin classes of $M$. We have
$\mathfrak{su}(3)\otimes\C=\mathfrak{sl}(3,\C)$, hence $T\otimes\C$
equals $\End_0(E)$, the bundle of trace-free complex endomorphisms.
The Chern character of $T\otimes\C$ equals (see for
instance~\cite{sa82})
$$
ch(T\otimes\C)=8+p_1+\frac{1}{12}(p_1^2-2p_2).
$$
On the other hand,
$$
ch(End(E))=ch(E\otimes\overline{E})=ch(End_0(E))+1.
$$
Now for a complex vector bundle with $c_1(E)=0$,
$$
ch(E)=3-c_2(E)+\frac{1}{2}c_3(E)+\frac{1}{12}c_2(E)^2
$$
and $c_i(E)=(-1)^ic_i(\overline{E})$ which implies
$$
ch(E\otimes\overline{E})=ch(E)\cup
ch(\overline{E})=9-6c_2(E)+\frac{3}{2}c_2(E)^2.
$$
As a consequence
\begin{equation}\label{pontrjagin}
p_1=p_1(\mf{su}(E)) = -6c_2(E),\quad p_2=p_2(\mf{su}(E))= 9c_2(E)^2.
\end{equation}
In particular, $p_1$ is divisible by 6 and we also rederive the
relation $4p_2=p_1^2$. Moreover $M$ is spinnable, hence the spin
index $\hat{A}\cup{\rm ch}(E)[M]$ is an integer. Since
$$
\hat{A}=1-p_1/24+(-4p_2+7p_1^2)/5760=1-p_1/24+p_1^2/960
$$
it follows
$$
\hat{A}\cup
ch(E)[M]=3\hat{A}[M]+p_1c_2(E)/24+c_2(E)^2/12=3\hat{A}[M]-p_1^2/216[M]\in\Z,
$$
which means $p_1^2[M]\in 216\Z$, proving the necessity of the
conditions.

\medskip

For the proof of the converse I am indebted to ideas of
M.~Crabb~\cite{cr02}. Let $B\subset M$ be an embedded open disc in
$M$ and consider the exact sequence
$$
K(M,M-B)\to K(M)\to K(M-B).
$$
We have $K(M,M-B)=\widetilde{K}(S^8)\cong\Z$ and the sequence is
split by the spin index
$$
x\in K(M)\mapsto\hat{A}\cup{\rm ch}(x)[M]\in\Z
$$
which therefore classifies the stable extensions over $M-B$ to $M$.
The first step consists in finding a stable complex vector bundle
$\xi$ over $M-B$ such that $c_1(\xi)=0$ and $\mf{su}(3)$ is stably
equivalent to $T_{|M-B}$. To that end, let
$[(M-B)_+,BSU(\infty)]\subset K(M-B)$ denote the set of pointed
homotopy classes, the subscript $+$ indicating a disjoint basepoint.
Let $(c_2,c_3)$ be the map which takes an equivalence class of
$[(M-B)_+,BSU(\infty)]$ to the second and third Chern class of the
associated bundle.

\begin{lem}\label{steenrod}
The image of the mapping
$$
(c_2,c_3):[(M-B)_+,BSU(\infty)]\to H^4(M,\Z)\oplus H^6(M,\Z)
$$
is the set $\{(u,v)\:|\:Sq^2\rho_2u=\rho_2v\}$.
\end{lem}

\begin{lemprf}
We first prove that for a complex vector bundle $\xi$ with
$c_1(\xi)=0$, we have
\begin{equation}\label{square}
Sq^2\rho_2c_2(\xi)=\rho_2c_3(\xi).
\end{equation}
Now if $W_i$ denote the Stiefel-Whitney classes of the {\em real}
vector bundle underlying $\xi$ this is equivalent to $Sq^2W_4=W_6$.
On the other hand, Wu's formula implies
$$
Sq^2W_4=W_2W_4+W_6
$$
and thus (\ref{square}) since $W_2=\rho_2c_1=0$. Next let
$i:F\hookrightarrow K(\Z,4)\times K(\Z,6)$ denote the homotopy fibre
of the induced map
$$
Sq^2\circ\rho_2+\rho_2:K(\Z,4)\times K(\Z,6)\to K(\Z_2,6).
$$
The relation (\ref{square}) implies that the map
$(c_2,c_3):BSU(\infty)\to K(\Z,4)\times K(\Z,6)$ is null-homotopic.
Consequently, $(c_2,c_3)$ lifts to a map $k:BSU(\infty)\to F$,
thereby inducing an isomorphism of homotopy groups
$\pi_i(BSU(\infty))\to\pi_i(F)$ for $i\leq 7$ and a surjection for
$i=8$. By the exact homotopy sequence for fibrations we conclude
that $\pi_4(F)=\Z$, $\pi_6(F)=2\Z$ and $\pi_i(F)=0$ for $i$
otherwise. On the other hand, the Chern class
$c_2:\widetilde{K}(S^4)=\Z\to H^4(S^4,\Z)=\Z$ is an isomorphism and
$c_3:\widetilde{K}(S^6)=\Z\to H^6(S^6,\Z)=\Z$ is multiplication
by~2. Since $M-B$ is at most 8-dimensional, it follows that the
induced map $k_*:[(M-B)_+,BSU(\infty)]\to[(M-B)_+,F]$ is surjective.
The horizontal row in the commutative diagram

\begin{center}

\begin{minipage}{11.5cm}

\unitlength0.4cm

\begin{picture}(2,5)
\put(11.15,0){$[(M-B)_+,BSU(\infty)]$}

\put(25.7,4.4){$H^6(M,\Z_2)$}

\put(20.4,5){$Sq^2\rho_2+\rho_2$}

\put(11.3,4.4){$H^4(M,\Z)\oplus H^6(M,\Z)$}

\put(0,4.4){$[(M-B)_+,F]$}

\put(7.9,5){$i_*$}

\put(16.2,2.2){$(c_2,c_3)$}

\put(7.6,2.2){$k_*$}

\put(10.8,0.2){\vector(-2,1){7.5}}

\put(5.7,4.65){\vector(1,0){5.2}}

\put(20.0,4.65){\vector(1,0){5.2}}

\put(15.50,0.8){\vector(0,1){3.1}}

\end{picture}

\end{minipage}

\end{center}

\vspace{0.5cm}

is exact, hence ${\rm im}\,(c_2,c_3)={\rm im}\,i_*=\ker\,
(Sq^2\circ\rho_2+\rho_2)$.
\end{lemprf}

By assumption $p_1\in H^4(M,\Z)$ is divisible by $6$ and therefore
we can write $p_1=-6u$ for $u\in H^4(M,\Z)$. On the other hand,
$p_1=2q_1$, where $q_1$ is the first spin characteristic class and
satisfies $\rho_2(q_1)=w_4$. Hence
$$
Sq^2\rho_2(u)=Sq^2w_4=w_2w_4+w_6=0,
$$
and the previous lemma implies the existence of a stable complex
vector bundle $\xi$ such that $c_1(\xi)=0$, $c_2(\xi)=u$ and
$c_3(\xi)=0$. From (\ref{pontrjagin}) it follows that
$p_1(\mf{su}(\xi))=p_1$ and since $w_2(\mf{su}(\xi))=0$,
$\mf{su}(\xi)$ and $T$ are stably equivalent over the 4-skeleton
$M^{(4)}$~\cite{wo82b}. Then $\mf{su}(\xi)$ and $T$ are stably
equivalent over $M-B$ as the restriction map $KO(M-B)\to
KO(M^{(4)})$ is injective. This follows from the exact sequence
$$
KO(M^{(i+1)},M^{(i)})\to KO(M^{(i+1)})\to KO(M^{(i)}).
$$
By definition
$KO(M^{(i+1)},M^{(i)})=\widetilde{KO}(M^{(i+1)}/M^{(i)})$ and
$M^{(i+1)}/M^{(i)}$ is a disjoint union of spheres $S^{i+1}$. But
$\widetilde{KO}(S^{i+1})=0$ for $i=4,\,5$ and $6$ and therefore the
map $KO(M^{(i+1)})\to KO(M^{(i)})$ is injective. Since $M=M^{(8)}$
is the disjoint union of $M^{(7)}$ and a finite number of open
embedded discs, the assertion follows. Next we extend $\xi$ over $B$
to a stable bundle on $M$. The condition to be represented by a
complex vector bundle $E$ of rank 3 is $c_4(\xi)=0$. As pointed out
above, such a bundle exists if
$$
\hat{A}\cup{\rm ch}(\xi)[M]=3\hat{A}[M]+p_1u/24+u^2/12
$$
is an integer and this holds by assumption. Next
$p_2(\mf{su}(\xi))=9u^2=p_2$ and as a consequence, $\mf{su}(\xi)$ is
stably isomorphic to $T$~\cite{wo82b}. Finally, two stably
isomorphic oriented real vector bundles of rank 8 are isomorphic as
$SO(8)$-bundles if they have the same Euler class. Since
$e(\mf{su}(E))=0$, we have $T=\mf{su}(E)$.
\end{prf}

\begin{cor}
If $M$ is closed and carries a $PSU(3)$-structure with vanishing
triality class, then
$$
\hat{A}[M]\in40\Z\mbox{ and }sgn(M)\in 640\Z.
$$
\end{cor}

\section{The topology of generalised exceptional structures}\label{topges}

First we introduce the (rather obvious) global notion of a
generalised metric.

\begin{definition}
A {\em generalised metric} over a manifold $M^n$ is an orthogonal
decomposition
$$
T\oplus T^*=V_+\oplus V_-
$$
into a positive and negative subbundle with respect to the natural
inner product $(\cdot\,,\cdot)$.
\end{definition}

As we have seen in Section~\ref{ges}, a generalised metric structure
corresponds to a unique pair $(g,b)$ of a Riemannian metric $g$ and
a 2-form $b$. In particular, any Riemannian manifold defines
trivially a generalised metric structure whose existence is
therefore unobstructed. The additional choice of an orientation then
gives rise to an $SO(n)\times SO(n)$ structure and in particular to
a globally defined $\Box$-operator.

Next, there is the natural notion of a generalised exceptional
structure.

\begin{definition}\hfill\newline
{\rm (i)} A {\em topological generalised $G_2$-structure} over a
7-manifold $M^7$ is defined by an even or odd form $\rho$ whose
stabiliser under the action of $Spin(7,7)$ is conjugate to
$G_2\times G_2$ at any point. We will denote this structure by the
pair $(M^7,\rho)$ and call $\rho$ the {\em structure form}.

{\rm (ii)} A {\em topological generalised $Spin(7)$-structure of
even or odd type} over an 8-manifold $M^8$ is defined by an even or
odd form $\rho$ whose stabiliser under the action of $Spin(8,8)$ is
conjugate to $Spin(7)\times Spin(7)$ at any point. We will denote
this structure by the pair $(M^8,\rho)$ and call $\rho$ the {\em
structure form}.

We refer to these structures also as {\em topological generalised
exceptional structures}.
\end{definition}

We will usually drop the adjective ``topological" and simply refer
to a generalised $G_2$- and $Spin(7)$-structures if there is no risk
of confusion. Propositions~\ref{linalg} and~\ref{lmap} assert in
particular the existence of a metric $g$ or equivalently an $SO(7)$-
or $SO(8)$-principal fibre bundle which admits a reduction to a
$G_2$- or $Spin(7)$. The inclusions $G_2\subset Spin(7)$ and
$Spin(7)\subset Spin(8)$ determine a unique spin structure for which
we can consider the associated spinor bundles $\Delta$ (over $M^7$)
or $\Delta_+$ and $\Delta_-$ (over $M^8$). From
Propositions~\ref{linalg},~\ref{lmap} and~\ref{selfdual} we deduce
the following statement.

\begin{thm}\label{topology}\hfill\newline
{\rm (i)} A topological generalised $G_2$-structure $(M^7,\rho)$ is
characterised by the following data,
\begin{itemize}
\item an orientation,
\item a metric $g$,
\item a 2-form $b$,
\item two unit spinors $\Psi_+,\Psi_-\in\Delta$ and a function $F$ such that
$$
e^{-F}L_b(\Psi_+\otimes\Psi_-)=\rho+\Box_{g,b}\rho.
$$
\end{itemize}

{\rm (ii)} An even or odd topological generalised
$Spin(7)$-structure $(M^8,\rho)$ for $\rho\in\Lambda^{ev}$ (even
type) or $\rho\in\Lambda^{od}$ (odd type) is characterised by the
following data,
\begin{itemize}
\item an orientation,
\item a metric $g$,
\item a 2-form $b$,
\item two unit half spinors $\Psi_+,\Psi_-\in\Delta$ of either equal (even type) or opposite chirality (odd type) and a function $F$ such
that
$$
e^{-F}L_b(\Psi_+\otimes\Psi_-)=\rho.
$$
\end{itemize}
\end{thm}

\begin{rmk}
In the case of a generalised $Spin(7)$-structure we fix the
orientation in such a way that $\Psi_+$ lies in $\Delta_+$ unless
otherwise stated.
\end{rmk}

Thus, a generalised exceptional structure may be regarded as a
Riemannian structure whose orthonormal frame bundle admits two
$G_2$- or $Spin(7)$-subbundles associated with the stabilisers of
$\Psi_{\pm}$ in $Spin(7)$ or $Spin(8)$. The $B$-field $b$
encapsulates the additional information in which copy of
$Spin(7)\times Spin(7)$ or $Spin(8)\times Spin(8)$ inside
$Spin(7,7)$ or $Spin(8,8)$ the stabiliser of $\rho$ is actually
sitting.

The most trivial example of a generalised structure would be a
Riemannian spin manifold which admits a nowhere vanishing spinor
$\Psi=\Psi_+=\Psi_-$ and where we put $b=0,\,F=0$. A generalised
structure defined by a single $G_2$- or $Spin(7)$-structure, possibly with a non-vanishing B-field $b$ and dilaton $F$ is said to be {\em straight} (in the
generalised $Spin(7)$-case this obviously makes only sense for
structures of even type). The existence of a nowhere vanishing
spinor field in dimension 7 or 8 is a classical result
(\cite{lami89} Subsection IV.10).

\begin{prp}\label{gentopob}\hfill\newline
{\rm (i)} A spinnable 7-fold $M^7$ always carries a nowhere
vanishing spinor and hence admits a topological $G_2$-structure.

{\rm (ii)} A differentiable 8-fold $M^8$ carries a nowhere vanishing
spinor if and only if $M^8$ is spin and for an appropriate choice of
orientation satisfies
$$
p_1(M)^2-4p_2(M)+8\chi(M)=0.
$$
\end{prp}

Conversely, any generalised structure gives rise to a Riemannian
spin manifold with a nowhere vanishing spinor and consequently, the
topological obstructions in Proposition~\ref{gentopob} are also
necessary. Since the spinor $\Psi_+$ induces a supersymmetric map
$$
X\in T\stackrel{\cong}{\mapsto}X\cdot\Psi_+\in\Delta_-
$$
we see that for a generalised $Spin(7)$-structure of odd type we
also need a nowhere vanishing vector field $X$ defined by
$X\cdot\Psi_+=\Psi_-$ which requires the vanishing of the Euler
class $\chi(M)$.

\begin{cor}\hfill\newline
{\rm (i)} A 7-fold $M$ carries a topological generalised
$G_2$-structure if and only if $M$ is spin.

{\rm (ii)} An 8-fold $M$ carries an even topological generalised
$Spin(7)$-structure if and only if $M$ is spin and
\begin{equation}\label{obgenspin71}
8\chi(M)+p_1(M)^2-4p_2(M)=0.
\end{equation}

{\rm (iii)} A differentiable 8-fold $M$ carries an odd topological
generalised $Spin(7)$-structure if and only if $M$ is spin and
\begin{equation}\label{obgenspin72}
\chi(M)=0,\quad p_1(M)^2-4p_2(M)=0.
\end{equation}
\end{cor}

\begin{ex}
The 7-sphere is spinnable and therefore admits a generalised
$G_2$-structure. The tangent bundle of the 8-sphere is stably
trivial and therefore all the Pontrjagin classes vanish. Since the
Euler class is non-trivial, there exists no generalised
$Spin(7)$-structure on an 8-sphere. However,
equations~(\ref{obgenspin71}) and (\ref{obgenspin72}) are
automatically satisfied for manifolds of the form $M^8=S^1\times
N^7$ with $N^7$ spinnable.
\end{ex}

An ordinary $G_2$- or $Spin(7)$-structure sitting inside the
orthonormal frame bundle $P$ induces trivially a generalised one
which settles the existence question. But do there exist genuine
generalised structures where the two $G_2$- or $Spin(7)$-subbundles
inside $P$ cannot be homotopically transformed into each other? And
if yes, can we classify them? In the case of a generalised
$Spin(7)$-structure the first question obviously only makes sense
for structures of even type and for the remainder of this section,
we shall only consider generalised $Spin(7)$-structures where the
two induced spinors live in $\Delta_+$.

To render our question more precise, we regard $G_2$- or (even)
$Spin(7)$-structures as being defined by a (continuous) section of
the sphere bundle $p_{\mathbb{S}}:\mathbb{S}\to M$ of $\Delta$. On
the space of sections $\Gamma(\mathbb{S})$ we then introduce the
following equivalence relation. Two spinors $\Psi_0$ and $\Psi_1$ in
$\Delta$ are considered to be equivalent (denoted
$\Psi_0\sim\Psi_1$) if and only if we can deform $\Psi_0$ into
$\Psi_1$ through sections which means that there exists a continuous
map $F:M\times I\to \mathbb{S}$ such that $F(x,0)=\Psi_0(x)$,
$F(x,1)=\Psi_1(x)$ and $p_{\mathbb{S}}\circ F(x,t)=x$. An
equivalence class will be denoted by $[\Psi]$. If two sections are
vertically homotopic, then the two corresponding $G_2$- or
$Spin(7)$-structures are isomorphic as principal $G_2$- or
$Spin(7)$-bundles over $M$. In particular the generalised structure
defined by the pair $(\Psi,\Psi_0)$ is equivalent to a straight
structure if and only if $\Psi\sim\Psi_0$. The set of generalised
structures with fixed $\Psi_+=\Psi$ is just
$Gen(M)=\Gamma(\mathbb{S})/\sim$ which is what we aim to determine.

The question if whether or not two sections are vertically homotopic
can be tackled by using obstruction theory (for a more detailed
account of obstruction theory than we need it here, see~\cite{st51}
or~\cite{wh78}). Assume more generally that we are given a fibre
bundle over a not necessarily compact $n$-fold $M^n$, and two
sections $s_1$ and $s_2$ which are vertically homotopy equivalent
over the $q$-skeleton $M^{[q]}$ of $M$. For technical simplicity, we
assume the fibre to be connected. The obstruction for extending the
vertical homotopy to the $q+1$-skeleton lies in
$H^{q+1}(M,\pi_{q+1}(F))$. In particular, there is the first
non-trivial obstruction $\delta(s_1,s_2)\in H^m(M,\pi_m(F))$ called
the {\em primary difference} of $s_1$ and $s_2$ for the least
integer $m$ such that $\pi_m(F)\not=0$. It is a homotopy invariant
of the sections $s_1$ and $s_2$ and enjoys the additivity property
\begin{equation}\label{additivity}
\delta(s_1,s_2)+\delta(s_2,s_3)=\delta(s_1,s_3).
\end{equation}
Coming back first to the generalised $G_2$-case we consider the
sphere bundle $\mathbb{S}$ over $M^7$ with fibre $S^7$.
Consequently, the primary difference between two sections lies in
$H^7(M,\Z)$ and is the only obstruction for two sections to be
vertically homotopy equivalent. Moreover, the additivity property
implies that $\delta(\Psi,\Psi_1)=\delta(\Psi,\Psi_2)$ if and only
if $\delta(\Psi_1,\Psi_2)=0$, that is, $\Psi_1\sim\Psi_2$, and for
any class $d\in H^7(M,\Z)$, there exists a section $\Psi_d$ such
that $d=\delta(\Psi,\Psi_d)$~\cite{st51},~\cite{wh78}. As a
consequence, we obtain the

\begin{prp}
The set of generalised $G_2$-structures can be identified with
$$
Gen(M)=H^7(M,\Z)=\left\{\begin{array}{ll} \Z, & \mbox{ if $M$ is
compact}\\ 0, & \mbox{ if $M$ is non-compact}\end{array}\right..
$$
\end{prp}

Generalised $G_2$-structures are therefore classified by an integer
invariant which over a compact $M^7$ has the natural interpretation
as the number of points (counted with an appropriate sign
convention) where the two $G_2$-structures coincide. To see this we
associate with every equivalence class $[\Phi]$ the intersection
class $\#(\Psi(M),\Phi(M))\in H_{14}(\mathbb{S},\Z)$ of the
7-dimensional oriented submanifolds $\Psi(M)$ and $\Phi(M)$ inside
$\mathbb{S}$. Since the total space of the sphere bundle is
14-dimensional, $\#(\Psi,\Phi)$ counts the number of points in $M$
where the two spinors $\Psi$ and $\Phi$ coincide. Taking the cup
product of the Poincar\'e duals of $\Psi(M)$ and $\Phi(M)$ then sets
up a map
$$
[\Phi]\in Gen(M)\mapsto PD(\Psi(M))\cup PD(\Phi(M))\in
H^{14}(\mathbb{S},\Z).
$$
On the other hand, the Gysin sequence implies that integration along
the fibre defines an isomorphism
$\pi_{\mathbb{S}*}:H^{14}(\mathbb{S},\Z)\to H^7(M,\Z)$ and
therefore, any generalised $G_2$-structure induced by the
equivalence class $[\Phi]$ over a compact 7-fold $M$ gives rise to a
well-defined cohomology class $d(\Psi,\Phi)\in H^7(M,\Z)$. The
following theorem shows this class to coincide with the primary
difference (for the proof I benefited from discussions with
W.~Sutherland and M.~Crabb).

\begin{thm}\label{intersection}
We have
$$
d(\Psi,\Phi)=\delta(\Psi,\Phi).
$$
In particular generalised $G_2$-structures are classified by the number of points where the two underlying $G_2$-structures coincide.
\end{thm}

\begin{prf}
We regard the spinor bundle $\Delta$ as an 8-dimensional oriented
real vector bundle over M and consider the two sections $\Psi$ and
$\Phi$ of the sphere bundle. The primary difference
$\delta(\Psi,\Phi)$ can be represented geometrically by the zero-set
(a finite set with signs) of the section $(m,x)\mapsto(1 - x)\Psi(m)
+ x\Phi(m)$ of the pullback of $\Delta$ to $M\times\R$, and deformed
to be transverse to the zero-section. In particular, if $\Psi$ and
$-\Phi$ never coincide then the primary difference is 0. Therefore
the intersection number, defined geometrically by making $\Psi(M)$
and $\Phi(M)$ transverse and taking the coincidence set, will be
$\delta(\Psi,-\Phi)$ (with appropriate sign conventions). By virtue
of~(\ref{additivity}), we have $\delta(\Psi,-\Phi) =
\delta(\Psi,\Phi)+\delta(\Phi,-\Phi)$. The difference class
$\delta(\Phi,-\Phi)$ corresponds to the self-intersection number
$\#(\Phi(M),\Phi(M))$ which is 0 since $M$ is 7-dimensional. It
follows that we can identify the intersection class $d(\Psi,\Phi)$
with the primary difference $\delta(\Psi,\Phi)$.
\end{prf}

The index we have just defined is a topological obstruction to
closed strong integrability as we will explain in the next chapter
(see Definition~\ref{intstruc} and Corollaries~\ref{vanishingg2}
and~\ref{vanishing}). Outside the singularity set, the structure
group reduces to $SU(3)$, the stabiliser of two orthogonal spinors.
There we can express the structure form $\rho$ in terms of
$SU(3)$-invariants, but as can be seen from the normal form given in
Proposition~\ref{normalform}, this description breaks down at a
singular point where $\sin(a)=0$, that is, both spinors are
parallel.

A generalised structure is said to be {\em exotic} if it is defined
by two inequivalent spinors. Here are two examples.

\begin{ex}\hfill\newline
(i) Choose a $G_2$-structure over $S^7$ represented by the spinor
$\Psi$. Since the tangent bundle of $S^7$ is trivial, so is the sphere bundle of $\Delta$, i.e. $\mathbb{S}=S^7\times S^7$. Consequently,
$$
Gen(S^7)=\Gamma(\mathbb{S})/\!\sim\;=[S^7,S^7]=\pi_7(S^7)=\Z
$$
and any map $\Phi:S^7\to S^7$ which is not homotopic to a constant gives rise to an exotic
generalised $G_2$-structure for $\Psi\equiv const$.

(ii) The tangent bundle of the manifold $M^8=S^1\times S^7$ is also
trivial. As observed above, $M$ carries a unit spinor and we can
again trivialise the sphere bundle $\mathbb{S}$, so that
$Gen(M)=[S^1\times S^7,S^7]$, containing the set
$[S^7,S^7]=\pi_7(S^7)=\Z$. Choosing a non-trivial homotopy class in
$\pi_7(S^7)$ which we extend trivially to $S^1\times S^7$ defines an
exotic generalised $Spin(7)$-structure.
\end{ex}

In the $Spin(7)$-case the sphere bundle is 15-dimensional with fibre
isomorphic to $S^7$ and an 8-dimensional base, so that two
transverse sections will intersect in a curve. We meet the first
obstruction for the existence of a vertical homotopy in $H^7(M,\Z)$
which by Poincar\'e duality trivially vanishes if $H_1(M,\Z)=0$
(e.g. if $M$ is simply connected). The second obstruction lies in
the top cohomology module. Since $\pi_8(S^7)=\Z_2$, we obtain

\begin{prp}
Let $M^8$ be an 8-manifold with $H_1(M,\Z)=0$ which admits a
generalised $Spin(7)$-structure. Then the set of generalised
$Spin(7)$-structures with a fixed $\Psi_+=\Psi$ is given by
$$
Gen(M)=H^8(M,\pi_8(S^7))\cong\left\{\begin{array}{ll} \Z_2, & \mbox{if $M$ compact}\\
0, & \mbox{if $M$ non-compact}\end{array}\right..
$$
\end{prp}

Again we see that over a non-compact $M^8$ with $H_1(M,\Z)=0$ any two sections of
the sphere bundle can be vertically transformed into each other
while this is clearly not the case for the general compact case as
shown by the example above.

\begin{rmk}
The stable homotopy group $\pi_{n+k}(S^n)$ is isomorphic to the
framed cobordism group of $k$-manifolds~\cite{sw02}. It is
conceivable that the $\Z_2$-class is the framed (or spin) cobordism
class of the 1-manifold where the two sections coincide.
\end{rmk}


\chapter{Geometry}\label{geometry}


The motivation for adopting the supersymmetric approach was to
provide a framework within which we can understand the algebra of
special metrics associated with stable forms. Stability also
provides us with a natural set of integrability conditions arising
from Hitchin's variational principle. These integrability conditions
define the geometry, following the general philosophy of this
thesis, through closed forms. Based on the papers~\cite{hi01}
and~\cite{hi02}, we will formulate in the first section various
kinds of this principle. It leads to the notions of integrable
$PSU(3)$- and to (closed) strongly or weakly integrable generalised
$G_2$-structures. The variational principle is perfectly general and
in particular, no metric is needed to set up these equations.

Starting on a Riemannian manifold however, we can adopt the
supersymmetric point of view and reformulate these integrability
conditions in terms of the corresponding spinorial invariants,
providing thus a {\em spinorial solution} to the variational
problem. As a result, we obtain some well-known equations of
mathematical physics, namely the Rarita-Schwinger equation
(Theorem~\ref{raritaschwinger}) and the supersymmetry equations of
supergravity of type~IIA/B with bosonic background fields
(Theorem~\ref{integrability}). We will also analyse related
geometries as introduced in Chapter~\ref{linearalgebra} where we
impose the integrability conditions from the variational or the
supersymmetric ansatz. Moreover, there is also a practical value to
the use of spinors as they are easier to manipulate than forms. In
particular, we will use these spinor field equations to derive
integrability conditions on the Ricci tensor (Sections
\ref{geomexclass} and~\ref{geompropgen}). A further important
consequence of the spinorial formulation is a no-go theorem for
closed strongly integrable generalised $G_2$- and
$Spin(7)$-structures in the sense that any solution over a compact
manifold is induced by two classical $G_2$- or $Spin(7)$-structures
with a parallel spinor (Corollary~\ref{vanishingg2}). In this sense,
the (unconstrained) variational principle leads to classical
geometries. However, we shall see that the notion of weak
integrability potentially yields an interesting new type of geometry
(cf. the discussion in Section~\ref{exgeneral}) for which
non-trivial compact examples could exist. The spinorial picture is
less useful for the construction of explicit examples and therefore,
we step back to work with forms again. A particularly striking
advantage of forms over spinors in this context is the device of
T-duality which easily yields non-trivial generalised solutions from
classical ones (Section~\ref{exgeneral}).

\section{The variational problem}\label{unc_var_prob}

\subsection{The unconstrained variational problem in the classical
case}

We recall the set-up of the variational problem as given
in~\cite{hi01}. Let $(M^n,\rho)$ be a closed and oriented manifold
which is endowed with a {\em stable} form $\rho$ of pure degree $p$.
This means that $\rho$ is a section of the open subset of
$\Omega^p(M)$ whose intersection with any fiber $\Lambda^pT_x^*M$
can be identified with an open orbit $U$ (cf. Section
\ref{stableforms}). This provides us with a natural volume form
$\phi(\rho)$ which asks to be integrated
$$
V(\rho)=\int\limits_M\phi(\rho).
$$
Since stability is a generic condition, we can differentiate the
volume functional and we shall consider its variation over a fixed
cohomology class in $H^p(M,\R)$.

\begin{thm}\label{variation}$\!\!\!${\rm~\cite{hi01}}\hspace{2pt}
A closed stable form $\rho\in\Omega^p(M)$ is a critical point in its
cohomology class if and only if $d\hat{\rho}=0$.
\end{thm}

\begin{prf}
The first variation of $V$ is
$$
\delta
V_{\rho}(\dot{\rho})=\int\limits_MDV_{\rho}(\dot{\rho})=\int\limits_M\hat{\rho}\wedge\dot{\rho}.
$$
Since we vary over a fixed cohomology class, $\dot{\rho}$ is exact,
i.e. $\dot{\rho}=d\alpha$. As a consequence of Stokes' theorem, this
variation vanishes for all $\alpha\in\Omega^{p-1}(M)$ if and only if
$$
\int\limits_M\hat{\rho}\wedge
d\alpha=\pm\int\limits_Md\hat{\rho}\wedge\alpha=0
$$
which holds precisely when $\hat{\rho}$ is closed.
\end{prf}

\begin{rmk}
Note that the condition $d\hat{\rho}=0$ is non-linear in $\rho$ as
the $\wedge$-operation also depends on $\rho$. In this way, we can
see the variational problem as performing a non-linear version of
Hodge theory~\cite{hi01}.
\end{rmk}

\subsection{The unconstrained variational problem in the generalised
case}\label{uvpgen}

Here we follow the approach given in~\cite{hi02} for generalised
Calabi-Yau manifolds. Assume $M^7$ to be closed, oriented and
provided with a stable form $\rho$ giving rise to a generalised
$G_2$-structure. Since the structure group $GL_+(7)$ lifts to
$Spin(7,7)$, $\rho$ is a section of a vector bundle associated with
$GL_+(7)$ and with typical fiber isomorphic (as a {\em vector
space}) to $\Lambda^{ev,od}\R^{7*}$. This becomes the spin
representation for $GL_+(7)\subset Spin(7,7)$ if we twist with the
square root $\Omega^7(M)^{1/2}$ (compare with the remark on
Page~\pageref{glvs.u}). In this subsection, we shall, for sake of
clarity, denote by $\lambda^{ev,od}$ the spin representation of
$GL_+(7)\subset Spin(7,7)$, while $GL_+(7)$ acts on
$\Lambda^{ev,od}$ as usual by $A^*\rho$.

Let us consider the volume functional $\phi$ as described in
Proposition~\ref{volformg2g2}. It is homogeneous of degree 2 and can
therefore be considered from a representation theoretic point of
view as a $GL(7)_+$-equivariant function
$$
\phi:U\subset\Lambda^{ev,od}\to\Lambda^7,
$$
since for an element $A\in GL_+(7)$ we have
$$
\phi(A^*\rho)=(\det A)^{-1}\cdot\phi(\sqrt{\det A}\cdot
A^*\rho)=(\det A)^{-1}\phi(A\bullet\rho)=(\det A)^{-1}\phi(\rho).
$$
The non-degenerate form $q$ now takes values in $\Lambda^7$ so that
$$
D\phi_{\rho}(\dot{\rho})=q(\hat{\rho},\dot{\rho})
$$
is also a volume form. As in the classical case stable forms are
sections of the open subset in $\Omega^{ev,od}(M)$ and $\phi$ is a
well-defined function on this set. We obtain the following analogue
of Theorem~\ref{variation}.

\begin{thm}\label{genvariation}
A closed stable form $\rho\in\Omega^{ev,od}(M)$ is a critical
point in its cohomology class if and only if $d\hat{\rho}=0$.
\end{thm}

\begin{prf}
The argument is a direct translation of the proof of Theorem
\ref{variation}. Take the first variation of $V(\rho)$ over $[\rho]$
$$
\delta
V(d\alpha)=\int_MD\phi(d\alpha)=\int_Mq(\hat{\rho},d\alpha)=\int_M\hat{\rho}\wedge\sigma(d\alpha).
$$
Moreover, it is straightforward to check that
$$
d\sigma(\tau^{ev,od})=\mp\sigma(d\tau^{ev,od})
$$
for any even or odd form $\tau\in\Omega^{ev,od}(M^7)$. Thus
$$
\delta V(\dot{\rho})=\int_Mq(\hat{\rho},d\alpha)=\pm\int_Mq(
d\hat{\rho},\alpha)
$$
and the variation vanishes for all $\alpha$ if and only if
$\hat{\rho}$ is closed.
\end{prf}

\subsection{The constrained variational problem}\label{cvpgen}

Next we generalise the constrained variational problem for classical
$G_2$-manifolds (there is no $PSU(3)$-analogue) described
in~\cite{hi01} to our setup. On a closed oriented manifold $M^7$,
the non-degenerate pairing $q$ induces a non-degenerate pairing
between the spaces of forms $\Omega^{ev}(M)$ and $\Omega^{od}(M)$ by
$$
\int_Mq(\sigma,\tau).
$$
If $\sigma=d\gamma$ is exact, then Stokes' theorem implies that
$$
\int_Mq(d\gamma,\tau)=\int_Mq(\gamma,d\tau),
$$
and this vanishes for all $\gamma$ if and only if $\tau$ is closed.
Hence we obtain a non-degenerate pairing which allows us to identify
$$
\Omega^{ev}_{\rm
exact}(M^7)^*\cong\Omega^{od}(M^7)/\Omega^{od}_{\rm closed}(M^7).
$$
The exterior differential maps the latter space isomorphically
onto $\Omega^{ev}_{\rm exact}(M^7)$ so that
$$
\Omega^{ev}_{\rm exact}(M^7)^*\cong\Omega^{ev}_{\rm exact}(M^7).
$$
Consequently, the non-degenerate pairing becomes a non-degenerate
quadratic form on $\Omega^{ev}_{exact}(M^7)$ given by
$$
Q(d\gamma)=\int_Mq(\gamma,d\gamma).
$$

The same conclusion holds for odd instead of even forms.

\begin{thm}\hfill\newline
A stable form $\rho\in\Omega^{ev,od}_{\rm exact}(M)$ is a critical
point subject to the constraint $Q(\rho)=const$ if and only if there
exists a real constant $\lambda$ with $d\hat{\rho}=\lambda\rho$.
\end{thm}

\begin{prf}
From the proof of Theorem~\ref{genvariation}, the first variation of
$V$ at $\rho=d\gamma$ is
$$
(\delta V)_{\rho}(d\alpha)=\int\limits_Mq(\hat{\rho},d\alpha)
$$
whereas the differential of $Q$ at $\rho$ is
$$
(\delta Q)_{\rho}(d\alpha)=2\int\limits_Mq(\alpha,\rho).
$$
By Lagrange's theorem, we see that for a critical point we have
$d\hat{\rho}=\lambda\rho$.
\end{prf}

\subsection{The twisted variational problem}\label{tvpgen}

To produce the most general form $e^{b/2}\Psi_+\otimes\Psi_-$ of a
spinor with stabiliser $G_2\times G_2$, we first considered the
``unperturbed" spinor $\rho_0=\Psi_+\otimes\Psi_-$ and then hit it
with a B-field. The equations for $\rho_0$ to define a critical
point are
$$
d\rho_0=0,\quad d\hat{\rho}_0=0
$$
which under the action of a $B$-field are equivalent to the
inhomogeneous set of equations
$$
d\rho_0+\frac{1}{2}db\wedge\rho_0=0,\quad
d\hat{\rho_0}+\frac{1}{2}db\wedge\hat{\rho_0}=0.
$$
Following~\cite{hi02} we can, more generally, consider the twisted
differential operator
$$
d_H\rho=d\rho+H\wedge\rho,
$$
where $H$ is now a closed (but not necessarily exact) 3-form. The
closedness guarantees that $d_H$ still defines a differential
complex. Hence it makes sense to speak of a $d_H$-cohomology class
and we can set up the variational problem to take place over such a
cohomology class. The following theorem is a direct translation from
Theorem 13 in~\cite{hi02}.

\begin{thm}
A $d_H$-closed stable form $\rho\in\Omega^{ev,od}(M)$ is a critical
point of $V(\rho)$ in its $d_H$-cohomology class if and only if
$d_H\hat{\rho}=0$.
\end{thm}

\begin{prf}
This time, the first variation is
$$
(\delta
V)_{\rho}(d\alpha)=\int\limits_Mq(\hat{\rho},d\alpha+H\wedge\alpha)=\pm\int\limits_Mq(
d\hat{\rho},\alpha)+\hat{\rho}\wedge\sigma(\alpha)\wedge\sigma(H),
$$
where $\rho$ and $\alpha$ are both either even or odd. Since
$$
\hat{\rho}\wedge\sigma(H\wedge\alpha)=\hat{\rho}\wedge\sigma(\alpha)\wedge\sigma(H)=H\wedge\hat{\rho}\wedge\sigma(\alpha)
$$
the first variation is
$$
(\delta V)_{\rho}(d\alpha)=\int\limits_Mq(
d\hat{\rho}+H\wedge\hat{\rho},\alpha).
$$
This vanishes for all $\alpha$ if and only if $d_H\hat{\rho}=0$.
\end{prf}

\section{Integrable variational structures and related geometries in the classical
case}\label{spinorial}

\subsection{The spinorial solution of the variational problem and related
geometries}\label{class_sol}

Next we will consider the geometries introduced in Section
\ref{susy} together with the following integrability conditions. For
geometries of type~$1_0$, $2_0$ or $3_0$, that is, the reduction of
the structure group of the tangent bundle to $PSU(3)$, $SO(3)\times
SO(3)\times SO(2)$ or $SO(3)\times SO(5)$ (cf. Section~\ref{relgeo})
is induced by an honest 3-form $\rho\in\Omega^3(M)$, we impose the
integrability condition coming from the variational principle
$d\rho=\dstar\rho=0$. If the geometry is defined by a supersymmetric
map $\gamma_{\pm}:\Lambda^1\to\Delta_{\pm}$, then this geometry is
said to be {\em integrable} if and only if $\gamma_{\pm}$ is
harmonic with respect to the twisted Dirac operator on
$T^*M\otimes\Delta$. In the case of $PSU(3)$, we will show these
conditions to coincide. The latter condition therefore provides what
we call a spinorial solution to the variational problem. We will
analyse this case in some detail before we briefly investigate the
related geometries of type~2 and 3.

\bigskip

\textbf{$PSU(3)$-geometry}

The formal solution to the variational problem in the classical case
suggests the following definition.

\begin{definition}
A $PSU(3)$-manifold $(M,\rho)$ is called {\em integrable} if and
only if
\begin{equation}\label{psu3intcond}
d\rho=0,\quad \dstar\rho=0.
\end{equation}
\end{definition}

Recall that in~\ref{susy} we proved that as a consequence of the
triality principle the $PSU(3)$-invariant 3-form $\rho$ gives rise
to two isometries $\gamma_{\pm}:\Lambda^1\to\Delta_{\pm}$. In this
section we want to derive an equivalent formulation of integrability
in terms of
$\gamma=\gamma_+\oplus\gamma_-\in\Lambda^1\otimes(\Delta_+\oplus\Delta_-)$.
This means that we start with a spinnable (not necessarily compact)
Riemannian manifold $(M^8,g,\gamma)$ whose metric is induced by the
$PSU(3)$-structure associated with $\gamma$, and ask for the
conditions on $\gamma$ which makes this structure integrable.

To this end we have to analyse the action of the covariant
derivative on $\rho$ and $\gamma$, but first a general outline of
the strategy is in order. The following discussion will also prove
useful for the generalised case so we start with a generic
Riemannian manifold $(M^n,g)$. Consider the Levi-Civita connection
form
$$
Z:TP_{SO(n)}\to\mf{so}(n),
$$
which is an $\mf{so}(n)$-valued 1-form on the orthonormal frame
bundle $P_{SO(n)}$. If we are given a reduction from this
$SO(n)$-bundle to a subbundle $P_G$ acted on by the subgroup
$G\leqslant SO(n)$, then the Lie algebra $\mf{so}(n)$ splits into
the subalgebra $\mf{g}$ and an orthogonal complement $\mf{m}$
$$
\mf{so}(n)=\mf{g}\oplus\mf{m}.
$$
We can decompose the restriction of the Levi-Civita form $Z$ to
$P_G$ accordingly into
$$
Z_{|TP_G}=\widetilde{Z}\oplus\Gamma
$$
where $\widetilde{Z}$ is a $\mf{g}$-valued connection form on $P_G$.
It defines a connection for the $G$-fibre bundle and in particular,
any $G$-invariant object will be parallel for the induced affine
connection $\widetilde{\nabla}$. Note however that this property is
in general not sufficient to guarantee the uniqueness of
$\widetilde{\nabla}$ (see, for instance Remark 2 in Section 4.2
of~\cite{br03}). $\Gamma$ is a tensorial 1-form of type $Ad$, i.e.
it can be seen as a form on $M$ taking values in the associated
bundle
$$
P_G\times_{\mathop{\rm Ad}}\mf{m}.
$$
At $x\in M^n$ we can interpret $\Gamma_x$ as an element of
$\R^n\otimes\mf{m}$ and in this sense $\Gamma$ measures the failure
of the Levi-Civita connection to reduce to $P_G$. It is also
referred to as the {\em torsion} of the structure. The decomposition
of the $G$-module $\R^n\otimes\mf{m}$ into irreducible modules gives
rise to the different geometrical {\em types} of $G$-structures
according to the non-trivial components of the torsion. Knowing the
type allows us to extract further geometrical information. For
instance, if we consider a $G_2$-structure on a 7-manifold, then
$\mf{m}=\Lambda^2_7$ is isomorphic to the vector representation of
$G_2$. Consequently, $\R^7\otimes\mf{m}$ decomposes into
$$
\R^7\otimes\mf{m}=\Lambda^0_1\oplus\Lambda^1_7\oplus\Lambda^2_{14}\oplus\Lambda^3_{27}.
$$
A $G_2$-structure whose associated $\Gamma$ takes values in say
$\Lambda^0_1=\R$ is then said to be of {\em type} $\Lambda^0_1$
(these $G_2$-structures are also called {\em nearly parallel
$G_2$-manifolds}, see~\cite{fkms97}). In particular, the metric
structure it defines is Einstein. Generally speaking, the type can
be detected by looking at the action of $Z$ on a $G$-invariant
object. For structures of type $\Lambda^0_1$ the covariant
derivative of the 3-form $\varphi$ and the corresponding spinor
$\Psi$ is given by
$$
\nabla_X\varphi=-2\lambda X\llcorner\star\varphi\mbox{ and }
\nabla_X\Psi=\lambda X\cdot\Psi.
$$
Here, as in the sequel, the covariant derivative operator $\nabla$
will {\em always} denote the Levi-Civita connection. We will revisit
this example later on when we consider generalised $G_2$-structures.

Coming back to the mainstream of the development, we write the
covariant derivative of the $PSU(3)$-invariant 3-form as
$$
\nabla\rho=T(\rho)
$$
with $T\in\Lambda^1\otimes\Lambda^2$. $T$ acts through its
$\Lambda^2$ factor which we view as the Lie algebra of $SO(8)$.
Since $\mf{su}(3)\leqslant\Lambda^2$ stabilises $\rho$, we can
assume that $T\in\Lambda^1\otimes\mf{su}(3)^{\perp}$ so that $T$
really encodes the action of the Lie algebra valued 1-form $\Gamma$
above. Consequently, we refer to $T$ as the {\em torsion} of the
$PSU(3)$-structure. Unless the space is symmetric or the holonomy
reduces to a proper subgroup of $PSU(3)$, inspection of Berger's
list~\cite{be55} implies that the torsion does not vanish. In this
case we say that the integrable $PSU(3)$-structure is {\em
non-trivial}. To make contact with the variational principle, we
remark that $d\rho$ and $d^{\star}\rho$ are the image of $T$ under
the $PSU(3)$-equivariant morphisms
$$
\tilde{d}:\Lambda^1\otimes\mf{su}(3)^{\perp}\to\Lambda^4
$$
and
$$
\widetilde{d^{\star}}:\Lambda^1\otimes\mf{su}(3)^{\perp}\to\Lambda^2
$$
which map the tensor $A(\rho)\in\Lambda^1\otimes\Lambda^3$ to
$\Lambda^4$ and $\Lambda^2$ by skew-symmetrisation and contraction
respectively. These last two operations also map $\nabla\rho$ to
$d\rho$ and $-d^{\star}\rho$ and therefore,
$$
\tilde{d}T=d\rho,\quad\widetilde{d^{\star}}T=-d^{\star}\rho.
$$
The integrability conditions~(\ref{psu3intcond}) are thus equivalent to
$$
T\in\ker\tilde{d}\cap\ker\widetilde{d^{\star}}.
$$
The next step consists in decomposing
$\Lambda^1\otimes\mf{su}(3)^{\perp}$ into irreducible modules. Let $\wedge:\Lambda^1\otimes\mf{su}(3)^{\perp}\to\Lambda^3$ denote the natural skewing map. Then $\Lambda^1\otimes\mf{su}(3)^{\perp}\cong\ker\wedge\oplus\rho^{\perp}$,
where $\rho^{\perp}=[1,1]\oplus\llbracket3,0\rrbracket\oplus[2,2]$ is the orthogonal complement of $\rho$ in $\Lambda^3$. Moreover, the natural contraction map $\llcorner:\ker\wedge\subset\Lambda^1\otimes\mf{su}(3)^{\perp}\to\Lambda^1$ splits $\ker\wedge$ into a direct sum isomorphic to $\ker\,\llcorner\oplus\Lambda^1$ where $\ker\,\llcorner\cong[2,2]\oplus\llbracket4,1\rrbracket$. On the other hand we obtain from Proposition~\ref{psu3decomp}
$$
\Lambda^2=[1,1]\oplus\llbracket3,0\rrbracket,\quad\Lambda^4=2[1,1]\oplus2[2,2].
$$
Schur's Lemma implies at once
$$
\llbracket4,1\rrbracket\leqslant\ker\tilde{d}\cap\ker\widetilde{d^{\star}}.
$$
In order to determine the kernel completely, it will be convenient
to complexify the modules. We find
$$
\begin{array}{ccccl}
& &\Lambda^1\otimes\Lambda^2_{10+} & = & (1,1)_+\oplus(0,3)\oplus(2,2)_+\oplus(1,4)\\
\Lambda^1\otimes\mf{su}(3)^{\perp}\otimes\C & = & \oplus & & \\
& &\Lambda^1\otimes\Lambda^2_{10-} & = & (1,1)_-\oplus(3,0)\oplus(2,2)_-\oplus(4,1).
\end{array}
$$
The modules $(1,1)_{\pm}$ and $(2,2)_{\pm}$ have non-trivial projections to both $\ker\wedge$ and $\rho^{\perp}$. In particular, they map non-trivially under $\wedge$.

\begin{prp}\label{fegr}
The $PSU(3)$-structure defined by $\rho$ is integrable if and only
if
$$
T\in\llbracket4,1\rrbracket.
$$
\end{prp}

\begin{rmk}
The main difficulty of the proof resides in the fact that the modules $[1,1]$ and $[2,2]$ appear twice in $\Lambda^4$. Even if both components in $\Lambda^1\otimes\mf{su}(3)^{\perp}$ map non-trivially under $\widetilde{d}$, we can still have a component of the kernel isomorphic to $[1,1]$ or $[2,2]$ which trivially intersects our fixed copies of $[1,1]$ or $[2,2]$. Therefore, the argument presented in~\cite{hi01} is inconclusive.
\end{rmk}

\begin{prf}
We only have to check that $d\rho=0$ and $d^*\rho=0$ implies $T\in\llbracket4,1\rrbracket$. Let us determine the kernel of $\widetilde{d^*}$ and recall that $a_{\pm}^*\rho=\pm\sqrt{3}i\star(a_{\pm}\wedge\rho)$ for any $a_{\pm}\in\Lambda^2_{10\pm}$ (Proposition~\ref{psu3decomp}). By complexifying, it follows that restricted to the $PSU(3)$-invariant modules $\Lambda^1\otimes\Lambda^2_{10\pm}$,
\begin{equation}\label{dstar_van}
\widetilde{d^*}(\sum e_j\otimes a_j^{\pm})=\mp\sqrt{3}i\sum e_j\llcorner\star(a^{\pm}_j\wedge\rho)=\pm\sqrt{3}i\sum \star(e_j\wedge a^{\pm}_j\wedge\rho).
\end{equation}
In virtue of the remarks preceding the theorem, the kernel of the skewing map $\Lambda^1\otimes\Lambda^2_{\pm}$ is isomorphic to $(1,4)$ and $(4,1)$, so this vanishes outside these modules if and only if $\sum e_i\wedge a_i^{\pm}$ lies in $\widetilde{1}\oplus[2,2]$, the kernel of the map which wedges 3-forms with $\rho$. Invoking Schur's Lemma, $\ker\widetilde{d^*}\cong[1,1]\oplus2[2,2]\oplus\llbracket4,1\rrbracket$, 
where the precise embedding of $[1,1]$ will be of no importance to us.

Next we consider the operator $\widetilde{d}$. If we can show that it is surjective, then $\ker\widetilde{d}\cong\llbracket3,0\rrbracket\oplus\llbracket4,1\rrbracket$ and consequently, the kernels of $\widetilde{d^*}$ and $\widetilde{d}$ intersect in $\llbracket4,1\rrbracket$. Let $\iota_{\rho^{\perp}}$ denote the injection of $\rho^{\perp}$ into $\Lambda^1\otimes\Lambda^2_{20}$ obtained by projecting the natural embedding of $\Lambda^3$ into $\Lambda^1\otimes\Lambda^2$. We first prove the relation
\begin{equation}\label{surjd}
b_3(\alpha)=\frac{1}{2}\widetilde{d}(\iota_{\rho^{\perp}}(\alpha)),\quad \alpha\in\rho^{\perp}\subset\Lambda^3
\end{equation}
which shows that $\ker b_4\subset\im\widetilde{d}$. By~(\ref{b_decomp}), the kernel of $b_3$ is isomorphic to $\mathbf{1}\oplus\llbracket1,2\rrbracket$, so the claim needs only to be checked for the module $[1,1]\oplus[2,2]$ in $\Lambda^3$. A sample vector is obtained by 
\begin{equation}\label{sample}
p_3(e_{128})=\alpha_8\oplus\alpha_{27}={\textstyle\frac{1}{8}}(5e_{128}+\sqrt{3}e_{345}+\sqrt{3}e_{367}-2e_{458}+2e_{678}),
\end{equation}
where $p_3=b^*_4b_3$. That both components $\alpha_8$ and $\alpha_{27}$ are non-trivial can be seen as follows. Restricting $p_3$ to $\Lambda^3_8$ and $\Lambda^3_{27}$ is multiplication by real scalars $x_1$ and $x_2$ since the modules are representations of real type. If one, say $x_1$, vanished, then $p^2_3(e_{128})=p_3(\alpha_{27})=x_2\cdot\alpha_{27}$. However
$$
p^2_3(e_{128})={\textstyle\frac{1}{64}}(39e_{128}+7\sqrt{3}e_{345}+7\sqrt{3}e_{367}-18e_{458}+18e_{678})
$$
which is not a multiple of~(\ref{sample}). Moreover, we have indeed 
\begin{eqnarray*}
b_3p_3(e_{128}) & = & {\textstyle\frac{1}{32}}(7\sqrt{3}e_{1245}+7\sqrt{3}e_{1267}-9e_{1468}-9e_{1578}+9e_{2478}-9e_{2568})\\
& = & {\textstyle\frac{1}{2}}\widetilde{d}(\iota_{\rho^{\perp}}p_3(e_{128}))
\end{eqnarray*} 
which proves~(\ref{surjd}). For the inclusion $\im b^*_5\subset\im\widetilde{d}$ we consider the vector $e_1\otimes e_{18}$ in $\ker\wedge$. Then $\widetilde{d}(e_1\otimes e_{18})=-e_{1238}/2-e_{1478}/4+e_{1568}/4$ takes values in both components of $\im b^*_5\subset\Lambda^4$ since
$b^*_4\widetilde{d}(e_1\otimes e_{18})=0$ and otherwise 
$$
b^*_5b_4\widetilde{d}(e_1\otimes e_{18})={\textstyle\frac{1}{32}}(-10e_{1238}-5e_{1478}+5e_{1568}+3e_{2468}+3e_{2578}+3e_{3458}-3e_{3678})
$$ 
would be a multiple of $\widetilde{d}(e_1\otimes e_{18})$. Hence $\widetilde{d}$ is surjective and the assertion follows.
\end{prf}

Let $\mc{D}$ denote the twisted Dirac operator on
$\Lambda^1\otimes\Delta$, so $\mc{D}$ can be locally written as
$$
\mc{D}(X\otimes\Psi)=\sum\limits_{i=1}^8\nabla_{e_i}X\otimes
e_i\cdot\Psi+X\otimes e_i\cdot\nabla_{e_i}\Psi.
$$
Consider the embedding of $\Delta$ into
$\Lambda^1\otimes\Delta=\Delta\oplus\Lambda^3\Delta_+\oplus\Lambda^3\Delta_-$
given by
$$
i(\psi)(X)=-\frac{1}{8}X\cdot\psi
$$
and the projection $p: \Lambda^1\otimes\Delta\to\Delta$
$$
p(\gamma)=\sum\limits_{i=1}^8e_i\cdot\gamma(e_i).
$$
With respect to the decomposition ${\rm Im}\,i\oplus\ker p$, the
Dirac operator $\mc{D}$ takes the form (see Prop. 2.7
in~\cite{wa91})
\begin{equation}\label{dodecomp}
\left(\begin{array}{cc} -\frac{3}{4}i\circ D\circ i^{-1} &
2i\circ\delta\\\frac{1}{4}P\circ i^{-1} & Q\end{array}\right),
\end{equation}
where $D:\Delta\to\Delta$ is just the usual Dirac-Operator,
$\delta:\Lambda^1\otimes\Delta\to\Delta$ the twisted co-differential, $P:\Delta\to
\Lambda^3\Delta_+\oplus\Lambda^3\Delta_-$ the Twistor operator
$P(\sigma)=\nabla\sigma-i\circ p(\nabla\sigma)$, and $Q$ the
so-called {\em Rarita-Schwinger} operator which arises in
supergravity and string theories. To introduce more physicist's
jargon still, we make the

\begin{definition}
We call a spin 3/2 field $\psi$, i.e. a section
$\psi\in\Gamma(\Lambda^3\Delta_-\oplus\Lambda^3\Delta_-)$, a {\em
Rarita-Schwinger field} if it satisfies
$$
Q(\psi)=0.
$$
\end{definition}

\begin{thm}\label{raritaschwinger}
The $PSU(3)$-structure is integrable if and only if
$\gamma=\gamma_+\oplus\gamma_-\in\Gamma(\Lambda^1\otimes\Delta)$ is
harmonic with respect to the twisted Dirac-operator $\mc{D}:\Gamma(\Lambda^1\otimes\Delta)\to\Gamma(\Lambda^1\otimes\Delta)$,
that is $\mc{D}\gamma=0$.
\end{thm}

\begin{prf}
The proof uses the same idea as Proposition~\ref{fegr}. Let us
consider the $PSU(3)$-equivariant morphism
$$
\widetilde{\mc{D}}:A\in\Lambda^1\otimes\mf{su}(3)^{\perp}\mapsto\mu(A(\gamma))\in\Lambda^1\otimes\Delta
$$
which maps $A(\gamma)\in\Lambda^1\otimes\Lambda^1\otimes\Delta$
through Clifford multiplication to $\mu(A(\gamma))$ in
$\Lambda^1\otimes\Delta$. 
Since the $PSU(3)$-invariant supersymmetric map $\gamma=\gamma_+\oplus\gamma_-$ has now components in both $\Delta_-\otimes\Lambda^1$ and $\Delta_+\otimes\Lambda^1$, we will split the Dirac operator $\widetilde{\mc{D}}=\widetilde{\mc{D}}_+\oplus\widetilde{\mc{D}}_-$ accordingly, i.e. $\widetilde{\mc{D}}_{\pm}(X\otimes a)=\mu_{\pm}(X\cdot a(\gamma_{\pm}))\in\Delta_{\pm}\otimes\Lambda^1$. We now have to show that
$$
\ker\widetilde{\mc{D}}_+\cap\ker\widetilde{\mc{D}}_-=\llbracket4,1\rrbracket.
$$
The intersection $\ker\widetilde{\mc{D}}_+\cap\ker\widetilde{\mc{D}}_-$ contains at least the module $\llbracket4,1\rrbracket$. First we show that $\llbracket3,0\rrbracket$ is not contained in this intersection by taking the vector 
\begin{eqnarray*}
\tau_{\llbracket3,0\rrbracket} & = & id\otimes\pi^2_{20}(\iota_{\rho^{\perp}}b_2(4e_{18}))\\
& = & -\sqrt{3}e_1\otimes e_{45}-\sqrt{3}e_1\otimes e_{67}+2e_2\otimes e_{38}-2e_3\otimes e_{28}+
\sqrt{3}e_4\otimes e_{15}+\\
& & e_4\otimes e_{78}-\sqrt{3}e_5\otimes e_{14}-e_5\otimes e_{68}+\sqrt{3}e_6\otimes e_{17}+e_6\otimes e_{58}-\sqrt{3}e_7\otimes e_{16}-\\
& & e_7\otimes e_{48}+2e_8\otimes e_{23}+e_8\otimes e_{47}-e_8\otimes e_{56}.
\end{eqnarray*}
A straightforward, if tedious, computation shows $\widetilde{\mc{D}}(\tau_{\llbracket3,0\rrbracket})\not=0$. For the remainder of the proof, it will again be convenient to complexify the torsion module $\Lambda^1\otimes\mf{su}(3)^{\perp}$ and to consider $(1,1)_{\pm}$ and $(2,2)_{\pm}$. The invariant 3-form $\rho$ induces equivariant maps $\rho_{\mp}:\Delta_{\pm}\to\Delta_{\mp}$ whose matrices with respect to the choices made in~(\ref{spin8_rep}) are given by~(\ref{arho1}) for $\rho_+$ and by its transpose for $\rho_-$. By Schur's Lemma, we have 
\begin{equation}\label{constant}
\widetilde{\mc{D}}_-((2,2)_+)=z\cdot\rho_-\otimes id\circ\widetilde{\mc{D}}_+((2,2)_+)
\end{equation}
for a complex scalar $z$. Since the operators $\widetilde{\mc{D}}_{\pm}$ are real and $(2,2)_-$ is the complex conjugate of $(2,2)_+$, the same relation holds for $(2,2)_-$ with $\bar{z}$. The vector $\tau_0=6(e_1\otimes e_{18}-e_2\otimes e_{28})$ is clearly in $\ker\,\llcorner\subset\ker\wedge$ and projecting the second factor to $\Lambda^2_{10+}$ yields 
\begin{eqnarray*}
id\otimes\pi^2_{10+}(\tau_0) & = & e_1\otimes(3e_{18}+i\sqrt{3}e_{23}-i\sqrt{3}e_{47}+i\sqrt{3}e_{56})+\\
& & e_2\otimes(i\sqrt{3}e_{13}-3e_{28}+i\sqrt{3}e_{46}+i\sqrt{3}e_{57}).
\end{eqnarray*}
Since any possible component in $(1,4)$ gets killed under $\widetilde{\mc{D}}$, we can plug this into~(\ref{constant}) to find $z=(1+i\sqrt{3})/8$ which shows that $(2,2)_{\pm}$ map non-trivially under $\widetilde{\mc{D}}$. On dimensional grounds, $\ker\widetilde{\mc{D}}_{\pm}$ therefore contains the module $(2,2)$ with multiplicity one. Their intersection, however, is trivial, for suppose otherwise. Let $(2,2)_0$ denote the corresponding copy in $\ker\widetilde{\mc{D}}_+$.
It is the graph of an isomorphism $M:(2,2)_+\to(2,2)_-$ since it intersects $(2,2)_{\pm}$ trivially. Now if $\tau=\tau_+\oplus M\tau_+\in(2,2)_0$ were in $\ker\widetilde{\mc{D}}_-$, then
\begin{eqnarray*}
\widetilde{\mc{D}}_-(\tau_+\oplus M\tau_+) & = & z\cdot\rho\otimes id\circ\widetilde{\mc{D}}_+(\tau_+)\oplus \bar{z}\cdot\rho\otimes id\circ\widetilde{\mc{D}}_+(M\tau_+)\\
& = & \rho\otimes id\circ\widetilde{\mc{D}}_+(z\cdot \tau_+\oplus \bar{z}\cdot M\tau_+)\\
& = & 0.
\end{eqnarray*}
Consequently, $z\cdot \tau_+\oplus \bar{z}\cdot M\tau_+\in\ker\widetilde{\mc{D}}_+$, that is, $\bar{z}\cdot M\tau_+=Mz\cdot \tau_+$ or $\bar{z}=z$ which is a contradiction. This shows that the kernels of $\widetilde{\mc{D}}_{\pm}$ intersect at most in $2(1,1)\oplus\llbracket2,3\rrbracket$ and furthermore, that the condition $\mc{D}(\gamma_+)=0$ or $\mc{D}(\gamma_-)=0$ on its own is not sufficient to guarantee the close- and cocloseness of $\rho$. This argument also applies to $(1,1)_{\pm}$. However, since $(1,1)$ appears twice in $\Delta_{\pm}\otimes\Lambda$, we first need to project onto $\Delta_{\mp}\cong(1,1)$ via Clifford multiplication before asserting the existence of a complex scalar $z$ such that
$$
\mu_+\circ\widetilde{\mc{D}}_-((1,1)_+)=z\cdot\rho_+(\mu_-\circ\widetilde{\mc{D}}_+((1,1)_+)).
$$
For the computation of $z$, we can use the vector 
$$
2\sqrt{3}ie_1\otimes\pi^2_{10+}(e_{18})= e_1\otimes(\sqrt{3}ie_{18}-e_{23}+e_{47}- e_{56})\in(1,1)_+\oplus(2,2)_+\oplus(1,4),
$$
as possible non-trivial components in $(2,2)_+\oplus(1,4)$ get killed under $\mu_{\mp}$. We find $z=2(1-\sqrt{3}i)$ which shows that $(1,1)$ occurs with multiplicity at most one in $\ker\widetilde{\mc{D}}_{\pm}$ and that it is not contained in their intersection. Consequently, $\ker\widetilde{\mc{D}}_+\cap\ker\widetilde{\mc{D}}_-=\llbracket4,1\rrbracket$, which proves the theorem in conjunction with Proposition~\ref{fegr}.
\end{prf}

\bigskip

\textbf{$SO(3)\times SO(3)\times SO(2)$- and $SO(3)\times
SO(5)$-geometry}

Next we briefly investigate geometries of type~2 and 3. Recall that
we divided these structures into the subtypes $0$, $+$ and $-$
according to the case where the structure group $SO(\Lambda^1)$,
$SO(\Delta_-)$ and $SO(\Delta_+)$ reduces to a group conjugated to
$SO(3)\times SO(3)\times SO(2)$ (type~2) and to $SO(3)\times SO(5)$
(type~3) (see~\ref{relgeo}). The discussion of integrable
$PSU(3)$-structures leaves us with natural integrability conditions
for geometries of type~2 and 3.

\begin{definition}\hfill\newline
{\rm (i)} We call a structure of type~$2_0$ or $3_0$ {\em
integrable} if and only if the associated 3-form $\rho_{2_0}$ or
$\rho_{3_0}$ satisfies
\begin{equation}\label{harmform}
d\rho_{2_0,3_0}=0,\quad\dstar\rho_{2_0,3_0}=0.
\end{equation}
{\rm (ii)} We call a structure of type~$2_{\pm}$ or $3_{\pm}$ {\em
integrable} if and only if the spinor-valued 1-form
$\gamma=\gamma_{\pm}\in\Gamma(T^*M\otimes\Delta_{\pm})$ is harmonic
with respect to the twisted Dirac operator on
$\Gamma(T^*M\otimes\Delta)$, i.e.
\begin{equation}\label{harmonicity}
\mc{D}\gamma=0.
\end{equation}
\end{definition}

The integrability conditions~(\ref{harmform}) and
(\ref{harmonicity}) can be analysed in terms of the resulting
torsion along the lines of the $PSU(3)$-case, although the situation
becomes more difficult as the decomposition of type 2-structures
into irreducible modules of $\Lambda^1\otimes\mf{g}^{\perp}$ breaks
up into small pieces of various multiplicities. We shall deal with
two natural questions. For integrable structures of subtype~$0$, it
makes sense to ask if the holonomy of the metric reduces to the
groups $SO(3)\times SO(3)\times SO(2)$ or $SO(3)\times SO(5)$. We
will show that this is not to be expected for a generic metric of
type~$2_0$ or $3_0$. An example illustrating this fact will be given
in the next section. On the other hand, if the structure is of
subtype~$2_{\pm}$ we saw in Section~\ref{relgeo} that we obtain an
almost complex structure on the tangent bundle. We will show that
under the integrability condition~(\ref{harmonicity}) not even
partial integrability (in the sense that some of the torsion
components of an almost complex structure as described on p.~39
in~\cite{sa89} vanish) can be expected. We start with the analysis
of integrable structures of subtype~$0$.

\medskip

{\em Integrable structures of type~$2_0$.}\hspace{0.2cm} The
structure group reduces from $SO(8)=SO(\Lambda^1)$ to $SO(3)\times
SO(3)\times SO(2)$. We find the decomposition
\begin{eqnarray*}
\Lambda^1\otimes(\mf{su}(2)\oplus\mf{su}(2)\oplus \mf{t}^1)^{\perp}
& =
&\phantom{\oplus}([2,0,0]\oplus[0,0,2]\oplus\llbracket0,0,2\rrbracket)\otimes([2,2,0]\oplus\llbracket2,0,2\rrbracket\oplus\llbracket0,2,2\rrbracket)\\
& = &\phantom{\oplus}\llbracket0,2,4\rrbracket
\oplus2\llbracket0,0,2\rrbracket\oplus\llbracket0,2,2\rrbracket\oplus3[0,2,0]\oplus\llbracket0,4,2\rrbracket\oplus\\
&
&\oplus\llbracket2,0,4\rrbracket\oplus\llbracket2,0,2\rrbracket\oplus3[2,0,0]\oplus3\llbracket2,2,2\rrbracket
\oplus2[2,2,0]\oplus\\
& & \oplus[2,4,0]\oplus\llbracket4,0,2\rrbracket\oplus[4,2,0].
\end{eqnarray*}
Moreover, we have
$$
\Lambda^2=\mathbf{1}\oplus[0,2,0]\oplus\llbracket0,2,2\rrbracket\oplus[2,0,0]\oplus\llbracket2,0,2\rrbracket\oplus[2,2,0]
$$
and
$$
\Lambda^4=2\llbracket0,0,2\rrbracket\oplus2[0,2,0]\oplus2[2,0,0]\oplus2[2,2,0]\oplus2\llbracket2,2,2\rrbracket.
$$
From this we can partially determine the kernel by an application of
Schur's lemma. For instance, the kernel of the map
$\widetilde{\delta}:\Lambda^1\otimes(\mf{su}(2)\oplus\mf{su}(2)\oplus
\mf{t}^1)^{\perp}\to\Lambda^2$ restricted to the sub-module
$3\llbracket2,2,2\rrbracket$ contains on dimensional grounds at
least a module isomorphic to $\llbracket2,2,2\rrbracket$ etc.. Using
also
$$
\dim\ker\widetilde{d}\cap\dim\ker\widetilde{\delta}=\dim\ker\widetilde{d}+\dim\ker\widetilde{\delta}-\dim(\ker\widetilde{d}+\ker\widetilde{\delta}),
$$
we immediately draw the

\begin{prp}
We have
$$
\llbracket0,2,4\rrbracket
\oplus\llbracket0,4,2\rrbracket\oplus\llbracket2,0,4\rrbracket\oplus\llbracket2,2,2\rrbracket\oplus\llbracket4,0,2\rrbracket\oplus[2,4,0]\oplus[4,2,0]\subset\ker
\tilde{d}\cap\ker\widetilde{d^{\star}}.
$$
\end{prp}

In particular, the holonomy of an integrable structure of type~$2_0$
should not reduce to $SO(3)\times SO(3)\times SO(2)$. This is
confirmed by the example constructed in Section~\ref{geomexclass}.

\medskip

{\em Integrable structures of type~$3_0$.}\hspace{0.2cm} Next, we
assume that the structure group reduces from $SO(8)$ to $SO(3)\times
SO(5)$. We obtain
\begin{eqnarray*}
\Lambda^1\otimes(\mf{so}(3)\oplus\mf{so}(5))^{\perp} & = & ([2,0,0]\oplus[0,1,0])\otimes[2,1,0]\\
& = &
\phantom{(}[0,1,0]\oplus[2,0,0]\oplus[2,0,2]\oplus[2,1,0]\oplus[2,2,0]\oplus[4,1,0]
\end{eqnarray*}
while
$$
\Lambda^2=[0,0,2]\oplus[2,0,0]\oplus[2,1,0]
$$
and
$$
\Lambda^4=2[0,1,0]+2[2,0,2]
$$
which implies the

\begin{prp}
We have
$$
[2,2,0]\oplus[4,1,0]\subset\ker\tilde{d}\cap\ker\widetilde{d^{\star}}.
$$
\end{prp}

Again we give an explicit example where the holonomy does not reduce
to $SO(3)\times SO(5)$ in Section~\ref{geomexclass}.
\medskip

{\em Integrable structures of type~$2_{\pm}$.}\hspace{0.2cm} Here,
the structure group $SO(8)=SO(\Delta_{\pm})$ of the $Spin(8)$-bundle
reduces to $SO(3)\times SO(3)\times SO(2)$ (cf. Table
\ref{g2action}). By choosing the orientation appropriately we may
assume that we deal with a structure of type~$2_+$ which means that
$\gamma=\gamma_+\in\Lambda^3\Delta_-$ defines an isometry
$\gamma_+:\Lambda^1\to\Delta_+$. As we saw in Section~\ref{relgeo},
this forces the structure group of the tangent bundle to reduce to a
subgroup of $U(4)$, that is, we obtain an almost complex structure
on the tangent bundle. To discuss the integrability of this almost
complex structure under the condition (\ref{harmonicity}), we
decompose
\begin{eqnarray}
\Lambda^1\otimes\mf{u}(4)^{\perp} & = &
\llbracket0,1,1,-3\rrbracket\oplus\llbracket1,1,0,-1\rrbracket\oplus\llbracket1,0,0,-3\rrbracket\oplus\llbracket1,0,0,1\rrbracket\nonumber\\
& = & \llbracket
V\rrbracket\oplus\llbracket\lambda^{1,2}_0\rrbracket\oplus\llbracket\lambda^{0,3}\rrbracket\oplus\llbracket\lambda^{0,1}\rrbracket,\label{u4_tors}
\end{eqnarray}
where the labeling of the irreducible $U(4)$-modules was chosen in
accordance with~\cite{sa89}. From the point of view of $U(4)$, we
have $\Lambda^1\otimes\C=\C^4\oplus\overline{\C}^4$ where $U(4)$
acts on its vector representation with weights
$$
-\alpha_1+\alpha_2+\alpha_4,-\alpha_2+\alpha_3+\alpha_4,-\alpha_3+\alpha_4,\alpha_1+\alpha_4.
$$
This is refined by the actual action of $SU(2)\times SU(2)\times
U(1)/\Z_2$ which acts on $\C^4=(1,1,1)$ with weights
$$
-\sigma_1-\sigma_2+\tau,-\sigma_1+\sigma_2+\tau,\sigma_1-\sigma_2+\tau,\sigma_1+\sigma_2+\tau.
$$
Substituting these weights and decomposing the $U(4)$-modules
appearing in~(\ref{u4_tors}) accordingly, we obtain
$$
\begin{array}{lclcl}
\llbracket V\rrbracket & = & \llbracket0,1,1,-3\rrbracket & = &
\llbracket1,1,3\rrbracket\oplus\llbracket1,3,3\rrbracket\oplus\llbracket3,1,3\rrbracket\\
\llbracket\lambda^{1,2}_0\rrbracket & = & \llbracket1,1,0,-1\rrbracket & = &\llbracket1,1,-1\rrbracket\oplus\llbracket1,3,-1\rrbracket\oplus\llbracket3,1,-1\rrbracket\\
\llbracket\lambda^{0,3}\rrbracket & = & \llbracket1,0,0,-3\rrbracket
& = &\llbracket1,1,-3\rrbracket
\\
\llbracket\lambda^{0,1}\rrbracket & = & \llbracket1,0,0,1\rrbracket
& = & \llbracket1,1,1\rrbracket.
\end{array}
$$
The actual torsion of the structure lives in
\begin{eqnarray*}
\Lambda^1\otimes(\mf{su}(2)\oplus\mf{su}(2)\oplus \mf{t}^1)^{\perp} & = & \phantom{\oplus2}\llbracket1,1,1\rrbracket\otimes([2,2,0]\oplus\llbracket2,0,2\rrbracket\oplus\llbracket0,2,2\rrbracket)\\
& = &
\phantom{\oplus}2\llbracket1,1,3\rrbracket\oplus\llbracket1,3,3\rrbracket\oplus3\llbracket1,1,1\rrbracket\oplus2\llbracket1,3,1\rrbracket\oplus\llbracket3,1,3\rrbracket\oplus\\
& & \oplus2\llbracket3,1,1\rrbracket\oplus\llbracket3,3,1\rrbracket
\end{eqnarray*}
and is given by the kernel of the map
$$
\widetilde{\mc{D}}:\Lambda^1\otimes(\mf{su}(2)\oplus\mf{su}(2)\oplus\mf{t}^1)^{\perp}\to\Lambda^1\otimes\Delta_-.
$$
Now
$$
\Lambda^1\otimes\Delta_-=\llbracket1,1,1\rrbracket\otimes(\mf{su}(2)\oplus\mf{su}(2)\oplus
E_2)=3\llbracket1,1,1\rrbracket\oplus\llbracket3,1,1\rrbracket\oplus\llbracket1,3,1\rrbracket\oplus\llbracket1,1,3\rrbracket
$$
which implies

\begin{prp}
We have
$$
\llbracket1,1,3\rrbracket\oplus\llbracket1,3,1\rrbracket\oplus\llbracket3,1,1\rrbracket\oplus\llbracket1,3,3\rrbracket\oplus\llbracket3,1,3\rrbracket\oplus\llbracket3,3,1\rrbracket\subset\ker\tilde{\mc{D}}.
$$
\end{prp}

In particular, a generic integrable structure of type~$2_{\pm}$ is
expected to have torsion in the components $\llbracket V\rrbracket$,
$\llbracket\lambda_0^{1,2}\rrbracket$ and
$\llbracket\lambda^{0,3}\rrbracket$, and the almost complex
structure is not integrable.

\medskip

{\em Integrable structures of type~$3_{\pm}$.}\hspace{0.2cm} We
finally look at structures of type~$3_{\pm}$ where the structure
group $SO(8)=SO(\Delta_{\pm})$ of the $Spin(8)$-bundle reduces to
$SO(3)\times SO(5)$ (cf. Table~\ref{g3action}). Again we assume to
deal with a structure of subtype~$+$. We have
$$
\Lambda^1\otimes(\mf{so}(3)\oplus\mf{so}(5))^{\perp}=[1,0,1]\otimes[2,1,0]=[1,0,1]\oplus[1,1,1]\oplus[3,0,1]\oplus[3,1,1].
$$
On the other hand,
$$
\Lambda^1\otimes\Delta_-=[1,0,1]\otimes([2,0,0]\oplus[0,1,0])=2[1,0,1]\oplus[3,0,1]\oplus[1,1,1]
$$
and therefore

\begin{prp}
We have
$$
[3,1,1]\subset\ker\tilde{\mc{D}}.
$$
\end{prp}

Recall that an almost quaternionic structure on $M^8$ is defined by
an $Sp(1)\cdot Sp(2)$-invariant 4-form
$\Omega$~\cite{sa82},~\cite{sa89}. This defines a quaternionic
K\"ahler structure if and only if
$$
\nabla\Omega=T(\Omega)=0.
$$
Since
$$
\Lambda^3\cong\Lambda^5=[1,0,1]\oplus[1,1,1]\oplus[3,0,1]
$$
and
$\tilde{d}:\Lambda^1\otimes(\mf{so}(3)\oplus\mf{so}(5))^{\perp}\to\Lambda^5$
is onto~\cite{sa89}, the condition $d\Omega=0$ implies $T\in[3,1,1]$
and therefore $D\gamma=0$.

\begin{cor}\label{quatex}
Any almost quaternionic 8-manifold $(M^8,\Omega)$ with $d\Omega=0$ defines an integrable structure of type $3_{\pm}$.
\end{cor}

\subsection{Geometrical properties and examples}\label{geomexclass}

As in the case of special holonomy, a solution to the spinor field
equation~(\ref{harmonicity}) puts constraints on the Ricci tensor of
the metric. To see this we recall that according to Proposition 2.8
in~\cite{wa91} (using the notation of~(\ref{dodecomp})), we have
\begin{equation}\label{bochner}
(D\circ\delta-\delta\circ D)(\gamma)=\frac{1}{2}p(\gamma\circ{\rm
Ric})
\end{equation}
for any $\gamma\in\Gamma(T^*M\otimes\Delta)$. Now write
$\gamma=\sum_ie_i\otimes\gamma_i$ and regard ${\rm Ric}$ as an
endomorphism of $T$ so that
$$
\gamma\circ{\rm Ric}=\sum_{i,j}{\rm Ric}_{ij}e_i\otimes\gamma_j.
$$
If the structure is integrable,~(\ref{bochner}) implies
$$
p(\gamma\circ{\rm Ric})=\sum\limits_{i,j}{\rm
Ric}_{ij}e_i\cdot\gamma_j=0
$$
(alternatively, see~\cite{hi01} for a direct derivation of this
equation). This means that ${\rm Ric}$ is in the kernel of the map
\begin{equation}\label{rickernel}
A\in\bigodot^2\Lambda^1\mapsto
p(A\circ\gamma)=\sum_{i,j}A_{ij}e_i\gamma_j\in\Delta,
\end{equation}
which is invariant under the stabiliser of $\gamma$. In the case of
a $PSU(3)$-structure, we have $\odot^2\wedge^1=\mathbf{1}\oplus
[1,1]\oplus [2,2]$ and $\Delta=2[1,1]$. Since this map is
non-trivial, ${\rm Ric}$ vanishes on the module $[1,1]$. Similarly,
we can decompose the symmetric tensors for structures of type
$2_{\pm}$ and $3_{\pm}$ into
$$
\bigodot^2\Lambda^1=\mathbf{1}\oplus[2,0,0]\oplus[0,2,0]\oplus[2,2,0]\oplus\llbracket0,0,2\rrbracket\oplus\llbracket2,2,2\rrbracket
$$
and
$$
\bigodot^2\Lambda^1=\mathbf{1}+[0,1,0]+[2,0,2]
$$
respectively. The map~(\ref{rickernel}) takes then values in
$\Delta_{\mp}$ which carries a Lie bracket. It is immediate that for
an integrable type 3-structure, ${\rm Ric}$ must vanish on the
5-dimensional module $[0,1,0]$ while for type~$2$-structures it
vanishes on the component isomorphic to $\Delta_{\mp}$ (take the
matrices~(\ref{matrix2}) and check with the base elements
$e_i\otimes e_j+e_j\otimes e_i$ to see that the image has maximal
rank).

\begin{prp}\label{ricclass}\hfill\newline
{\rm (i)}{\rm~\cite{hi01}}\hspace{2pt} If $g$ is an integrable
$PSU(3)$-metric, then ${\rm Ric}$ vanishes on the component $[1,1]$
inside $\odot^2\Lambda^1$.

{\rm (ii)} If $g$ is an integrable metric of type~$2_{\pm}$, then
${\rm Ric}$ vanishes on the components
$\Delta_{\mp}\cong[2,0,0]\oplus[0,2,0]\oplus\llbracket0,0,2\rrbracket$
inside $\odot^2\Lambda^1$.

{\rm (iii)} If $g$ is an integrable metric of type~$3_{\pm}$, then
${\rm Ric}$ vanishes on the component $[0,1,0]$ inside
$\odot^2\Lambda^1$. In particular, this holds for any almost
quaternionic structure whose $Sp(1)\cdot Sp(2)$-invariant 4-form is
closed.
\end{prp}

\begin{rmk}
A quaternionic K\"ahler manifold is Einstein~\cite{sa89}, that is,
the Ricci-tensor vanishes on $[0,1,0]\oplus[2,0,2]$ if $\dim=8$. The
weaker condition $d\Omega=0$ still guarantees the vanishing on the
smaller component $[0,1,0]$.
\end{rmk}

\bigskip

\textbf{Examples}

(i) ({\em local description of type~$3_0$-geometries}) Consider a
3-form $\rho$ which induces a supersymmetric map
$\Delta_+\otimes\Delta_-$ of spin 3/2. Let ${\rm
Ann}(\rho)\subset\Gamma(T)$ denote the annihilator of $\rho$, that
is the distribution whose sections contracted with $\rho$ yield 0.
If the form $\rho$ is closed, then ${\rm Ann}(\rho)$ is involutive
since
$$
[X,Y]\llcorner\rho=\mc{L}_x(Y\llcorner\rho)-Y\llcorner\mc{L}_X\rho=-Y\llcorner
X\llcorner d\rho=0.
$$
If, in addition, the form is co-closed, then ${\rm Ann(\star\rho)}$
also defines an involutive distribution. In the case of a
$PSU(3)$-structure the annihilator of $\rho$ and $\star\rho$ is
$\{0\}$ as a consequence of the genericity of the form, but for
forms of type~2 and 3, we get something non-trivial. In particular,
we can apply Frobenius' theorem for type~3 structures to assert the
existence of local coordinates $(x_1,\ldots,x_8)$ such that ${\rm
Ann}(\star\rho)=\langle\partial_{x_1},\partial_{x_2},\partial_{x_3}\rangle$
and ${\rm
Ann}(\rho)=\langle\partial_{x_4},\ldots,\partial_{x_8}\rangle$.
Since $\star\rho$ and $\rho$ are the volume forms for ${\rm
Ann}(\rho)$ and ${\rm Ann}(\star\rho)$ respectively, the condition
to be of unit norm is just
$$
\det_{i,j=1,2,3}g_{ij}=\det_{i,j=1,\ldots,5}g_{ij}=1.
$$

{\rm (ii)} ({\em integrable $PSU(3)$-structures of positive,
negative and zero scalar curvature with vanishing torsion}) As
mentioned above, any integrable $PSU(3)$-metric with vanishing
torsion $T$ is either flat or symmetric. For example, consider the
3-form $ \rho(X,Y,Z)=B(X,[Y,Z])$ on the Lie algebra $\mf{su}(3)$,
the symmetric space $SU(3)=SU(3)\times SU(3)/SU(3)$ or its
non-compact dual $SL(3,\C)/SU(3)$. Since $\rho$ and $\star\rho$ are
${\rm Ad}$-invariant forms, they are closed and hence they define
integrable $PSU(3)$-structures. These three cases correspond to
zero, positive and negative scalar curvature respectively. In
particular, vanishing torsion does not imply Ricci-flatness.

(iii) ({\em non-trivial local examples of integrable
$PSU(3)$-structures and type~$2_0$- and $3_0$-structures which are
Ricci-flat}) This example will be built out of a hyperk\"ahler
4-manifold $M^4$ with a triholomorphic vector field. Let $U\equiv
U(x,y,z)$ be a positive harmonic function defined on some domain
$D\subset\R^3$ and $\theta$ a 1-form on $\R^3$ for which
$$
dU=\star d\theta
$$
holds. By the Gibbons-Hawking ansatz, the metric on $D\times \R$
\begin{equation}\label{hypmet}
g=U(dx^2+dy^2+dz^2)+\frac{1}{U}(dt+\theta)^2
\end{equation}
is hyperk\"ahler with associated K\"ahler forms given by
\begin{eqnarray*}
\omega_1 & = & Udy\wedge dz + dx\wedge (dt+\theta)\\
\omega_2 & = & Udx\wedge dy + dz\wedge (dt+\theta)\\
\omega_3 & = & Udx\wedge dz - dy\wedge (dt+\theta).
\end{eqnarray*}
The vector field $X=\frac{\partial}{\partial t}$ is triholomorphic,
that is it defines an infinitesimal transformation which preserves
any of the three complex structures induced by $\omega_1$,
$\omega_2$ or $\omega_3$. Conversely, a hyperk\"ahler metric on a
4-dimensional manifold which admits a triholomorphic vector field is
locally of the form~(\ref{hypmet}).

Following~\cite{ar02}, we define the 2-form $\widetilde{\omega}_3$
by changing the sign in $\omega_3$, that is
$$ \widetilde{\omega}_3=Udx\wedge dz + dy\wedge(dt+\theta).$$

This 2-form is closed if and only if
$$
U\equiv U(x,z)
$$
for $d\omega_3=0$ implies
\begin{equation}\label{de45}
d(Udx\wedge dz)=d(dy\wedge(dt+\theta)),
\end{equation}
so that
$$
d\widetilde{\omega}_3=2d(Udx\wedge dz)=2\frac{\partial U}{\partial
y}dy\wedge dx\wedge dz.
$$
Pick such a $U$ and take the standard coordinates $x_1,\ldots,x_4$
of the Euclidean space $(\R^4,g_0)$. Put
$$
\begin{array}{llll} e^1=dx_1, & e^2=dx_2, & e^3=dx_3 & e^8=dx_4 \\
e^4=\sqrt{U}dy, & e^5=-\frac{1}{\sqrt{U}}(dt+\theta), &
e^6=-\sqrt{U}dx, & e^7=\sqrt{U}dz
\end{array}
$$
which we take as an orthonormal basis on $M^4\times \R^4$. Endowed
with the orientation defined by $(e_4,\ldots,e_7)$, the forms
$\omega_i$ are anti-self-dual on $M^4$, while the forms
$\widetilde{\omega}_1=Udy\wedge dz - dx\wedge (dt+\theta)$,
$\widetilde{\omega}_2=Udx\wedge dy - dz\wedge (dt+\theta)$ and
$\widetilde{\omega}_3$ are self-dual. An example of an integrable
$PSU(3)$-structure is then provided by (cf.~\ref{rho1}
and~\ref{starrho1})
\begin{eqnarray*}
\rho_1 & = & dx_1\wedge dx_2\wedge
dx_3+\frac{1}{2}(dx_1\wedge\omega_1+dx_2\wedge\omega_2+dx_3\wedge\omega_3)+\frac{\sqrt{3}}{2}dx_4\wedge\widetilde{\omega}_3\\
& = &
e_{123}+\frac{1}{2}e_1\wedge\omega_1+\frac{1}{2}e_2\wedge\omega_2+\frac{1}{2}e_3\wedge\omega_3+\frac{\sqrt{3}}{2}e_8\wedge\widetilde{\omega}_3.
\end{eqnarray*}
Obviously, the equality $d\rho_1=0$ holds. Moreover, we have
\begin{eqnarray*}
\star\rho_1 & = & \phantom{+}Udx\wedge dz\wedge
dy\wedge(dt+\theta)\wedge dx_4-\frac{1}{2}\omega_1\wedge dx_2\wedge
dx_3\wedge dx_4\\ & & +\frac{1}{2}\omega_2\wedge dx_1\wedge
dx_3\wedge dx_4-\frac{1}{2}\omega_3\wedge dx_1\wedge dx_2\wedge
dx_4+\frac{\sqrt{3}}{2}\widetilde{\omega}_3\wedge dx_1\wedge
dx_2\wedge
dx_3\\
& = & \phantom{+}e_{45678}-\frac{1}{2}\omega_1\wedge
e_{238}+\frac{1}{2}\omega_2\wedge e_{138}-\frac{1}{2}\omega_3\wedge
e_{128}+\frac{\sqrt{3}}{2}\widetilde{\omega}_3\wedge e_{123},
\end{eqnarray*}
and one immediately checks that $\star\rho_1$ is also closed.

This ansatz yields then a structure of type~$2_0$ and $3_0$ by
defining
$$
\rho_2=\frac{1}{\sqrt{2}}e_{145}+\frac{1}{\sqrt{2}}e_{367}\mbox{ and
} \rho_3=e_{145}.
$$
Then
$$
\star\rho_2=\frac{1}{\sqrt{2}}e_{23678}+\frac{1}{\sqrt{2}}e_{12458}\mbox{
and }\star\rho_3=e_{23678},
$$
and all these forms are closed since~(\ref{de45}) implies
$de_{45}=de_{67}=0$.

To check that the holonomy is not contained in the respective
stabiliser groups $PSU(3)$, $SO(3)\times SO(3)\times SO(2)$ and
$SO(3)\times SO(5)$, we consider the specific example defined by
$$
U(x,y,z)=x\mbox{ on } \{x>0\}\mbox{ and } \theta=ydz,
$$
and show that $\nabla\rho_{i}\not=0$. The metric $g$ on
$M^4\times\R^4$ is given by
$$
g=dx_1^2+dx_2^2+dx_3^2+dx_4^2+xdx^2+xdy^2+(x+\frac{y^2}{x})dz^2+\frac{1}{x}dt^2+2\frac{y}{x}dzdt.
$$
with orthonormal basis
\begin{equation*}
\begin{array}{llll}
e_1=\partial_{x_1}, & e_2=\partial_{x_2}, & e_3=\partial_{x_3}, &
e_8=\partial_{x_4},\\ e_4=\frac{1}{\sqrt{x}}\partial_y, &
e_5=-\sqrt{x}\partial_t, & e_6=-\frac{1}{\sqrt{x}}\partial_x, &
e_7=\frac{1}{\sqrt{x}}(\partial_z-y\partial_t).
\end{array}
\end{equation*}
The only non-trivial brackets are
$$
\begin{array}{ll}
[e_4,e_6]=-\frac{1}{2\sqrt{x}^3}e_4 &
[e_5,e_6]=\frac{1}{2\sqrt{x}^3}e_5\\
\lbrack e_4,e_7\rbrack=\frac{1}{\sqrt{x}^3}e_5 &
[e_6,e_7]=\frac{1}{2\sqrt{x}^3}e_7.
\end{array}
$$
Since the anti-self-dual 2-forms $\omega_1$, $\omega_2$ and
$\omega_3$ are the associated K\"ahler forms of the hyperk\"ahler
structure on $M$, we have $\nabla\omega_i=0$. In particular,
$$
\nabla(e^6\wedge e^7)=\nabla(e^4\wedge
e^5)=-\frac{1}{12}\cdot\frac{1}{\sqrt{x^3}}(e_4\otimes\widetilde{\omega_1}+e_5\otimes\widetilde{\omega_2})
$$
holds. Thus
\begin{eqnarray*}
\nabla\rho_1 & = &
-\frac{1}{4\sqrt{3}}\cdot\frac{1}{\sqrt{x^3}}(e_4\otimes\widetilde{\omega}_1\wedge
e_8+e_5\otimes\widetilde{\omega}_2\wedge e_8)\\
\nabla\rho_2 & = &
-\frac{1}{12\sqrt{2x^3}}(e_4\otimes\widetilde{\omega}_1\wedge
(e_1+e_3)+e_5\otimes\widetilde{\omega}_2\wedge (e_1+e_3))\\
\nabla\rho_3 & = &
-\frac{1}{12\sqrt{x^3}}(e_4\otimes\widetilde{\omega}_1\wedge
e_1+e_5\otimes\widetilde{\omega}_2\wedge e_1),
\end{eqnarray*}
which shows that the holonomy does not reduce to the corresponding
stabiliser subgroup. This example also shows that Ricci-flatness
does not imply the vanishing of the torsion for an integrable
structure of any type.

(iv) ({\em non-trivial compact $PSU(3)$-manifold which is not
Einstein}) Recall that by defining a basis of (anti)-self-dual forms
$\omega_{1\pm}=e^{47}\pm e^{56}$, $\omega_{2\pm}=e^{46}\mp e^{57}$
and $\omega_{3\pm}=e^{45}\pm e^{67}$ we can write
\begin{equation}\label{rho1s}
\rho=
e_{123}+\frac{1}{2}e_1\wedge\omega_{1-}+\frac{1}{2}e_2\wedge\omega_{2-}+\frac{1}{2}e_3\wedge\omega_{3-}
+\frac{\sqrt{3}}{2}e_8\wedge\omega_{3+}.
\end{equation}
and
\begin{equation}\label{starrho1s}
\star\rho=e_{45678}-\frac{1}{2} e_{238}\wedge\omega_1+\frac{1}{2}
e_{138}\wedge\omega_2-\frac{1}{2}e_{128}\wedge\omega_3
+\frac{\sqrt{3}}{2}e_{123}\wedge\omega_{3+}.
\end{equation}
Note that
\begin{equation}\label{wedge}
\omega_{i\pm}\wedge\omega_{j\mp}=0,\quad
\omega_{i\pm}\wedge\omega_{j\pm}=\pm2\delta_{ij}e_{4567}.
\end{equation}
Take $N$ to be a 6-dimensional nilmanifold associated with the Lie
algebra $\mf{g}=\langle e_2,\ldots,e_8\rangle$ whose structure
constants are determined by
$$
de_i=\left\{\begin{array}{rl}0,&i=2,\ldots,7\\\omega_{1+}=e_{47}+e_{56},&i=8\end{array}\right.
$$
(the labeling of the vectors follows our conventional use of a
$PSU(3)$-frame), that is the only non-trivial structure constants
are
$$
c_{478}=-c_{748}=c_{568}=-c_{658}=1.
$$
Let $G$ be the associated simply-connected Lie group. The
rationality of the structure constants guarantees the existence of a
lattice $\Gamma$ for which $N=\Gamma\backslash G$ is
compact~\cite{ma51}. We let $M=T^2\times N$ with $e_i=dt_i,\,i=1,2$
on the torus, hence $de_i=0$. We take the basis $e_1,\ldots,e_8$ to
be orthonormal on $M$ and let $g$ denote the corresponding metric.
Then~(\ref{rho1s}), (\ref{starrho1s}),~(\ref{wedge}) imply
$$
d\rho =
\frac{\sqrt{3}}{2}de_8\wedge\omega_{3+}=\frac{\sqrt{3}}{2}\omega_{1+}\wedge\omega_{3+}=0
$$
and similarly, we have
\begin{eqnarray*}
d\star\rho & = & e_{4567}\wedge de_8-\frac{1}{2}\omega_{1-}\wedge
e_{23}\wedge de_8+\frac{1}{2}\omega_{2-}\wedge e_{13}\wedge
de_8-\frac{1}{2}\omega_{3-}\wedge e_{12}\wedge de_8\\
& = & e_{4567}\wedge \omega_{1+}-\frac{1}{2}\omega_{1-}\wedge
e_{23}\wedge \omega_{1+}+\frac{1}{2}\omega_{2-}\wedge e_{13}\wedge
\omega_{1+}-\frac{1}{2}\omega_{3-}\wedge e_{12}\wedge \omega_{1+}\\
& = & 0.
\end{eqnarray*}
Next we compute the covariant derivatives $\nabla_{e_i}e_j$. By
Koszul's formula, we have
$$
2g(\nabla_{e_i}e_j,e_k)=g([e_i,e_j],e_k)+g([e_k,e_i],e_j)+g(e_i,[e_k,e_j])=c_{ijk}+c_{kij}+c_{kji}.
$$
It follows that
$$
\nabla
e_i=\left\{\begin{array}{rl}0,&i=1,\,2,\,3\\-\frac{1}{2}(e_7\otimes
e_8+e_8\otimes e_7),&i=4\\-\frac{1}{2}(e_6\otimes e_8+e_8\otimes
e_6),&i=5\\\frac{1}{2}(e_5\otimes e_8+e_8\otimes
e_5),&i=6\\\frac{1}{2}(e_4\otimes e_8+e_8\otimes
e_4),&i=7\\\frac{1}{2}(-e_4\otimes e_7+e_7\otimes e_4-e_5\otimes
e_6+e_6\otimes e_5),&i=8\\.\end{array}\right.
$$
Now
$$
\nabla_{e_4}(e_8\wedge\omega_{3+})=e_{457}
$$
and since the coefficient of $e_8\wedge\omega_{3+}$ is irrational
while all the remaining ones are rational, we deduce
$$
\nabla_{e_4}\rho\not=0
$$
and therefore the metric is non-symmetric. Moreover, a
straightforward computation shows the diagonal of the Ricci-tensor
$\ric_{ij}$ to be given by
$$
\ric_{ii}=\left\{\begin{array}{rl}0,&i=1,\,2,\,3,\,8\\\frac{1}{2},&i=4,\,5,\,6,\,7\end{array}\right.
$$
In particular, it follows that $(M,g)$ is not Einstein, that is,
$\ric$ will usually have a non-trivial $[2,2]$-component (cf.
Proposition~\ref{ricclass}).

{\rm (v)} ({\em non-trivial compact example of a type
$3_+$-structure}) In~\cite{sa01}, Salamon constructed a compact
almost quaternionic 8-manifold $M$ whose structure form $\Omega$ is
closed, but not parallel. The example is of the form $M=N^6\times
T^2$, where $N^6$ is a compact nilmanifold associated with the Lie
algebra given by
$$
de_i=\left\{\begin{array}{rl}0,&i=1,\,2,\,3,\,5\\e_{15},&i=4\\e_{13},&i=6\end{array}\right.
$$
According to Corollary~\ref{quatex}, this implies the existence of a
compact non-trivial structure of type $3_{\pm}$.

\begin{prp}\hfill\newline
{\rm (i)} Compact integrable $PSU(3)$-manifolds with non-vanishing
torsion and which are not Einstein exist. On the other hand,
vanishing torsion does not imply Ricci-flatness. Moreover, there
exist non-trivial integrable local Ricci-flat examples.

{\rm (ii)} Compact integrable structures of type $3_+$ with
non-vanishing torsion exist. Moreover, there exist non-trivial local
and Ricci-flat examples of type $2_0$- and $3_0$-structures.
\end{prp}

\section{Integrable variational structures and related geometries in the generalised
case}\label{gengeom}

\subsection{The spinorial solution of the variational problem and related
geometries}\label{spinsolgencase}

Having analysed the variational problem for classical structures in
terms of a spinorial field equation we want to approach the
generalised case in a similar vein. First, we give a somewhat
technical definition which has the merit of encapsulating at once
the equations coming from the unconstrained, constrained and the
twisted variational problem.

\begin{definition}\label{intstruc}
Let $H$ be a (closed) 3-form and $\lambda$ be a real, non-zero
constant.

{\rm (i)} A topological generalised $G_2$-structure $(M,\rho)$ is
said to be {\em (closed) strongly integrable} with respect to $H$ if
and only if
$$
d_H\rho=0,\quad d_H\hat{\rho}=0.
$$

{\rm (ii)} A topological generalised $G_2$-structure $(M,\rho)$ is
said to be {\em (closed) weakly integrable with respect to $H$ and
with Killing number $\lambda$} if and only if
$$
d_H\rho=\lambda\hat{\rho}.
$$
We call such a structure {\em even} or {\em odd} according to the
parity of the form $\rho$. If we do not wish to distinguish the
type, we will refer to both structures as {\em (closed) weakly
integrable}.

Similarly, the structures in (i) and (ii) will also be referred to
as {\em integrable} if a condition applies to both (closed) weakly
and strongly integrable structures.

{\rm (iii)} A topological generalised $Spin(7)$-structure $(M,\rho)$
will be called {\em (closed) integrable with respect to $H$} if and
only if
$$
d_H\rho=0.
$$
\end{definition}

\begin{rmk}
As we shall see in Corollary~\ref{torsioncomp}, the constant
$\lambda$ describes the $\Lambda^4_0$-component of $d\varphi_{\pm}$,
where $\varphi_{\pm}$ are the underlying stable 3-forms of a weakly
integrable structure. A classical $G_2$-structure having only
torsion components in this module is said to be {\em
nearly-parallel} and can be equivalently characterised by the
existence of a Killing spinor in the sense of~\cite{bfgk91} (see
also the discussion in Section~\ref{class_sol}). The constant
$\lambda$ is usually referred to as the {\em Killing number} which
explains our terminology.
\end{rmk}

We want to interpret the equations in Definition~\ref{intstruc} in
terms of the data of Theorem~\ref{topology}. The key will be again
the interplay between the twisted Dirac operators $\mc{D}$ and
$\mc{\widehat{D}}$ on $\Delta\otimes\Delta$, given in a local basis
by
$$
\mc{D}(\Psi\otimes\Phi)=\sum
s_i\cdot\nabla_{s_i}\Psi\otimes\Phi+s_i\cdot\Psi\otimes\nabla_{s_i}\Phi
$$
and
$$
\widehat{\mc{D}}(\Psi\otimes\Phi)=\sum \nabla_{s_i}\Psi\otimes
s_i\cdot\Phi+\Psi\otimes s_i\cdot\nabla_{s_i}\Phi,
$$
and the corresponding Dirac operators on $\Lambda^p$ (under the
usual identification of $\Delta\otimes\Delta$ with forms)
$d+d^{\star}$ and $(-1)^p(d-d^{\star})$. For the generalised version
however, we have to use the twisted isomorphism $L_b$ (as defined in
Section~\ref{ges}) between $\Delta\otimes\Delta$ and
$\Lambda^{ev,od}$. In order to describe its effect on the level of
differential operators we define for a $p$-form $\alpha$ the
operator
$$
d^{\Box}\alpha^p=(-1)^{n(p+1)+1}\Box\, d\,\Box\alpha^p.
$$
This definition ensures that in absence of a B-field, the identity
$$
d^{\Box}=d^{\star}
$$
holds for $n=7$ or $8$ (see Lemma~\ref{hat}).

\begin{prp}\label{dirac}(Cf. the notation of page~\pageref{Lnot})

{\rm (i)} For $n=7$ we have
$$
L^{ev,od}_b\circ\mc{D}=dL^{od,ev}_b+
d^{\Box}L_b^{od,ev}+\frac{1}{2}e^{b/2}\wedge(db\llcorner
L^{od,ev}-db\wedge L^{od,ev})
$$
and
$$
L^{ev,od}_b\circ\widehat{\mc{D}}=
\pm(dL_b^{od,ev}-d^{\Box}L^{od,ev}_b)\mp\frac{1}{2}e^{b/2}\wedge(db\llcorner
L^{od,ev}+ db\wedge L^{od,ev}).
$$

{\rm (ii)} For $n=8$ we have
$$
\begin{array}{ccc}
L^{ev}_{b,\pm}\circ\mc{D} & = &
dL^{od,ev}_{b,\pm}+d^{\Box}L^{od,ev}_{b,\pm}+\frac{1}{2}e^{b/2}\wedge(db\llcorner
L^{od,ev}_{\pm}- db\wedge L^{od,ev}_{\pm})\\
L^{od}_{b,\pm}\circ\mc{D} & = &
dL^{od,ev}_{b,\pm}+d^{\Box}L^{od,ev}_{b,\pm}+\frac{1}{2}e^{b/2}\wedge(db\llcorner
L^{od,ev}_{\pm}- db\wedge L^{od,ev}_{\pm})
\end{array}
$$
and
\begin{eqnarray*}
L^{ev}_{b,\pm}\circ\widehat{\mc{D}} & = &
\phantom{-}dL^{od}_{b,\pm}-d^{\Box}L_{b,\pm}^{od}-\frac{1}{2}e^{b/2}\wedge(db\llcorner
L^{od}_{\pm}+ db\wedge L^{od}_{\pm})\\
L^{od}_{b,\pm}\circ\widehat{\mc{D}} & = &
-dL^{ev}_{b,\pm}+d^{\Box}L_{b,\pm}^{ev}+\frac{1}{2}e^{b/2}\wedge(db\llcorner
L^{ev}_{\pm}+ db\wedge L^{ev}_{\pm})
\end{eqnarray*}
\end{prp}

\begin{prf}
For $b=0$ the assertion
$$
L^{ev,od}\circ\mc{D}=dL^{od,ev}+d^*L^{od,ev}
$$
and
$$
L^{ev,od}\circ\widehat{\mc{D}}=\pm(dL^{od,ev}-d^*L^{od,ev}).
$$
is classical (see, for instance,~\cite{be87} Subsection 1.I; compare
also Section~\ref{class_sol}).

Now consider the case of an arbitrary B-field $b$ which gives rise
to the twisted maps
$L^{ev,od}_b:\Delta\otimes\Delta\to\Lambda^{ev,od}T^*$. Then
$$
L^{ev,od}_b\circ\mc{D}=e^{b/2}\wedge
L^{ev,od}\circ\mc{D}=e^{b/2}\wedge\left(dL^{od,ev}+d^{\star}L^{od,ev}\right).
$$
The first term on the right hand side equals
$$
e^{b/2}\wedge dL^{od,ev}=dL^{od,ev}_b-\frac{1}{2}e^{b/2}\wedge
db\wedge L^{od,ev}.
$$
For the second term, we have to treat the cases $n=7$ and $n=8$
separately. If $n=7$, then
$$
d^{\star}L^{od,ev}=\mp\!\star\dstar L^{od,ev}.
$$
Proposition~\ref{selfdual} and Lemma~\ref{hat} give
\begin{eqnarray*}
e^{b/2}\wedge d^*L^{od,ev} & = & \mp e^{b/2}\wedge \star\dstar L^{od,ev}\\
& = & \mp e^{b/2}\wedge \star
d\sigma(L^{ev,od})\\
& = & \mp e^{b/2}\wedge\star
d(e^{b/2}\wedge\sigma(L^{ev,od}_b))\\
& = & \mp e^{b/2}\wedge \star(e^{b/2}\wedge(\frac{1}{2}db\wedge\sigma(L^{ev,od}_b)+ d\sigma(L^{ev,od}_b)))\\
& = & \mp\frac{1}{2} e^{b/2}\wedge\star
(db\wedge \sigma(L^{ev,od}))\mp\Box_{g,b}dL^{ev,od}_b\\
& = & \phantom{\pm}\frac{1}{2} e^{b/2}\wedge(db\llcorner
L^{od,ev})+d^{\Box}L^{od,ev}_b,
\end{eqnarray*}
hence
$$
L^{ev,od}_b\circ\mc{D}=dL^{od,ev}_b+
d^{\Box}L_b^{od,ev}+\frac{1}{2}e^{b/2}\wedge(db\llcorner
L^{od,ev}-db\wedge L^{od,ev}).
$$
Similarly, we obtain
\begin{eqnarray*}
L^{ev,od}_b\circ\widehat{\mc{D}} & = & \pm
e^{b/2}\wedge(dL^{od,ev}\pm\star \dstar L^{od,ev})\\
& = &
\pm(dL_b^{od,ev}-d^{\Box}L^{od,ev}_b)\mp\frac{1}{2}e^{b/2}\wedge(db\llcorner
L^{od,ev}+ db\wedge L^{od,ev}).
\end{eqnarray*}
Next we treat the case $n=8$. Here we have to take into account the
chirality of the spinors. Recall our convention that
$L^{ev,od}_{b,\pm}$ denotes the restriction to
$\Delta_*\otimes\Delta_{+,-}$ where $*=+$ for $L^{ev}_b$ and $*=-$
for $L^{od}_b$. So far we know that
$$
L^{ev,od}_{b,\pm}\circ\mc{D}=e^{b/2}\wedge (dL^{od,ev}_{\pm}-\star
\dstar L^{od,ev}_{\pm}).
$$
Lemma~\ref{hat} and Proposition~\ref{selfdual} imply
$$
\star L^{ev}_{\pm}=\pm\sigma(L^{ev}_{\pm})\mbox{ and}\star
L^{od}_{\pm}=\mp\sigma(L^{od}_{\pm}).
$$
Consequently, we have to deal with both cases separately. Now
\begin{eqnarray*}
e^{b/2}\wedge d^*L^{ev}_{\pm} & = & -e^{b/2}\wedge\star \dstar
L^{ev}_{\pm}\\
& = & \mp e^{b/2}\wedge\star d\sigma(L^{ev}_{\pm})\\
& = & \mp e^{b/2}\wedge\star d(e^{b/2}\wedge
\sigma(L^{ev}_{b,\pm}))\\
& = & \mp\frac{1}{2}
e^{b/2}\wedge\star(db\wedge\sigma(L^{ev}_{\pm}))\pm
e^{b/2}\wedge\star(e^{b/2}\wedge[dL^{ev}_{b,\pm}]^{\wedge})\\
& = & \pm\frac{1}{2}e^{b/2}\wedge(db\llcorner
*\sigma(L^{ev}_{\pm}))\pm\Box dL^{ev}_{b,\pm}\\
& = & \phantom{\pm}\frac{1}{2}e^{b/2}\wedge(db\llcorner
L^{ev}_{\pm})+d^{\Box}L^{ev}_{b,\pm}.
\end{eqnarray*}
Similarly, we find
$$
-e^{b/2}\wedge\star \dstar
L^{od}_{\pm}=\frac{1}{2}e^{b/2}\wedge(db\llcorner
L^{od}_{\pm})+d^{\Box}L^{od}_{b,\pm}.
$$
Putting these pieces together yields
$$
L^{ev,od}_{b,\pm}\circ\mc{D}=dL^{od,ev}_{b,\pm}+d^{\Box}L^{ev,od}_{b,\pm}+\frac{1}{2}e^{b/2}(db\llcorner
L^{ev,od}_{\pm}- db\wedge L^{od,ev}_{\pm}).
$$
A similar expression holds for
$L^{ev,od}_{b,\pm}\circ\widehat{\mc{D}}$.
\end{prf}

Theorem~\ref{topology} states that a generalised structure $\rho$ on
$M$ can be represented by $e^{-F}L_b(\Psi_+\otimes\Psi_-)$. The
integrability conditions in Definition~\ref{intstruc} translate as
follows.

\begin{thm}\label{integrability}\hfill\newline
{\rm (i)} A generalised $G_2$-structure $(M^7,\rho)$ is weakly
integrable with Killing number $\lambda$ if and only if
$e^{-F}L_b(\Psi_+\otimes\Psi_-)=\rho+\Box_{g,b}\rho$ satisfies (with
$T=db/2+H$)
\begin{enumerate}
\item for all vector fields $X$,
$$
\nabla_X\Psi_{\pm}\pm\frac{1}{4}(X\llcorner T)\cdot\Psi_+=0,
$$
\item
$$
(dF\pm\frac{1}{2}T\pm\lambda)\cdot\Psi_{\pm}=0
$$
in case of an even structure and
$$
(dF\pm\frac{1}{2}T+\lambda)\cdot\Psi_{\pm}=0
$$
in case of an odd structure.
\end{enumerate}

The structure $e^{-F}L_b(\Psi_+\otimes\Psi_-)=\rho+\Box_{g,b}\rho$
is strongly integrable if and only if these equations hold for
$\lambda=0$.

{\rm (ii)} A generalised $Spin(7)$-structure $(M^8,\rho)$ is
integrable if and only if $e^{-F}L_b(\Psi_+\otimes\Psi_-)=\rho$
satisfies (with $T=db/2+H$)
\begin{enumerate}
\item for all vector fields $X$,
$$
\nabla_X\Psi_{\pm}\pm\frac{1}{4}(X\llcorner T)\cdot\Psi_{\pm}=0
$$
\item
$$
(dF\pm\frac{1}{2}T)\cdot\Psi_{\pm}=0
$$
\end{enumerate}
We refer to the equation involving the covariant derivative of the
spinor as the {\em (generalised) Killing equation} and to the
equation involving the differential of $F$ as the {\em dilaton
equation}.
\end{thm}

\begin{rmk}
The generalised Killing equation basically states that we have two
metric connections $\nabla^{\pm}$ preserving the underlying
$G_{2\pm}$- or $Spin(7)_{\pm}$-structures whose torsion tensor (as
it is to be defined in~\ref{geompropgen}) is {\em skew}. The dilaton
equation then serves to identify the components of the torsion with
respect to the decomposition into irreducible $G_{2\pm}$- or
$Spin(7)_{\pm}$-modules with the additional data $dF$ and $\lambda$.
\end{rmk}

\begin{prf}
We first consider the strong integrability condition $d_H\rho=0$,
$d_H\hat{\rho}=0$ which is equivalent to
\begin{equation}\label{intcondH}
d\rho=-H\wedge\rho,\quad d\Box\rho=-H\wedge\Box\rho.
\end{equation}
We write $\rho=e^{-F} L_b^{ev}(\Psi_+\otimes\Psi_-)$,
$e^{-F}L_{b,+}^{ev}(\Psi_+\otimes\Psi_-)$ or $e^{-F}
L_{b,-}^{od}(\Psi_+\otimes\Psi_-)$ according to the cases $G_2$,
$Spin(7)$ even or $Spin(7)$ odd. Taking $\Box$ of the second
equation in~(\ref{intcondH}) gives
$$
\begin{array}{lcll}
d^{\Box} e^{-F}L^{ev}_b & = & e^{b/2}\wedge
(H\llcorner e^{-F}L^{ev}), & \quad(G_2)\\
d^{\Box} e^{-F}L^{ev}_{b,+} & = & e^{b/2}\wedge
(H\llcorner e^{-F}L^{ev}_+),& \quad(Spin(7),\:even)\\
d^{\Box} e^{-F}L^{od}_{b,-} & = & e^{b/2}\wedge
(H\llcorner e^{-F}L_-^{od}),& \quad(Spin(7),\:odd)\\
\end{array}
$$
For example,
\begin{eqnarray*}
d^{\Box}e^{-F}L^{ev}_{b,+} & = & -\Box(H\wedge\Box e^{-F}L^{ev}_{b,+})\\
& = & -e^{b/2}\wedge\star(e^{b/2}\wedge\sigma(H\wedge
e^{-F}L^{ev}_{b,_+}))\\
& = & -e^{b/2}\wedge\star(H\wedge\sigma(e^{-F}L^{ev}_+))\\
& = & \phantom{-}e^{b/2}\wedge(H\llcorner e^{-F}L^{ev}_+).
\end{eqnarray*}
From Proposition~\ref{dirac} we obtain the following expressions for
the Dirac operator $\mc{D}$ (where we put $ T=db/2+H$).
$$
\begin{array}{lcl}
L^{od}\circ\mc{D}(e^{-F}\Psi_+\otimes\Psi_-) & = &  T\llcorner
L^{ev}(e^{-F}\Psi_+\otimes\Psi_-)- T\wedge
L^{ev}(e^{-F}\Psi_+\otimes\Psi_-)\\
L^{ev}_+\circ\mc{D}(e^{-F}\Psi_+\otimes\Psi_-) & = & T\llcorner
L^{od}_+(e^{-F}\Psi_+\otimes\Psi_-)- T\wedge
L^{od}_+(e^{-F}\Psi_+\otimes\Psi_-)\\
L^{od}_-\circ\mc{D}(e^{-F}\Psi_+\otimes\Psi_-) & = & T\llcorner
L^{ev}_-(e^{-F}\Psi_+\otimes\Psi_-)- T\wedge
L^{ev}_-(e^{-F}\Psi_+\otimes\Psi_-).\\
\end{array}
$$
Applying Corollary~\ref{baction} to these equations entails
\begin{eqnarray}\label{inteq1}
\mc{D}(e^{-F}\Psi_+\otimes\Psi_-) & =& D(e^{-F}\Psi_+)\otimes\Psi_-+e^{-F}\sum\limits_i s_i\cdot\Psi_+\otimes\nabla_{s_i}\Psi_-\\
& = & \frac{1}{4}e^{-F}(-
T\cdot\Psi_+\otimes\Psi_-+\sum\limits_is_i\cdot\Psi_+\otimes
(s_i\llcorner  T)\cdot\Psi_-),\nonumber
\end{eqnarray}
where $D$ denotes the usual Dirac operator associated with the
Clifford bundle $(\Delta,\,q_{\Delta})$. Since
$q_{\Delta}(s_i\cdot\Psi_+,s_j\cdot\Psi_+)=\delta_{ij}$, the
contraction of this equation with $q_{\Delta}(\cdot,s_m\cdot\Psi_+)$
yields
$$
q_{\Delta}(D(e^{-F}\Psi_+),s_m\cdot\Psi_+)\Psi_-+e^{-F}\nabla_{s_m}\Psi_-=
\frac{1}{4}e^{-F}(-q_{\Delta}(T\cdot\Psi_+,s_m\cdot\Psi_+)\Psi_-+(s_m\llcorner
 T)\cdot\Psi_-),
$$
implying
$$
e^{-F}\nabla_{s_m}\Psi_-=-\frac{1}{4}q_{\Delta}(4D(e^{-F}\Psi_+)+e^{-F}
T\cdot\Psi_+,s_m\cdot\Psi_+)\Psi_-+\frac{1}{4}e^{-F}(s_m\llcorner
 T)\cdot\Psi_-.
$$
However, from this expression of the covariant derivative we
deduce
$$
s_m.q_{\Delta}(\Psi_-,\Psi_-) =
2q_{\Delta}(\nabla_{s_m}\Psi_-,\Psi_-)=-\frac{1}{2}q_{\Delta}(4e^FD(e^{-F}\Psi_+)+
T\cdot\Psi_+,s_m\cdot\Psi_+)=0
$$
as $s_m\llcorner  T\in\mf{so}(T,g)$ and consequently
$q_{\Delta}((s_m\llcorner  T)\cdot\Psi_-,\Psi_-)=0$. It follows
\begin{equation}\label{cov-}
\nabla_{s_m}\Psi_-=\frac{1}{4}(s_m\llcorner  T)\cdot\Psi_-.
\end{equation}
In order to derive the corresponding expression for the spinor
$\Psi_+$, we consider the Dirac operator $\widehat{\mc{D}}$.
Firstly, we find
$$
\begin{array}{lcl}
L^{od}\circ\widehat{\mc{D}}(\Psi_+\otimes e^{-F}\Psi_-) & = &
\phantom{-}T\llcorner e^{-F}L^{ev}(\Psi_+\otimes \Psi_-)+ T\wedge
e^{-F}L^{ev}(\Psi_+\otimes\Psi_-)\\
L^{ev}_+\circ\widehat{\mc{D}}(\Psi_+\otimes e^{-F}\Psi_-) & = & -
T\llcorner e^{-F}L^{od}_+(\Psi_+\otimes\Psi_-)- T\wedge
e^{-F}L^{od}_+(\Psi_+\otimes\Psi_-)\\
L^{od}_-\circ\widehat{\mc{D}}(\Psi_+\otimes e^{-F}\Psi_-) & = &
\phantom{-}T\llcorner e^{-F}L^{ev}_-(\Psi_+\otimes\Psi_-)+ T\wedge
e^{-F}L^{ev}_-(\Psi_+\otimes\Psi_-).\\
\end{array}
$$
We appeal again to Lemma~\ref{baction} to deduce
\begin{eqnarray}\label{inteq2}
\widehat{\mc{D}}(\Psi_+\otimes e^{-F}\Psi_-) & = & \Psi_+\otimes D(e^{-F}\Psi_-)+ e^{-F}\sum\limits_i\nabla_{s_i}\cdot\Psi_+\otimes s_i\cdot\Psi_-\\
& = & \frac{1}{4}e^{-F}(\Psi_+\otimes
T\cdot\Psi_--\sum\limits_i(s_i\llcorner  T)\cdot\Psi_+\otimes
s_i\cdot\Psi_-).\nonumber
\end{eqnarray}
Now contraction with $q_{\Delta}(\cdot,s_m\cdot\Psi_-)$ gets
$$
e^{-F}\nabla_{s_m}\Psi_+=\frac{1}{4}q_{\Delta}(-4D(e^{-F}\Psi_-)+
e^{-F} T\cdot
\Psi_-,s_m\cdot\Psi_-)\Psi_+-\frac{1}{4}e^{-F}(s_m\llcorner
 T)\cdot\Psi_+.
$$
We derive the norm of $\Psi_+$ and obtain
$$
s_m.q_{\Delta}(\Psi_+,\Psi_+) =
2q_{\Delta}(\nabla_{s_m}\Psi_+,\Psi_+)=\frac{1}{2}q_{\Delta}(-4e^FD(e^{-F}\Psi_-)+
T\cdot \Psi_-,s_m\cdot\Psi_-)=0.
$$
Hence in all three cases ($G_2$, $Spin(7)$ even, $Spin(7)$ odd) we
find
\begin{equation}\label{cov+}
\nabla_{s_m}\Psi_+=-\frac{1}{4}(s_m\llcorner  T)\cdot\Psi_+.
\end{equation}
If we now consider the inhomogeneous equation
$d_H\rho=\lambda\hat{\rho}$, then we pick up an extra term in
(\ref{inteq1}) and~(\ref{inteq2}). For example,~(\ref{inteq2})
becomes
\begin{equation}\label{diraceven}
\widehat{\mc{D}}(\Psi_+\otimes e^{-F}\Psi_-) =-\lambda
e^F\Psi_+\otimes\Psi_-+\frac{1}{4}e^{-F}(\Psi_+\otimes
T\cdot\Psi_--\sum\limits_i(s_i\llcorner T)\cdot\Psi_+\otimes
s_i\cdot\Psi_-),
\end{equation}
for weak integrability of even type and
\begin{equation}\label{diracodd}
\widehat{\mc{D}}(\Psi_+\otimes e^{-F}\Psi_-) =\lambda
e^F\Psi_+\otimes\Psi_-+\frac{1}{4}e^{-F}(\Psi_+\otimes
T\cdot\Psi_--\sum\limits_i(s_i\llcorner T)\cdot\Psi_+\otimes
s_i\cdot\Psi_-)
\end{equation}
for weak integrability of odd type. But these extra terms cancel if
we contract by $q_{\Delta}(\cdot,s_m\cdot\Psi_-)$ and similarly for
the expression involving $\mc{D}$. Consequently we derive the same
condition on the covariant derivative of $\Psi_{\pm}$. The
difference between these integrability conditions is only visible in
the dilaton equation to which we turn next.

Here, we contract (\ref{diraceven}) and (\ref{diracodd}) with
$q_{\Delta}(\cdot,\Psi_+)$ and obtain
\begin{eqnarray*}
De^{-F}\Psi_- & = & -\lambda
e^{-F}\Psi_-+\frac{1}{4}e^{-F}T\cdot\Psi_-,\quad\mbox{\rm
even case}\\
De^{-F}\Psi_- & = & \phantom{-}\lambda
e^{-F}\Psi_-+\frac{1}{4}e^{-F}T\cdot\Psi_-,\quad\mbox{\rm
odd case.}\\
\end{eqnarray*}
On the other hand, equation~(\ref{cov-}) entails
$$
De^{-F}\Psi_-=e^{-F}(-dF\cdot\Psi_-+\frac{3}{4} T\cdot\Psi_-),
$$
which gives the dilaton equation for $\Psi_-$. Using $\mc{D}$ in
conjunction with (\ref{cov+}) implies the dilaton equation for
$\Psi_+$.

Conversely assume that the two spinors $\Psi_{\pm}$ satisfy the
generalised Killing equations. We consider the $G_2$-case first.
Since
$$
d_H(e^{-F}L^{ev,od}_b)=e^{-F}e^{b/2}\wedge( T\wedge L^{ev,od}-
dF\wedge L^{ev,od}_b+dL^{ev,od})
$$
it suffices to show that
\begin{equation}\label{closed}
dL^{ev,od}-dF\wedge L^{ev,od}=\lambda L^{od,ev}-T\wedge L^{ev,od}.
\end{equation}
Let us consider a weak structure of even type (the odd case, again,
is analogous). Then
\begin{eqnarray*}
(dL^{ev}-dF\wedge L^{ev})(\Psi_+\otimes\Psi_-) & = &
\phantom{+}\frac{1}{2}L^{od}(-dF\cdot\Psi_+\otimes\Psi_-+\Psi_+\otimes
dF\cdot\Psi_-)+\\
& &+\sum\limits_i s_i\wedge\nabla_{s_i}L^{ev}(\Psi_+\otimes\Psi_-)\\
& = & \phantom{+}\lambda
L^{od}(\Psi_+\otimes\Psi_-)+\frac{1}{4}L^{od}(
T\cdot\Psi_+\otimes\Psi_-+\Psi_+\otimes
T\cdot\Psi_-)+\\
& &
+\frac{1}{2}\sum\limits_iL^{od}(s_i\cdot\nabla_{s_i}\Psi_+\otimes\Psi_--
\nabla_{s_i}\Psi_+\otimes s_i\cdot\Psi_-+\\
& & +s_i\cdot\Psi_+\otimes\nabla_{s_i}\Psi_--\Psi_+\otimes s_i\cdot\nabla_{s_i}\Psi_-)\\
& = &
\phantom{+}\lambda L^{od}(\Psi_+\otimes\Psi_-)-\frac{1}{8}L^{od}( T\cdot\Psi_+\otimes\Psi_-+\Psi_+\otimes  T\cdot\Psi_-\\
& & -\sum\limits_i (s_i\llcorner  T)\cdot\Psi_+\otimes
s_i\cdot\Psi_--\sum\limits_i s_i\cdot\Psi_+\otimes(s_i\llcorner
T)\cdot\Psi_-)\\
& = & \phantom{-}(\lambda L^{od} -T\wedge
L^{ev})(\Psi_+\otimes\Psi_-)
\end{eqnarray*}
which proves~(\ref{closed}).

It remains to establish the $Spin(7)$-case. Here, we have to show
that the forms $e^{-F}L^{ev}_{b+}$ and $e^{-F}L^{od}_{b-}$ are
closed under $d_H$. Equivalently, we can show that
$$
dL^{ev}_+-dF\wedge L^{ev}_+=- T\wedge L^{ev}_+
$$
and
$$
dL^{od}_--dF\wedge L^{od}_-=- T\wedge L^{od}_-,
$$
which is proven in the same way as in the $G_2$-case.
\end{prf}

\begin{rmk}\hfill\newline
{\rm (i)} Similarly, we can introduce a 1-form $\alpha$ and consider
the twisted differential operator $d_{\alpha}$. This would lead to a
substitution of $dF$ by $dF+\alpha$ in the dilaton equation. If
$\alpha$ is in addition closed, such a geometry would naturally
appear as the twisted variational problem with a 1-form. However for
later applications the twisted problem with a 3-form is more natural
as it fits into the natural setup of T-duality.

{\rm (ii)} The generalised Killing and the dilaton equation with
closed $T$ occur in physics as supersymmetric solutions of
type~IIA/B supergravity in the NS-NS sector (see, for
instance,~\cite{gmpw02}).
\end{rmk}

\subsection{Geometrical properties}\label{geompropgen}

As in the classical case, the spinorial formulation is useful to
deduce global information about the geometrical properties of
integrable structures. In particular we want to compute the Ricci
tensor and to that end, we need to understand the previous theorem
again from a rather representation theoretic point of view.

Recall that the {\em torsion tensor} measures the difference between
an arbitrary metric connection $\widetilde{\nabla}$ and the
Levi-Civita connection $\nabla$ and is defined by~\cite{friv02}
$$
g(\widetilde{\nabla}_XY,Z)=g(\nabla_XY,Z)+\frac{1}{2}{\rm
Tor}(X,Y,Z).
$$
In the situation of Theorem~\ref{integrability}, we have two
principal subbundles induced by $\Psi_+$ and $\Psi_-$ inside the
spin bundle which carry metric connections $\nabla^+$ and $\nabla^-$
such that
$$
\nabla^{\pm}_X\Psi_{\pm}=\nabla_X\Psi_{\pm}\pm\frac{1}{4}(X\llcorner
T)\cdot\Psi_{\pm}.
$$
Hence these connections are induced by the metric linear connections
$$
\nabla_X^{\pm}Y=\nabla_XY\pm\frac{1}{2}T(X,Y,\cdot)
$$
and the theorem entails the torsion ${\rm Tor}_{\pm}$ of
$\nabla^{\pm}$ to be
$$
{\rm Tor}_{\pm}=\pm T.
$$
Starting with spinors, it makes therefore sense to consider a
broader class of geometries, namely geometries with two linear
metric connections whose torsion tensors are skew and sum up to 0.
In the form picture, we then ask for closed-ness of the $B$-field
free form $\rho_0=L(\Psi_+\otimes\Psi_-)$ with respect to the
twisted differential operator $d_T$. The untwisted variational
problem singles out those geometries with exact torsion as
$d_{db}\rho_0=0$ holds if and only if $d(e^b\rho_0)=0$. Hence these
geometries encapsulate all the structures we introduced and
conveniently avoid distinguishing between ``internal" torsion $db$
coming from the ubiquitous $B$-field and ``external" torsion $H$
which might or might not be present. Consequently, we refer to $T$
as the {\em torsion} of the generalised exceptional structure and
regard it as a part of the intrinsic data of an integrable
structure. Connections with skew-symmetric torsion have been
systematically studied in the recent mathematical literature,
e.g.~\cite{agfr04},~\cite{friv02},~\cite{friv03} and~\cite{iv01}.
The case of a $G_2$-connection with closed skew-symmetric torsion
has been considered in~\cite{chsw04}, where such $G_2$-structures
are called {\em strong} (conflicting with our notion of strong
integrability which is why we refer to these structures as {\em
closed}).

Again, we call such a geometry {\em non-trivial} if it has
non-trivial torsion. A generalised geometry which is trivial can be
completely understood in terms of classical geometries. If the
torsion vanishes, then both spinors $\Psi_+$ and $\Psi_-$ are
parallel with respect to the Levi-Civita connection. Consequently,
the holonomy of the underlying Riemannian manifold reduces to $G_2$
or $Spin(7)$. In the case of a generalised $G_2$-structure, it
follows in particular that the underlying topological generalised
$G_2$-structure cannot be exotic (cf. Section~\ref{topges}). If the
spinors were to be linearly dependent at one point, covariant
constancy would imply global linear dependency and we would have an
ordinary manifold of holonomy $G_2$. If the two spinors are linearly
independent at some (and hence at any) point, then the holonomy
reduces to an $SU(3)$-principal fibre bundle which is the
intersection of the two $G_2$-structures. In this case, $M$ is
locally isometric to $CY^3\times N^1$ where $CY^3$ is a Calabi-Yau
3-fold. Similarly, if the two spinors defining a generalised
$Spin(7)$-structure of even type are linearly dependent at one
point, then the holonomy of the underlying Riemannian manifold is
contained in $Spin(7)$. Otherwise, the holonomy reduces to a
subgroup of $Spin(6)=SU(4)$. For structures of odd type, vanishing
torsion means that the holonomy reduces to $G_2$ and thus $M^8$ is
locally isometric to $N^7\times N^1$ where $N^7$ has holonomy $G_2$.

For clarity of exposition, we will deal with the $G_2$- and the
$Spin(7)$-cases separately.

\bigskip

\textbf{Closed integrable generalised $\mathbf{G_2}$-structures}

In order to fix the notation, we state the

\begin{prp}\label{t_comp}$\!\!\!${\rm~\cite{friv03}}\hspace{2pt}
For any $G_2$-structure with stable form $\varphi$ there exist
unique differential forms $\lambda\in\Omega^0(M)$,
$\theta\in\Omega^1(M)$, $\xi\in\Omega^2_{14}(M,\varphi)$, and
$\tau\in\Omega^3_{27}(M,\varphi)$ so that the differentials of
$\varphi$ and $\star\varphi$ are given by
$$
\begin{array}{ccl}
d\varphi & = & -\lambda\star\varphi+\frac{3}{4}\theta\wedge\varphi+\star\tau\\
\dstar\varphi & = & \theta\wedge\star\varphi+\xi\wedge\varphi.
\end{array}
$$
\end{prp}

To specify the torsion tensor of a connection is in general not
sufficient to guarantee its uniqueness. However, this is true for
$G_2$-connections with skew-torsion. Using the notation of the
previous proposition, we can assert the following result.

\begin{prp}\label{torsion}$\!\!\!${\rm~\cite{friv02},~\cite{friv03}}\hspace{2pt} For a $G_2$-structure with stable form $\varphi$ the following
statements are equivalent.

{\rm (i)} The $G_2$-structure is {\em integrable}, i.e. $\xi=0$.

{\rm (ii)} There exists a unique linear connection
$\widetilde{\nabla}$ whose torsion tensor ${\rm Tor}$ is skew and
which preserves the $G_2$-structure, i.e.
$$
\widetilde{\nabla}\varphi=0.
$$

The torsion can be expressed by
\begin{eqnarray}
{\rm Tor} & = & -\star
d\varphi-\frac{7}{6}\lambda\cdot\varphi+\star(\theta\wedge\varphi)\nonumber\\
& = &
-\frac{1}{6}\lambda\cdot\varphi+\frac{1}{4}\star(\theta\wedge\varphi)-\star\tau\label{tor}.
\end{eqnarray}

Moreover, the Clifford action of the torsion 3-form on the induced
spinor $\Psi$ is
\begin{equation}\label{comp}
{\rm Tor}\cdot\Psi=\frac{7}{6}\lambda\Psi-\theta\cdot\Psi.
\end{equation}
\end{prp}

Using the labeling of Proposition~\ref{t_comp} with additional
subscripts $\pm$ to indicate the torsion forms of $\nabla^{\pm}$,
equations~(\ref{tor}) and~(\ref{comp}) read
$$
{\rm Tor}_{\pm}=\pm
T=-\frac{1}{6}\lambda_{\pm}\cdot\varphi_{\pm}+\frac{1}{4}\star(\theta_{\pm}\wedge\varphi_{\pm})-\star\tau_{\pm}
$$
and
$$
{\rm Tor}_{\pm}\cdot\Psi_{\pm}=\pm
T\cdot\Psi_{\pm}=\frac{7}{6}\lambda_{\pm}\Psi_{\pm}-\theta_{\pm}\cdot\Psi_{\pm}.
$$
Note also the general $G_2$-identity
$(X\llcorner\star\varphi)\cdot\Psi=4X\cdot\Psi$. In view of Theorem
\ref{integrability} we can use the dilaton equation to relate the
torsion components to the additional parameters $dF$ and $\lambda$.
We have
$$
\pm T\cdot\Psi_{\pm}=\mp2\lambda\Psi_{\pm}-2dF\cdot\Psi_{\pm}
$$
if the structure is even and
$$
\pm T\cdot\Psi_{\pm}=-2\lambda\Psi_{\pm}-2dF\cdot\Psi_{\pm}
$$
if the structure is odd.

\begin{cor}\label{torsioncomp}
If the generalised $G_2$-structure is weakly integrable, then there
exist two linear connections $\nabla^{\pm}$ preserving the
$G_{2\pm}$-structure with skew torsion $\pm T$. These connections
are uniquely determined. The differentials of the stable forms
$\varphi_{\pm}$ and $\star\varphi_{\pm}$ are given as follows. If
the structure is weakly integrable and of even type, then
$$
\begin{array}{ccl}
d\varphi_+ & = & \frac{12}{7}\lambda\star\varphi_++\frac{3}{2}dF\wedge\varphi_+-\star T_{27+}\\
\dstar\varphi_+ & = & 2dF\wedge\star\varphi_+
\end{array}
$$
and
$$
\begin{array}{ccl}
d\varphi_- & = & -\frac{12}{7}\lambda\star\varphi_-+\frac{3}{2}dF\wedge\varphi_-+\star T_{27-}\\
\dstar\varphi_- & = & 2dF\wedge\star\varphi_-,
\end{array}
$$
where $T_{27\pm}$ denotes the projection of $T$ onto
$\Omega^3_{27}(M,\varphi_{\pm})$. Moreover, the torsion can be
expressed by the formula
\begin{equation}\label{tor1}
{\rm Tor}_{\pm}=\pm T=-e^{2F}\star
de^{-2F}\varphi_{\pm}\pm2\lambda\cdot\varphi_{\pm}.
\end{equation}
If the structure is weakly integrable and of odd type, then
$$
\begin{array}{ccl}
d\varphi_+ & = & \frac{12}{7}\lambda\star\varphi_++\frac{3}{2}dF\wedge\varphi_+-\star T_{27+}\\
\dstar\varphi_+ & = & 2dF\wedge\star\varphi_+
\end{array}
$$
and
$$
\begin{array}{ccl}
d\varphi_- & = & \frac{12}{7}\lambda\star\varphi_-+\frac{3}{2}dF\wedge\varphi_-+\star T_{27-}\\
\dstar\varphi_- & = & 2dF\wedge\star\varphi_-.
\end{array}
$$
The torsion can be expressed by the formula
\begin{equation}\label{tor2}
{\rm Tor}_{\pm}=\pm T=-e^{2F}\star
de^{-2F}\varphi_{\pm}+2\lambda\cdot\varphi_{\pm}.
\end{equation}

The strongly integrable case follows if we set $\lambda=0$.

Conversely, if we are given two $G_2$-structures defined by the
stable forms $\varphi_+$ and $\varphi_-$ inducing the same metric, a
constant $\lambda$ and a function $F$ such that~(\ref{tor1}) or
(\ref{tor2}) defines a (closed) 3-form $T$, then the corresponding
spinors $\Psi_{\pm}$ satisfy the integrability condition of Theorem
\ref{integrability} and hence define a (closed) integrable
generalised $G_2$-structure of even or odd type.
\end{cor}

Depending on the torsion components which are turned on, we will
also say that the $G_{2\pm}$-structure has torsion in the components
$\mathbf{1}_{\pm}$, $\mathbf{7}_{\pm}$ or $\mathbf{27}_{\pm}$.

Next we investigate closed structures only. The
expressions~(\ref{tor1}) and~(\ref{tor2}) for the torsion form $T$
have an interesting consequence (see also~\cite{gmpw02} and
\cite{ivpa01}). Assume $M$ to be compact and endowed with a closed
weakly integrable structure of even type. Then~(\ref{tor1}) and
Stokes' Theorem imply
$$
\int_Me^{-2F}T\wedge \star T=\mp\int_MT\wedge
d(e^{-2F}\varphi_{\pm})+2\lambda\int_MT\wedge
e^{-2F}\star\varphi_{\pm}=\frac{4}{7}\lambda^2\int_Me^{-2F}vol_M.
$$
Here we have used that $dT=0$ and the fact that the projection of
$T$ on $\varphi_{\pm}$ is given by
$T_{1{\pm}}=2\lambda\varphi_{\pm}$. Since the left hand side is
strictly positive (unless $T\equiv 0$), we need $\lambda\not=0$. If
we start with an odd structure instead, we also find
$$
\int_Me^{-2F}T\wedge \star T=\mp\int_MT\wedge
d(e^{-2F}\varphi_{\pm})\pm2\lambda\int_MT\wedge
e^{-2F}\star\varphi_{\pm}=\frac{4}{7}\lambda^2\int_Me^{-2F}vol_M,
$$
since now $T_{1{\pm}}=\pm2\lambda\varphi_{\pm}$.

Consequently, we obtain the following no-go theorem.

\begin{cor}\label{vanishingg2}
If $M$ is compact and carries a closed integrable generalised
$G_2$-structure, then $T=0$ if and only if $\lambda=0$. In that case
the spinors $\Psi_{\pm}$ are parallel with respect to the
Levi-Civita connection.
\end{cor}

Next we compute the Ricci tensor. As a corollary, we will be able to
exclude the existence of interesting closed strongly integrable
examples on homogeneous spaces.

To begin with, let $\ric$ and $\ric^{\pm}$ denote the Ricci tensor
associated with the Levi-Civita connection and the connections
$\nabla^{\pm}$. The relationship between $\ric$ and $\ric^{\pm}$ was
analysed in~\cite{friv02}. Generally speaking, if we have a
$G$-structure with a $G$-preserving, metric linear connection
$\widetilde{\nabla}$ that has closed, skew torsion, and a
$G$-invariant spinor $\Psi$, then the following identities hold.

\begin{prp}\label{skew_prop}$\!\!\!${\rm~\cite{friv02}}\hspace{2pt}
The Ricci tensor associated with $\widetilde{\nabla}$ is determined
by the equation
$$
\widetilde{\ric}(X)\cdot\Psi=(\widetilde{\nabla}_X{\rm
Tor}+\frac{1}{2}X\llcorner dT)\cdot\Psi
$$
and relates to the metric Ricci tensor through
$$
\ric(X,Y)=\widetilde{\ric}(X,Y)+\frac{1}{2}d^*{\rm
Tor}(X,Y)+\frac{1}{4}g(X\llcorner {\rm Tor},Y\llcorner {\rm Tor}).
$$
\end{prp}

\begin{thm}\label{riccitensor}
The Ricci-tensor of a closed integrable generalised $G_2$-structure
is given by
$$
\ric(X,Y)=-\frac{7}{2}H^F(X,Y)+\frac{1}{4}g(X\llcorner T, Y\llcorner
T),
$$
where $H^F(X,Y)=X.Y.F-\nabla_XY.F$ is the Hessian of the dilaton
$F$.
\end{thm}

\begin{prf}
According to the previous proposition we obtain
$$
\ric(X,Y)=\frac{1}{2}(\ric^+(X,Y)+\ric^-(X,Y))+\frac{1}{4}g(X\llcorner
T,Y\llcorner T).
$$
Consequently, it remains to show that
\begin{equation}\label{ricnull}
\ric^+(X,Y)+\ric^-(X,Y)=-7H^F(X,Y).
\end{equation}
Since $\nabla^{\pm}$ preserves the $G_{2\pm}$-structures and $dT=0$,
only $T_{7\pm}$, the $\Lambda^3_{7\pm}$ component of $T$, impacts on
$\ric^{\pm}$. Hence we can write
$$
\gamma_{X\pm}\llcorner\star\varphi_{\pm}=\nabla^{\pm}_XT_{7\pm}=\pm\frac{1}{2}\nabla^{\pm}_X\star(dF\wedge\varphi_{\pm})
$$
and thus
$$
\ric^{\pm}(X)\cdot\Psi_{\pm}=\pm\nabla_X^{\pm}T\cdot\Psi_{\pm}=\pm(\gamma_{X\pm}\llcorner\star\varphi_{\pm})\cdot\Psi_{\pm}=\pm4\gamma_{X\pm}\cdot\Psi_{\pm}.
$$
Moreover, for any $G_2$-structure $(\varphi,g)$ we have
$g(X\wedge\varphi,X\wedge\varphi)=4g(X,X)$ and therefore,
\begin{eqnarray}
\ric^{\pm}(X,Y) & = & \pm4g(\gamma_{X\pm},Y)\nonumber\\
& = & \pm
g(\star(\gamma_{X\pm}\wedge\varphi_{\pm}),\star(Y\wedge\varphi_{\pm}))\nonumber\\
& = & \mp g(\gamma_{X\pm}\llcorner\star\varphi_{\pm},\star(Y\wedge\varphi_{\pm}))\nonumber\\
& = &
-\frac{1}{2}g(\nabla^{\pm}_X\star(dF\wedge\varphi_{\pm}),\star(Y\wedge\varphi_{\pm}))\nonumber\\
& = &
\phantom{\pm}\frac{1}{2}g(Y\wedge\nabla^{\pm}_X\star(dF\wedge\varphi_{\pm}),\star\varphi_{\pm}).\label{lastterm}
\end{eqnarray}
We shall now derive~(\ref{ricnull}) by a pointwise computation in a
frame that satisfies $\nabla_{e_i}e_j=0$, or equivalently,
$\nabla_{e_i}^{\pm}e_j=\pm\frac{1}{2}\sum_kT_{ijk}e_k$ in a fixed
point $x$. Continuing with~(\ref{lastterm}),
\begin{eqnarray*}
\ric^{\pm}(e_i,e_j) & = &
\phantom{-}\frac{1}{2}g(e_j\wedge\nabla^{\pm}_{e_i}\star(dF\wedge\varphi_{\pm}),\star\varphi_{\pm})\nonumber\\
& = & \phantom{-}\frac{1}{2}g(\nabla^{\pm}_{e_i}(e_j\wedge
\star(dF\wedge\varphi_{\pm}))-\nabla^{\pm}_{e_i}e_j\wedge
\star(dF\wedge\varphi_{\pm}),\star\varphi_{\pm})\nonumber\\
& = & \phantom{-}\frac{1}{2}e_i.g(e_j\wedge
\star(dF\wedge\varphi_{\pm}),\star\varphi_{\pm})-\frac{1}{2}g(\nabla^{\pm}_{e_i}e_j\wedge
\star(dF\wedge\varphi_{\pm}),\star\varphi_{\pm})\nonumber\\
& = & -\frac{7}{2}e_i.e_j.F\pm\frac{7}{4}\sum_kT_{ijk}e_k.F.
\end{eqnarray*}
On the other hand, the Hessian at the point $x$ evaluated in this
basis is
$$
H^F(e_i,e_j)=e_i.e_j.F-\nabla_{e_i}e_j.F=e_i.e_j.F,
$$
which yields the assertion.
\end{prf}

\begin{cor}
For a closed integrable generalised $G_2$-structure, the scalar
curvature of $\nabla^{\pm}$ is
$$
{\rm Scal}^{\pm}=\frac{7}{2}\Delta F,
$$
where $\Delta(\cdot)=-{\rm Tr_g}H^{(\cdot)}$ is the Laplacian
associated with the metric $g$. Hence, the scalar curvature of the
metric $g$ is given by
$$
{\rm Scal}={\rm Tr}(\ric)=\frac{7}{2}\Delta F+\frac{1}{4}\norm{T}^2.
$$
\end{cor}

Since $T=0$ implies the covariant constancy of the spinors
$\Psi_{\pm}$ with respect to the Levi-Civita connection, an
integrable structure with vanishing torsion is necessarily
Ricci-flat. It is interesting to note however, that already the
vanishing of the dilaton $F$ is sufficient.

\begin{prp}\label{ricciflat}
The metric $g$ of a closed strongly integrable generalised
$G_2$-structure is
$$
\mbox{\rm Ricci-flat}\Leftrightarrow dF=0 \Leftrightarrow T=0.
$$
\end{prp}

\begin{prf}
The only non-trivial statement to prove is that a constant dilaton
causes Ricci-flatness. But by Corollary~\ref{torsioncomp} the
condition $dF=0$ implies that the underlying $G_{2\pm}$-structures
are co-calibrated which means $\dstar\varphi_{\pm}=0$. As we can see
earlier in the proof of Theorem~\ref{riccitensor}, the Ricci tensors
$\ric^{\pm}$ also vanish. Consequently, we can appeal to Theorem 5.4
of~\cite{friv02} which asserts that under these circumstances the
underlying $G_{2\pm}$-structures actually have holonomy contained in
$G_2$ if $T_{1\pm}$ or $T_{27\pm}$ vanishes. Since we assumed strong
integrability, $T_{1\pm}=0$ and the assertion follows.
\end{prf}

\begin{rmk}

(i) The proposition fails if the torsion is not closed as an example
in Section~\ref{exgeneral} will show.

(ii) We can always get rid of the dilaton $F$ by the conformal
transformation $\tilde{g}=e^{-F}g$,
$\tilde{\varphi}_{\pm}=e^{-\frac{3}{2}F}\varphi_{\pm}$. The
transformed torsion 3-form $\widetilde{\rm Tor}_{\pm}$ is given by
$$
\widetilde{\rm
Tor}_{\pm}=\pm\widetilde{T}=e^{-2F}\cdot(T-\frac{1}{2}\star
dF\wedge\varphi_{\pm})
$$
(\cite{friv03}, Corollary 4.2). Since the resulting
$G_{2\pm}$-structures are co-calibrated, they have non-positive
scalar curvature. However, the previous proposition does not apply,
because $\widetilde{T}$ is no longer closed.
\end{rmk}

Since for a homogeneous generalised $G_2$-structure the dilaton $F$
is necessarily constant, the previous proposition implies the
following

\begin{cor}
Any homogeneous closed strongly integrable generalised
$G_2$-structure is Ricci-flat and hence both underlying topological
$G_2$-structures have actually holonomy $G_2$.
\end{cor}

\bigskip

\textbf{Closed integrable $\mathbf{Spin(7)}$-structures}

We now turn to integrable generalised $Spin(7)$-structures which we
discuss along the lines of the generalised $G_2$-case. We summarise
the results we need in the next proposition.

\begin{prp}$\!\!\!${\rm~\cite{iv01}}\hspace{2pt}
Let $(M^8,\Omega)$ be a $Spin(7)$-manifold. There
exist unique forms $\theta\in\Omega^1(M)$ and
$\tau\in\Omega^5_{48}(M)$ such that
$$
d\Omega=\theta\wedge\Omega+\tau.
$$
In particular,
$$
\theta=\frac{1}{7}\star(d^*\Omega\wedge\Omega).
$$
There exists a unique linear metric connection $\widetilde{\nabla}$
with skew-symmetric torsion which preserves the $Spin(7)$-structure,
i.e. $\widetilde{\nabla}\Omega=0$. It is determined through
\begin{equation}\label{torspin7}
{\rm Tor}=-d^*\Omega-\frac{7}{6}\star(\theta\wedge\Omega)
\end{equation}
and the action of ${\rm Tor}$ on the corresponding spinor $\Psi$ is
given by
$$
{\rm Tor}\cdot\Psi=-\frac{7}{6}\theta\cdot\Psi.
$$
Finally, the following expressions hold for the Riemannian scalar
curvature ${\rm Scal}$ and the scalar curvature $\widetilde{\rm
Scal}$ of $\widetilde{\nabla}$
\begin{equation}\label{scalspin7}
{\rm Scal}=\frac{49}{18}\norm{\theta}^2-\frac{1}{12}\norm{{\rm
Tor}}^2+\frac{7}{2}d^*\theta,\quad \widetilde{\rm Scal}=
\frac{49}{18}\norm{\theta}^2-\frac{1}{3}\norm{{\rm
Tor}}^2+\frac{7}{2}d^*\theta.
\end{equation}
\end{prp}

For an integrable generalised $Spin(7)$-structure we have again,
according to Theorem~\ref{integrability}, two metric connections
$\nabla^{\pm}$ preserving the $Spin(7)_{\pm}$-structure, and whose
torsion is ${\rm Tor}_{\pm}=\pm T$. The dilaton equation then
implies
$$
(dF\pm\frac{1}{2}T)\cdot\Psi_{\pm}=(dF-\frac{7}{12}\theta_{\pm})\cdot\Psi_{\pm}=0,
$$
and hence
$$
\theta_{\pm}=\frac{12}{7}dF.
$$
As a consequence, formula~(\ref{torspin7}) can be written in the
more succinct form as given in

\begin{cor}
If the generalised $Spin(7)$-structure is integrable, then
\begin{equation}\label{torspin7cor}
\pm T=e^{2F}\star de^{-2F}\Omega_{\pm}
\end{equation}
and
\begin{equation}\label{domegaspin7}
d\Omega_{\pm}=\frac{12}{7}dF\wedge\Omega_{\pm}\pm\star T_{48\pm},
\end{equation}
where $T_{48\pm}$ denotes the projection of $T$ onto
$\Omega^3_{48}(M,\Omega_{\pm})$.

Conversely, if we are given two $Spin(7)_{\pm}$-invariant forms
inducing the same metric $g$, a function $F$ and a closed 3-form $T$
such that~(\ref{torspin7cor}) and~(\ref{domegaspin7}) hold, then the
corresponding spinors $\Psi_{\pm}$ satisfy Theorem
\ref{integrability} and hence define an integrable generalised
$Spin(7)$-structure.
\end{cor}

This implies again a vanishing theorem over compact manifolds if $T$
is closed since
$$
\int_Me^{-2F}T\wedge\star T=\pm\int_MT\wedge d(e^{-2F}\Omega_{\pm}).
$$

\begin{cor}\label{vanishing}
If $M$ is compact and carries a closed integrable generalised
$Spin(7)$-structure, then $T=0$ and consequently, the spinors
$\Psi_{\pm}$ are parallel with respect to the Levi-Civita
connection.
\end{cor}

Next we compute the Ricci tensor.

\begin{thm}\label{riccitensorspin7}
The Ricci-tensor of an integrable generalised $Spin(7)$-structure is
given by
$$
\ric(X,Y)=-2H^F(X,Y)+\frac{1}{4}g(X\llcorner T, Y\llcorner T),
$$
where $H^F(X,Y)=X.Y.F-\nabla_XY.F$ is the Hessian of the dilaton
$F$.
\end{thm}

\begin{prf}
The proof is essentially the same as for Theorem~\ref{riccitensor},
that is we need to compute the sum $\ric^++\ric^-$. Again we fix a
local frame $e_1,\ldots,e_8$ that satisfies $\nabla_{e_i}e_j=0$ at a
point. By proposition 7.2. of~\cite{iv01},
$$
\ric^{\pm}(X)=\mp\star(\nabla_X^{\pm}T\wedge\Omega_{\pm}).
$$
Since the connections $\nabla^{\pm}$ preserve the $Spin(7)_{\pm}$
structure,
\begin{eqnarray*}
\ric^{\pm}(e_i,e_j) & = & \pm g(\nabla^{\pm}_{e_i}T\wedge\Omega_{\pm},\star e_j)\\
& = & \pm \frac{1}{7}g(\nabla^{\pm}_{e_i}T\wedge\Omega_{\pm},(e_j\llcorner\Omega_{\pm})\wedge\Omega_{\pm})\\
& = & \pm g(\nabla^{\pm}_{e_i}T,e_j\llcorner\Omega_{\pm})\\
& = & -\frac{2}{7}g(e_j\wedge\nabla^{\pm}_{e_i}\star(dF\wedge\Omega_{\pm}),\Omega_{\pm})\\
& = & -\frac{2}{7}e_i.g(e_j\wedge
\star(dF\wedge\Omega_{\pm}),\Omega_{\pm})+\frac{2}{7}g(\nabla^{\pm}_{e_i}e_j\wedge
\star(dF\wedge\Omega_{\pm}),\Omega_{\pm})\\
& = &
-\frac{2}{7}e_i.g(dF\wedge\Omega_{\pm},e_j\wedge\Omega_{\pm})+\frac{2}{7}\sum_kT_{ijk}g(dF\wedge\Omega_{\pm},e_k\wedge\Omega_{\pm})\\
& = & -2e_i.e_j.F\pm\sum_kT_{ijk}e_k.F,
\end{eqnarray*}
and the result follows as in the $G_2$-case.
\end{prf}

\begin{cor}
For a closed integrable generalised $Spin(7)$-structure, the scalar
curvature of $\nabla^{\pm}$ is
$$
{\rm Scal}^{\pm}=2\Delta F=\frac{7}{6}d^*\theta_{\pm},
$$
where $\Delta(\cdot)=-{\rm Tr}_gH^{(\cdot)}$ is the Laplacian on
functions associated with the metric $g$. Hence, the scalar
curvature of $g$ is given by
$$
{\rm Scal}={\rm Tr}(\ric)=2\Delta F+\frac{1}{4}\norm{T}^2.
$$
\end{cor}

Finally we can again exclude the possibility of any interesting
homogeneous structures.

\begin{prp}\label{ricciflatspin7}
A closed integrable generalised $Spin(7)$-structure is
$$
\mbox{\rm Ricci-flat}\Leftrightarrow dF=0 \Leftrightarrow T=0.
$$
\end{prp}

\begin{prf}
The previous corollary shows that if $dF=0$, then the scalar
curvature ${\rm Scal}^{\pm}$ vanishes and ${\rm
Scal}=\norm{T}^2/4\geq0$. But~(\ref{scalspin7}) implies ${\rm
Scal}=-\norm{T}^2/12\leq0$, hence $\norm{T}=0$. The remaining
implications are obvious.
\end{prf}

\begin{rmk}
Note that as in the $G_2$-case, we can scale $dF$ away by the
conformal transformation $\tilde{g}=e^{-\frac{3}{7}F}g$,
$\widetilde{\Omega}_{\pm}=e^{-\frac{12}{7}F}$, but we obtain again a
torsion form which is not closed~\cite{iv01}.
\end{rmk}

\begin{cor}
Any homogeneous closed integrable generalised $Spin(7)$-structure is
Ricci-flat and hence both underlying topological
$Spin(7)$-structures have actually holonomy $Spin(7)$.
\end{cor}

\subsection{Examples}\label{exgeneral}

\textbf{Straight structures}

If we want to construct explicit examples of the various geometries
we introduced, then it is most natural to adopt the form point of
view and to write the structure forms in terms of the underlying
$SU(3)$- or $G_2$-structures as at the end of Section~\ref{ges}.
Imposing the integrability condition on the structure form yields
differential conditions on the $SU(3)$- or $G_2$-invariants which
can then be discussed in the terms of the respective representation
theory.

\begin{ex}({\em strongly integrable straight structures})\hfill\newline
We first look for strongly integrable generalised structures which
are {\em straight}, that is structures defined by one spinor
$\Psi=\Psi_+=\Psi_-$ (cf. Section~\ref{topges}).

(i) ({\em straight generalised $G_2$-structures}) Let us consider a
generalised $G_2$-structure which is induced by a classical
$G_2$-structure $(M^7,\varphi)$. Given a (closed) 3-form $T$ and a
function $F$, we want to solve the equations of strong integrability
$$
d_Te^F(1-\star\varphi)=0,\quad d_Te^F(-\varphi+vol_g)=0.
$$
Expanding and ordering by degree, these equations hold if and only
if
\begin{eqnarray*}
dF+T-dF\wedge\star\varphi-T\wedge\star\varphi & = & 0\\
dF\wedge\varphi+d\varphi+T\wedge\varphi & = & 0
\end{eqnarray*}
which is equivalent to
$$
dF=0,\: T=0,\:d\varphi=0\mbox{ and } \dstar\varphi=0.
$$
In particular, the holonomy group of $(M,g_{\varphi})$ is contained
in $G_2$~\cite{fegr82}.

(ii) ({\em straight generalised $Spin(7)$-structures, even type})
Similarly, we can consider a generalised $Spin(7)$-structure of even
type coming from a $Spin(7)$-manifold $(M^8,\Omega)$ where strong
integrability
$$
d_Te^F(1-\Omega+vol_g)=0
$$
is equivalent to
$$
dF=0,\:T=0\mbox{ and }d\Omega=0.
$$
Again this causes the holonomy of $(M,g_{\Omega})$ to reduce to
$Spin(7)$~\cite{br87}.

(iii) ({\em straight generalised $Spin(7)$-structures, odd type}) A
trivial example of odd type is provided by product manifolds of the
form $N^7\times\R$, where $N^7$ carries a metric of holonomy $G_2$.
According to Corollary~\ref{normalformcor1}, the form
$dt\wedge(-1+\star\varphi)-\varphi+vol_N$ defines a generalised
$Spin(7)$-structure of odd type which is clearly closed. Multiplying
with a dilaton and considering the twisted differential $d_T$
implies once more $dF=0$ and $T=0$.
\end{ex}

As a result, we obtain the

\begin{prp}
Any straight generalised $G_2$- or $Spin(7)$-structure is strongly
integrable if and only if the holonomy of the underlying Riemannian
manifold is contained in $G_2$ or $Spin(7)$. Such a structure is
necessarily trivial (i.e. $T=0$) and therefore closed.
\end{prp}

\bigskip

\textbf{Compact examples of strongly integrable generalised
structures}

Here is a non-trivial compact example with constant dilaton.

\begin{ex}({\em compact strongly integrable generalised
$G_2$-structures}) Consider the 6-dimensional nilmanifold $N$
associated with the Lie algebra $\mf{g}$ spanned by the orthonormal
basis $e_2,\ldots,e_7$ determined by the relations
$$
de_i=\left\{\begin{array}{cl}\phantom{-}e_{37},&\,i=4\\-e_{35},&\,i=6\\\phantom{-}0,&\,\mbox{else.}\end{array}\right.
$$
We then define $M=N\times S^1$ which we endow with the product
metric $g=g_N+dt\otimes dt$. On $N$ we choose the $SU(3)$-structure
coming from
$$
\omega=-e_{23}-e_{45}+e_{67},\quad\psi_+=e_{356}-e_{347}-e_{257}-e_{246}
$$
(cf. Section~\ref{su3}) which induces a generalised $G_2$-structure
on $M$ with $\alpha=e_1=dt$. We let $T=e_{167}+e_{145}$ and consider
the equations of strong integrability
$$
d_T\rho=0,\quad d_T\hat{\rho}=0,
$$
where $\rho=\omega+\psi_+\wedge\alpha-\omega^3/6$ and
$\hat{\rho}=\alpha-\psi_--\omega^2\wedge\alpha/2$ (using the
notation of Section~\ref{ges}). These are equivalent to
$$
d\omega=0,\quad d\psi_+\wedge\alpha=-T\wedge\omega,\quad
T\wedge\alpha=d\psi_-,\quad T\wedge\psi_-=0.
$$
By design, $T\wedge\alpha=0$ and $T\wedge\psi_-=0$. Moreover,
$$
d\omega=-de_4\wedge e_5+de_6\wedge e_7=0
$$
and
$$
d\psi_-=-e_3\wedge de_4\wedge e_6+e_{34}\wedge de_6+e_{25}\wedge
de_6+e_2\wedge de_4\wedge e_7=0.
$$
Finally, we have
$$
d\psi_+\wedge\alpha=-e_{12367}-e_{12345}=-T\wedge\omega.
$$
Note that $dT=2e_{1357}$ in accordance with
Corollary~\ref{vanishingg2}.
\end{ex}

From the next lemma and the previous example we immediately infer
the existence of a compact integrable generalised
$Spin(7)$-structures of even type.

\begin{lem}\label{g2spin7}
If $(M,\rho_0,T)$ is a (closed) strongly integrable generalised
$G_2$-manifold with $\rho_0$ even, then $(M^7\times
S^1,dt\wedge\hat{\rho}_0+\rho_0)$ is a (closed) integrable
generalised $Spin(7)$-manifold of even type.
\end{lem}

\begin{prf}
According to Corollary~\ref{normalformcor1},
$\rho=dt\wedge\hat{\rho}_0+\rho_0$ defines a generalised
$Spin(7)$-structure and it is immediate to check that $d_T\rho=0$
and $d_T\hat{\rho}=0$ imply $d_T\rho=0$.
\end{prf}

In order to construct an odd example we start with a compact
calibrated $G_2$-manifold $(M^7,\varphi)$, i.e. $d\varphi=0$. An
example can be easily obtained by taking the nilmanifold $N$ used
for the compact generalised $G_2$-example. We swap orientations so
that $d\psi_+=0$ and $d\psi_-\wedge dt\not=0$, hence
$$
d\varphi=d\omega\wedge dt+d\psi_+=0,\quad
d\star\varphi=d\psi_-\wedge dt+d\omega\wedge\omega=d\psi_-\wedge dt.
$$
Write $d\star\varphi=\xi\wedge\varphi$ and put $T=-dt\wedge\xi$.
Corollary~\ref{normalformcor1} implies that
$$
\rho=dt\wedge(-1+\star_M\varphi)-\varphi+vol_M
$$
defines an odd generalised $Spin(7)$-structure on $M^7\times S^1$
equipped with the product metric $g=g_M+dt\otimes dt$. Now
\begin{eqnarray*}
d_T\rho & = & -dt\wedge d\star_M\varphi-d\varphi+T\wedge\rho\\
& = & -dt\wedge\xi\wedge\varphi-T\wedge\varphi\\
& = & 0.
\end{eqnarray*}
In particular, it follows that the torsion of a compact calibrated
$G_2$-manifold can never be closed.

\begin{prp}
There exist compact strongly integrable generalised $G_2$- and
$Spin(7)$-manifolds of either type. In particular, we see that in
Corollaries~\ref{vanishingg2},~\ref{vanishing} and
Propositions~\ref{ricciflat},~\ref{ricciflatspin7} it is necessary
to assume $T$ to be closed.
\end{prp}

\bigskip

\textbf{Weakly integrable structures}

Next we look for weakly integrable $G_2$-structures. The natural
idea would be to use nearly parallel $G_2$-manifolds to produce weak
analogues of straight strongly integrable structures. However, any
such construction results to be integrable if and only if it is
strongly integrable. In this sense weak structures really define a
new type of geometry without a classical, i.e. straight counterpart.

To see this, assume we were given a straight weakly integrable
structure of even type, i.e. $\varphi_+=\varphi_-$ and $\pm
T=-e^{2F}\star d(e^{-2F}\varphi)\pm2\lambda\varphi$
(Corollary~\ref{torsioncomp}). In particular, we have $\star
de^{-2F}\varphi=0$ which is equivalent to
$d\varphi=2dF\wedge\varphi$. It follows $\lambda=0$ in contradiction
to our initial assumption. For a straight structure of odd type,
Corollary~\ref{torsioncomp} implies at once $T=0$ and hence
$\lambda=0$.

The absence of a straight example already hints at the fact that
examples, if they exist, are quite hard to find. It is instructive
to see where the difficulties arise.

\begin{ex}({\em ansatz for weakly integrable structures with constant dilaton})\hfill\newline
(i) ({\em even type, homogeneous case}) First we consider the even
case, that is $\rho$ is even and satisfies
$d_T\rho=\lambda\hat{\rho}$ for some closed 3-form $T$. Moreover, we
require the generalised $G_2$-structure to be homogeneous. Since the
only invariant 0-form is a constant, the 1-component of $\hat{\rho}$
has to vanish which in turn implies that $\sin(a)=0$ (cf.
Proposition~\ref{normalform}). The equation for weak integrability
becomes
$$
d_T(1-\star\varphi)=\lambda(-\varphi+vol_g),
$$
implying $T=-\lambda\varphi$ and
$-\lambda\varphi\wedge\star\varphi=\lambda vol_g$. Since
$\varphi\wedge\star\varphi=7vol_g$, we have $\lambda=0$ and
therefore $T=0$ and $\dstar\varphi=0$. Hence we obtain a classical
co-calibrated $G_2$-structure with $\lambda=0$ which does not define
a weakly integrable structure.

(ii) ({\em odd type}) For any odd structure $\rho\in\Omega^{od}(M)$
we have $d_T\rho=\lambda\hat{\rho}$ and it follows immediately that
the 0-form of $\hat{\rho}$ vanishes. Hence $\cos(a)\equiv 0$ and the
two spinors $\Psi_+$ and $\Psi_-$ are orthogonal. Written in terms
of the underlying $SU(3)$-structure, the corresponding stable
3-forms $\varphi_{\pm}$ are given by
$$
\varphi_{\pm}=\omega\wedge\alpha\pm\psi_+,
$$
cf.~(\ref{varphi+-}). The odd and even forms $\rho$ and $\hat{\rho}$
are now $\rho=\alpha-\psi_--\omega^2\wedge\alpha/2$ and
$\hat{\rho}=\omega+\psi_+\wedge\alpha-\omega^3/6$. The weak
integrability condition becomes
$$
d\alpha-d\psi_--\frac{1}{2}d(\omega^2\wedge\alpha)+T\wedge(\alpha-\psi_--\frac{1}{2}\omega^2\wedge\alpha)=\lambda
(\omega+\psi_+\wedge\alpha-\frac{1}{6}\omega^3)
$$
which is equivalent to
\begin{eqnarray*}
d\alpha_{\phantom{-}} & = & \phantom{-}\lambda\omega\\
T\wedge\alpha-d\psi_- & = &
\lambda\psi_+\wedge\alpha\\
T\wedge\psi_-+\frac{1}{2}d(\omega^2\wedge\alpha)_{\phantom{-}} & = &
\frac{\lambda}{6}\omega^3.
\end{eqnarray*}
Using $d\alpha=\lambda\omega$ finally yields
\begin{eqnarray}
d\alpha & = & \phantom{-}\lambda\omega\nonumber\\
T\wedge\alpha_{\phantom{-}}& = &
\phantom{-}\lambda\psi_+\wedge\alpha+d\psi_- \nonumber\\
T\wedge\psi_- & = & -\frac{1}{3}\lambda\omega^3\label{weaksol}.
\end{eqnarray}
The shape of these equations suggest to start with a compact
Calabi-Yau manifold with integral K\"ahler class in order to
construct a weakly integrable structure over the associated
$S^1$-principal fibre bundle $P_{\omega}$ endowed with a connection
1-form $\theta=\alpha/\lambda$ such that $d\theta=\omega$. However,
the second equation of~(\ref{weaksol}) implies $T=\lambda\psi_+$
(modulo possible components in the kernel of $\wedge\alpha$), but
this fails to satisfy the last equation as
$\psi_+\wedge\psi_-=2\omega^3/3$. Again, the algebraic constraints
imply $\lambda=0$.
\end{ex}

The examples we have found so far are, except for the compact one,
rather trivial. However, they are still valuable as the device of
T-duality permits us to go from a trivial straight solution existing
over a non-trivial principal $S^1$-bundle to a non-trivial solution
over a trivial $S^1$-bundle. This tool will occupy us next.

\bigskip

\textbf{T-duality}

T-duality is the name of a procedure in string theory which consists
in changing the topology of a model which is an $S^1$-fibration (or
more generally, a $T^n$-fibration -- the ``T" therefore referring to
a ``torus") by replacing this fibre by a different one of same type,
but without destroying integrability. It thus interchanges type IIA
with type IIB string theory. A basic instance of T-duality is the
$R\to 1/R$ invariance in string theory compactified on a circle of
radius $R$. In this way it relates strongly coupled string theories
to weakly coupled ones and opens the way for perturbative methods in
the investigation of strongly coupled models. Mathematically
speaking, T-duality transforms the data $(g,b,F)$, that is, a
generalised metric consisting of an ordinary (usually Lorentzian)
metric $g$ and a 2-form $b$ as well as the dilaton $F$, all living
on a principal $S^1$-bundle $P\to M$ with connection form $\theta$,
into a generalised metric $(g^t,b^t)$ and a dilaton $F^t$ living on
a {\em new} principal $S^1$-fiber bundle $P^t\to M$ with connection
form $\theta^t$. The coordinate description of this dualising
procedure is well-known in physics under the name of {\em Buscher
rules} (see, for instance Appendix A in~\cite{kstt02}). Topological
issues of T-duality are discussed in~\cite{bem03} which will be our
standard reference. 

We will briefly introduce the formalism and set up the notation as
far as we will need it here. Consider a principal $S^1$-bundle
$p:P\to M$ which carries a connection form $\theta$. We denote by
$X$ the vertical vector field dual to $\theta$, so
$X\llcorner\theta=1$, and by $\mc{F}$ the curvature 2-form which we
regard as a 2-form on $M$, that is $d\theta=p^*\mc{F}$. Moreover, we
assume to be given an $S^1$-invariant closed integral 3-form $T$
such that the 2-form $\mc{F}^t$ defined by
$$
p^*\mc{F}^t=-X\llcorner T
$$
is also closed and integral. In practice we shall often assume
$T=0$. Integrality of $\mc{F}^t$ ensures the existence of another
principal $S^1$-bundle $P^t$, the {\em T-dual} of $P$ defined by the
choice of a connection form $\theta^t$ with $d\theta^t=p^*\mc{F}^t$.
Here and from now on, we ease notation and drop the pull-back $p^*$
bearing in mind that the forms $\mc{F}$ and $\mc{F}^t$ etc. live on
$M$. Writing $T=-\theta\wedge \mc{F}^t+\mc{T}$ for a 3-form
$\mc{T}\in\Omega^3(M)$, we define the T-dual of $T$ by
$$
T^t=-\theta^t\wedge\mc{F}+\mc{T}.
$$
Note that $T^t$ is closed and integral if and only if $T$ is closed
and integral.

To make contact with our situation, consider an $S^1$-invariant
exceptional generalised structure defined by $\rho$ which we write
$$
\rho=\theta\wedge\rho_0+\rho_1.
$$
The {\em T-dual of} $\rho$ is defined to be
$$
\rho^t=\theta^t\wedge\rho_1+\rho_0.
$$
In particular, T-duality reverses the parity of forms and maps even
to odd and odd to even forms. It is enacted by multiplication with
the element $X\oplus\theta\in Pin(n,n)$ on $\rho$ (giving the form
$(X\oplus\theta)\bullet\rho=\theta\wedge\rho_1+\rho_0$) followed by
the substitution $\theta\to\theta^t$.

The crucial feature of T-duality we shall need in the sequel is that
it preserves the $Spin(n,n)$-orbit structure on $\Lambda^{ev,od}$.
To see this, we decompose
$$
TP\oplus T^*P\cong T\oplus\R X\oplus T^*M\oplus\R\theta,\quad
TP\oplus T^*P\cong T\oplus\R X^t\oplus T^*M\oplus\R\theta^t.
$$
Then consider the map $\tau:TP\oplus T^*P\to TP^t\oplus
T^*P^t$~\cite{ca05} defined with respect to this splitting by
$$
\tau(V+uX\oplus\xi +v\theta)=-V+vX^t\oplus-\xi+u\theta^t.
$$
It satisfies
$$
(a\bullet\rho)^t=\tau(a)\bullet\rho^t
$$
for any $a\in TP\oplus T^*P$ and in particular, $\tau(a)^2=-(a,a)$:
\begin{eqnarray*}
((V+uX\oplus\xi +v\theta)\bullet\rho)^t & = &
(-\theta\wedge(V\llcorner\rho_0)+V\llcorner\rho_1+u\rho_0-\theta\wedge\xi\wedge\rho_0+\xi\wedge\rho_1+v\theta\wedge\rho_1)^t\\
& = & (\theta\wedge(-V\llcorner\rho_0-\xi\wedge\rho_0+v\rho_1)+V\llcorner\rho_1+u\rho_0+\xi\wedge\rho_1)^t\\
& = & \phantom{(}\theta^t\wedge(V\llcorner\rho_1+u\rho_0+\xi\wedge\rho_1)-V\llcorner\rho_0-\xi\wedge\rho_0+v\rho_1\\
& = & (-V+vX^t\oplus-\xi+u\theta^t)\bullet\rho^t.
\end{eqnarray*}
Hence this map extends to an isomorphism $\cliff(TP\oplus
T^*P)\cong\cliff(TP^t\oplus T^*P^t)$ and any orbit of the form
$Spin(TP\oplus T^*P)/G$ gets mapped to an equivalent orbit
$Spin(TP^t\oplus T^*P^t)/G^t$ where $G$ and $G^t$ are isomorphic as
abstract groups.

As an illustration of this, consider a generalised $G_2$-structure
over $P$ with structure form $\rho=\theta\wedge\rho_0+\rho_1$ and
companion $\hat{\rho}=\theta\wedge\hat{\rho_0}+\hat{\rho}_1$ (note
the abuse of notation: while $\hat{\rho}$ denotes the
$\wedge$-operation on $\Omega^*(P)$ as introduced in Section
\ref{stableforms}, $\hat{\rho}_i$ is a mere notation for the forms
pulled back from the basis $M$). These have T-duals $\rho^t$ and
$\hat{\rho}^t$. Since $\rho$ and $\hat{\rho}$ have the same
stabiliser $G$ inside $Spin(TP\oplus T^*P)\cong Spin(7,7)$, it
follows that $\rho^t$ and $\hat{\rho}^t$ are stabilised by the same
$G^t$ inside $Spin(TP^t\oplus T^*P^t)\cong Spin(7,7)$ which is hence
isomorphic to $G_2\times G_2$. By invariance, $\hat{\rho}^t$ and
$\hat{\rho^t}$ coincide up to a constant (which we henceforth
ignore). The integrability condition transforms as follows.

\begin{lem}
$$
d_T\rho=\lambda\hat{\rho}\mbox{ if and only if
}d_T^t\rho^t=-\lambda\widehat{\rho^t}
$$
\end{lem}

\begin{prf}
Since $T=-\theta\wedge \mc{F}^t+\mc{T}$, expansion of the left hand
side yields
\begin{eqnarray*}
d_T\rho & = & d\theta\wedge\rho_0-\theta\wedge d\rho_0+d\rho_1-\theta\wedge\mc{F}^t\wedge\rho_1+\mc{T}\wedge\theta\wedge\rho_0+\mc{T}\wedge\rho_1\\
& = &
\mc{F}\wedge\rho_0+d\rho_1+\mc{T}\wedge\rho_1+\theta\wedge(-\mc{F}^t\wedge\rho_1-d\rho_0-\mc{T}\wedge\rho_0)\\
& = & \theta\wedge\lambda\hat{\rho}_0+\lambda\hat{\rho}_1,
\end{eqnarray*}
and this is equivalent to
\begin{equation}\label{weaktdual}
\lambda\hat{\rho}_0=-\mc{F}^t\wedge\rho_1-d\rho_0-\mc{T}\wedge\rho_0\mbox{
and
}\lambda\hat{\rho}_1=\mc{F}\wedge\rho_0+d\rho_1+\mc{T}\wedge\rho_1.
\end{equation}
On the other hand, we have
\begin{eqnarray*}
d_{T^t}\rho^t & = & \phantom{-}d\theta^t\wedge\rho_1-\theta^t\wedge
d\rho_1+d\rho_0-\theta^t\wedge \mc{F}\wedge\rho_0+\mc{T}\wedge
\theta^t\wedge\rho_1+\mc{T}\wedge\rho_0\\
& = &
\phantom{-}\mc{F}^t\wedge\rho_1+d\rho_0+\mc{T}\wedge\rho_0+\theta^t\wedge(-\mc{F}\wedge\rho_0-d\rho_1-\mc{T}\wedge\rho_1)\\
& = & -\theta^t\wedge\lambda\hat{\rho}_1-\lambda\hat{\rho}_0
\end{eqnarray*}
which gives precisely~(\ref{weaktdual}).
\end{prf}

\begin{cor}
Assume that $p:P\to M$ is a principal $S^1$-bundle with connection
form $\theta$, $X$ the vertical vector field with $\theta(X)=1$, and
that $T$ is a closed integral 3-form such that $X\llcorner T$ is
also closed and integral.

{\rm (i)} A generalised $G_2$-structure $\rho$ is closed strongly
integrable with torsion $T$ if and only if $\rho^t$ is a closed
strongly integrable $G_2$-structure with torsion $T^t$,
$$
\begin{array}{c}\mbox{closed strongly int. gen. $G_2$ with torsion $T$ on
$P$}\\\stackrel{t}{\longleftrightarrow}\\\mbox{closed strongly int.
gen. $G_2$ with torsion $T^t$ on $P^t$}.
\end{array}
$$
{\rm (ii)} A generalised $G_2$-structure $\rho$ is closed weakly
integrable of even or odd type with torsion $T$ and Killing number
$\lambda$ if and only if $\rho^t$ is a closed weakly integrable
$G_2$-structure of odd or even type with torsion $T^t$ and Killing
number $-\lambda$,
$$
\begin{array}{c}\mbox{closed even/odd weakly int. gen. $G_2$ with torsion $T$,
Killing number $\lambda$ on
$P$}\\\stackrel{t}{\longleftrightarrow}\\\mbox{closed odd/even
weakly int. gen. $G_2$ with torsion $T^t$, Killing number $-\lambda$
on $P^t$}.
\end{array}
$$
{\rm (iii)} A generalised $Spin(7)$-structure $\rho$ is closed
integrable of even or odd type with torsion $T$ if and only if
$\rho^t$ is a closed integrable $Spin(7)$-structure of odd or even
type with torsion $T^t$,
$$
\begin{array}{c}\mbox{closed even/odd int. gen. $Spin(7)$ with torsion $T$ on
$P$}\\\stackrel{t}{\longleftrightarrow}\\\mbox{closed odd/even int.
gen. $Spin(7)$ with torsion $T^t$ on $P^t$}.
\end{array}
$$
\end{cor}

We put the previous corollary into action as follows. Start with a
non-trivial principal $S^1$-fibre bundle $P$ with connection form
$\theta$ and which admits a metric of holonomy $G_2$ or $Spin(7)$.
We put $T=0$ so that $\mc{F}^t$ -- the curvature of the T-dual
bundle $P^t$ -- vanishes. This implies the triviality of $P^t$. The
resulting straight generalised structure is strongly integrable and
so is its dual, but according to the T-duality rules, we acquire
non-trivial torsion given by $T^t=\theta^t\wedge\mc{F}$. Local
examples of such $G_2$-structures exist in abundance~\cite{apsa04}.
Applying Lemma~\ref{g2spin7} yields local examples of closed
integrable generalised $Spin(7)$-structures of even type over a
trivial $S^1$-fibration $P=M^7\times S^1$. The torsion $T$ is
integral and contracting with the dual of a connection form on $P$
yields 0. According to the T-duality rules we obtain a closed
integrable generalised $Spin(7)$-manifold of odd type whose torsion
equals $T^t=T$.

\begin{prp}
There exist local examples of non-trivial closed strongly integrable
generalised $G_2$- and $Spin(7)$-manifolds of either type.
\end{prp}

\end{document}

%% file: fig1.pstex_t
\begin{picture}(0,0)%
\includegraphics{fig1.pstex}%
\end{picture}%
\setlength{\unitlength}{2763sp}%
\begingroup\makeatletter\ifx\SetFigFont\undefined%
\gdef\SetFigFont#1#2#3#4#5{%
  \reset@font\fontsize{#1}{#2pt}%
  \fontfamily{#3}\fontseries{#4}\fontshape{#5}%
  \selectfont}%
\fi\endgroup%
\begin{picture}(3400,3240)(3268,-5194)
\put(6076,-2536){\makebox(0,0)[lb]{\smash{{\SetFigFont{8}{9.6}{\rmdefault}{\mddefault}{\updefault}$V_+$}}}}
\put(6376,-3736){\makebox(0,0)[lb]{\smash{{\SetFigFont{8}{9.6}{\rmdefault}{\mddefault}{\updefault}$T$}}}}
\put(6076,-4936){\makebox(0,0)[lb]{\smash{{\SetFigFont{8}{9.6}{\rmdefault}{\mddefault}{\updefault}$V_-$}}}}
\put(5476,-3811){\makebox(0,0)[lb]{\smash{{\SetFigFont{8}{9.6}{\rmdefault}{\mddefault}{\updefault}$v$}}}}
\put(4426,-2986){\makebox(0,0)[lb]{\smash{{\SetFigFont{8}{9.6}{\rmdefault}{\mddefault}{\updefault}$Pv$}}}}
\put(4726,-2086){\makebox(0,0)[lb]{\smash{{\SetFigFont{8}{9.6}{\rmdefault}{\mddefault}{\updefault}$T^*$}}}}
\end{picture}%